\documentclass{amsart}
\pdfoutput=1
\usepackage[encapsulated]{CJK}
\usepackage{ucs}
\usepackage[utf8x]{inputenc}
\usepackage{amssymb,latexsym,graphicx}
\usepackage[section]{placeins}
\usepackage{enumerate}
\usepackage[T1]{fontenc}
\usepackage{mathrsfs}
\usepackage{algorithmic}
\usepackage{algorithm}
\usepackage{pgf,tikz}
\usetikzlibrary{arrows,shapes,positioning,calc}
\usepackage{cite}
\usepackage{longtable}
\usepackage{endnotes}

\theoremstyle{plain}
\newtheorem{theorem}{Theorem}

\theoremstyle{definition}

\theoremstyle{remark}

\newtheorem*{remark*}{Remark}

\newtheorem*{example*}{Example}

\begin{document}\begin{CJK*}{UTF8}{min}
\title[$11$ crossing knot bridges]{Minimal bridge projections for $11$-crossing prime knots}
\author{Chad Musick}

\begin{abstract}We give the bridge indices for $11$-crossing prime knots and give a minimal bridge projection for each of these knots. The results on the indices may be easily summarized: all of these knots that are not rational knots or Montesinos knots have bridge index three.
\end{abstract}

\maketitle

\section{introduction}

A knot is the image of a smooth proper embedding $S^1 \hookrightarrow S^3$. We view a knot diagram of a knot $K$ as a projection in general position of $K$ to $S^2$ such that there are finitely many double points and each has a transverse crossing as a pre-image. At each crossing, we designate one of the arcs in the neighborhood as lying above and the other as lying below. We adopt one of the standard definitions of the \emph{bridge index} of a knot $K$ as being the smallest number of over arcs in any knot diagram of $K$. These over arcs are the \emph{bridges} of the diagram.

To describe knot diagrams, we adopt the DT code of Dowker and Thistlewaite \cite{dtcode} for reasons of clarity and convenience. The DT code of a knot diagram is given by choosing an arbitrary single-point on the diagram and, proceeding along the orientation of the knot, labeling each crossing with consecutive integers, beginning at $1$, as they are encountered. If an even label is applied as a result of passing under, it is given a negative sign. This will result in each crossing being labeled with a pair of integers, one even and one odd. The DT code is then the sequence of even integers obtained from listing the even integer associated to the odd integers $1, 3, 5, \dots$ in ascending order. The diagram resulting from a valid DT code is unique up to isotopy. We note here that this is not true for diagrams in $\mathbb{R}^2$, which is our motivation for viewing diagrams as lying in $S^2$.

We are interested in the bridge index of the set of prime knots that have minimal crossing number $11$. Of those knots, each has a bridge index of two, three, or four. For knots with bridge index two, the presentation in two bridge form is sufficient to establish the bridge index. For knots with bridge index three, we note that none of these are rational knots, a fact which may be checked for each but is not shown in this paper, and exhibit a three bridge presentation. 

For knots with bridge index four, we use the following theorem, from Boileau and Zieschang \cite{boileau}:
\begin{theorem}[Th\'eor\`eme 1.1, {\cite{boileau}}]\label{T:main}
Soit $L = m(0 | e; (\alpha_1, \beta_1), \dots, (\alpha_r, \beta_r))$ un entrelacs de Montesinos \`a $r \geq 3$ branches, alors $w(L) = b(L) = r$.
\end{theorem}

This theorem says that if a Montesinos knot $K$ has $r \geq 3$ rational tangles, excluding integer tangles, then $K$ has bridge index $r$. We show that all of the knots we identify as four bridge knots are Montesinos knots with $r = 4$. Hence, by the theorem, their bridge indices are four.

We identify knots by their number as given in Knot Atlas \cite{knotatlas} and KnotInfo \cite{knotinfo}.

The following knots have bridge index two: $11a13$, $11a59$, $11a65$, $11a75$, $11a77$, $11a84$, $11a85$, $11a89$, $11a90$, $11a91$, $11a93$, $11a95$, $11a96$, $11a98$, $11a110$, $11a111$, $11a117$, $11a119$, $11a120$, $11a121$, $11a140$, $11a144$, $11a145$, $11a154$, $11a159$, $11a166$, $11a174$, $11a175$, $11a176$, $11a177$, $11a178$, $11a179$, $11a180$, $11a182$, $11a183$, $11a184$, $11a185$, $11a186$, $11a188$, $11a190$, $11a191$, $11a192$, $11a193$, $11a195$, $11a203$, $11a204$, $11a205$, $11a206$, $11a207$, $11a208$, $11a210$, $11a211$, $11a220$, $11a224$, $11a225$, $11a226$, $11a229$, $11a230$, $11a234$, $11a235$, $11a236$, $11a238$, $11a242$, $11a243$, $11a246$, $11a247$, $11a306$, $11a307$, $11a308$, $11a309$, $11a310$, $11a311$, $11a333$, $11a334$, $11a335$, $11a336$, $11a337$, $11a339$, $11a341$, $11a342$, $11a343$, $11a355$, $11a356$, $11a357$, $11a358$, $11a359$, $11a360$, $11a363$, $11a364$, $11a365$, and $11a367$.

The following knots have bridge index four: $11a43$, $11a44$, $11a47$, $11a57$, $11a231$, $11a263$, $11n71$, $11n72$, $11n73$, $11n74$, $11n75$, $11n76$, $11n77$, $11n78$, and $11n81$.

The remaining $11$-crossing prime knots have bridge index three.

The author is grateful to Yeonhee Jang for pointing out that the first few knots identified as $4$-bridge are Montesinos knots; this simplified the proof of the bridge index for these knots. In fact, all of the $4$-bridge knots among the $11$-crossing prime knots are Montesinos knots.

\section{The four bridge knots}

\begin{figure}[ht]
\begin{tikzpicture}[cm={0,1,-1,0,(0,0)}]
\draw (0, -.5) -- (0, .5) -- (1, .5) -- (1, -0.5) -- cycle;
\draw (.25, .5) -- (.25, 6) -- (2.75, 6) -- (2.75, 5) (.75, .5) -- (.75, 5.5) -- (2.25, 5.5) -- (2.25, 5);
\draw (2, 5) -- (3, 5) -- (3, 4) -- (2, 4) -- cycle;
\draw (2.25, 4) -- (2.25, 3.5) (2.75, 4) -- (2.75, 3.5);
\draw (2, 3.5) -- (3, 3.5) -- (3, 2.5) -- (2, 2.5) -- cycle;
\draw (2.25, 2.5) -- (2.25, 2) (2.75, 2.5) -- (2.75, 2);
\draw (2, -.5) -- (2, .5) -- (3, .5) -- (3, -.5) -- cycle;
\draw (2.25, .5) -- (2.25, 1) (2.75, .5) -- (2.75, 1);
\draw[dotted] (2.25, 1) -- (2.25, 2.5) (2.75, 1) -- (2.75, 2.5);
\draw (.25, -.5) -- (.25, -1.5) -- (2.75, -1.5) -- (2.75, -.5);
\draw (.75, -.5) -- (.75, -1) -- (2.25, -1) -- (2.25, -.5);
\node[anchor=center] at (.5, 0) {$\alpha_0$};
\node at (2.5, 4.5) {$\alpha_1 / \beta_1$};
\node at (2.5, 3) {$\alpha_2 / \beta_2$};
\node at (2.5, 0) {$\alpha_r / \beta_r$};
\end{tikzpicture}
\caption{A Montesinos link. Each box is a rational tangle.}\label{F:montesinos}
\end{figure}
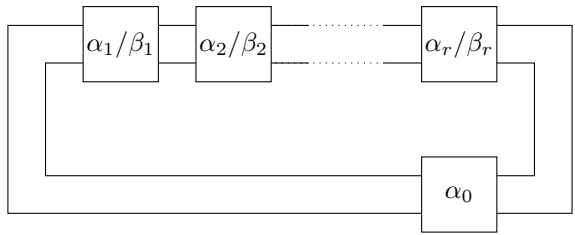

Montesinos knots are first defined in Montesinos \cite{montesinos}, and we describe them here. A Montesinos knot is a knot that is of the form shown in Fig. \ref{F:montesinos}, where each box represents a rational tangle, the box $\alpha_0$ is an integer tangle, and none of the other tangle boxes are integer tangles. For the definition of rational tangles we refer to Goldman and Kauffman \cite{tangles}.

In the previous section we claimed that knots $11a43$, $11a44$, $11a47$, $11a57$, $11a231$, $11n71\sim 11n78$ and $11n81$ are all Montesinos knots with $r = 4$. We may write any Montesinos knot as $(\alpha_0; \alpha_1 / \beta_1, \alpha_2 / \beta_2, \dots, \alpha_r / \beta_r)$. We do this and see that:
\begin{align*}
11a43 &= (0; -1/2, -2/3, -2/3, -2/3),\\
11a44 &= (0; -1/2, -1/3, -2/3, -2/3),\\
11a47 &= (0; -1/2, -2/3, -1/3, -2/3),\\
11a57 &= (0; -1/2, -1/3, -1/3, -2/3),\\
11a231 &= (0; -1/3, -1/2, -1/3, -2/3),\\
11a263 &= (0; 1/3, 1/2, 1/3, 1/3),\\
11n71 &= (0; -1/2, -2/3, 2/3, -2/3),\\
11n72 &= (0; 1/2, -2/3, -2/3, -2/3),\\
11n73 &= (0; 1/2, 2/3, -2/3, -2/3),\\
11n74 &= (0; 1/2, -2/3, 2/3, -2/3),\\
11n75 &= (0; 1/2, 2/3, 2/3, -2/3),\\
11n76 &= (0; -1/2, -1/3, 2/3, -2/3),\\
11n77 &= (0; 1/2, 1/3, -2/3, -2/3),\\
11n78 &= (0; 1/2, 1/3, 2/3, -2/3),\\
11n81 &= (0; 1/2, 1/3, 1/3, -2/3)
\end{align*}

and so by Theorem \ref{T:main} these knots have bridge index four.

\section{Minimal bridge presentations}

We give here the DT code for a minimal bridge presentation of each of the knots listed earlier. The format is as follows:
\[ \mbox{knot identifier} : \mbox{bridge index}\ |\ \mbox{knot DT code}\ |\ \mbox{DT code of bridge presentation}. \]

We do not prove here that each of the provided bridge positions represents the named knot. That each does so may be easily verified using the SnapPy \cite{snappy} tool, except for $11a367$, which is not hyperbolic. That knot may be verified using the continued fraction from rational tangles instead.\\

\noindent$11a1$:  $3$ | $(4$, $8$, $10$, $14$, $2$, $16$, $20$, $6$, $22$, $12$, $18)$ | $(12$, $16$, $58$, $60$, $14$, $-92$, $-90$, $-94$, $-32$, $-40$, $120$, $102$, $108$, $112$, $98$, $116$, $124$, $106$, $104$, $122$, $118$, $100$, $110$, $-26$, $-34$, $-38$, $-22$, $-42$, $-30$, $-96$, $-28$, $-44$, $-24$, $-36$, $50$, $18$, $56$, $62$, $46$, $64$, $54$, $20$, $52$, $66$, $48$, $128$, $126$, $114$, $-6$, $-74$, $-80$, $-86$, $-68$, $-88$, $-78$, $-76$, $-8$, $-4$, $-72$, $-82$, $-84$, $-70$, $-2$, $-10)$\\

\noindent$11a2$ : $3$ | $(4$, $8$, $10$, $14$, $2$, $18$, $6$, $20$, $12$, $22$, $16)$ | $(60$, $78$, $46$, $74$, $64$, $56$, $124$, $54$, $66$, $72$, $48$, $118$, $114$, $120$, $50$, $70$, $68$, $52$, $122$, $58$, $62$, $76$, $-94$, $-132$, $-110$, $-140$, $-104$, $-126$, $-100$, $-88$, $-90$, $-98$, $-128$, $-106$, $-142$, $-108$, $-130$, $-96$, $-92$, $-86$, $-26$, $-24$, $116$, $152$, $182$, $150$, $166$, $168$, $148$, $180$, $154$, $158$, $176$, $144$, $172$, $162$, $82$, $-136$, $-134$, $-112$, $-138$, $-102$, $156$, $178$, $146$, $170$, $164$, $84$, $80$, $160$, $174$, $-32$, $-18$, $-8$, $-42$, $-2$, $-38$, $-12$, $-14$, $-36$, $-28$, $-22$, $-4$, $-44$, $-6$, $-20$, $-30$, $-34$, $-16$, $-10$, $-40)$\\

\noindent$11a3$ : $3$ | $(4$, $8$, $10$, $14$, $2$, $18$, $6$, $20$, $22$, $12$, $16)$ | $(22$, $74$, $76$, $20$, $24$, $32$, $-50$, $-62$, $-68$, $-56$, $-44$, $-46$, $-54$, $-66$, $-64$, $-52$, $-48$, $-60$, $-70$, $-58$, $-4$, $102$, $120$, $100$, $110$, $112$, $98$, $118$, $104$, $72$, $108$, $114$, $96$, $116$, $106$, $-90$, $-78$, $-42$, $40$, $34$, $14$, $30$, $26$, $18$, $38$, $36$, $16$, $28$, $-84$, $-10$, $-2$, $-6$, $-88$, $-92$, $-80$, $-12$, $-82$, $-94$, $-86$, $-8)$\\

\noindent$11a4$ : $3$ | $(4$, $8$, $10$, $14$, $2$, $18$, $6$, $22$, $20$, $12$, $16)$ | $(108$, $52$, $144$, $168$, $178$, $154$, $150$, $174$, $172$, $148$, $156$, $180$, $166$, $142$, $54$, $106$, $110$, $50$, $146$, $170$, $176$, $152$, $-12$, $-20$, $-134$, $-78$, $-84$, $-140$, $104$, $112$, $48$, $122$, $98$, $118$, $44$, $116$, $100$, $124$, $132$, $160$, $184$, $162$, $130$, $126$, $102$, $114$, $46$, $120$, $-16$, $-40$, $-8$, $-24$, $-138$, $-82$, $-80$, $-136$, $-22$, $-10$, $-42$, $-14$, $-18$, $-38$, $-6$, $-26$, $-28$, $-4$, $-36$, $158$, $182$, $164$, $128$, $-86$, $-76$, $-62$, $-94$, $-68$, $-70$, $-92$, $-60$, $-34$, $-2$, $-30$, $-56$, $-88$, $-74$, $-64$, $-96$, $-66$, $-72$, $-90$, $-58$, $-32)$\\

\noindent$11a5$ : $3$ | $(4$, $8$, $10$, $14$, $2$, $20$, $16$, $6$, $22$, $12$, $18)$  | $(42$, $52$, $32$, $56$, $38$, $46$, $48$, $36$, $92$, $90$, $34$, $50$, $44$, $40$, $54$, $-74$, $-84$, $-64$, $-100$, $-70$, $-78$, $-80$, $-68$, $-102$, $-66$, $-82$, $-76$, $-72$, $-86$, $-88$, $-18$, $114$, $124$, $104$, $128$, $110$, $118$, $120$, $108$, $130$, $106$, $122$, $116$, $94$, $-96$, $-62$, $-98$, $58$, $60$, $112$, $126$, $-12$, $-24$, $-2$, $-28$, $-8$, $-16$, $-20$, $-6$, $-30$, $-4$, $-22$, $-14$, $-10$, $-26)$\\

\noindent$11a6$ : $3$ | $(4$, $8$, $10$, $16$, $2$, $18$, $20$, $6$, $14$, $22$, $12)$ | $(32$, $36$, $54$, $48$, $42$, $60$, $40$, $50$, $52$, $38$, $58$, $44$, $46$, $56$, $124$, $-70$, $-72$, $-74$, $-68$, $-84$, $-62$, $-80$, $-112$, $-78$, $-64$, $-86$, $-66$, $-76$, $-114$, $-82$, $96$, $108$, $88$, $104$, $122$, $126$, $120$, $102$, $90$, $110$, $94$, $98$, $106$, $-16$, $-6$, $-26$, $-22$, $-10$, $-12$, $-20$, $-28$, $-4$, $-18$, $-14$, $-8$, $-24$, $92$, $100$, $30$, $34$, $-2$, $-116$, $-118)$\\

\noindent$11a7$ : $3$ | $(4$, $8$, $10$, $16$, $2$, $18$, $20$, $6$, $22$, $12$, $14)$ | $(50$, $40$, $32$, $28$, $36$, $44$, $46$, $54$, $76$, $52$, $48$, $42$, $34$, $-86$, $-100$, $-98$, $-88$, $-84$, $-102$, $-96$, $-90$, $-82$, $-104$, $-80$, $-92$, $-94$, $-78$, $114$, $38$, $30$, $-6$, $-26$, $-8$, $-58$, $-4$, $-24$, $-10$, $-56$, $112$, $116$, $106$, $72$, $64$, $62$, $70$, $108$, $118$, $110$, $68$, $60$, $66$, $74$, $-12$, $-22$, $-2$, $-18$, $-16$, $-14$, $-20)$\\

\noindent$11a8$ : $3$ | $(4$, $8$, $10$, $16$, $2$, $18$, $20$, $6$, $22$, $14$, $12)$ | $(84$, $146$, $158$, $72$, $166$, $138$, $92$, $98$, $250$, $100$, $90$, $140$, $164$, $70$, $160$, $144$, $86$, $82$, $148$, $156$, $74$, $168$, $136$, $94$, $96$, $134$, $170$, $76$, $154$, $150$, $80$, $88$, $142$, $162$, $-198$, $-240$, $-226$, $-184$, $-212$, $-206$, $-190$, $-232$, $-234$, $-192$, $-204$, $-246$, $-220$, $-178$, $-218$, $-248$, $-126$, $-52$, $-44$, $-10$, $-24$, $-66$, $-30$, $-4$, $-38$, $-58$, $152$, $78$, $172$, $132$, $-16$, $-18$, $-216$, $-180$, $-222$, $-244$, $-202$, $-194$, $-236$, $-230$, $-188$, $-208$, $-210$, $-186$, $-228$, $-238$, $-196$, $-200$, $-242$, $-224$, $-182$, $-214$, $-20$, $-14$, $-48$, $254$, $104$, $288$, $274$, $118$, $268$, $294$, $262$, $112$, $280$, $282$, $110$, $260$, $292$, $270$, $120$, $272$, $290$, $102$, $252$, $176$, $256$, $106$, $286$, $276$, $116$, $266$, $296$, $264$, $114$, $278$, $284$, $108$, $258$, $174$, $130$, $-128$, $-50$, $-46$, $-12$, $-22$, $-64$, $-32$, $-2$, $-36$, $-60$, $-122$, $-56$, $-40$, $-6$, $-28$, $-68$, $-26$, $-8$, $-42$, $-54$, $-124$, $-62$, $-34)$\\

\noindent$11a9$ : $3$ | $(4$, $8$, $10$, $16$, $2$, $18$, $20$, $22$, $6$, $12$, $14)$ | $(18$, $26$, $24$, $20$, $16$, $14$, $-58$, $-56$, $-36$, $-38$, $-30$, $-32$, $-40$, $-34$, $44$, $50$, $42$, $46$, $64$, $48$, $-4$, $-52$, $-54$, $-2$, $-6$, $22$, $28$, $62$, $60$, $-12$, $-10$, $-8)$\\

\noindent$11a10$ : $3$ | $(4$, $8$, $10$, $16$, $2$, $18$, $20$, $22$, $6$, $14$, $12)$ | $(34$, $48$, $46$, $36$, $56$, $38$, $44$, $50$, $112$, $116$, $54$, $40$, $42$, $52$, $114$, $-18$, $-64$, $-66$, $-58$, $-106$, $-94$, $-104$, $-60$, $-68$, $-62$, $-102$, $-96$, $-108$, $78$, $90$, $70$, $86$, $118$, $88$, $-4$, $-14$, $-26$, $-22$, $-10$, $-8$, $-20$, $-28$, $-16$, $-2$, $-6$, $-12$, $-24$, $74$, $82$, $30$, $84$, $72$, $92$, $76$, $80$, $32$, $-100$, $-98$, $-110)$\\

\noindent$11a11$ : $3$ | $(4$, $8$, $10$, $16$, $2$, $20$, $18$, $6$, $14$, $22$, $12)$ | $(36$, $28$, $40$, $32$, $92$, $94$, $100$, $88$, $98$, $96$, $90$, $34$, $38$, $-76$, $-70$, $-66$, $-80$, $-72$, $-74$, $-78$, $-68$, $30$, $56$, $60$, $-18$, $-12$, $-8$, $-22$, $-14$, $-16$, $-20$, $-10$, $52$, $110$, $108$, $102$, $114$, $104$, $106$, $112$, $54$, $62$, $50$, $58$, $-46$, $-82$, $-64$, $-84$, $-48$, $-86$, $-44$, $-2$, $-26$, $-4$, $-42$, $-6$, $-24)$\\

\noindent$11a12$ : $3$ | $(4$, $8$, $10$, $16$, $2$, $20$, $22$, $18$, $6$, $14$, $12)$ | $(62$, $108$, $106$, $102$, $112$, $-20$, $-2$, $-6$, $104$, $110$, $114$, $100$, $124$, $118$, $96$, $120$, $122$, $98$, $116$, $126$, $14$, $-4$, $-18$, $-8$, $-22$, $-36$, $-28$, $-30$, $-34$, $-24$, $-40$, $-38$, $-26$, $-32$, $56$, $52$, $64$, $60$, $48$, $10$, $46$, $16$, $44$, $12$, $50$, $58$, $66$, $54$, $-70$, $-90$, $-76$, $-82$, $-84$, $-42$, $-86$, $-80$, $-78$, $-88$, $-72$, $-94$, $-68$, $-92$, $-74)$\\

\noindent$11a13$ : $2$ | $(4$, $8$, $10$, $16$, $2$, $22$, $20$, $18$, $6$, $14$, $12)$ | $(84$, $110$, $62$, $106$, $88$, $80$, $114$, $66$, $102$, $92$, $76$, $118$, $70$, $98$, $96$, $72$, $120$, $74$, $94$, $100$, $68$, $116$, $78$, $90$, $104$, $64$, $112$, $82$, $86$, $108$, $-148$, $-158$, $-138$, $-168$, $-128$, $-178$, $-122$, $-174$, $-132$, $-164$, $-142$, $-154$, $-152$, $-144$, $-162$, $-134$, $-172$, $-124$, $-180$, $-126$, $-170$, $-136$, $-160$, $-146$, $-150$, $-156$, $-140$, $-166$, $-130$, $-176$, $204$, $230$, $182$, $226$, $208$, $200$, $234$, $186$, $222$, $212$, $196$, $238$, $190$, $218$, $216$, $192$, $240$, $194$, $214$, $220$, $188$, $236$, $198$, $210$, $224$, $184$, $232$, $202$, $206$, $228$, $-28$, $-38$, $-18$, $-48$, $-8$, $-58$, $-2$, $-54$, $-12$, $-44$, $-22$, $-34$, $-32$, $-24$, $-42$, $-14$, $-52$, $-4$, $-60$, $-6$, $-50$, $-16$, $-40$, $-26$, $-30$, $-36$, $-20$, $-46$, $-10$, $-56)$\\

\noindent$11a14$ : $3$ | $(4$, $8$, $12$, $2$, $14$, $18$, $6$, $20$, $10$, $22$, $16)$ | $(36$, $72$, $66$, $80$, $68$, $70$, $122$, $124$, $112$, $104$, $102$, $114$, $126$, $120$, $-12$, $-16$, $-24$, $-88$, $-86$, $-26$, $-14$, $30$, $78$, $64$, $74$, $34$, $38$, $28$, $40$, $32$, $76$, $-20$, $-92$, $-82$, $-96$, $-84$, $-90$, $-22$, $-18$, $-94$, $108$, $98$, $118$, $128$, $116$, $100$, $106$, $110$, $-10$, $-4$, $-48$, $-62$, $-46$, $-6$, $-8$, $-44$, $-60$, $-50$, $-2$, $-54$, $-56$, $-42$, $-58$, $-52)$\\

\noindent$11a15$ : $3$ | $(4$, $8$, $12$, $2$, $14$, $18$, $6$, $20$, $22$, $10$, $16)$ | $(44$, $154$, $144$, $34$, $140$, $158$, $48$, $160$, $138$, $36$, $146$, $152$, $42$, $46$, $156$, $142$, $-68$, $-54$, $-82$, $-86$, $-76$, $-60$, $-62$, $-74$, $122$, $98$, $102$, $108$, $92$, $116$, $168$, $118$, $94$, $106$, $104$, $96$, $120$, $166$, $114$, $90$, $110$, $162$, $164$, $112$, $-130$, $-22$, $-4$, $-30$, $-10$, $-16$, $-18$, $-8$, $-32$, $-6$, $-20$, $-128$, $-132$, $-24$, $-2$, $-28$, $-12$, $-14$, $-84$, $38$, $148$, $150$, $40$, $124$, $100$, $-52$, $-70$, $-66$, $-56$, $-80$, $-88$, $-78$, $-58$, $-64$, $-72$, $-50$, $-136$, $-126$, $-134$, $-26)$\\

\noindent$11a16$ : $3$ | $(4$, $8$, $12$, $2$, $14$, $18$, $6$, $22$, $20$, $10$, $16)$ | $(42$, $150$, $134$, $142$, $158$, $162$, $164$, $156$, $140$, $136$, $152$, $40$, $44$, $148$, $132$, $144$, $160$, $-10$, $-14$, $-80$, $-58$, $-60$, $-78$, $-12$, $36$, $102$, $118$, $108$, $92$, $168$, $94$, $110$, $116$, $100$, $34$, $98$, $114$, $112$, $96$, $166$, $90$, $106$, $120$, $104$, $-124$, $-22$, $-2$, $-26$, $-128$, $-32$, $-8$, $-16$, $-122$, $-20$, $-4$, $-28$, $-130$, $-30$, $-6$, $-18$, $138$, $154$, $38$, $46$, $146$, $-74$, $-64$, $-54$, $-84$, $-86$, $-52$, $-66$, $-72$, $-76$, $-62$, $-56$, $-82$, $-88$, $-50$, $-68$, $-70$, $-48$, $-126$, $-24)$\\

\noindent$11a17$ : $3$ | $(4$, $8$, $12$, $2$, $16$, $6$, $20$, $10$, $22$, $14$, $18)$ | $(102$, $128$, $106$, $34$, $66$, $36$, $46$, $56$, $62$, $40$, $42$, $60$, $58$, $44$, $38$, $64$, $-82$, $-68$, $-86$, $-116$, $-94$, $-114$, $-88$, $-122$, $-2$, $-126$, $-4$, $-120$, $-90$, $-112$, $-92$, $-118$, $-84$, $142$, $148$, $52$, $50$, $150$, $48$, $54$, $146$, $144$, $140$, $98$, $132$, $110$, $136$, $-20$, $-28$, $-12$, $-76$, $-74$, $-70$, $-80$, $-16$, $-24$, $134$, $96$, $138$, $108$, $130$, $100$, $104$, $-72$, $-78$, $-14$, $-26$, $-22$, $-18$, $-30$, $-10$, $-8$, $-32$, $-6$, $-124)$\\

\noindent$11a18$ : $3$ | $(4$, $8$, $12$, $2$, $16$, $6$, $20$, $18$, $10$, $22$, $14)$ | $(36$, $88$, $98$, $96$, $90$, $104$, $108$, $110$, $102$, $92$, $94$, $100$, $112$, $106$, $-12$, $-16$, $-24$, $-84$, $-82$, $-26$, $-14$, $30$, $72$, $64$, $76$, $34$, $38$, $28$, $40$, $32$, $74$, $-20$, $-8$, $-4$, $-2$, $-10$, $-18$, $-22$, $-6$, $66$, $70$, $114$, $68$, $-86$, $-80$, $-46$, $-62$, $-48$, $-78$, $-50$, $-60$, $-44$, $-54$, $-56$, $-42$, $-58$, $-52)$\\

\noindent$11a19$ : $3$ | $(4$, $8$, $12$, $2$, $16$, $18$, $6$, $20$, $22$, $14$, $10)$ | $(84$, $34$, $42$, $128$, $130$, $112$, $132$, $126$, $122$, $136$, $116$, $66$, $114$, $134$, $124$, $-16$, $-104$, $-98$, $-22$, $-24$, $-96$, $-106$, $32$, $86$, $82$, $36$, $40$, $78$, $90$, $30$, $88$, $80$, $38$, $-52$, $-4$, $-44$, $-6$, $-10$, $-12$, $-26$, $-20$, $-100$, $-102$, $-18$, $-28$, $-14$, $-8$, $74$, $68$, $118$, $138$, $120$, $70$, $72$, $92$, $76$, $-54$, $-64$, $-50$, $-2$, $-46$, $-60$, $-58$, $-108$, $-94$, $-110$, $-56$, $-62$, $-48)$\\

\noindent$11a20$ : $3$ | $(4$, $8$, $12$, $2$, $18$, $6$, $20$, $22$, $10$, $16$, $14)$ | $(38$, $70$, $104$, $72$, $40$, $36$, $-96$, $-78$, $-86$, $-88$, $-76$, $-4$, $-6$, $-74$, $-90$, $-84$, $-80$, $-94$, $-98$, $-100$, $-92$, $-82$, $32$, $16$, $20$, $28$, $26$, $22$, $14$, $34$, $42$, $30$, $18$, $-48$, $-56$, $-102$, $-52$, $106$, $68$, $112$, $124$, $122$, $110$, $66$, $108$, $120$, $126$, $114$, $116$, $128$, $118$, $24$, $-54$, $-50$, $-64$, $-46$, $-58$, $-12$, $-2$, $-8$, $-62$, $-44$, $-60$, $-10$)\\

\noindent$11a21$ : $3$ | $(4$, $8$, $12$, $2$, $18$, $6$, $22$, $20$, $10$, $16$, $14)$ | $(26$, $82$, $100$, $86$, $96$, $90$, $92$, $94$, $88$, $98$, $84$, $102$, $-50$, $-52$, $-22$, $-6$, $-18$, $-10$, $-14$, $64$, $34$, $68$, $30$, $72$, $74$, $104$, $76$, $70$, $32$, $66$, $36$, $-12$, $-16$, $-8$, $-20$, $-4$, $-80$, $-2$, $28$, $24$, $-48$, $-54$, $-44$, $-58$, $-40$, $-62$, $-38$, $-60$, $-42$, $-56$, $-46$, $-78)$\\

\noindent$11a22$ : $3$ | $(4$, $8$, $12$, $2$, $18$, $14$, $6$, $20$, $22$, $10$, $16)$ | $(40$, $26$, $44$, $36$, $30$, $48$, $32$, $34$, $46$, $28$, $38$, $42$, $-62$, $-56$, $-80$, $-50$, $-76$, $-52$, $-66$, $-58$, $-60$, $-64$, $-54$, $-78$, $98$, $74$, $102$, $94$, $88$, $106$, $90$, $92$, $104$, $86$, $84$, $-68$, $-8$, $-14$, $100$, $96$, $72$, $70$, $-82$, $-20$, $-2$, $-24$, $-4$, $-18$, $-10$, $-12$, $-16$, $-6$, $-22)$\\

\noindent$11a23$ : $3$ | $(4$, $8$, $12$, $2$, $18$, $14$, $6$, $22$, $20$, $10$, $16)$ | $(54$, $24$, $28$, $26$, $22$, $82$, $84$, $74$, $86$, $80$, $-68$, $-36$, $-30$, $-38$, $40$, $78$, $88$, $76$, $42$, $-6$, $-4$, $-8$, $-20$, $-10$, $-14$, $-16$, $-34$, $-32$, $-18$, $-12$, $48$, $52$, $56$, $44$, $46$, $58$, $50$, $-62$, $-2$, $-66$, $-70$, $-60$, $-72$, $-64)$\\

\noindent$11a24$ : $3$ | $(4$, $8$, $12$, $2$, $18$, $16$, $6$, $20$, $22$, $14$, $10)$ | $(86$, $94$, $74$, $26$, $78$, $98$, $82$, $30$, $32$, $84$, $96$, $76$, $-122$, $-58$, $-66$, $-130$, $-68$, $-56$, $-72$, $-4$, $-22$, $-10$, $-100$, $-12$, $-20$, $-2$, $-70$, $90$, $134$, $140$, $44$, $46$, $142$, $132$, $88$, $92$, $-110$, $-120$, $-124$, $-60$, $-64$, $-128$, $-116$, $-114$, $-106$, $-108$, $-112$, $-118$, $-126$, $-62$, $80$, $28$, $34$, $54$, $36$, $38$, $52$, $146$, $50$, $40$, $136$, $138$, $42$, $48$, $144$, $-16$, $-104$, $-6$, $-24$, $-8$, $-102$, $-14$, $-18)$\\

\noindent$11a25$ : $3$ | $(4$, $8$, $12$, $2$, $18$, $20$, $6$, $10$, $22$, $14$, $16)$ | $(116$, $104$, $94$, $86$, $158$, $160$, $88$, $92$, $106$, $114$, $118$, $102$, $96$, $84$, $98$, $100$, $120$, $112$, $108$, $90$, $-14$, $-6$, $-26$, $-28$, $-8$, $-12$, $-38$, $-16$, $-4$, $-24$, $-30$, $156$, $162$, $146$, $52$, $40$, $54$, $148$, $164$, $154$, $60$, $46$, $-140$, $-126$, $-70$, $-72$, $-128$, $-138$, $-82$, $-62$, $-80$, $-136$, $-130$, $-74$, $-68$, $-66$, $-76$, $-132$, $-134$, $-78$, $-64$, $-10$, $50$, $42$, $56$, $150$, $166$, $152$, $58$, $44$, $48$, $122$, $110$, $-36$, $-18$, $-2$, $-22$, $-32$, $-142$, $-124$, $-144$, $-34$, $-20)$\\

\noindent$11a26$ : $3$ | $(4$, $8$, $12$, $2$, $20$, $18$, $6$, $10$, $22$, $14$, $16)$ | $(114$, $124$, $102$, $122$, $116$, $48$, $40$, $54$, $70$, $68$, $56$, $42$, $46$, $118$, $120$, $104$, $126$, $112$, $166$, $-94$, $-134$, $-84$, $-82$, $-136$, $-96$, $-36$, $-92$, $-132$, $-86$, $-30$, $-6$, $-4$, $-32$, $-88$, $-130$, $-90$, $-34$, $-2$, $-8$, $-28$, $174$, $160$, $172$, $106$, $128$, $110$, $168$, $164$, $178$, $76$, $62$, $-142$, $-156$, $-20$, $-16$, $-152$, $-146$, $-98$, $-138$, $-80$, $-140$, $-100$, $-144$, $-154$, $-18$, $108$, $170$, $162$, $176$, $78$, $60$, $64$, $74$, $50$, $38$, $52$, $72$, $66$, $58$, $44$, $-24$, $-12$, $-148$, $-150$, $-14$, $-22$, $-158$, $-26$, $-10)$\\

\noindent$11a27$ : $3$ | $(4$, $8$, $12$, $14$, $2$, $18$, $22$, $20$, $10$, $16$, $6)$ | $(158$, $186$, $92$, $174$, $170$, $88$, $190$, $154$, $208$, $242$, $64$, $224$, $226$, $62$, $240$, $210$, $152$, $192$, $86$, $168$, $176$, $94$, $184$, $160$, $156$, $188$, $90$, $172$, $-8$, $-40$, $-198$, $-100$, $-204$, $-46$, $-2$, $-50$, $-14$, $-34$, $-32$, $-16$, $-18$, $-30$, $-148$, $-110$, $-132$, $-126$, $-116$, $-142$, $-24$, $180$, $164$, $82$, $78$, $74$, $214$, $236$, $58$, $230$, $220$, $68$, $246$, $70$, $218$, $232$, $56$, $234$, $216$, $72$, $244$, $66$, $222$, $228$, $60$, $238$, $212$, $76$, $-106$, $-136$, $-122$, $-120$, $-138$, $-20$, $-28$, $-146$, $-112$, $-130$, $-128$, $-114$, $-144$, $-26$, $-22$, $-140$, $-118$, $-124$, $-134$, $-108$, $-150$, $84$, $166$, $178$, $96$, $182$, $162$, $80$, $-104$, $-194$, $-36$, $-12$, $-52$, $-4$, $-44$, $-202$, $-98$, $-200$, $-42$, $-6$, $-54$, $-10$, $-38$, $-196$, $-102$, $-206$, $-48)$\\

\noindent$11a28$ : $3$ | $(4$, $8$, $12$, $16$, $2$, $18$, $20$, $22$, $10$, $6$, $14)$ | $(46$, $40$, $84$, $82$, $42$, $44$, $80$, $86$, $90$, $76$, $94$, $128$, $96$, $78$, $88$, $-16$, $-72$, $-22$, $-24$, $-70$, $-64$, $-102$, $-66$, $-68$, $-26$, $-20$, $-74$, $-18$, $-28$, $-14$, $-6$, $132$, $100$, $136$, $98$, $130$, $92$, $-110$, $-124$, $-118$, $-104$, $-62$, $-106$, $-120$, $-122$, $-108$, $-112$, $-126$, $-116$, $-2$, $-10$, $134$, $60$, $32$, $54$, $52$, $34$, $38$, $48$, $58$, $30$, $56$, $50$, $36$, $-114$, $-4$, $-8$, $-12)$\\

\noindent$11a29$ : $3$ | $(4$, $8$, $12$, $16$, $2$, $18$, $22$, $20$, $10$, $6$, $14)$ | $(104$, $110$, $78$, $126$, $88$, $234$, $98$, $116$, $72$, $120$, $94$, $238$, $92$, $122$, $74$, $114$, $100$, $232$, $86$, $128$, $80$, $108$, $106$, $82$, $130$, $84$, $230$, $102$, $112$, $76$, $124$, $90$, $236$, $96$, $118$, $-152$, $-244$, $-164$, $-140$, $-256$, $-176$, $-132$, $-172$, $-252$, $-144$, $-160$, $-240$, $-156$, $-148$, $-248$, $-168$, $-136$, $-260$, $-178$, $-258$, $-138$, $-166$, $-246$, $-150$, $-154$, $-242$, $-162$, $-142$, $-254$, $-174$, $202$, $262$, $282$, $182$, $222$, $208$, $196$, $268$, $276$, $188$, $216$, $214$, $190$, $274$, $270$, $194$, $210$, $220$, $184$, $280$, $264$, $200$, $204$, $226$, $-44$, $-2$, $-56$, $-32$, $-14$, $-68$, $-20$, $-26$, $-62$, $-8$, $-38$, $-50$, $-48$, $-40$, $-6$, $-60$, $-28$, $-18$, $-70$, $-16$, $-30$, $-58$, $-4$, $-42$, $-46$, $-134$, $-170$, $-250$, $-146$, $-158$, $272$, $192$, $212$, $218$, $186$, $278$, $266$, $198$, $206$, $224$, $180$, $228$, $-52$, $-36$, $-10$, $-64$, $-24$, $-22$, $-66$, $-12$, $-34$, $-54)$\\

\noindent$11a30$ : $3$ | $(4$, $8$, $12$, $18$, $2$, $16$, $20$, $6$, $10$, $22$, $14)$ | $(128$, $148$, $108$, $160$, $166$, $114$, $142$, $134$, $122$, $62$, $256$, $254$, $60$, $238$, $272$, $220$, $270$, $240$, $58$, $252$, $258$, $232$, $278$, $226$, $264$, $246$, $-198$, $-8$, $-182$, $-214$, $-100$, $234$, $276$, $224$, $266$, $244$, $54$, $248$, $262$, $228$, $280$, $230$, $260$, $250$, $56$, $242$, $268$, $222$, $274$, $236$, $154$, $-22$, $-40$, $-86$, $-76$, $-50$, $-70$, $-92$, $-34$, $-28$, $-98$, $-64$, $-44$, $-82$, $-80$, $-46$, $-66$, $-96$, $-30$, $-32$, $-94$, $-68$, $-48$, $-78$, $-84$, $-42$, $-20$, $-24$, $-38$, $-88$, $-74$, $-52$, $-72$, $-90$, $-36$, $-26$, $120$, $136$, $140$, $116$, $168$, $158$, $106$, $150$, $126$, $130$, $146$, $110$, $162$, $164$, $112$, $144$, $132$, $124$, $152$, $104$, $156$, $170$, $118$, $138$, $-176$, $-14$, $-204$, $-192$, $-2$, $-188$, $-208$, $-18$, $-172$, $-216$, $-180$, $-10$, $-200$, $-196$, $-6$, $-184$, $-212$, $-102$, $-210$, $-186$, $-4$, $-194$, $-202$, $-12$, $-178$, $-218$, $-174$, $-16$, $-206$, $-190)$\\

\noindent$11a31$ : $3$ | $(4$, $8$, $14$, $2$, $16$, $18$, $6$, $22$, $20$, $10$, $12)$  | $(26$, $14$, $64$, $16$, $24$, $28$, $-56$, $-62$, $-52$, $-40$, $-42$, $-50$, $-60$, $-58$, $-84$, $-88$, $-76$, $-86$, $30$, $106$, $92$, $108$, $32$, $34$, $110$, $94$, $104$, $100$, $98$, $112$, $96$, $102$, $-54$, $-38$, $-44$, $-48$, $-36$, $-46$, $70$, $22$, $18$, $66$, $72$, $74$, $68$, $20$, $-80$, $-8$, $-2$, $-12$, $-4$, $-6$, $-82$, $-90$, $-78$, $-10)$\\

\noindent$11a32$ : $3$ | $(4$, $8$, $14$, $2$, $16$, $18$, $6$, $22$, $20$, $12$, $10)$ | $(108$, $128$, $88$, $130$, $106$, $110$, $126$, $90$, $-48$, $-66$, $-68$, $-46$, $-76$, $-58$, $-56$, $-8$, $-2$, $-12$, $-52$, $-62$, $-72$, $118$, $98$, $174$, $92$, $124$, $112$, $104$, $132$, $102$, $114$, $122$, $94$, $172$, $96$, $120$, $116$, $100$, $176$, $180$, $34$, $198$, $32$, $178$, $-78$, $-154$, $-140$, $-170$, $-144$, $-150$, $-164$, $-134$, $-160$, $-84$, $-82$, $-158$, $-136$, $-166$, $-148$, $-146$, $-168$, $-138$, $-156$, $-80$, $-86$, $-162$, $28$, $194$, $38$, $184$, $18$, $22$, $188$, $42$, $190$, $24$, $16$, $182$, $36$, $196$, $30$, $26$, $192$, $40$, $186$, $20$, $-142$, $-152$, $-6$, $-4$, $-14$, $-50$, $-64$, $-70$, $-44$, $-74$, $-60$, $-54$, $-10)$\\

\noindent$11a33$ : $3$ | $(4$, $8$, $14$, $2$, $16$, $18$, $20$, $6$, $10$, $22$, $12)$ | $(122$, $160$, $112$, $132$, $150$, $102$, $142$, $140$, $104$, $152$, $130$, $114$, $56$, $116$, $128$, $154$, $106$, $138$, $144$, $100$, $148$, $134$, $110$, $158$, $124$, $120$, $-162$, $-190$, $226$, $260$, $212$, $236$, $250$, $202$, $246$, $240$, $208$, $256$, $230$, $218$, $266$, $220$, $54$, $118$, $126$, $156$, $108$, $136$, $146$, $-176$, $-70$, $-198$, $-64$, $-182$, $-170$, $-76$, $-192$, $-58$, $-188$, $-164$, $-82$, $-80$, $-166$, $-186$, $-60$, $-194$, $-74$, $-172$, $-180$, $-66$, $-200$, $-68$, $-178$, $-174$, $-72$, $-196$, $-62$, $-184$, $-168$, $-78$, $-84$, $222$, $264$, $216$, $232$, $254$, $206$, $242$, $244$, $204$, $252$, $234$, $214$, $262$, $224$, $228$, $258$, $210$, $238$, $248$, $-14$, $-36$, $-96$, $-42$, $-8$, $-20$, $-30$, $-90$, $-48$, $-2$, $-52$, $-86$, $-26$, $-24$, $-4$, $-46$, $-92$, $-32$, $-18$, $-10$, $-40$, $-98$, $-38$, $-12$, $-16$, $-34$, $-94$, $-44$, $-6$, $-22$, $-28$, $-88$, $-50)$\\

\noindent$11a34$ : $3$ | $(4$, $8$, $14$, $2$, $16$, $18$, $20$, $6$, $22$, $12$, $10)$ | $(28$, $16$, $32$, $24$, $64$, $26$, $30$, $-58$, $-94$, $-82$, $-78$, $-90$, $-62$, $-54$, $-56$, $-60$, $-92$, $-80$, $20$, $68$, $74$, $66$, $22$, $34$, $18$, $70$, $72$, $128$, $114$, $116$, $126$, $104$, $-88$, $-76$, $-84$, $-96$, $-86$, $100$, $108$, $122$, $120$, $110$, $98$, $102$, $106$, $124$, $118$, $112$, $-52$, $-40$, $-42$, $-10$, $-8$, $-44$, $-38$, $-50$, $-2$, $-14$, $-4$, $-48$, $-36$, $-46$, $-6$, $-12)$\\

\noindent$11a35$ : $3$ | $(4$, $8$, $14$, $2$, $16$, $18$, $22$, $6$, $20$, $12$, $10)$ | $(66$, $82$, $56$, $182$, $188$, $50$, $76$, $72$, $60$, $86$, $62$, $70$, $78$, $52$, $186$, $184$, $54$, $80$, $68$, $64$, $84$, $58$, $180$, $190$, $-116$, $-94$, $-132$, $-100$, $-110$, $-122$, $-88$, $-126$, $-106$, $-104$, $-128$, $-90$, $-120$, $-118$, $-92$, $-130$, $-102$, $-108$, $-124$, $164$, $194$, $174$, $154$, $148$, $136$, $144$, $158$, $170$, $198$, $168$, $160$, $142$, $138$, $150$, $152$, $176$, $192$, $162$, $166$, $196$, $172$, $156$, $146$, $-30$, $-8$, $-46$, $-44$, $-6$, $-32$, $-28$, $-10$, $-48$, $-12$, $-26$, $-34$, $-4$, $-42$, $-18$, $-20$, $-40$, $-2$, $-36$, $-24$, $-14$, $74$, $140$, $-112$, $-98$, $-134$, $-96$, $-114$, $-178$, $-16$, $-22$, $-38)$\\

\noindent$11a36$ : $3$ | $(4$, $8$, $14$, $2$, $16$, $20$, $6$, $12$, $22$, $10$, $18)$ | $(72$, $32$, $36$, $26$, $40$, $28$, $110$, $104$, $106$, $108$, $30$, $38$, $-56$, $-42$, $-60$, $-52$, $-46$, $-44$, $-54$, $-58$, $92$, $98$, $86$, $100$, $112$, $102$, $88$, $96$, $94$, $90$, $68$, $76$, $-16$, $-20$, $-12$, $-48$, $-50$, $-14$, $-18$, $66$, $74$, $70$, $34$, $-4$, $-22$, $-10$, $-8$, $-24$, $-6$, $-84$, $-2$, $-80$, $-64$, $-78$, $-62$, $-82)$\\

\noindent$11a37$ : $3$ | $(4$, $8$, $14$, $2$, $16$, $20$, $18$, $6$, $22$, $12$, $10)$ | $(58$, $98$, $88$, $96$, $106$, $-44$, $-36$, $-52$, $-2$, $-8$, $-46$, $-34$, $-50$, $-4$, $-6$, $-48$, $92$, $102$, $110$, $114$, $116$, $108$, $104$, $94$, $90$, $100$, $112$, $-38$, $-42$, $-40$, $56$, $12$, $22$, $32$, $24$, $14$, $54$, $16$, $26$, $30$, $20$, $10$, $18$, $28$, $-78$, $-68$, $-66$, $-80$, $-86$, $-76$, $-70$, $-64$, $-82$, $-84$, $-62$, $-72$, $-74$, $-60)$\\

\noindent$11a38$ : $3$ | $(4$, $8$, $14$, $2$, $18$, $16$, $20$, $6$, $10$, $22$, $12)$ | $(58$, $46$, $50$, $54$, $62$, $66$, $60$, $56$, $48$, $-4$, $-34$, $-42$, $-36$, $-2$, $-16$, $-6$, $52$, $44$, $94$, $90$, $96$, $-86$, $-80$, $-74$, $-72$, $-82$, $-84$, $-70$, $-76$, $-78$, $-68$, $-40$, $-38$, $92$, $98$, $30$, $18$, $28$, $100$, $26$, $20$, $32$, $22$, $24$, $64$, $-10$, $-12$, $-88$, $-8$, $-14)$\\

\noindent$11a39$ : $3$ | $(4$, $8$, $14$, $2$, $18$, $20$, $6$, $22$, $12$, $10$, $16)$ | $(50$, $72$, $84$, $58$, $64$, $90$, $66$, $56$, $44$, $42$, $48$, $52$, $70$, $86$, $60$, $62$, $88$, $68$, $54$, $46$, $-120$, $-76$, $-118$, $-122$, $-100$, $-102$, $-124$, $-116$, $-94$, $-108$, $-130$, $-110$, $-92$, $-114$, $-126$, $-104$, $-98$, $82$, $-18$, $-40$, $-20$, $-96$, $-106$, $-128$, $-112$, $152$, $134$, $160$, $144$, $142$, $166$, $140$, $146$, $158$, $132$, $154$, $150$, $136$, $162$, $78$, $80$, $164$, $138$, $148$, $156$, $-30$, $-8$, $-14$, $-36$, $-24$, $-2$, $-74$, $-4$, $-26$, $-34$, $-12$, $-10$, $-32$, $-28$, $-6$, $-16$, $-38$, $-22)$\\

\noindent$11a40$ : $3$ | $(4$, $8$, $14$, $2$, $18$, $20$, $16$, $6$, $22$, $10$, $12)$ | $(54$, $28$, $56$, $52$, $26$, $24$, $88$, $90$, $76$, $92$, $86$, $-70$, $-40$, $-34$, $82$, $96$, $80$, $84$, $94$, $78$, $-10$, $-22$, $-12$, $-16$, $-18$, $-38$, $-32$, $-30$, $-36$, $-20$, $-14$, $46$, $50$, $58$, $42$, $44$, $60$, $48$, $-64$, $-8$, $-2$, $-4$, $-68$, $-72$, $-62$, $-74$, $-66$, $-6)$\\

\noindent$11a41$ : $3$ | $(4$, $8$, $14$, $2$, $18$, $20$, $16$, $6$, $22$, $12$, $10)$ | $(62$, $30$, $28$, $60$, $64$, $32$, $26$, $96$, $98$, $84$, $100$, $94$, $-78$, $-38$, $-46$, $-40$, $92$, $102$, $86$, $48$, $90$, $104$, $88$, $-4$, $-12$, $-24$, $-14$, $-18$, $-20$, $-36$, $-44$, $-42$, $-34$, $-22$, $-16$, $54$, $58$, $66$, $50$, $52$, $68$, $56$, $-72$, $-6$, $-2$, $-10$, $-76$, $-80$, $-70$, $-82$, $-74$, $-8)$\\

\noindent$11a42$ : $3$ | $(4$, $8$, $14$, $2$, $20$, $16$, $6$, $12$, $22$, $10$, $18)$ | $(56$, $28$, $30$, $26$, $94$, $100$, $84$, $98$, $96$, $86$, $102$, $92$, $-72$, $-32$, $-38$, $40$, $88$, $104$, $90$, $-4$, $-6$, $-24$, $-8$, $-18$, $-14$, $-12$, $-20$, $-34$, $-36$, $-22$, $-10$, $-16$, $48$, $52$, $60$, $44$, $42$, $58$, $54$, $46$, $62$, $50$, $-66$, $-78$, $-2$, $-74$, $-70$, $-80$, $-64$, $-82$, $-68$, $-76)$\\

\noindent$11a43$ : $4$ | $(4$, $8$, $14$, $2$, $20$, $16$, $6$, $18$, $12$, $22$, $10)$ | $(92$, $88$, $18$, $28$, $-60$, $-46$, $-22$, $-50$, $-56$, $40$, $34$, $76$, $-48$, $-58$, $-62$, $-44$, $-66$, $-94$, $-112$, $-4$, $-64$, $84$, $68$, $74$, $78$, $36$, $114$, $38$, $80$, $72$, $70$, $82$, $86$, $-100$, $-106$, $-8$, $-108$, $-98$, $-96$, $-110$, $-6$, $-104$, $-102$, $-20$, $-54$, $-52$, $116$, $14$, $24$, $12$, $32$, $42$, $30$, $10$, $26$, $16$, $90$, $-2)$\\

\noindent$11a44$ : $4$ | $(4$, $8$, $14$, $2$, $20$, $16$, $18$, $6$, $12$, $22$, $10)$ | $(72$, $76$, $60$, $80$, $68$, $66$, $22$, $62$, $108$, $64$, $-104$, $-26$, $82$, $46$, $-10$, $-18$, $-16$, $-4$, $-58$, $-6$, $-14$, $-20$, $-12$, $-8$, $-56$, $-2$, $74$, $70$, $78$, $-96$, $-86$, $-106$, $-102$, $-28$, $-100$, $-92$, $-90$, $-88$, $-94$, $-98$, $-84$, $24$, $40$, $34$, $52$, $110$, $50$, $36$, $38$, $48$, $30$, $44$, $42$, $32$, $-54)$\\

\noindent$11a45$ : $3$ | $(4$, $8$, $14$, $2$, $20$, $18$, $6$, $22$, $12$, $10$, $16)$ | $(22$, $16$, $60$, $54$, $58$, $64$, $-84$, $-74$, $-68$, $-88$, $-70$, $-72$, $-86$, $-66$, $-10$, $-38$, $-6$, $-8$, $56$, $62$, $14$, $24$, $20$, $18$, $26$, $28$, $34$, $-78$, $-80$, $-52$, $-76$, $-82$, $-50$, $100$, $90$, $104$, $92$, $98$, $32$, $36$, $30$, $96$, $94$, $102$, $-46$, $-2$, $-42$, $-12$, $-40$, $-4$, $-48$, $-44)$\\

\noindent$11a46$ : $3$ | $(4$, $8$, $14$, $2$, $20$, $18$, $16$, $6$, $12$, $22$, $10)$ | $(98$, $82$, $78$, $102$, $76$, $84$, $96$, $100$, $80$, $-4$, $-10$, $-58$, $-64$, $-52$, $-68$, $-54$, $-62$, $-60$, $-8$, $-6$, $-16$, $-2$, $-12$, $-56$, $-66$, $90$, $70$, $94$, $86$, $74$, $140$, $72$, $88$, $92$, $-132$, $-126$, $-138$, $-122$, $-104$, $-114$, $-112$, $-106$, $-120$, $-136$, $-128$, $-130$, $-134$, $-118$, $-108$, $-110$, $-116$, $38$, $20$, $44$, $144$, $42$, $18$, $40$, $142$, $46$, $22$, $36$, $32$, $26$, $50$, $28$, $30$, $48$, $24$, $34$, $-124$, $-14)$\\

\noindent$11a47$ : $4$ | $(4$, $8$, $14$, $2$, $20$, $22$, $18$, $6$, $12$, $16$, $10)$ | $(16$, $36$, $28$, $44$, $46$, $30$, $34$, $-98$, $-22$, $-86$, $84$, $40$, $-58$, $-50$, $-102$, $-94$, $-96$, $-100$, $-52$, $-56$, $-60$, $-48$, $-104$, $106$, $74$, $68$, $20$, $18$, $70$, $72$, $108$, $32$, $-12$, $-4$, $-88$, $-24$, $-90$, $-6$, $-10$, $-62$, $-14$, $-2$, $-54$, $38$, $26$, $42$, $78$, $64$, $80$, $82$, $66$, $76$, $-8$, $-92)$\\

\noindent$11a48$ : $3$ | $(4$, $8$, $14$, $10$, $2$, $16$, $20$, $6$, $22$, $12$, $18)$ | $(38$, $122$, $128$, $134$, $116$, $42$, $118$, $132$, $130$, $120$, $40$, $114$, $136$, $126$, $124$, $138$, $112$, $140$, $-54$, $-52$, $-64$, $82$, $94$, $70$, $98$, $78$, $142$, $76$, $100$, $72$, $92$, $84$, $80$, $96$, $-16$, $-6$, $-26$, $-30$, $-2$, $-20$, $-12$, $-10$, $-22$, $-34$, $-32$, $-24$, $-8$, $-14$, $-18$, $-4$, $-28$, $74$, $90$, $86$, $36$, $88$, $-106$, $-62$, $-44$, $-66$, $-50$, $-56$, $-102$, $-58$, $-48$, $-68$, $-46$, $-60$, $-104$, $-110$, $-108)$\\

\noindent$11a49$ : $3$ | $(4$, $8$, $14$, $10$, $2$, $18$, $6$, $20$, $22$, $12$, $16)$ | $(68$, $248$, $254$, $274$, $228$, $60$, $240$, $262$, $266$, $236$, $56$, $232$, $270$, $258$, $244$, $64$, $224$, $278$, $250$, $252$, $276$, $226$, $62$, $242$, $260$, $268$, $234$, $-170$, $-212$, $-202$, $-180$, $-160$, $-222$, $-192$, $-190$, $66$, $246$, $256$, $272$, $230$, $58$, $238$, $264$, $-16$, $-50$, $-26$, $-76$, $-6$, $-40$, $-36$, $-2$, $-72$, $-30$, $-46$, $-12$, $-82$, $-20$, $-54$, $-22$, $-80$, $-10$, $-44$, $-32$, $-34$, $-42$, $-8$, $-78$, $-24$, $-52$, $-18$, $-84$, $-14$, $-48$, $-28$, $-74$, $-4$, $-38$, $98$, $128$, $124$, $106$, $146$, $90$, $136$, $116$, $114$, $138$, $88$, $144$, $108$, $122$, $130$, $96$, $152$, $100$, $280$, $102$, $150$, $94$, $132$, $120$, $110$, $142$, $86$, $140$, $112$, $118$, $134$, $92$, $148$, $104$, $126$, $-158$, $-182$, $-200$, $-214$, $-168$, $-172$, $-210$, $-204$, $-178$, $-162$, $-220$, $-194$, $-188$, $-154$, $-186$, $-196$, $-218$, $-164$, $-176$, $-206$, $-208$, $-174$, $-166$, $-216$, $-198$, $-184$, $-156$, $-70)$\\

\noindent$11a50$ : $3$ | $(4$, $8$, $14$, $10$, $2$, $18$, $6$, $22$, $20$, $12$, $16)$ | $(38$, $130$, $150$, $144$, $136$, $42$, $134$, $146$, $148$, $132$, $40$, $90$, $92$, $76$, $106$, $74$, $94$, $88$, $-54$, $-52$, $-64$, $80$, $102$, $70$, $98$, $84$, $158$, $86$, $96$, $72$, $104$, $78$, $82$, $100$, $-16$, $-6$, $-118$, $-108$, $-122$, $-10$, $-12$, $-20$, $-2$, $-114$, $-112$, $-124$, $-110$, $-116$, $-4$, $-18$, $-14$, $-8$, $-120$, $140$, $154$, $126$, $36$, $128$, $152$, $142$, $138$, $156$, $-24$, $-34$, $-56$, $-50$, $-66$, $-44$, $-62$, $-28$, $-30$, $-60$, $-46$, $-68$, $-48$, $-58$, $-32$, $-26$, $-22)$\\

\noindent$11a51$ : $3$ | $(4$, $8$, $14$, $10$, $2$, $20$, $16$, $6$, $22$, $12$, $18)$ | $(40$, $58$, $32$, $62$, $36$, $44$, $54$, $82$, $50$, $48$, $80$, $56$, $42$, $38$, $60$, $-104$, $-94$, $-90$, $-108$, $-100$, $-98$, $-110$, $-76$, $-68$, $-116$, $-70$, $-74$, $-112$, $-96$, $-102$, $-106$, $-92$, $34$, $148$, $118$, $144$, $142$, $120$, $150$, $-66$, $-114$, $-72$, $52$, $46$, $78$, $130$, $132$, $138$, $124$, $154$, $126$, $136$, $134$, $128$, $152$, $122$, $140$, $146$, $-18$, $-86$, $-12$, $-24$, $-2$, $-28$, $-8$, $-64$, $-6$, $-30$, $-4$, $-22$, $-14$, $-84$, $-16$, $-20$, $-88$, $-10$, $-26)$\\

\noindent$11a52$ : $3$ | $(4$, $8$, $14$, $12$, $2$, $18$, $22$, $20$, $10$, $16$, $6)$ | $(24$, $20$, $28$, $34$, $32$, $30$, $36$, $88$, $-40$, $-50$, $-44$, $-46$, $-48$, $-42$, $-52$, $-66$, $-54$, $-68$, $-74$, $90$, $86$, $94$, $82$, $96$, $84$, $92$, $64$, $-8$, $-12$, $-38$, $-14$, $-6$, $-10$, $56$, $62$, $16$, $60$, $58$, $18$, $26$, $22$, $-2$, $-78$, $-70$, $-72$, $-76$, $-4$, $-80)$\\

\noindent$11a53$ : $3$ | $(4$, $8$, $14$, $16$, $2$, $18$, $20$, $22$, $10$, $12$, $6)$ | $(84$, $96$, $72$, $74$, $94$, $86$, $82$, $98$, $70$, $76$, $92$, $88$, $80$, $100$, $-102$, $-104$, $-18$, $-6$, $124$, $112$, $134$, $114$, $122$, $34$, $126$, $110$, $132$, $116$, $120$, $32$, $128$, $108$, $130$, $30$, $-58$, $-46$, $-36$, $-48$, $-60$, $-66$, $-54$, $-42$, $-40$, $-52$, $-64$, $-62$, $-50$, $-38$, $-44$, $-56$, $118$, $68$, $78$, $90$, $-22$, $-10$, $-2$, $-14$, $-26$, $-28$, $-16$, $-4$, $-8$, $-20$, $-106$, $-24$, $-12)$\\

\noindent$11a54$ : $3$ | $(4$, $8$, $14$, $16$, $2$, $18$, $20$, $22$, $12$, $10$, $6)$ | $(138$, $222$, $122$, $154$, $100$, $70$, $258$, $84$, $86$, $256$, $68$, $102$, $62$, $106$, $64$, $252$, $90$, $80$, $262$, $74$, $96$, $156$, $98$, $72$, $260$, $82$, $88$, $254$, $66$, $104$, $-174$, $-108$, $-170$, $-178$, $-210$, $-188$, $-160$, $-118$, $-2$, $-114$, $-164$, $-184$, $-214$, $-182$, $-166$, $-112$, $-4$, $-120$, $-158$, $-190$, $-208$, $-176$, $-172$, $236$, $152$, $124$, $220$, $136$, $140$, $224$, $-194$, $-204$, $-32$, $-54$, $-18$, $-14$, $-50$, $-36$, $-200$, $-198$, $-38$, $-48$, $-12$, $-20$, $-56$, $-30$, $-206$, $-192$, $246$, $226$, $142$, $134$, $218$, $126$, $150$, $234$, $238$, $240$, $232$, $148$, $128$, $216$, $132$, $144$, $228$, $244$, $248$, $94$, $76$, $264$, $78$, $92$, $250$, $242$, $230$, $146$, $130$, $-16$, $-52$, $-34$, $-202$, $-196$, $-40$, $-46$, $-10$, $-22$, $-58$, $-28$, $-26$, $-60$, $-24$, $-8$, $-44$, $-42$, $-6$, $-110$, $-168$, $-180$, $-212$, $-186$, $-162$, $-116)$\\

\noindent$11a55$ : $3$ | $(4$, $8$, $16$, $2$, $18$, $20$, $22$, $6$, $14$, $10$, $12)$ | $(68$, $74$, $82$, $90$, $84$, $76$, $66$, $78$, $86$, $88$, $80$, $102$, $-54$, $-2$, $-12$, $-14$, $-22$, $-4$, $-10$, $-16$, $-20$, $-6$, $-8$, $-18$, $40$, $32$, $70$, $72$, $34$, $42$, $46$, $38$, $30$, $-96$, $-98$, $-24$, $-100$, $-94$, $-64$, $-50$, $-58$, $-92$, $-62$, $-48$, $-60$, $44$, $36$, $28$, $104$, $26$, $-56$, $-52)$\\

\noindent$11a56$ : $3$ | $(4$, $8$, $16$, $2$, $18$, $20$, $22$, $6$, $14$, $12$, $10)$ | $(60$, $80$, $68$, $52$, $50$, $-48$, $-8$, $136$, $116$, $126$, $144$, $124$, $118$, $138$, $132$, $112$, $130$, $140$, $120$, $122$, $142$, $128$, $114$, $134$, $-94$, $-14$, $-28$, $-34$, $-42$, $-20$, $-22$, $-40$, $-36$, $-26$, $-16$, $-46$, $-30$, $-32$, $-44$, $-18$, $-24$, $-38$, $64$, $56$, $76$, $72$, $12$, $10$, $70$, $78$, $58$, $62$, $82$, $66$, $54$, $74$, $-90$, $-98$, $-4$, $-106$, $-84$, $-104$, $-2$, $-100$, $-88$, $-110$, $-92$, $-96$, $-6$, $-108$, $-86$, $-102)$\\

\noindent$11a57$ : $4$ | $(4$, $8$, $16$, $2$, $20$, $22$, $18$, $6$, $12$, $14$, $10)$ | $(22$, $54$, $44$, $60$, $50$, $48$, $62$, $46$, $52$, $104$, $-74$, $-76$, $20$, $56$, $42$, $58$, $-26$, $-2$, $-18$, $32$, $92$, $-70$, $-80$, $-64$, $-84$, $-66$, $-78$, $-72$, $-40$, $-68$, $-82$, $98$, $86$, $102$, $94$, $90$, $34$, $106$, $36$, $88$, $96$, $100$, $-10$, $-16$, $-4$, $-28$, $-30$, $-6$, $-14$, $-12$, $-8$, $-38$, $-24)$\\

\noindent$11a58$ : $3$ | $(4$, $8$, $16$, $2$, $20$, $22$, $18$, $6$, $14$, $12$, $10)$ | $(132$, $162$, $144$, $124$, $154$, $152$, $122$, $146$, $160$, $130$, $50$, $80$, $82$, $52$, $128$, $158$, $148$, $120$, $150$, $156$, $126$, $142$, $164$, $134$, $-106$, $-170$, $-112$, $-42$, $-6$, $-14$, $-34$, $-36$, $-12$, $-8$, $-40$, $-30$, $-18$, $-2$, $-46$, $-108$, $-168$, $-110$, $-44$, $-4$, $-16$, $-32$, $-38$, $-10$, $66$, $58$, $88$, $74$, $202$, $136$, $166$, $140$, $198$, $70$, $92$, $62$, $-178$, $-98$, $-184$, $-114$, $-172$, $-104$, $-190$, $-192$, $-22$, $-26$, $-196$, $-186$, $-100$, $-176$, $-118$, $-180$, $-96$, $-182$, $-116$, $-174$, $-102$, $-188$, $-194$, $-24$, $138$, $200$, $72$, $90$, $60$, $64$, $94$, $68$, $56$, $86$, $76$, $48$, $78$, $84$, $54$, $-28$, $-20)$\\

\noindent$11a59$ : $2$ | $(4$, $8$, $16$, $2$, $22$, $20$, $18$, $6$, $14$, $12$, $10)$ | $(70$, $44$, $74$, $66$, $48$, $78$, $62$, $52$, $82$, $58$, $56$, $84$, $54$, $60$, $80$, $50$, $64$, $76$, $46$, $68$, $72$, $-104$, $-110$, $-98$, $-116$, $-92$, $-122$, $-86$, $-126$, $-88$, $-120$, $-94$, $-114$, $-100$, $-108$, $-106$, $-102$, $-112$, $-96$, $-118$, $-90$, $-124$, $154$, $128$, $158$, $150$, $132$, $162$, $146$, $136$, $166$, $142$, $140$, $168$, $138$, $144$, $164$, $134$, $148$, $160$, $130$, $152$, $156$, $-20$, $-26$, $-14$, $-32$, $-8$, $-38$, $-2$, $-42$, $-4$, $-36$, $-10$, $-30$, $-16$, $-24$, $-22$, $-18$, $-28$, $-12$, $-34$, $-6$, $-40)$\\

\noindent$11a60$ : $3$ | $(4$, $8$, $16$, $10$, $2$, $18$, $20$, $6$, $22$, $12$, $14)$ | $(30$, $22$, $44$, $46$, $20$, $28$, $32$, $24$, $-58$, $-50$, $-66$, $-60$, $-56$, $-52$, $-64$, $-62$, $-54$, $26$, $18$, $70$, $72$, $-38$, $-68$, $-40$, $42$, $80$, $88$, $86$, $78$, $76$, $84$, $90$, $82$, $74$, $-6$, $-48$, $-8$, $-4$, $-16$, $-36$, $-10$, $-2$, $-14$, $-34$, $-12)$\\

\noindent$11a61$ : $3$ | $(4$, $8$, $16$, $10$, $2$, $18$, $20$, $6$, $22$, $14$, $12)$ | $(104$, $134$, $128$, $98$, $110$, $140$, $112$, $96$, $126$, $52$, $118$, $90$, $120$, $54$, $124$, $94$, $114$, $138$, $108$, $100$, $130$, $132$, $102$, $106$, $136$, $-148$, $-86$, $-156$, $122$, $92$, $116$, $50$, $196$, $226$, $198$, $172$, $202$, $222$, $192$, $188$, $218$, $206$, $176$, $182$, $212$, $-152$, $-160$, $-82$, $-144$, $-74$, $-168$, $-68$, $-66$, $-166$, $-76$, $-142$, $-80$, $-162$, $-150$, $-88$, $-154$, $-158$, $-84$, $-146$, $-72$, $-170$, $-70$, $-64$, $-164$, $-78$, $190$, $220$, $204$, $174$, $184$, $214$, $210$, $180$, $178$, $208$, $216$, $186$, $194$, $224$, $200$, $-42$, $-4$, $-18$, $-28$, $-56$, $-26$, $-20$, $-62$, $-34$, $-12$, $-10$, $-36$, $-48$, $-2$, $-44$, $-40$, $-6$, $-16$, $-30$, $-58$, $-24$, $-22$, $-60$, $-32$, $-14$, $-8$, $-38$, $-46)$\\

\noindent$11a62$ : $3$ | $(4$, $8$, $16$, $10$, $2$, $18$, $20$, $22$, $6$, $12$, $14)$ | $(58$, $42$, $26$, $34$, $50$, $48$, $32$, $28$, $44$, $54$, $38$, $-66$, $-84$, $-102$, $-74$, $-76$, $-100$, $-82$, $-68$, $-64$, $-86$, $-104$, $-72$, $-78$, $-98$, $-80$, $-70$, $-62$, $-88$, $-90$, $152$, $146$, $130$, $128$, $144$, $94$, $136$, $122$, $138$, $96$, $142$, $126$, $132$, $148$, $154$, $150$, $-18$, $-118$, $-10$, $140$, $124$, $134$, $92$, $60$, $56$, $40$, $24$, $36$, $52$, $46$, $30$, $-14$, $-6$, $-114$, $-22$, $-112$, $-4$, $-16$, $-120$, $-12$, $-8$, $-116$, $-20$, $-110$, $-2$, $-106$, $-108)$\\

\noindent$11a63$ : $3$ | $(4$, $8$, $16$, $10$, $2$, $18$, $20$, $22$, $6$, $14$, $12)$ | $(182$, $210$, $154$, $160$, $204$, $188$, $176$, $216$, $148$, $166$, $198$, $194$, $170$, $174$, $190$, $202$, $162$, $152$, $212$, $180$, $184$, $208$, $156$, $158$, $206$, $186$, $178$, $214$, $150$, $164$, $200$, $192$, $172$, $-26$, $-218$, $-18$, $-58$, $-34$, $-6$, $-46$, $252$, $232$, $274$, $238$, $246$, $74$, $258$, $226$, $268$, $66$, $266$, $224$, $260$, $72$, $244$, $240$, $272$, $230$, $254$, $78$, $250$, $234$, $276$, $236$, $248$, $76$, $256$, $228$, $270$, $68$, $264$, $222$, $262$, $70$, $-108$, $-90$, $-130$, $-120$, $-80$, $-118$, $-132$, $-92$, $-106$, $-146$, $-140$, $-100$, $-98$, $-138$, $-112$, $-86$, $-126$, $-124$, $-84$, $-114$, $-136$, $-96$, $-102$, $-142$, $-144$, $-104$, $-94$, $-134$, $-116$, $-82$, $-122$, $-128$, $-88$, $-110$, $242$, $168$, $196$, $-22$, $-62$, $-30$, $-10$, $-50$, $-42$, $-2$, $-38$, $-54$, $-14$, $-16$, $-56$, $-36$, $-4$, $-44$, $-48$, $-8$, $-32$, $-60$, $-20$, $-220$, $-24$, $-64$, $-28$, $-12$, $-52$, $-40)$\\

\noindent$11a64$ : $3$ | $(4$, $8$, $16$, $10$, $2$, $20$, $22$, $18$, $6$, $14$, $12)$ | $(56$, $34$, $52$, $60$, $38$, $48$, $70$, $66$, $44$, $42$, $64$, $72$, $50$, $36$, $58$, $54$, $-112$, $-106$, $-118$, $-86$, $-76$, $-126$, $-78$, $-84$, $-120$, $-104$, $-114$, $-110$, $-108$, $-116$, $-88$, $-74$, $-124$, $-80$, $-82$, $-122$, $90$, $146$, $128$, $150$, $94$, $152$, $130$, $144$, $-102$, $-22$, $-14$, $68$, $46$, $40$, $62$, $156$, $134$, $140$, $160$, $138$, $136$, $158$, $142$, $132$, $154$, $92$, $148$, $-18$, $-98$, $-10$, $-26$, $-4$, $-32$, $-2$, $-28$, $-8$, $-100$, $-20$, $-16$, $-96$, $-12$, $-24$, $-6$, $-30)$\\

\noindent$11a65$ : $2$ | $(4$, $8$, $16$, $10$, $2$, $22$, $20$, $18$, $6$, $14$, $12)$ | $(90$, $80$, $100$, $70$, $110$, $60$, $114$, $66$, $104$, $76$, $94$, $86$, $84$, $96$, $74$, $106$, $64$, $116$, $62$, $108$, $72$, $98$, $82$, $88$, $92$, $78$, $102$, $68$, $112$, $-140$, $-162$, $-118$, $-166$, $-136$, $-144$, $-158$, $-122$, $-170$, $-132$, $-148$, $-154$, $-126$, $-174$, $-128$, $-152$, $-150$, $-130$, $-172$, $-124$, $-156$, $-146$, $-134$, $-168$, $-120$, $-160$, $-142$, $-138$, $-164$, $206$, $196$, $216$, $186$, $226$, $176$, $230$, $182$, $220$, $192$, $210$, $202$, $200$, $212$, $190$, $222$, $180$, $232$, $178$, $224$, $188$, $214$, $198$, $204$, $208$, $194$, $218$, $184$, $228$, $-24$, $-46$, $-2$, $-50$, $-20$, $-28$, $-42$, $-6$, $-54$, $-16$, $-32$, $-38$, $-10$, $-58$, $-12$, $-36$, $-34$, $-14$, $-56$, $-8$, $-40$, $-30$, $-18$, $-52$, $-4$, $-44$, $-26$, $-22$, $-48)$\\

\noindent$11a66$ : $3$ | $(4$, $8$, $16$, $12$, $2$, $18$, $20$, $22$, $10$, $6$, $14)$ | $(132$, $344$, $366$, $110$, $104$, $372$, $338$, $138$, $126$, $350$, $360$, $116$, $98$, $378$, $332$, $144$, $120$, $356$, $354$, $122$, $142$, $334$, $376$, $100$, $114$, $362$, $348$, $128$, $136$, $340$, $370$, $106$, $108$, $368$, $342$, $134$, $130$, $346$, $364$, $112$, $102$, $374$, $336$, $140$, $124$, $352$, $358$, $118$, $-150$, $-200$, $-312$, $-166$, $-184$, $-296$, $-182$, $-168$, $-310$, $-198$, $-152$, $-326$, $-208$, $-320$, $-158$, $-192$, $-304$, $-174$, $-176$, $-302$, $-190$, $-160$, $-318$, $-206$, $236$, $234$, $284$, $286$, $264$, $214$, $256$, $386$, $242$, $228$, $278$, $292$, $270$, $220$, $250$, $392$, $248$, $222$, $272$, $294$, $276$, $226$, $244$, $388$, $254$, $216$, $266$, $288$, $282$, $232$, $238$, $382$, $260$, $-36$, $-92$, $-20$, $-52$, $-76$, $-4$, $-68$, $-60$, $-12$, $-84$, $-44$, $-28$, $-328$, $-30$, $-42$, $-86$, $-14$, $-58$, $-70$, $-2$, $-74$, $-54$, $-18$, $-90$, $-38$, $-34$, $-94$, $-22$, $-50$, $-78$, $-6$, $-66$, $-62$, $-10$, $-82$, $-46$, $-26$, $-96$, $-24$, $-48$, $-80$, $-8$, $-64$, $274$, $224$, $246$, $390$, $252$, $218$, $268$, $290$, $280$, $230$, $240$, $384$, $258$, $212$, $262$, $380$, $330$, $-146$, $-204$, $-316$, $-162$, $-188$, $-300$, $-178$, $-172$, $-306$, $-194$, $-156$, $-322$, $-210$, $-324$, $-154$, $-196$, $-308$, $-170$, $-180$, $-298$, $-186$, $-164$, $-314$, $-202$, $-148$, $-32$, $-40$, $-88$, $-16$, $-56$, $-72)$\\

\noindent$11a67$ : $3$ | $(4$, $8$, $16$, $12$, $2$, $18$, $22$, $20$, $10$, $6$, $14)$ | $(18$, $92$, $94$, $96$, $104$, $-44$, $-32$, $-40$, $-52$, $-50$, $-38$, $-34$, $-46$, $-8$, $-42$, $100$, $108$, $58$, $60$, $110$, $98$, $102$, $106$, $56$, $62$, $112$, $64$, $-4$, $-78$, $-88$, $-76$, $-2$, $-6$, $26$, $10$, $28$, $16$, $20$, $90$, $24$, $12$, $30$, $14$, $22$, $-36$, $-48$, $-54$, $-66$, $-72$, $-84$, $-82$, $-70$, $-68$, $-80$, $-86$, $-74)$\\

\noindent$11a68$ : $3$ | $(4$, $8$, $16$, $14$, $2$, $18$, $20$, $22$, $10$, $12$, $6)$ | $(32$, $42$, $44$, $30$, $48$, $38$, $36$, $50$, $56$, $58$, $52$, $34$, $40$, $46$, $-96$, $-64$, $-102$, $-90$, $-112$, $-92$, $-100$, $-66$, $-98$, $-94$, $-110$, $-88$, $-104$, $-108$, $-86$, $-106$, $54$, $124$, $118$, $-4$, $-26$, $-10$, $-14$, $-22$, $-62$, $-16$, $-20$, $-60$, $-18$, $80$, $74$, $114$, $70$, $122$, $120$, $68$, $116$, $76$, $82$, $84$, $78$, $72$, $-12$, $-24$, $-2$, $-6$, $-28$, $-8)$\\

\noindent$11a69$ : $3$ | $(4$, $8$, $16$, $14$, $2$, $18$, $20$, $22$, $12$, $10$, $6)$ | $(136$, $160$, $112$, $158$, $138$, $54$, $234$, $248$, $68$, $152$, $144$, $60$, $240$, $242$, $62$, $146$, $150$, $66$, $246$, $236$, $56$, $140$, $156$, $114$, $162$, $134$, $-106$, $-186$, $-76$, $-94$, $-174$, $-88$, $-82$, $-180$, $-100$, $258$, $224$, $272$, $128$, $168$, $120$, $264$, $218$, $266$, $122$, $170$, $126$, $270$, $222$, $260$, $116$, $164$, $132$, $276$, $228$, $254$, $-192$, $-202$, $-16$, $-36$, $-214$, $-28$, $-24$, $-210$, $-40$, $-12$, $-198$, $-196$, $-108$, $-188$, $-74$, $-96$, $-176$, $-86$, $-84$, $-178$, $-98$, $-72$, $-190$, $-110$, $-194$, $-200$, $-14$, $-38$, $-212$, $-26$, $124$, $268$, $220$, $262$, $118$, $166$, $130$, $274$, $226$, $256$, $252$, $230$, $52$, $232$, $250$, $70$, $154$, $142$, $58$, $238$, $244$, $64$, $148$, $-32$, $-20$, $-206$, $-44$, $-8$, $-2$, $-50$, $-104$, $-184$, $-78$, $-92$, $-172$, $-90$, $-80$, $-182$, $-102$, $-48$, $-4$, $-6$, $-46$, $-204$, $-18$, $-34$, $-216$, $-30$, $-22$, $-208$, $-42$, $-10)$\\

\noindent$11a70$ : $3$ | $(4$, $8$, $18$, $12$, $2$, $16$, $20$, $6$, $10$, $22$, $14)$ | $(172$, $236$, $186$, $198$, $224$, $326$, $218$, $204$, $180$, $242$, $178$, $206$, $216$, $328$, $226$, $196$, $188$, $234$, $170$, $104$, $-130$, $-156$, $-24$, $-14$, $-146$, $-140$, $-8$, $-30$, $-162$, $-124$, $-166$, $-34$, $-4$, $-136$, $-150$, $-18$, $-20$, $-152$, $-134$, $-2$, $-36$, $-168$, $-126$, $-160$, $-28$, $-10$, $-142$, $-144$, $-12$, $-26$, $-158$, $-128$, $-132$, $-154$, $-22$, $-16$, $-148$, $-138$, $-6$, $-32$, $-164$, $192$, $230$, $332$, $212$, $210$, $174$, $238$, $184$, $200$, $222$, $324$, $220$, $202$, $182$, $240$, $176$, $208$, $214$, $330$, $228$, $194$, $190$, $232$, $-298$, $-260$, $-262$, $-314$, $-282$, $-244$, $-278$, $-318$, $-266$, $-256$, $-294$, $-302$, $-304$, $-292$, $-254$, $-268$, $-320$, $-276$, $-246$, $-284$, $-312$, $-40$, $-310$, $-286$, $-248$, $-274$, $-322$, $-270$, $-252$, $-290$, $-306$, $-300$, $-296$, $-258$, $-264$, $-316$, $-280$, $72$, $110$, $46$, $98$, $90$, $54$, $118$, $64$, $80$, $334$, $78$, $66$, $116$, $52$, $92$, $96$, $48$, $112$, $70$, $74$, $108$, $44$, $100$, $88$, $56$, $120$, $62$, $82$, $84$, $60$, $122$, $58$, $86$, $102$, $42$, $106$, $76$, $68$, $114$, $50$, $94$, $-272$, $-250$, $-288$, $-308$, $-38)$\\

\noindent$11a71$ : $3$ | $(4$, $10$, $12$, $14$, $2$, $22$, $18$, $20$, $6$, $8$, $16)$ | $(74$, $88$, $128$, $132$, $92$, $70$, $78$, $84$, $126$, $86$, $76$, $72$, $90$, $130$, $-6$, $-16$, $-18$, $-4$, $-26$, $-8$, $-94$, $-10$, $-24$, $-2$, $-20$, $-14$, $-98$, $136$, $40$, $42$, $138$, $52$, $134$, $-56$, $-60$, $-120$, $-110$, $-108$, $-122$, $-62$, $-54$, $-64$, $-124$, $-106$, $-112$, $-118$, $-58$, $68$, $80$, $82$, $66$, $30$, $34$, $48$, $142$, $46$, $36$, $28$, $38$, $44$, $140$, $50$, $32$, $-104$, $-114$, $-116$, $-102$, $-100$, $-96$, $-12$, $-22)$\\

\noindent$11a72$ : $3$ |  $(4$, $10$, $12$, $14$, $2$, $22$, $18$, $20$, $8$, $6$, $16)$ | $(154$, $112$, $186$, $122$, $140$, $168$, $98$, $172$, $136$, $126$, $182$, $108$, $158$, $150$, $148$, $160$, $106$, $180$, $128$, $134$, $174$, $100$, $166$, $142$, $120$, $188$, $114$, $228$, $116$, $190$, $118$, $144$, $164$, $102$, $176$, $132$, $130$, $178$, $104$, $162$, $146$, $152$, $156$, $110$, $184$, $124$, $138$, $170$, $-270$, $-208$, $-240$, $-290$, $-250$, $-218$, $-260$, $-280$, $-198$, $-230$, $-194$, $-234$, $-202$, $-276$, $-264$, $-214$, $-246$, $-294$, $-244$, $-212$, $-266$, $-274$, $-204$, $-236$, $-286$, $-254$, $-222$, $-256$, $-284$, $-282$, $-258$, $-220$, $-252$, $-288$, $-238$, $-206$, $-272$, $-268$, $-210$, $-242$, $-292$, $-248$, $-216$, $-262$, $-278$, $-200$, $-232$, $-192$, $350$, $308$, $312$, $346$, $354$, $304$, $374$, $326$, $332$, $368$, $298$, $360$, $340$, $318$, $382$, $-28$, $-82$, $-196$, $310$, $348$, $352$, $306$, $376$, $324$, $334$, $366$, $296$, $362$, $338$, $320$, $380$, $224$, $384$, $316$, $342$, $358$, $300$, $370$, $330$, $328$, $372$, $302$, $356$, $344$, $314$, $386$, $226$, $378$, $322$, $336$, $364$, $-72$, $-38$, $-18$, $-92$, $-8$, $-48$, $-62$, $-56$, $-54$, $-2$, $-86$, $-24$, $-32$, $-78$, $-66$, $-44$, $-12$, $-96$, $-14$, $-42$, $-68$, $-76$, $-34$, $-22$, $-88$, $-4$, $-52$, $-58$, $-60$, $-50$, $-6$, $-90$, $-20$, $-36$, $-74$, $-70$, $-40$, $-16$, $-94$, $-10$, $-46$, $-64$, $-80$, $-30$, $-26$, $-84)$\\

\noindent$11a73$ : $3$ | $(4$, $10$, $12$, $14$, $16$, $2$, $18$, $20$, $22$, $6$, $8)$ | $(122$, $154$, $106$, $136$, $140$, $102$, $150$, $126$, $44$, $236$, $238$, $46$, $128$, $148$, $100$, $142$, $134$, $108$, $156$, $120$, $-86$, $-68$, $-164$, $-62$, $-92$, $266$, $218$, $256$, $252$, $222$, $270$, $112$, $160$, $116$, $274$, $226$, $248$, $260$, $214$, $262$, $246$, $228$, $276$, $118$, $158$, $110$, $268$, $220$, $254$, $-56$, $-170$, $-74$, $-80$, $-176$, $-186$, $-200$, $-22$, $-18$, $-204$, $-182$, $-180$, $-84$, $-70$, $-166$, $-60$, $-94$, $-52$, $-174$, $-78$, $-76$, $-172$, $-54$, $-96$, $-58$, $-168$, $-72$, $-82$, $-178$, $-184$, $-202$, $-20$, $114$, $272$, $224$, $250$, $258$, $216$, $264$, $244$, $230$, $40$, $232$, $242$, $50$, $132$, $144$, $98$, $146$, $130$, $48$, $240$, $234$, $42$, $124$, $152$, $104$, $138$, $-8$, $-32$, $-190$, $-196$, $-26$, $-14$, $-208$, $-2$, $-38$, $-88$, $-66$, $-162$, $-64$, $-90$, $-36$, $-4$, $-210$, $-12$, $-28$, $-194$, $-192$, $-30$, $-10$, $-212$, $-6$, $-34$, $-188$, $-198$, $-24$, $-16$, $-206)$\\

\noindent$11a74$ : $3$ | $(4$, $10$, $12$, $14$, $18$, $2$, $6$, $20$, $22$, $8$, $16)$ | $(60$, $158$, $164$, $54$, $66$, $152$, $170$, $174$, $130$, $176$, $168$, $154$, $64$, $56$, $162$, $160$, $58$, $62$, $156$, $166$, $178$, $132$, $172$, $-14$, $-20$, $-40$, $-148$, $-120$, $-88$, $-114$, $-112$, $-90$, $-122$, $-150$, $-42$, $-18$, $-16$, $-44$, $-12$, $-22$, $-38$, $-4$, $-30$, $144$, $184$, $138$, $80$, $50$, $70$, $74$, $46$, $76$, $68$, $52$, $82$, $140$, $186$, $142$, $84$, $146$, $182$, $136$, $78$, $48$, $72$, $-98$, $-104$, $-10$, $-24$, $-36$, $-2$, $-32$, $-28$, $-6$, $180$, $134$, $-124$, $-92$, $-110$, $-116$, $-86$, $-118$, $-108$, $-94$, $-126$, $-102$, $-100$, $-128$, $-96$, $-106$, $-8$, $-26$, $-34)$\\

\noindent$11a75$ : $2$ | $(4$, $10$, $12$, $14$, $18$, $2$, $6$, $22$, $20$, $8$, $16)$ | $(112$, $158$, $98$, $126$, $144$, $84$, $140$, $130$, $94$, $154$, $116$, $108$, $162$, $102$, $122$, $148$, $88$, $136$, $134$, $90$, $150$, $120$, $104$, $164$, $106$, $118$, $152$, $92$, $132$, $138$, $86$, $146$, $124$, $100$, $160$, $110$, $114$, $156$, $96$, $128$, $142$, $-200$, $-222$, $-178$, $-244$, $-172$, $-228$, $-194$, $-206$, $-216$, $-184$, $-238$, $-166$, $-234$, $-188$, $-212$, $-210$, $-190$, $-232$, $-168$, $-240$, $-182$, $-218$, $-204$, $-196$, $-226$, $-174$, $-246$, $-176$, $-224$, $-198$, $-202$, $-220$, $-180$, $-242$, $-170$, $-230$, $-192$, $-208$, $-214$, $-186$, $-236$, $276$, $322$, $262$, $290$, $308$, $248$, $304$, $294$, $258$, $318$, $280$, $272$, $326$, $266$, $286$, $312$, $252$, $300$, $298$, $254$, $314$, $284$, $268$, $328$, $270$, $282$, $316$, $256$, $296$, $302$, $250$, $310$, $288$, $264$, $324$, $274$, $278$, $320$, $260$, $292$, $306$, $-36$, $-58$, $-14$, $-80$, $-8$, $-64$, $-30$, $-42$, $-52$, $-20$, $-74$, $-2$, $-70$, $-24$, $-48$, $-46$, $-26$, $-68$, $-4$, $-76$, $-18$, $-54$, $-40$, $-32$, $-62$, $-10$, $-82$, $-12$, $-60$, $-34$, $-38$, $-56$, $-16$, $-78$, $-6$, $-66$, $-28$, $-44$, $-50$, $-22$, $-72)$\\

\noindent$11a76$ : $3$ | $(4$, $10$, $12$, $14$, $18$, $2$, $20$, $8$, $22$, $6$, $16)$ | $(120$, $162$, $176$, $184$, $154$, $186$, $174$, $164$, $118$, $122$, $114$, $168$, $170$, $190$, $158$, $180$, $-136$, $-6$, $-12$, $-130$, $-26$, $-22$, $-126$, $-16$, $-2$, $-140$, $-28$, $-132$, $-10$, $-8$, $-134$, $-30$, $-138$, $-4$, $-14$, $-128$, $-24$, $42$, $52$, $64$, $32$, $62$, $54$, $40$, $72$, $44$, $196$, $198$, $46$, $70$, $38$, $56$, $60$, $34$, $66$, $50$, $-142$, $-92$, $-96$, $-146$, $-144$, $-94$, $116$, $166$, $172$, $188$, $156$, $182$, $178$, $160$, $192$, $194$, $200$, $48$, $68$, $36$, $58$, $-104$, $-84$, $-78$, $-110$, $-148$, $-98$, $-90$, $-74$, $-88$, $-100$, $-150$, $-108$, $-80$, $-82$, $-106$, $-152$, $-102$, $-86$, $-76$, $-112$, $-20$, $-124$, $-18)$\\

\noindent$11a77$ : $2$ | $(4$, $10$, $12$, $14$, $18$, $2$, $22$, $20$, $8$, $16$, $6)$ | $(206$, $164$, $248$, $138$, $232$, $180$, $190$, $222$, $148$, $258$, $154$, $216$, $196$, $174$, $238$, $132$, $242$, $170$, $200$, $212$, $158$, $254$, $144$, $226$, $186$, $184$, $228$, $142$, $252$, $160$, $210$, $202$, $168$, $244$, $134$, $236$, $176$, $194$, $218$, $152$, $260$, $150$, $220$, $192$, $178$, $234$, $136$, $246$, $166$, $204$, $208$, $162$, $250$, $140$, $230$, $182$, $188$, $224$, $146$, $256$, $156$, $214$, $198$, $172$, $240$, $-310$, $-372$, $-272$, $-348$, $-334$, $-286$, $-386$, $-296$, $-324$, $-358$, $-262$, $-362$, $-320$, $-300$, $-382$, $-282$, $-338$, $-344$, $-276$, $-376$, $-306$, $-314$, $-368$, $-268$, $-352$, $-330$, $-290$, $-390$, $-292$, $-328$, $-354$, $-266$, $-366$, $-316$, $-304$, $-378$, $-278$, $-342$, $-340$, $-280$, $-380$, $-302$, $-318$, $-364$, $-264$, $-356$, $-326$, $-294$, $-388$, $-288$, $-332$, $-350$, $-270$, $-370$, $-312$, $-308$, $-374$, $-274$, $-346$, $-336$, $-284$, $-384$, $-298$, $-322$, $-360$, $466$, $424$, $508$, $398$, $492$, $440$, $450$, $482$, $408$, $518$, $414$, $476$, $456$, $434$, $498$, $392$, $502$, $430$, $460$, $472$, $418$, $514$, $404$, $486$, $446$, $444$, $488$, $402$, $512$, $420$, $470$, $462$, $428$, $504$, $394$, $496$, $436$, $454$, $478$, $412$, $520$, $410$, $480$, $452$, $438$, $494$, $396$, $506$, $426$, $464$, $468$, $422$, $510$, $400$, $490$, $442$, $448$, $484$, $406$, $516$, $416$, $474$, $458$, $432$, $500$, $-50$, $-112$, $-12$, $-88$, $-74$, $-26$, $-126$, $-36$, $-64$, $-98$, $-2$, $-102$, $-60$, $-40$, $-122$, $-22$, $-78$, $-84$, $-16$, $-116$, $-46$, $-54$, $-108$, $-8$, $-92$, $-70$, $-30$, $-130$, $-32$, $-68$, $-94$, $-6$, $-106$, $-56$, $-44$, $-118$, $-18$, $-82$, $-80$, $-20$, $-120$, $-42$, $-58$, $-104$, $-4$, $-96$, $-66$, $-34$, $-128$, $-28$, $-72$, $-90$, $-10$, $-110$, $-52$, $-48$, $-114$, $-14$, $-86$, $-76$, $-24$, $-124$, $-38$, $-62$, $-100)$\\

\noindent$11a78$ : $3$ | $(4$, $10$, $12$, $14$, $20$, $2$, $16$, $6$, $22$, $8$, $18)$ | $(10$, $96$, $98$, $12$, $-66$, $-64$, $-68$, $-80$, $48$, $58$, $40$, $54$, $52$, $42$, $60$, $46$, $14$, $50$, $56$, $-74$, $-86$, $-6$, $-2$, $-82$, $-78$, $-70$, $-88$, $-72$, $-76$, $-84$, $-4$, $44$, $116$, $114$, $108$, $90$, $104$, $102$, $92$, $110$, $112$, $94$, $100$, $106$, $-24$, $-36$, $-18$, $-30$, $-62$, $-28$, $-20$, $-38$, $-22$, $-26$, $-34$, $-16$, $-32$, $-8)$\\

\noindent$11a79$ : $3$ | $(4$, $10$, $12$, $14$, $20$, $2$, $18$, $8$, $22$, $6$, $16)$ | $(78$, $66$, $68$, $80$, $76$, $-94$, $-88$, $-104$, $-90$, $-92$, $-102$, $-86$, $-96$, $-4$, $112$, $122$, $110$, $26$, $114$, $120$, $20$, $16$, $14$, $22$, $118$, $116$, $24$, $12$, $18$, $-42$, $-46$, $-54$, $-34$, $-28$, $-36$, $-52$, $-48$, $-40$, $-32$, $-30$, $-38$, $-50$, $72$, $62$, $106$, $58$, $108$, $64$, $70$, $82$, $74$, $60$, $-44$, $-56$, $-10$, $-2$, $-6$, $-98$, $-84$, $-100$, $-8)$\\

\noindent$11a80$ : $3$ | $(4$, $10$, $12$, $14$, $22$, $2$, $18$, $20$, $6$, $8$, $16)$ | $(116$, $162$, $102$, $130$, $144$, $88$, $148$, $198$, $154$, $94$, $138$, $136$, $96$, $156$, $122$, $110$, $168$, $108$, $124$, $126$, $106$, $166$, $112$, $120$, $158$, $98$, $134$, $140$, $92$, $152$, $196$, $150$, $90$, $142$, $132$, $100$, $160$, $118$, $114$, $164$, $104$, $128$, $146$, $-250$, $-228$, $-192$, $-220$, $-258$, $-214$, $-186$, $-234$, $-244$, $-176$, $-204$, $-170$, $-208$, $-180$, $-240$, $-238$, $-182$, $-210$, $-172$, $-202$, $-174$, $-246$, $-232$, $-188$, $-216$, $-260$, $-218$, $-190$, $-230$, $-248$, $-252$, $-226$, $-194$, $-222$, $-256$, $-212$, $-184$, $-236$, $-242$, $-178$, $-206$, $292$, $338$, $200$, $334$, $296$, $288$, $342$, $282$, $302$, $328$, $264$, $320$, $310$, $274$, $-224$, $-254$, $-38$, $336$, $294$, $290$, $340$, $280$, $304$, $326$, $262$, $322$, $308$, $276$, $272$, $312$, $318$, $266$, $330$, $300$, $284$, $344$, $286$, $298$, $332$, $268$, $316$, $314$, $270$, $278$, $306$, $324$, $-22$, $-54$, $-68$, $-8$, $-84$, $-14$, $-62$, $-60$, $-16$, $-28$, $-48$, $-74$, $-2$, $-78$, $-44$, $-32$, $-34$, $-42$, $-80$, $-4$, $-72$, $-50$, $-26$, $-18$, $-58$, $-64$, $-12$, $-86$, $-10$, $-66$, $-56$, $-20$, $-24$, $-52$, $-70$, $-6$, $-82$, $-40$, $-36$, $-30$, $-46$, $-76)$\\

\noindent$11a81$ : $3$ | $(4$, $10$, $12$, $14$, $22$, $2$, $18$, $20$, $8$, $6$, $16)$ | $(144$, $220$, $120$, $228$, $244$, $78$, $88$, $234$, $238$, $84$, $82$, $240$, $232$, $90$, $76$, $246$, $226$, $118$, $222$, $146$, $142$, $218$, $122$, $230$, $242$, $80$, $86$, $236$, $-16$, $-24$, $-48$, $-8$, $-32$, $-114$, $-38$, $-2$, $-42$, $-110$, $-200$, $-178$, $-168$, $-208$, $-170$, $-176$, $-202$, $-108$, $-44$, $-4$, $-36$, $-116$, $-34$, $-6$, $-46$, $124$, $216$, $140$, $148$, $224$, $-152$, $-194$, $-184$, $-162$, $-26$, $-14$, $-54$, $-18$, $-22$, $-50$, $-10$, $-30$, $-158$, $-188$, $-190$, $-156$, $66$, $100$, $98$, $68$, $248$, $74$, $92$, $106$, $60$, $134$, $210$, $130$, $56$, $128$, $212$, $136$, $62$, $104$, $94$, $72$, $250$, $70$, $96$, $102$, $64$, $138$, $214$, $126$, $58$, $132$, $-20$, $-52$, $-12$, $-28$, $-160$, $-186$, $-192$, $-154$, $-150$, $-196$, $-182$, $-164$, $-204$, $-174$, $-172$, $-206$, $-166$, $-180$, $-198$, $-112$, $-40)$\\

\noindent$11a82$ : $3$ | $(4$, $10$, $12$, $16$, $2$, $22$, $18$, $20$, $6$, $14$, $8)$ | $(88$, $60$, $108$, $68$, $80$, $96$, $52$, $100$, $76$, $72$, $104$, $56$, $92$, $84$, $64$, $112$, $114$, $62$, $86$, $90$, $58$, $106$, $70$, $78$, $98$, $-164$, $-134$, $-124$, $-154$, $-174$, $-180$, $-148$, $-118$, $-140$, $-170$, $-158$, $-128$, $-130$, $-160$, $-168$, $-138$, $-120$, $-150$, $-178$, $-176$, $-152$, $-122$, $-136$, $-166$, $-162$, $-132$, $-126$, $-156$, $-172$, $-142$, $-116$, $-146$, $144$, $234$, $204$, $252$, $214$, $224$, $242$, $196$, $244$, $222$, $216$, $250$, $202$, $236$, $-182$, $-32$, $232$, $206$, $254$, $212$, $226$, $240$, $198$, $246$, $220$, $218$, $248$, $200$, $238$, $228$, $210$, $256$, $208$, $230$, $110$, $66$, $82$, $94$, $54$, $102$, $74$, $-18$, $-46$, $-8$, $-186$, $-28$, $-36$, $-2$, $-40$, $-24$, $-190$, $-12$, $-50$, $-14$, $-192$, $-22$, $-42$, $-4$, $-34$, $-30$, $-184$, $-6$, $-44$, $-20$, $-194$, $-16$, $-48$, $-10$, $-188$, $-26$, $-38)$\\

\noindent$11a83$ : $3$ | $(4$, $10$, $12$, $16$, $2$, $22$, $18$, $20$, $8$, $14$, $6)$ | $(74$, $62$, $82$, $88$, $86$, $122$, $104$, $124$, $114$, $112$, $126$, $106$, $120$, $-10$, $-6$, $-90$, $-8$, $28$, $116$, $110$, $128$, $108$, $118$, $30$, $32$, $26$, $66$, $70$, $78$, $-18$, $-38$, $-44$, $-24$, $-12$, $-14$, $-22$, $-42$, $-40$, $-20$, $-16$, $-36$, $-46$, $-50$, $-34$, $-48$, $84$, $64$, $72$, $76$, $60$, $80$, $68$, $-54$, $-94$, $-2$, $-98$, $-58$, $-100$, $-4$, $-92$, $-52$, $-102$, $-56$, $-96)$\\

\noindent$11a84$ : $2$ | $(4$, $10$, $12$, $16$, $2$, $22$, $20$, $18$, $6$, $14$, $8)$ | $(162$, $116$, $194$, $130$, $148$, $176$, $102$, $180$, $144$, $134$, $190$, $112$, $166$, $158$, $120$, $198$, $126$, $152$, $172$, $106$, $184$, $140$, $138$, $186$, $108$, $170$, $154$, $124$, $200$, $122$, $156$, $168$, $110$, $188$, $136$, $142$, $182$, $104$, $174$, $150$, $128$, $196$, $118$, $160$, $164$, $114$, $192$, $132$, $146$, $178$, $-244$, $-270$, $-218$, $-296$, $-208$, $-280$, $-234$, $-254$, $-260$, $-228$, $-286$, $-202$, $-290$, $-224$, $-264$, $-250$, $-238$, $-276$, $-212$, $-300$, $-214$, $-274$, $-240$, $-248$, $-266$, $-222$, $-292$, $-204$, $-284$, $-230$, $-258$, $-256$, $-232$, $-282$, $-206$, $-294$, $-220$, $-268$, $-246$, $-242$, $-272$, $-216$, $-298$, $-210$, $-278$, $-236$, $-252$, $-262$, $-226$, $-288$, $362$, $316$, $394$, $330$, $348$, $376$, $302$, $380$, $344$, $334$, $390$, $312$, $366$, $358$, $320$, $398$, $326$, $352$, $372$, $306$, $384$, $340$, $338$, $386$, $308$, $370$, $354$, $324$, $400$, $322$, $356$, $368$, $310$, $388$, $336$, $342$, $382$, $304$, $374$, $350$, $328$, $396$, $318$, $360$, $364$, $314$, $392$, $332$, $346$, $378$, $-44$, $-70$, $-18$, $-96$, $-8$, $-80$, $-34$, $-54$, $-60$, $-28$, $-86$, $-2$, $-90$, $-24$, $-64$, $-50$, $-38$, $-76$, $-12$, $-100$, $-14$, $-74$, $-40$, $-48$, $-66$, $-22$, $-92$, $-4$, $-84$, $-30$, $-58$, $-56$, $-32$, $-82$, $-6$, $-94$, $-20$, $-68$, $-46$, $-42$, $-72$, $-16$, $-98$, $-10$, $-78$, $-36$, $-52$, $-62$, $-26$, $-88)$\\

\noindent$11a85$ : $2$ | $(4$, $10$, $12$, $16$, $2$, $22$, $20$, $18$, $8$, $14$, $6)$ | $(172$, $122$, $204$, $140$, $154$, $190$, $108$, $186$, $158$, $136$, $208$, $126$, $168$, $176$, $118$, $200$, $144$, $150$, $194$, $112$, $182$, $162$, $132$, $212$, $130$, $164$, $180$, $114$, $196$, $148$, $146$, $198$, $116$, $178$, $166$, $128$, $210$, $134$, $160$, $184$, $110$, $192$, $152$, $142$, $202$, $120$, $174$, $170$, $124$, $206$, $138$, $156$, $188$, $-272$, $-246$, $-298$, $-220$, $-314$, $-230$, $-288$, $-256$, $-262$, $-282$, $-236$, $-308$, $-214$, $-304$, $-240$, $-278$, $-266$, $-252$, $-292$, $-226$, $-318$, $-224$, $-294$, $-250$, $-268$, $-276$, $-242$, $-302$, $-216$, $-310$, $-234$, $-284$, $-260$, $-258$, $-286$, $-232$, $-312$, $-218$, $-300$, $-244$, $-274$, $-270$, $-248$, $-296$, $-222$, $-316$, $-228$, $-290$, $-254$, $-264$, $-280$, $-238$, $-306$, $384$, $334$, $416$, $352$, $366$, $402$, $320$, $398$, $370$, $348$, $420$, $338$, $380$, $388$, $330$, $412$, $356$, $362$, $406$, $324$, $394$, $374$, $344$, $424$, $342$, $376$, $392$, $326$, $408$, $360$, $358$, $410$, $328$, $390$, $378$, $340$, $422$, $346$, $372$, $396$, $322$, $404$, $364$, $354$, $414$, $332$, $386$, $382$, $336$, $418$, $350$, $368$, $400$, $-60$, $-34$, $-86$, $-8$, $-102$, $-18$, $-76$, $-44$, $-50$, $-70$, $-24$, $-96$, $-2$, $-92$, $-28$, $-66$, $-54$, $-40$, $-80$, $-14$, $-106$, $-12$, $-82$, $-38$, $-56$, $-64$, $-30$, $-90$, $-4$, $-98$, $-22$, $-72$, $-48$, $-46$, $-74$, $-20$, $-100$, $-6$, $-88$, $-32$, $-62$, $-58$, $-36$, $-84$, $-10$, $-104$, $-16$, $-78$, $-42$, $-52$, $-68$, $-26$, $-94)$\\

\noindent$11a86$ : $3$ | $(4$, $10$, $12$, $16$, $18$, $2$, $20$, $6$, $8$, $22$, $14)$ | $(94$, $108$, $96$, $38$, $62$, $36$, $98$, $110$, $92$, $116$, $124$, $128$, $112$, $90$, $114$, $126$, $-22$, $-100$, $-64$, $-104$, $-18$, $-28$, $-26$, $-16$, $-2$, $-4$, $-14$, $-24$, $-30$, $-20$, $-102$, $120$, $132$, $52$, $48$, $40$, $60$, $34$, $54$, $46$, $42$, $58$, $32$, $56$, $44$, $-80$, $-70$, $-106$, $-66$, $130$, $122$, $118$, $134$, $50$, $-68$, $-78$, $-88$, $-82$, $-72$, $-8$, $-10$, $-74$, $-84$, $-86$, $-76$, $-12$, $-6)$\\

\noindent$11a87$ : $3$ | $(4$, $10$, $12$, $16$, $18$, $2$, $20$, $8$, $6$, $22$, $14)$ | $(54$, $36$, $26$, $40$, $32$, $30$, $42$, $28$, $34$, $52$, $56$, $38$, $-90$, $-80$, $-44$, $-76$, $-46$, $-82$, $-88$, $-92$, $-78$, $64$, $74$, $60$, $58$, $-96$, $-84$, $-86$, $-94$, $-98$, $-8$, $-12$, $-14$, $-6$, $-24$, $-4$, $-16$, $-10$, $62$, $66$, $72$, $106$, $112$, $104$, $70$, $68$, $102$, $110$, $108$, $100$, $-48$, $-22$, $-2$, $-18$, $-50$, $-20)$\\

\noindent$11a88$ : $3$ | $(4$, $10$, $12$, $16$, $18$, $2$, $20$, $22$, $8$, $6$, $14)$ | $(54$, $34$, $40$, $92$, $166$, $146$, $150$, $170$, $96$, $36$, $38$, $94$, $168$, $148$, $-12$, $-74$, $-84$, $-122$, $-138$, $-124$, $-86$, $-72$, $-10$, $-26$, $-14$, $-76$, $-82$, $-80$, $-78$, $-16$, $-24$, $-8$, $-70$, $-88$, $90$, $164$, $144$, $152$, $172$, $98$, $174$, $154$, $142$, $162$, $-108$, $-110$, $-126$, $-136$, $-120$, $56$, $52$, $32$, $42$, $66$, $68$, $64$, $44$, $30$, $50$, $58$, $156$, $140$, $160$, $62$, $46$, $28$, $48$, $60$, $158$, $-20$, $-4$, $-104$, $-114$, $-130$, $-132$, $-116$, $-102$, $-2$, $-18$, $-22$, $-6$, $-106$, $-112$, $-128$, $-134$, $-118$, $-100)$\\

\noindent$11a89$ : $2$ | $(4$, $10$, $12$, $16$, $18$, $2$, $22$, $6$, $20$, $8$, $14)$| $(162$, $224$, $136$, $188$, $198$, $126$, $214$, $172$, $152$, $234$, $146$, $178$, $208$, $120$, $204$, $182$, $142$, $230$, $156$, $168$, $218$, $130$, $194$, $192$, $132$, $220$, $166$, $158$, $228$, $140$, $184$, $202$, $122$, $210$, $176$, $148$, $236$, $150$, $174$, $212$, $124$, $200$, $186$, $138$, $226$, $160$, $164$, $222$, $134$, $190$, $196$, $128$, $216$, $170$, $154$, $232$, $144$, $180$, $206$, $-282$, $-336$, $-244$, $-320$, $-298$, $-266$, $-352$, $-260$, $-304$, $-314$, $-250$, $-342$, $-276$, $-288$, $-330$, $-238$, $-326$, $-292$, $-272$, $-346$, $-254$, $-310$, $-308$, $-256$, $-348$, $-270$, $-294$, $-324$, $-240$, $-332$, $-286$, $-278$, $-340$, $-248$, $-316$, $-302$, $-262$, $-354$, $-264$, $-300$, $-318$, $-246$, $-338$, $-280$, $-284$, $-334$, $-242$, $-322$, $-296$, $-268$, $-350$, $-258$, $-306$, $-312$, $-252$, $-344$, $-274$, $-290$, $-328$, $398$, $460$, $372$, $424$, $434$, $362$, $450$, $408$, $388$, $470$, $382$, $414$, $444$, $356$, $440$, $418$, $378$, $466$, $392$, $404$, $454$, $366$, $430$, $428$, $368$, $456$, $402$, $394$, $464$, $376$, $420$, $438$, $358$, $446$, $412$, $384$, $472$, $386$, $410$, $448$, $360$, $436$, $422$, $374$, $462$, $396$, $400$, $458$, $370$, $426$, $432$, $364$, $452$, $406$, $390$, $468$, $380$, $416$, $442$, $-46$, $-100$, $-8$, $-84$, $-62$, $-30$, $-116$, $-24$, $-68$, $-78$, $-14$, $-106$, $-40$, $-52$, $-94$, $-2$, $-90$, $-56$, $-36$, $-110$, $-18$, $-74$, $-72$, $-20$, $-112$, $-34$, $-58$, $-88$, $-4$, $-96$, $-50$, $-42$, $-104$, $-12$, $-80$, $-66$, $-26$, $-118$, $-28$, $-64$, $-82$, $-10$, $-102$, $-44$, $-48$, $-98$, $-6$, $-86$, $-60$, $-32$, $-114$, $-22$, $-70$, $-76$, $-16$, $-108$, $-38$, $-54$, $-92)$\\

\noindent$11a90$ : $2$ | $(4$, $10$, $12$, $16$, $18$, $2$, $22$, $20$, $8$, $6$, $14)$ | $(120$, $158$, $90$, $150$, $128$, $112$, $166$, $98$, $142$, $136$, $104$, $172$, $106$, $134$, $144$, $96$, $164$, $114$, $126$, $152$, $88$, $156$, $122$, $118$, $160$, $92$, $148$, $130$, $110$, $168$, $100$, $140$, $138$, $102$, $170$, $108$, $132$, $146$, $94$, $162$, $116$, $124$, $154$, $-236$, $-190$, $-200$, $-246$, $-226$, $-180$, $-210$, $-256$, $-216$, $-174$, $-220$, $-252$, $-206$, $-184$, $-230$, $-242$, $-196$, $-194$, $-240$, $-232$, $-186$, $-204$, $-250$, $-222$, $-176$, $-214$, $-258$, $-212$, $-178$, $-224$, $-248$, $-202$, $-188$, $-234$, $-238$, $-192$, $-198$, $-244$, $-228$, $-182$, $-208$, $-254$, $-218$, $292$, $330$, $262$, $322$, $300$, $284$, $338$, $270$, $314$, $308$, $276$, $344$, $278$, $306$, $316$, $268$, $336$, $286$, $298$, $324$, $260$, $328$, $294$, $290$, $332$, $264$, $320$, $302$, $282$, $340$, $272$, $312$, $310$, $274$, $342$, $280$, $304$, $318$, $266$, $334$, $288$, $296$, $326$, $-64$, $-18$, $-28$, $-74$, $-54$, $-8$, $-38$, $-84$, $-44$, $-2$, $-48$, $-80$, $-34$, $-12$, $-58$, $-70$, $-24$, $-22$, $-68$, $-60$, $-14$, $-32$, $-78$, $-50$, $-4$, $-42$, $-86$, $-40$, $-6$, $-52$, $-76$, $-30$, $-16$, $-62$, $-66$, $-20$, $-26$, $-72$, $-56$, $-10$, $-36$, $-82$, $-46)$\\

\noindent$11a91$ : $2$ | $(4$, $10$, $12$, $16$, $20$, $2$, $22$, $18$, $8$, $6$, $14)$ | $(178$, $236$, $136$, $220$, $194$, $162$, $252$, $152$, $204$, $210$, $146$, $246$, $168$, $188$, $226$, $130$, $230$, $184$, $172$, $242$, $142$, $214$, $200$, $156$, $256$, $158$, $198$, $216$, $140$, $240$, $174$, $182$, $232$, $132$, $224$, $190$, $166$, $248$, $148$, $208$, $206$, $150$, $250$, $164$, $192$, $222$, $134$, $234$, $180$, $176$, $238$, $138$, $218$, $196$, $160$, $254$, $154$, $202$, $212$, $144$, $244$, $170$, $186$, $228$, $-336$, $-274$, $-372$, $-300$, $-310$, $-362$, $-264$, $-346$, $-326$, $-284$, $-382$, $-290$, $-320$, $-352$, $-258$, $-356$, $-316$, $-294$, $-378$, $-280$, $-330$, $-342$, $-268$, $-366$, $-306$, $-304$, $-368$, $-270$, $-340$, $-332$, $-278$, $-376$, $-296$, $-314$, $-358$, $-260$, $-350$, $-322$, $-288$, $-384$, $-286$, $-324$, $-348$, $-262$, $-360$, $-312$, $-298$, $-374$, $-276$, $-334$, $-338$, $-272$, $-370$, $-302$, $-308$, $-364$, $-266$, $-344$, $-328$, $-282$, $-380$, $-292$, $-318$, $-354$, $434$, $492$, $392$, $476$, $450$, $418$, $508$, $408$, $460$, $466$, $402$, $502$, $424$, $444$, $482$, $386$, $486$, $440$, $428$, $498$, $398$, $470$, $456$, $412$, $512$, $414$, $454$, $472$, $396$, $496$, $430$, $438$, $488$, $388$, $480$, $446$, $422$, $504$, $404$, $464$, $462$, $406$, $506$, $420$, $448$, $478$, $390$, $490$, $436$, $432$, $494$, $394$, $474$, $452$, $416$, $510$, $410$, $458$, $468$, $400$, $500$, $426$, $442$, $484$, $-80$, $-18$, $-116$, $-44$, $-54$, $-106$, $-8$, $-90$, $-70$, $-28$, $-126$, $-34$, $-64$, $-96$, $-2$, $-100$, $-60$, $-38$, $-122$, $-24$, $-74$, $-86$, $-12$, $-110$, $-50$, $-48$, $-112$, $-14$, $-84$, $-76$, $-22$, $-120$, $-40$, $-58$, $-102$, $-4$, $-94$, $-66$, $-32$, $-128$, $-30$, $-68$, $-92$, $-6$, $-104$, $-56$, $-42$, $-118$, $-20$, $-78$, $-82$, $-16$, $-114$, $-46$, $-52$, $-108$, $-10$, $-88$, $-72$, $-26$, $-124$, $-36$, $-62$, $-98)$\\

\noindent$11a92$ : $3$ | $(4$, $10$, $12$, $16$, $22$, $2$, $18$, $20$, $8$, $14$, $6)$ | $(264$, $254$, $336$, $284$, $116$, $334$, $286$, $282$, $338$, $256$, $266$, $302$, $318$, $244$, $324$, $296$, $272$, $344$, $276$, $292$, $328$, $248$, $314$, $306$, $400$, $308$, $312$, $250$, $330$, $290$, $278$, $342$, $270$, $298$, $322$, $242$, $320$, $300$, $268$, $258$, $340$, $280$, $288$, $332$, $252$, $310$, $-52$, $-46$, $-356$, $-24$, $-74$, $-68$, $-30$, $-350$, $-40$, $-58$, $-86$, $-12$, $-6$, $-16$, $-82$, $-62$, $-36$, $-346$, $-34$, $-64$, $-80$, $-18$, $-4$, $-88$, $-56$, $-42$, $-352$, $-28$, $-70$, $-72$, $-26$, $-354$, $-44$, $-54$, $-90$, $-50$, $-48$, $-358$, $-22$, $-76$, $-66$, $-32$, $-348$, $-38$, $-60$, $-84$, $-14$, $118$, $136$, $262$, $168$, $98$, $154$, $146$, $106$, $176$, $128$, $126$, $178$, $108$, $144$, $156$, $96$, $166$, $162$, $92$, $160$, $140$, $112$, $182$, $122$, $132$, $172$, $102$, $150$, $-378$, $-198$, $-214$, $-394$, $-226$, $-360$, $-190$, $-370$, $-78$, $-20$, $-2$, $-192$, $-372$, $-236$, $-384$, $-204$, $-208$, $-388$, $-232$, $-366$, $-186$, $-8$, $-364$, $-230$, $-390$, $-210$, $-202$, $-382$, $-238$, $-374$, $-194$, $-218$, $-398$, $-222$, $-224$, $-396$, $-216$, $-196$, $-376$, $-240$, $-380$, $-200$, $-212$, $-392$, $-228$, $-362$, $-10$, $-188$, $-368$, $-234$, $-386$, $-206$, $274$, $294$, $326$, $246$, $316$, $304$, $402$, $138$, $114$, $184$, $120$, $134$, $260$, $170$, $100$, $152$, $148$, $104$, $174$, $130$, $124$, $180$, $110$, $142$, $158$, $94$, $164$, $-220)$\\

\noindent$11a93$ : $2$ | $(4$, $10$, $12$, $16$, $22$, $2$, $20$, $18$, $8$, $14$, $6)$ | $(126$, $176$, $108$, $144$, $158$, $94$, $162$, $140$, $112$, $180$, $122$, $130$, $172$, $104$, $148$, $154$, $98$, $166$, $136$, $116$, $184$, $118$, $134$, $168$, $100$, $152$, $150$, $102$, $170$, $132$, $120$, $182$, $114$, $138$, $164$, $96$, $156$, $146$, $106$, $174$, $128$, $124$, $178$, $110$, $142$, $160$, $-236$, $-214$, $-258$, $-192$, $-274$, $-198$, $-252$, $-220$, $-230$, $-242$, $-208$, $-264$, $-186$, $-268$, $-204$, $-246$, $-226$, $-224$, $-248$, $-202$, $-270$, $-188$, $-262$, $-210$, $-240$, $-232$, $-218$, $-254$, $-196$, $-276$, $-194$, $-256$, $-216$, $-234$, $-238$, $-212$, $-260$, $-190$, $-272$, $-200$, $-250$, $-222$, $-228$, $-244$, $-206$, $-266$, $310$, $360$, $292$, $328$, $342$, $278$, $346$, $324$, $296$, $364$, $306$, $314$, $356$, $288$, $332$, $338$, $282$, $350$, $320$, $300$, $368$, $302$, $318$, $352$, $284$, $336$, $334$, $286$, $354$, $316$, $304$, $366$, $298$, $322$, $348$, $280$, $340$, $330$, $290$, $358$, $312$, $308$, $362$, $294$, $326$, $344$, $-52$, $-30$, $-74$, $-8$, $-90$, $-14$, $-68$, $-36$, $-46$, $-58$, $-24$, $-80$, $-2$, $-84$, $-20$, $-62$, $-42$, $-40$, $-64$, $-18$, $-86$, $-4$, $-78$, $-26$, $-56$, $-48$, $-34$, $-70$, $-12$, $-92$, $-10$, $-72$, $-32$, $-50$, $-54$, $-28$, $-76$, $-6$, $-88$, $-16$, $-66$, $-38$, $-44$, $-60$, $-22$, $-82)$\\

\noindent$11a94$ : $3$ | $(4$, $10$, $12$, $18$, $2$, $8$, $20$, $22$, $6$, $16$, $14)$ | $(30$, $92$, $98$, $34$, $96$, $94$, $32$, $100$, $90$, $106$, $104$, $88$, $102$, $108$, $-42$, $-40$, $-48$, $68$, $52$, $64$, $110$, $62$, $54$, $70$, $66$, $-10$, $-4$, $-16$, $-24$, $-22$, $-18$, $-2$, $-12$, $-8$, $-6$, $-14$, $-26$, $-20$, $60$, $56$, $72$, $28$, $74$, $58$, $-78$, $-44$, $-38$, $-50$, $-36$, $-46$, $-80$, $-86$, $-76$, $-84$, $-82)$\\

\noindent$11a95$ : $2$ | $(4$, $10$, $12$, $18$, $2$, $8$, $22$, $20$, $6$, $16$, $14)$ | $(112$, $98$, $126$, $84$, $140$, $74$, $136$, $88$, $122$, $102$, $108$, $116$, $94$, $130$, $80$, $144$, $78$, $132$, $92$, $118$, $106$, $104$, $120$, $90$, $134$, $76$, $142$, $82$, $128$, $96$, $114$, $110$, $100$, $124$, $86$, $138$, $-186$, $-164$, $-208$, $-146$, $-204$, $-168$, $-182$, $-190$, $-160$, $-212$, $-150$, $-200$, $-172$, $-178$, $-194$, $-156$, $-216$, $-154$, $-196$, $-176$, $-174$, $-198$, $-152$, $-214$, $-158$, $-192$, $-180$, $-170$, $-202$, $-148$, $-210$, $-162$, $-188$, $-184$, $-166$, $-206$, $256$, $242$, $270$, $228$, $284$, $218$, $280$, $232$, $266$, $246$, $252$, $260$, $238$, $274$, $224$, $288$, $222$, $276$, $236$, $262$, $250$, $248$, $264$, $234$, $278$, $220$, $286$, $226$, $272$, $240$, $258$, $254$, $244$, $268$, $230$, $282$, $-42$, $-20$, $-64$, $-2$, $-60$, $-24$, $-38$, $-46$, $-16$, $-68$, $-6$, $-56$, $-28$, $-34$, $-50$, $-12$, $-72$, $-10$, $-52$, $-32$, $-30$, $-54$, $-8$, $-70$, $-14$, $-48$, $-36$, $-26$, $-58$, $-4$, $-66$, $-18$, $-44$, $-40$, $-22$, $-62)$\\

\noindent$11a96$ : $2$ | $(4$, $10$, $12$, $18$, $14$, $2$, $22$, $8$, $20$, $6$, $16)$ | $(170$, $212$, $128$, $228$, $154$, $186$, $196$, $144$, $238$, $138$, $202$, $180$, $160$, $222$, $122$, $218$, $164$, $176$, $206$, $134$, $234$, $148$, $192$, $190$, $150$, $232$, $132$, $208$, $174$, $166$, $216$, $124$, $224$, $158$, $182$, $200$, $140$, $240$, $142$, $198$, $184$, $156$, $226$, $126$, $214$, $168$, $172$, $210$, $130$, $230$, $152$, $188$, $194$, $146$, $236$, $136$, $204$, $178$, $162$, $220$, $-286$, $-344$, $-252$, $-320$, $-310$, $-262$, $-354$, $-276$, $-296$, $-334$, $-242$, $-330$, $-300$, $-272$, $-358$, $-266$, $-306$, $-324$, $-248$, $-340$, $-290$, $-282$, $-348$, $-256$, $-316$, $-314$, $-258$, $-350$, $-280$, $-292$, $-338$, $-246$, $-326$, $-304$, $-268$, $-360$, $-270$, $-302$, $-328$, $-244$, $-336$, $-294$, $-278$, $-352$, $-260$, $-312$, $-318$, $-254$, $-346$, $-284$, $-288$, $-342$, $-250$, $-322$, $-308$, $-264$, $-356$, $-274$, $-298$, $-332$, $410$, $452$, $368$, $468$, $394$, $426$, $436$, $384$, $478$, $378$, $442$, $420$, $400$, $462$, $362$, $458$, $404$, $416$, $446$, $374$, $474$, $388$, $432$, $430$, $390$, $472$, $372$, $448$, $414$, $406$, $456$, $364$, $464$, $398$, $422$, $440$, $380$, $480$, $382$, $438$, $424$, $396$, $466$, $366$, $454$, $408$, $412$, $450$, $370$, $470$, $392$, $428$, $434$, $386$, $476$, $376$, $444$, $418$, $402$, $460$, $-46$, $-104$, $-12$, $-80$, $-70$, $-22$, $-114$, $-36$, $-56$, $-94$, $-2$, $-90$, $-60$, $-32$, $-118$, $-26$, $-66$, $-84$, $-8$, $-100$, $-50$, $-42$, $-108$, $-16$, $-76$, $-74$, $-18$, $-110$, $-40$, $-52$, $-98$, $-6$, $-86$, $-64$, $-28$, $-120$, $-30$, $-62$, $-88$, $-4$, $-96$, $-54$, $-38$, $-112$, $-20$, $-72$, $-78$, $-14$, $-106$, $-44$, $-48$, $-102$, $-10$, $-82$, $-68$, $-24$, $-116$, $-34$, $-58$, $-92)$\\

\noindent$11a97$ : $3$ | $(4$, $10$, $12$, $18$, $16$, $2$, $20$, $22$, $8$, $6$, $14)$ | $(16$, $26$, $12$, $10$, $-46$, $-48$, $-58$, $-54$, $-52$, $-60$, $-66$, $-62$, $-50$, $-56$, $14$, $18$, $24$, $42$, $20$, $22$, $-64$, $-6$, $-8$, $44$, $76$, $70$, $84$, $72$, $74$, $82$, $68$, $78$, $80$, $-38$, $-32$, $-2$, $-28$, $-4$, $-34$, $-40$, $-36$, $-30)$\\

\noindent$11a98$ : $2$ | $(4$, $10$, $12$, $18$, $16$, $2$, $22$, $20$, $8$, $6$, $14)$ | $(106$, $140$, $80$, $132$, $114$, $98$, $148$, $88$, $124$, $122$, $90$, $150$, $96$, $116$, $130$, $82$, $142$, $104$, $108$, $138$, $78$, $134$, $112$, $100$, $146$, $86$, $126$, $120$, $92$, $152$, $94$, $118$, $128$, $84$, $144$, $102$, $110$, $136$, $-170$, $-206$, $-216$, $-180$, $-160$, $-196$, $-226$, $-190$, $-154$, $-186$, $-222$, $-200$, $-164$, $-176$, $-212$, $-210$, $-174$, $-166$, $-202$, $-220$, $-184$, $-156$, $-192$, $-228$, $-194$, $-158$, $-182$, $-218$, $-204$, $-168$, $-172$, $-208$, $-214$, $-178$, $-162$, $-198$, $-224$, $-188$, $258$, $292$, $232$, $284$, $266$, $250$, $300$, $240$, $276$, $274$, $242$, $302$, $248$, $268$, $282$, $234$, $294$, $256$, $260$, $290$, $230$, $286$, $264$, $252$, $298$, $238$, $278$, $272$, $244$, $304$, $246$, $270$, $280$, $236$, $296$, $254$, $262$, $288$, $-18$, $-54$, $-64$, $-28$, $-8$, $-44$, $-74$, $-38$, $-2$, $-34$, $-70$, $-48$, $-12$, $-24$, $-60$, $-58$, $-22$, $-14$, $-50$, $-68$, $-32$, $-4$, $-40$, $-76$, $-42$, $-6$, $-30$, $-66$, $-52$, $-16$, $-20$, $-56$, $-62$, $-26$, $-10$, $-46$, $-72$, $-36)$\\

\noindent$11a99$ : $3$ | $(4$, $10$, $12$, $20$, $16$, $2$, $18$, $8$, $22$, $6$, $14)$ | $(42$, $76$, $48$, $64$, $68$, $52$, $72$, $38$, $46$, $78$, $44$, $40$, $74$, $50$, $66$, $-8$, $-56$, $-60$, $-58$, $-82$, $-96$, $-98$, $-80$, $-102$, $-92$, $-86$, $-108$, $70$, $62$, $114$, $-112$, $-104$, $-90$, $-88$, $-106$, $-110$, $-84$, $-94$, $-100$, $122$, $140$, $116$, $128$, $134$, $30$, $136$, $126$, $118$, $142$, $120$, $124$, $138$, $32$, $132$, $130$, $34$, $36$, $-14$, $-24$, $-2$, $-20$, $-18$, $-4$, $-26$, $-12$, $-54$, $-10$, $-28$, $-6$, $-16$, $-22)$\\

\noindent$11a100$ : $3$ | $(4$, $10$, $12$, $20$, $16$, $2$, $22$, $8$, $14$, $6$, $18)$ | $(62$, $92$, $82$, $118$, $76$, $98$, $102$, $72$, $114$, $86$, $88$, $112$, $70$, $104$, $96$, $78$, $120$, $80$, $94$, $106$, $68$, $110$, $90$, $84$, $116$, $74$, $100$, $-12$, $-28$, $-38$, $-64$, $222$, $108$, $-154$, $-128$, $-170$, $-138$, $-144$, $-164$, $-122$, $-160$, $-148$, $-134$, $-174$, $-132$, $-150$, $-158$, $-124$, $-166$, $-142$, $-140$, $-168$, $-126$, $-156$, $-152$, $-130$, $-172$, $-136$, $-146$, $-162$, $208$, $182$, $58$, $192$, $198$, $218$, $176$, $214$, $202$, $188$, $54$, $186$, $204$, $212$, $178$, $220$, $196$, $194$, $60$, $180$, $210$, $206$, $184$, $56$, $190$, $200$, $216$, $-20$, $-46$, $-4$, $-36$, $-30$, $-10$, $-52$, $-14$, $-26$, $-40$, $-66$, $-42$, $-24$, $-16$, $-50$, $-8$, $-32$, $-34$, $-6$, $-48$, $-18$, $-22$, $-44$, $-2)$\\

\noindent$11a101$ : $3$ | $(4$, $10$, $14$, $12$, $2$, $20$, $18$, $22$, $6$, $8$, $16)$ | $(152$, $146$, $130$, $168$, $172$, $162$, $136$, $140$, $158$, $64$, $62$, $156$, $142$, $134$, $164$, $174$, $166$, $132$, $144$, $154$, $60$, $66$, $160$, $138$, $-14$, $-26$, $-36$, $-4$, $-202$, $-198$, $-184$, $-122$, $-186$, $-196$, $-204$, $-6$, $-34$, $-28$, $-12$, $-46$, $-16$, $-24$, $-38$, $-2$, $-200$, $58$, $68$, $206$, $70$, $56$, $88$, $222$, $86$, $54$, $72$, $208$, $212$, $76$, $50$, $82$, $218$, $-20$, $-42$, $-92$, $-98$, $-108$, $-178$, $-116$, $-192$, $-190$, $-118$, $-180$, $-106$, $-100$, $-90$, $-102$, $-104$, $-182$, $-120$, $-188$, $-194$, $-114$, $-176$, $-110$, $-96$, $-94$, $-112$, $210$, $74$, $52$, $84$, $220$, $124$, $216$, $80$, $48$, $78$, $214$, $126$, $150$, $148$, $128$, $170$, $-8$, $-32$, $-30$, $-10$, $-44$, $-18$, $-22$, $-40)$\\

\noindent$11a102$ : $3$ | $(4$, $10$, $14$, $12$, $18$, $2$, $6$, $20$, $22$, $8$, $16)$ | $(92$, $110$, $86$, $78$, $84$, $108$, $94$, $90$, $112$, $88$, $46$, $-62$, $-50$, $-58$, $-12$, $-2$, $-16$, $-18$, $-4$, $-10$, $-56$, $-52$, $-64$, $-60$, $98$, $104$, $80$, $82$, $106$, $96$, $100$, $102$, $138$, $36$, $32$, $142$, $30$, $38$, $136$, $-122$, $-132$, $-120$, $-76$, $-66$, $-114$, $-70$, $-72$, $-116$, $-128$, $-126$, $-22$, $-124$, $-130$, $-118$, $-74$, $-68$, $34$, $140$, $28$, $40$, $134$, $44$, $24$, $48$, $26$, $42$, $-54$, $-8$, $-6$, $-20$, $-14)$\\

\noindent$11a103$ : $3$ | $(4$, $10$, $14$, $12$, $18$, $2$, $6$, $22$, $20$, $8$, $16)$ | $(66$, $74$, $58$, $132$, $136$, $-88$, $-6$, $130$, $138$, $142$, $126$, $150$, $118$, $152$, $124$, $144$, $10$, $86$, $12$, $146$, $122$, $154$, $120$, $148$, $128$, $140$, $-18$, $-32$, $-38$, $-48$, $-22$, $-28$, $-42$, $-44$, $-26$, $-24$, $-46$, $-40$, $-30$, $-20$, $-50$, $-16$, $-8$, $134$, $56$, $76$, $64$, $68$, $72$, $60$, $80$, $52$, $82$, $14$, $84$, $54$, $78$, $62$, $70$, $-96$, $-2$, $-92$, $-114$, $-100$, $-106$, $-108$, $-36$, $-34$, $-110$, $-104$, $-102$, $-112$, $-90$, $-4$, $-98$, $-116$, $-94)$\\

\noindent$11a104$ : $3$ | $(4$, $10$, $14$, $12$, $20$, $2$, $16$, $6$, $22$, $8$, $18)$ | $(56$, $38$, $68$, $44$, $114$, $42$, $70$, $40$, $54$, $58$, $36$, $66$, $46$, $136$, $50$, $62$, $-90$, $-78$, $-102$, $-108$, $-72$, $-112$, $-98$, $-82$, $-86$, $-94$, $-74$, $-106$, $-104$, $-76$, $-92$, $-88$, $-80$, $-100$, $-110$, $160$, $142$, $172$, $148$, $154$, $166$, $120$, $118$, $168$, $152$, $150$, $170$, $116$, $122$, $164$, $156$, $146$, $174$, $144$, $158$, $162$, $140$, $-124$, $-28$, $-134$, $-26$, $138$, $52$, $60$, $34$, $64$, $48$, $-84$, $-96$, $-8$, $-10$, $-16$, $-2$, $-20$, $-128$, $-32$, $-130$, $-22$, $-4$, $-14$, $-12$, $-6$, $-24$, $-132$, $-30$, $-126$, $-18)$\\

\noindent$11a105$ : $3$ | $(4$, $10$, $14$, $16$, $2$, $20$, $8$, $6$, $22$, $12$, $18)$ | $(52$, $34$, $118$, $38$, $56$, $48$, $30$, $44$, $60$, $42$, $98$, $40$, $58$, $46$, $-90$, $-64$, $-80$, $-74$, $-70$, $-84$, $-96$, $-62$, $-92$, $-88$, $-66$, $-78$, $-76$, $-68$, $-86$, $-94$, $142$, $124$, $130$, $148$, $136$, $104$, $102$, $134$, $150$, $132$, $100$, $106$, $138$, $146$, $128$, $126$, $144$, $140$, $122$, $-108$, $-4$, $-116$, $-6$, $120$, $32$, $50$, $54$, $36$, $-72$, $-82$, $-22$, $-20$, $-14$, $-28$, $-12$, $-110$, $-2$, $-114$, $-8$, $-24$, $-18$, $-16$, $-26$, $-10$, $-112)$\\

\noindent$11a106$ : $3$ | $(4$, $10$, $14$, $16$, $2$, $20$, $18$, $6$, $8$, $22$, $12)$ | $(54$, $52$, $86$, $-74$, $-66$, $-62$, $-78$, $-60$, $-68$, $-72$, $-4$, $-76$, $-64$, $10$, $16$, $50$, $20$, $6$, $22$, $14$, $12$, $24$, $8$, $18$, $-70$, $-58$, $-56$, $104$, $102$, $96$, $80$, $92$, $90$, $82$, $98$, $100$, $84$, $88$, $94$, $-40$, $-28$, $-46$, $-34$, $-36$, $-44$, $-26$, $-42$, $-38$, $-30$, $-48$, $-32$, $-2)$\\

\noindent$11a107$ : $3$ | $(4$, $10$, $14$, $16$, $2$, $20$, $18$, $8$, $6$, $22$, $12)$ | $(132$, $154$, $124$, $140$, $146$, $118$, $148$, $138$, $126$, $156$, $130$, $134$, $152$, $122$, $142$, $144$, $120$, $150$, $136$, $162$, $-38$, $-16$, $-6$, $-28$, $-32$, $-10$, $-12$, $-34$, $-26$, $-4$, $-18$, $-158$, $-20$, $-2$, $-24$, $-36$, $-14$, $-8$, $-30$, $48$, $56$, $70$, $164$, $64$, $62$, $42$, $72$, $54$, $50$, $76$, $46$, $58$, $68$, $166$, $66$, $60$, $44$, $74$, $52$, $-96$, $-80$, $-102$, $-112$, $-90$, $-86$, $-108$, $-106$, $-84$, $-92$, $-114$, $-100$, $-78$, $-98$, $-116$, $-94$, $-82$, $-104$, $-110$, $-88$, $128$, $40$, $-160$, $-22)$\\

\noindent$11a108$ : $3$ | $(4$, $10$, $14$, $16$, $2$, $20$, $22$, $18$, $6$, $8$, $12)$ | $(48$, $78$, $82$, $44$, $50$, $46$, $80$, $-4$, $-32$, $-38$, $-40$, $-34$, $-2$, $-12$, $-6$, $84$, $76$, $90$, $86$, $74$, $88$, $-36$, $-52$, $-58$, $-60$, $92$, $28$, $14$, $26$, $94$, $24$, $16$, $30$, $18$, $22$, $42$, $20$, $-68$, $-62$, $-56$, $-54$, $-64$, $-66$, $-72$, $-70$, $-8$, $-10)$\\

\noindent$11a109$ : $3$ | $(4$, $10$, $14$, $16$, $2$, $20$, $22$, $18$, $8$, $6$, $12)$ | $(54$, $40$, $34$, $48$, $58$, $50$, $36$, $38$, $52$, $56$, $42$, $32$, $46$, $112$, $44$, $-108$, $-72$, $-62$, $-80$, $-64$, $-70$, $-106$, $-110$, $-74$, $-60$, $-78$, $-66$, $-68$, $-76$, $116$, $92$, $82$, $96$, $120$, $98$, $84$, $90$, $114$, $118$, $94$, $-14$, $-4$, $-22$, $-30$, $-24$, $-6$, $-12$, $-16$, $-2$, $-20$, $-28$, $-26$, $100$, $86$, $88$, $102$, $-104$, $-8$, $-10$, $-18)$\\

\noindent$11a110$ : $2$ | $(4$, $10$, $14$, $16$, $2$, $22$, $20$, $18$, $6$, $8$, $12)$ | $(158$, $104$, $174$, $142$, $120$, $190$, $126$, $136$, $180$, $110$, $152$, $164$, $98$, $168$, $148$, $114$, $184$, $132$, $130$, $186$, $116$, $146$, $170$, $100$, $162$, $154$, $108$, $178$, $138$, $124$, $192$, $122$, $140$, $176$, $106$, $156$, $160$, $102$, $172$, $144$, $118$, $188$, $128$, $134$, $182$, $112$, $150$, $166$, $-228$, $-278$, $-206$, $-250$, $-256$, $-200$, $-272$, $-234$, $-222$, $-284$, $-212$, $-244$, $-262$, $-194$, $-266$, $-240$, $-216$, $-288$, $-218$, $-238$, $-268$, $-196$, $-260$, $-246$, $-210$, $-282$, $-224$, $-232$, $-274$, $-202$, $-254$, $-252$, $-204$, $-276$, $-230$, $-226$, $-280$, $-208$, $-248$, $-258$, $-198$, $-270$, $-236$, $-220$, $-286$, $-214$, $-242$, $-264$, $350$, $296$, $366$, $334$, $312$, $382$, $318$, $328$, $372$, $302$, $344$, $356$, $290$, $360$, $340$, $306$, $376$, $324$, $322$, $378$, $308$, $338$, $362$, $292$, $354$, $346$, $300$, $370$, $330$, $316$, $384$, $314$, $332$, $368$, $298$, $348$, $352$, $294$, $364$, $336$, $310$, $380$, $320$, $326$, $374$, $304$, $342$, $358$, $-36$, $-86$, $-14$, $-58$, $-64$, $-8$, $-80$, $-42$, $-30$, $-92$, $-20$, $-52$, $-70$, $-2$, $-74$, $-48$, $-24$, $-96$, $-26$, $-46$, $-76$, $-4$, $-68$, $-54$, $-18$, $-90$, $-32$, $-40$, $-82$, $-10$, $-62$, $-60$, $-12$, $-84$, $-38$, $-34$, $-88$, $-16$, $-56$, $-66$, $-6$, $-78$, $-44$, $-28$, $-94$, $-22$, $-50$, $-72)$\\

\noindent$11a111$ : $2$ | $(4$, $10$, $14$, $16$, $2$, $22$, $20$, $18$, $8$, $6$, $12)$ | $(168$, $110$, $184$, $152$, $126$, $200$, $136$, $142$, $194$, $120$, $158$, $178$, $104$, $174$, $162$, $116$, $190$, $146$, $132$, $204$, $130$, $148$, $188$, $114$, $164$, $172$, $106$, $180$, $156$, $122$, $196$, $140$, $138$, $198$, $124$, $154$, $182$, $108$, $170$, $166$, $112$, $186$, $150$, $128$, $202$, $134$, $144$, $192$, $118$, $160$, $176$, $-268$, $-218$, $-296$, $-240$, $-246$, $-290$, $-212$, $-274$, $-262$, $-224$, $-302$, $-234$, $-252$, $-284$, $-206$, $-280$, $-256$, $-230$, $-306$, $-228$, $-258$, $-278$, $-208$, $-286$, $-250$, $-236$, $-300$, $-222$, $-264$, $-272$, $-214$, $-292$, $-244$, $-242$, $-294$, $-216$, $-270$, $-266$, $-220$, $-298$, $-238$, $-248$, $-288$, $-210$, $-276$, $-260$, $-226$, $-304$, $-232$, $-254$, $-282$, $372$, $314$, $388$, $356$, $330$, $404$, $340$, $346$, $398$, $324$, $362$, $382$, $308$, $378$, $366$, $320$, $394$, $350$, $336$, $408$, $334$, $352$, $392$, $318$, $368$, $376$, $310$, $384$, $360$, $326$, $400$, $344$, $342$, $402$, $328$, $358$, $386$, $312$, $374$, $370$, $316$, $390$, $354$, $332$, $406$, $338$, $348$, $396$, $322$, $364$, $380$, $-64$, $-14$, $-92$, $-36$, $-42$, $-86$, $-8$, $-70$, $-58$, $-20$, $-98$, $-30$, $-48$, $-80$, $-2$, $-76$, $-52$, $-26$, $-102$, $-24$, $-54$, $-74$, $-4$, $-82$, $-46$, $-32$, $-96$, $-18$, $-60$, $-68$, $-10$, $-88$, $-40$, $-38$, $-90$, $-12$, $-66$, $-62$, $-16$, $-94$, $-34$, $-44$, $-84$, $-6$, $-72$, $-56$, $-22$, $-100$, $-28$, $-50$, $-78)$\\

\noindent$11a112$ : $3$ | $(4$, $10$, $14$, $16$, $18$, $2$, $20$, $6$, $22$, $12$, $8)$ | $(66$, $170$, $196$, $186$, $180$, $202$, $124$, $204$, $178$, $188$, $194$, $172$, $64$, $68$, $168$, $198$, $184$, $182$, $200$, $122$, $206$, $112$, $86$, $-158$, $-20$, $-6$, $-32$, $-36$, $-10$, $-16$, $-42$, $-102$, $-136$, $-138$, $-100$, $-40$, $-14$, $-12$, $-38$, $-30$, $-4$, $-22$, $-160$, $-156$, $-18$, $-8$, $-34$, $54$, $78$, $120$, $208$, $114$, $84$, $48$, $88$, $-154$, $-162$, $-24$, $-2$, $-28$, $-166$, $-150$, $74$, $58$, $90$, $50$, $82$, $116$, $210$, $118$, $80$, $52$, $92$, $56$, $76$, $72$, $60$, $176$, $190$, $192$, $174$, $62$, $70$, $-98$, $-140$, $-134$, $-104$, $-44$, $-110$, $-128$, $-146$, $-94$, $-144$, $-130$, $-108$, $-46$, $-106$, $-132$, $-142$, $-96$, $-148$, $-126$, $-152$, $-164$, $-26)$\\

\noindent$11a113$ : $3$ | $(4$, $10$, $14$, $16$, $18$, $2$, $20$, $22$, $8$, $12$, $6)$ | $(50$, $72$, $44$, $56$, $66$, $138$, $64$, $58$, $42$, $70$, $52$, $48$, $74$, $46$, $54$, $68$, $40$, $60$, $62$, $-88$, $-82$, $-132$, $-122$, $-76$, $-126$, $-128$, $-78$, $-120$, $-134$, $-84$, $-90$, $-86$, $-136$, $-118$, $-80$, $-130$, $-124$, $104$, $150$, $98$, $110$, $144$, $92$, $140$, $142$, $-36$, $-14$, $-8$, $-30$, $-20$, $-2$, $-24$, $-26$, $-4$, $-18$, $-32$, $-10$, $-12$, $112$, $96$, $148$, $106$, $102$, $152$, $100$, $108$, $146$, $94$, $114$, $-116$, $-38$, $-34$, $-16$, $-6$, $-28$, $-22)$\\

\noindent$11a114$ : $3$ | $(4$, $10$, $14$, $18$, $2$, $16$, $6$, $22$, $20$, $8$, $12)$ | $(182$, $216$, $196$, $258$, $250$, $204$, $208$, $246$, $262$, $192$, $220$, $186$, $188$, $222$, $190$, $264$, $244$, $210$, $202$, $252$, $256$, $198$, $214$, $180$, $184$, $218$, $194$, $260$, $248$, $206$, $-12$, $-36$, $-164$, $-44$, $-4$, $-52$, $-20$, $-28$, $-26$, $-22$, $-50$, $-2$, $-46$, $-162$, $-34$, $-14$, $-58$, $-10$, $-38$, $-166$, $-42$, $-6$, $-54$, $-18$, $-30$, $-24$, $78$, $74$, $108$, $176$, $90$, $62$, $96$, $170$, $102$, $68$, $84$, $268$, $82$, $70$, $104$, $172$, $94$, $60$, $92$, $174$, $106$, $72$, $80$, $266$, $242$, $212$, $200$, $254$, $-40$, $-8$, $-56$, $-16$, $-32$, $-160$, $-134$, $-136$, $-158$, $-112$, $-224$, $-116$, $-154$, $-140$, $-130$, $-238$, $-232$, $-124$, $-146$, $-148$, $-122$, $-230$, $-240$, $-132$, $-138$, $-156$, $-114$, $76$, $110$, $178$, $88$, $64$, $98$, $168$, $100$, $66$, $86$, $-228$, $-120$, $-150$, $-144$, $-126$, $-234$, $-236$, $-128$, $-142$, $-152$, $-118$, $-226$, $-48)$\\

\noindent$11a115$ : $3$ | $(4$, $10$, $14$, $18$, $2$, $20$, $6$, $22$, $12$, $8$, $16)$ | $(86$, $270$, $118$, $74$, $98$, $258$, $106$, $66$, $110$, $262$, $94$, $78$, $122$, $80$, $92$, $264$, $112$, $68$, $104$, $256$, $100$, $72$, $116$, $268$, $88$, $84$, $272$, $120$, $76$, $96$, $260$, $108$, $-162$, $-132$, $-192$, $-184$, $-140$, $-154$, $-170$, $-124$, $-200$, $-176$, $-148$, $-146$, $-178$, $-198$, $-126$, $-168$, $-156$, $-138$, $-186$, $-190$, $-134$, $-160$, $-164$, $-130$, $-194$, $-182$, $-142$, $-152$, $-172$, $276$, $306$, $212$, $230$, $288$, $294$, $224$, $218$, $300$, $282$, $236$, $206$, $312$, $204$, $238$, $280$, $302$, $216$, $226$, $292$, $290$, $228$, $214$, $304$, $278$, $240$, $202$, $310$, $208$, $234$, $284$, $298$, $220$, $222$, $296$, $286$, $232$, $210$, $308$, $-242$, $-52$, $-2$, $-48$, $-246$, $-30$, $-20$, $-60$, $-10$, $-40$, $-254$, $-38$, $-12$, $-62$, $-18$, $-32$, $-248$, $-46$, $-4$, $-54$, $274$, $82$, $90$, $266$, $114$, $70$, $102$, $-188$, $-136$, $-158$, $-166$, $-128$, $-196$, $-180$, $-144$, $-150$, $-174$, $-26$, $-24$, $-56$, $-6$, $-44$, $-250$, $-34$, $-16$, $-64$, $-14$, $-36$, $-252$, $-42$, $-8$, $-58$, $-22$, $-28$, $-244$, $-50)$\\

\noindent$11a116$ : $3$ | $(4$, $10$, $14$, $18$, $2$, $20$, $8$, $22$, $6$, $12$, $16)$ | $(60$, $38$, $86$, $44$, $54$, $94$, $52$, $46$, $88$, $36$, $62$, $58$, $40$, $84$, $42$, $56$, $64$, $-76$, $-102$, $-108$, $-82$, $-112$, $-98$, $-120$, $-126$, $-118$, $-96$, $-114$, $-80$, $-106$, $-104$, $-78$, $-74$, $90$, $48$, $50$, $92$, $142$, $140$, $164$, $134$, $148$, $-110$, $-100$, $-122$, $-124$, $-116$, $160$, $130$, $152$, $144$, $138$, $166$, $136$, $146$, $150$, $132$, $162$, $158$, $128$, $154$, $66$, $156$, $-70$, $-14$, $-8$, $-30$, $-24$, $-2$, $-20$, $-34$, $-18$, $-4$, $-26$, $-28$, $-6$, $-16$, $-68$, $-72$, $-12$, $-10$, $-32$, $-22)$\\

\noindent$11a117$ : $2$ | $(4$, $10$, $14$, $18$, $2$, $22$, $8$, $20$, $6$, $16$, $12)$ | $(190$, $128$, $214$, $166$, $152$, $228$, $142$, $176$, $204$, $118$, $200$, $180$, $138$, $224$, $156$, $162$, $218$, $132$, $186$, $194$, $124$, $210$, $170$, $148$, $232$, $146$, $172$, $208$, $122$, $196$, $184$, $134$, $220$, $160$, $158$, $222$, $136$, $182$, $198$, $120$, $206$, $174$, $144$, $230$, $150$, $168$, $212$, $126$, $192$, $188$, $130$, $216$, $164$, $154$, $226$, $140$, $178$, $202$, $-300$, $-262$, $-338$, $-240$, $-322$, $-278$, $-284$, $-316$, $-246$, $-344$, $-256$, $-306$, $-294$, $-268$, $-332$, $-234$, $-328$, $-272$, $-290$, $-310$, $-252$, $-348$, $-250$, $-312$, $-288$, $-274$, $-326$, $-236$, $-334$, $-266$, $-296$, $-304$, $-258$, $-342$, $-244$, $-318$, $-282$, $-280$, $-320$, $-242$, $-340$, $-260$, $-302$, $-298$, $-264$, $-336$, $-238$, $-324$, $-276$, $-286$, $-314$, $-248$, $-346$, $-254$, $-308$, $-292$, $-270$, $-330$, $422$, $360$, $446$, $398$, $384$, $460$, $374$, $408$, $436$, $350$, $432$, $412$, $370$, $456$, $388$, $394$, $450$, $364$, $418$, $426$, $356$, $442$, $402$, $380$, $464$, $378$, $404$, $440$, $354$, $428$, $416$, $366$, $452$, $392$, $390$, $454$, $368$, $414$, $430$, $352$, $438$, $406$, $376$, $462$, $382$, $400$, $444$, $358$, $424$, $420$, $362$, $448$, $396$, $386$, $458$, $372$, $410$, $434$, $-68$, $-30$, $-106$, $-8$, $-90$, $-46$, $-52$, $-84$, $-14$, $-112$, $-24$, $-74$, $-62$, $-36$, $-100$, $-2$, $-96$, $-40$, $-58$, $-78$, $-20$, $-116$, $-18$, $-80$, $-56$, $-42$, $-94$, $-4$, $-102$, $-34$, $-64$, $-72$, $-26$, $-110$, $-12$, $-86$, $-50$, $-48$, $-88$, $-10$, $-108$, $-28$, $-70$, $-66$, $-32$, $-104$, $-6$, $-92$, $-44$, $-54$, $-82$, $-16$, $-114$, $-22$, $-76$, $-60$, $-38$, $-98)$\\

\noindent$11a118$ : $3$ | $(4$, $10$, $14$, $18$, $12$, $2$, $8$, $20$, $22$, $6$, $16)$ | $(6$, $138$, $-136$, $-34$, $-40$, $-24$, $-50$, $70$, $90$, $56$, $86$, $74$, $66$, $94$, $60$, $82$, $78$, $62$, $96$, $64$, $76$, $84$, $58$, $92$, $68$, $72$, $88$, $-114$, $-124$, $-104$, $-134$, $-98$, $-130$, $-108$, $-120$, $-118$, $-110$, $-128$, $-100$, $-2$, $-4$, $-102$, $-126$, $-112$, $-116$, $-122$, $-106$, $-132$, $154$, $8$, $140$, $168$, $158$, $150$, $12$, $144$, $164$, $162$, $146$, $14$, $148$, $160$, $166$, $142$, $10$, $152$, $156$, $170$, $80$, $-32$, $-42$, $-22$, $-52$, $-16$, $-48$, $-26$, $-38$, $-36$, $-28$, $-46$, $-18$, $-54$, $-20$, $-44$, $-30)$\\

\noindent$11a119$ : $2$ | $(4$, $10$, $14$, $18$, $12$, $2$, $8$, $22$, $20$, $6$, $16)$ | $(110$, $128$, $92$, $146$, $78$, $142$, $96$, $124$, $114$, $106$, $132$, $88$, $150$, $82$, $138$, $100$, $120$, $118$, $102$, $136$, $84$, $152$, $86$, $134$, $104$, $116$, $122$, $98$, $140$, $80$, $148$, $90$, $130$, $108$, $112$, $126$, $94$, $144$, $-186$, $-204$, $-168$, $-222$, $-154$, $-218$, $-172$, $-200$, $-190$, $-182$, $-208$, $-164$, $-226$, $-158$, $-214$, $-176$, $-196$, $-194$, $-178$, $-212$, $-160$, $-228$, $-162$, $-210$, $-180$, $-192$, $-198$, $-174$, $-216$, $-156$, $-224$, $-166$, $-206$, $-184$, $-188$, $-202$, $-170$, $-220$, $262$, $280$, $244$, $298$, $230$, $294$, $248$, $276$, $266$, $258$, $284$, $240$, $302$, $234$, $290$, $252$, $272$, $270$, $254$, $288$, $236$, $304$, $238$, $286$, $256$, $268$, $274$, $250$, $292$, $232$, $300$, $242$, $282$, $260$, $264$, $278$, $246$, $296$, $-34$, $-52$, $-16$, $-70$, $-2$, $-66$, $-20$, $-48$, $-38$, $-30$, $-56$, $-12$, $-74$, $-6$, $-62$, $-24$, $-44$, $-42$, $-26$, $-60$, $-8$, $-76$, $-10$, $-58$, $-28$, $-40$, $-46$, $-22$, $-64$, $-4$, $-72$, $-14$, $-54$, $-32$, $-36$, $-50$, $-18$, $-68)$\\

\noindent$11a120$ : $2$ | $(4$, $10$, $14$, $18$, $12$, $2$, $22$, $6$, $20$, $8$, $16)$ | $(172$, $134$, $210$, $120$, $186$, $158$, $148$, $196$, $110$, $200$, $144$, $162$, $182$, $124$, $214$, $130$, $176$, $168$, $138$, $206$, $116$, $190$, $154$, $152$, $192$, $114$, $204$, $140$, $166$, $178$, $128$, $216$, $126$, $180$, $164$, $142$, $202$, $112$, $194$, $150$, $156$, $188$, $118$, $208$, $136$, $170$, $174$, $132$, $212$, $122$, $184$, $160$, $146$, $198$, $-262$, $-296$, $-228$, $-320$, $-238$, $-286$, $-272$, $-252$, $-306$, $-218$, $-310$, $-248$, $-276$, $-282$, $-242$, $-316$, $-224$, $-300$, $-258$, $-266$, $-292$, $-232$, $-324$, $-234$, $-290$, $-268$, $-256$, $-302$, $-222$, $-314$, $-244$, $-280$, $-278$, $-246$, $-312$, $-220$, $-304$, $-254$, $-270$, $-288$, $-236$, $-322$, $-230$, $-294$, $-264$, $-260$, $-298$, $-226$, $-318$, $-240$, $-284$, $-274$, $-250$, $-308$, $388$, $350$, $426$, $336$, $402$, $374$, $364$, $412$, $326$, $416$, $360$, $378$, $398$, $340$, $430$, $346$, $392$, $384$, $354$, $422$, $332$, $406$, $370$, $368$, $408$, $330$, $420$, $356$, $382$, $394$, $344$, $432$, $342$, $396$, $380$, $358$, $418$, $328$, $410$, $366$, $372$, $404$, $334$, $424$, $352$, $386$, $390$, $348$, $428$, $338$, $400$, $376$, $362$, $414$, $-46$, $-80$, $-12$, $-104$, $-22$, $-70$, $-56$, $-36$, $-90$, $-2$, $-94$, $-32$, $-60$, $-66$, $-26$, $-100$, $-8$, $-84$, $-42$, $-50$, $-76$, $-16$, $-108$, $-18$, $-74$, $-52$, $-40$, $-86$, $-6$, $-98$, $-28$, $-64$, $-62$, $-30$, $-96$, $-4$, $-88$, $-38$, $-54$, $-72$, $-20$, $-106$, $-14$, $-78$, $-48$, $-44$, $-82$, $-10$, $-102$, $-24$, $-68$, $-58$, $-34$, $-92)$\\

\noindent$11a121$ : $2$ | $(4$, $10$, $14$, $18$, $12$, $2$, $22$, $8$, $20$, $6$, $16)$ | $(168$, $206$, $130$, $230$, $144$, $192$, $182$, $154$, $220$, $120$, $216$, $158$, $178$, $196$, $140$, $234$, $134$, $202$, $172$, $164$, $210$, $126$, $226$, $148$, $188$, $186$, $150$, $224$, $124$, $212$, $162$, $174$, $200$, $136$, $236$, $138$, $198$, $176$, $160$, $214$, $122$, $222$, $152$, $184$, $190$, $146$, $228$, $128$, $208$, $166$, $170$, $204$, $132$, $232$, $142$, $194$, $180$, $156$, $218$, $-286$, $-324$, $-248$, $-348$, $-262$, $-310$, $-300$, $-272$, $-338$, $-238$, $-334$, $-276$, $-296$, $-314$, $-258$, $-352$, $-252$, $-320$, $-290$, $-282$, $-328$, $-244$, $-344$, $-266$, $-306$, $-304$, $-268$, $-342$, $-242$, $-330$, $-280$, $-292$, $-318$, $-254$, $-354$, $-256$, $-316$, $-294$, $-278$, $-332$, $-240$, $-340$, $-270$, $-302$, $-308$, $-264$, $-346$, $-246$, $-326$, $-284$, $-288$, $-322$, $-250$, $-350$, $-260$, $-312$, $-298$, $-274$, $-336$, $404$, $442$, $366$, $466$, $380$, $428$, $418$, $390$, $456$, $356$, $452$, $394$, $414$, $432$, $376$, $470$, $370$, $438$, $408$, $400$, $446$, $362$, $462$, $384$, $424$, $422$, $386$, $460$, $360$, $448$, $398$, $410$, $436$, $372$, $472$, $374$, $434$, $412$, $396$, $450$, $358$, $458$, $388$, $420$, $426$, $382$, $464$, $364$, $444$, $402$, $406$, $440$, $368$, $468$, $378$, $430$, $416$, $392$, $454$, $-50$, $-88$, $-12$, $-112$, $-26$, $-74$, $-64$, $-36$, $-102$, $-2$, $-98$, $-40$, $-60$, $-78$, $-22$, $-116$, $-16$, $-84$, $-54$, $-46$, $-92$, $-8$, $-108$, $-30$, $-70$, $-68$, $-32$, $-106$, $-6$, $-94$, $-44$, $-56$, $-82$, $-18$, $-118$, $-20$, $-80$, $-58$, $-42$, $-96$, $-4$, $-104$, $-34$, $-66$, $-72$, $-28$, $-110$, $-10$, $-90$, $-48$, $-52$, $-86$, $-14$, $-114$, $-24$, $-76$, $-62$, $-38$, $-100)$\\

\noindent$11a122$ : $3$ | $(4$, $10$, $14$, $18$, $16$, $2$, $20$, $22$, $8$, $12$, $6)$ | $(98$, $92$, $132$, $120$, $134$, $90$, $100$, $50$, $64$, $52$, $38$, $40$, $54$, $62$, $48$, $-78$, $-2$, $-12$, $-68$, $-86$, $-70$, $-10$, $-4$, $-76$, $-80$, $-66$, $-84$, $-72$, $-8$, $-6$, $-74$, $-82$, $126$, $140$, $88$, $136$, $122$, $130$, $94$, $96$, $128$, $124$, $138$, $-22$, $-102$, $-116$, $-110$, $-30$, $-112$, $-114$, $36$, $42$, $56$, $60$, $46$, $32$, $142$, $34$, $44$, $58$, $-106$, $-26$, $-18$, $-16$, $-28$, $-108$, $-118$, $-104$, $-24$, $-20$, $-14)$\\

\noindent$11a123$ : $3$ | $(4$, $10$, $14$, $20$, $2$, $8$, $18$, $22$, $6$, $12$, $16)$ | $(24$, $72$, $70$, $68$, $74$, $62$, $80$, $60$, $76$, $66$, $110$, $-52$, $-50$, $-20$, $-54$, $-10$, $-16$, $-58$, $-14$, $-12$, $-56$, $-18$, $-8$, $-4$, $106$, $112$, $108$, $64$, $78$, $-96$, $-82$, $-100$, $-92$, $-86$, $-104$, $-48$, $-102$, $-84$, $-94$, $-98$, $38$, $32$, $44$, $26$, $22$, $28$, $42$, $34$, $36$, $40$, $30$, $46$, $-6$, $-2$, $-90$, $-88)$\\

\noindent$11a124$ : $3$ | $(4$, $10$, $14$, $20$, $2$, $16$, $8$, $22$, $12$, $6$, $18)$ | $(22$, $60$, $56$, $70$, $54$, $62$, $-34$, $-48$, $-28$, $-44$, $-38$, $-36$, $-46$, $120$, $68$, $52$, $64$, $116$, $124$, $114$, $98$, $112$, $122$, $118$, $66$, $-88$, $-78$, $-72$, $-42$, $-40$, $-74$, $-76$, $-90$, $-86$, $-80$, $110$, $100$, $12$, $102$, $108$, $18$, $26$, $16$, $106$, $104$, $14$, $24$, $20$, $58$, $-4$, $-10$, $-32$, $-50$, $-30$, $-8$, $-6$, $-96$, $-2$, $-92$, $-84$, $-82$, $-94)$\\

\noindent$11a125$ : $3$ | $(4$, $10$, $14$, $20$, $2$, $16$, $18$, $8$, $22$, $6$, $12)$ | $(126$, $184$, $158$, $152$, $190$, $142$, $168$, $174$, $136$, $92$, $132$, $178$, $164$, $146$, $194$, $148$, $162$, $180$, $130$, $90$, $138$, $172$, $170$, $140$, $88$, $86$, $282$, $-118$, $-6$, $-106$, $-32$, $-26$, $-100$, $-12$, $-46$, $-42$, $-16$, $-96$, $-22$, $-36$, $-110$, $-2$, $-114$, $-122$, $-262$, $-204$, $-242$, $134$, $176$, $166$, $144$, $192$, $150$, $160$, $182$, $128$, $124$, $186$, $156$, $154$, $188$, $-50$, $-222$, $-224$, $-260$, $-202$, $-244$, $-240$, $-206$, $-264$, $-212$, $-234$, $-250$, $-196$, $-254$, $-230$, $-216$, $-54$, $-218$, $-228$, $-256$, $-198$, $-248$, $-236$, $-210$, $-266$, $-208$, $-238$, $-246$, $-200$, $-258$, $-226$, $-220$, $-52$, $-214$, $-232$, $-252$, $290$, $62$, $274$, $78$, $306$, $300$, $72$, $268$, $68$, $296$, $284$, $56$, $280$, $84$, $310$, $82$, $278$, $58$, $286$, $294$, $66$, $270$, $74$, $302$, $304$, $76$, $272$, $64$, $292$, $288$, $60$, $276$, $80$, $308$, $298$, $70$, $-44$, $-14$, $-98$, $-24$, $-34$, $-108$, $-4$, $-116$, $-120$, $-8$, $-104$, $-30$, $-28$, $-102$, $-10$, $-48$, $-40$, $-18$, $-94$, $-20$, $-38$, $-112)$\\

\noindent$11a126$ : $3$ | $(4$, $10$, $14$, $20$, $2$, $16$, $22$, $8$, $12$, $6$, $18)$ | $(38$, $30$, $66$, $28$, $70$, $80$, $72$, $122$, $106$, $120$, $74$, $78$, $68$, $-96$, $-64$, $-58$, $-18$, $-16$, $-60$, $-62$, $-14$, $-20$, $-56$, $-24$, $-10$, $-12$, $-22$, $76$, $118$, $108$, $124$, $110$, $116$, $-94$, $-98$, $-82$, $-102$, $-90$, $-92$, $-100$, $50$, $40$, $36$, $32$, $44$, $54$, $46$, $114$, $112$, $48$, $52$, $42$, $34$, $-86$, $-2$, $-8$, $-26$, $-6$, $-4$, $-88$, $-104$, $-84)$\\

\noindent$11a127$ : $3$ | $(4$, $10$, $14$, $20$, $2$, $22$, $18$, $8$, $12$, $6$, $16)$ | $(68$, $196$, $198$, $66$, $70$, $194$, $200$, $64$, $216$, $178$, $188$, $206$, $58$, $210$, $184$, $182$, $212$, $60$, $204$, $190$, $176$, $192$, $202$, $62$, $214$, $180$, $186$, $208$, $-160$, $-138$, $-170$, $-150$, $-100$, $-92$, $-94$, $-18$, $-50$, $-32$, $-28$, $-54$, $-22$, $-38$, $-44$, $-12$, $-4$, $120$, $218$, $122$, $88$, $118$, $-16$, $-48$, $-34$, $-26$, $-56$, $-24$, $-36$, $-46$, $-14$, $-6$, $-2$, $-10$, $-42$, $-40$, $-20$, $-52$, $-30$, $76$, $106$, $112$, $82$, $128$, $124$, $86$, $116$, $102$, $72$, $130$, $80$, $110$, $108$, $78$, $132$, $74$, $104$, $114$, $84$, $126$, $-144$, $-154$, $-166$, $-134$, $-164$, $-156$, $-142$, $-174$, $-146$, $-96$, $-90$, $-98$, $-148$, $-172$, $-140$, $-158$, $-162$, $-136$, $-168$, $-152$, $-8)$\\

\noindent$11a128$ : $3$ | $(4$, $10$, $14$, $20$, $16$, $2$, $6$, $12$, $22$, $8$, $18)$ | $(42$, $52$, $30$, $56$, $38$, $46$, $48$, $36$, $58$, $32$, $50$, $44$, $40$, $54$, $-74$, $-84$, $-120$, $-124$, $-88$, $-70$, $-78$, $-80$, $-68$, $-90$, $-66$, $-82$, $-76$, $-72$, $-86$, $-122$, $34$, $126$, $128$, $92$, $132$, $98$, $108$, $110$, $96$, $134$, $94$, $112$, $106$, $100$, $130$, $-12$, $-22$, $-2$, $-26$, $-8$, $-16$, $-60$, $-18$, $-6$, $-28$, $-4$, $-20$, $-62$, $-64$, $114$, $104$, $102$, $116$, $-118$, $-14$, $-10$, $-24)$\\

\noindent$11a129$ : $3$ | $(4$, $10$, $14$, $20$, $16$, $2$, $18$, $22$, $12$, $6$, $8)$ | $(60$, $70$, $24$, $68$, $62$, $58$, $72$, $54$, $-76$, $-84$, $-86$, $-78$, $102$, $92$, $110$, $94$, $100$, $16$, $98$, $96$, $18$, $22$, $66$, $64$, $20$, $-42$, $-34$, $-36$, $-40$, $-30$, $-28$, $-82$, $-88$, $-80$, $-26$, $-32$, $-38$, $56$, $52$, $104$, $90$, $108$, $50$, $106$, $-46$, $-8$, $-2$, $-12$, $-74$, $-14$, $-4$, $-6$, $-44$, $-48$, $-10)$\\

\noindent$11a130$ :$3$  | $(4$, $10$, $14$, $20$, $16$, $2$, $18$, $22$, $12$, $8$, $6)$ | $(20$, $44$, $56$, $42$, $54$, $46$, $-64$, $-72$, $-68$, $-60$, $-62$, $-70$, $16$, $74$, $14$, $22$, $18$, $76$, $12$, $78$, $82$, $-58$, $-38$, $-36$, $-28$, $-26$, $-34$, $-40$, $80$, $84$, $92$, $86$, $48$, $52$, $90$, $88$, $50$, $-66$, $-10$, $-4$, $-6$, $-8$, $-2$, $-30$, $-24$, $-32)$\\

\noindent$11a131$ : $3$ | $(4$, $10$, $14$, $22$, $18$, $2$, $20$, $6$, $8$, $12$, $16)$ | $(88$, $94$, $82$, $46$, $44$, $80$, $36$, $52$, $38$, $78$, $42$, $48$, $84$, $92$, $90$, $86$, $96$, $-112$, $-54$, $-116$, $-120$, $-58$, $-108$, $-60$, $-122$, $-114$, $100$, $144$, $138$, $106$, $134$, $148$, $128$, $152$, $130$, $34$, $50$, $40$, $-118$, $-56$, $-110$, $-62$, $-124$, $-68$, $-66$, $-126$, $-64$, $-70$, $-72$, $-14$, $-16$, $-74$, $-22$, $-8$, $136$, $146$, $98$, $102$, $142$, $140$, $104$, $132$, $150$, $-28$, $-2$, $-32$, $-24$, $-6$, $-10$, $-20$, $-76$, $-18$, $-12$, $-4$, $-26$, $-30)$\\

\noindent$11a132$ : $3$ | $(4$, $10$, $16$, $12$, $18$, $2$, $20$, $6$, $22$, $8$, $14)$ | $(56$, $18$, $16$, $22$, $102$, $104$, $20$, $-44$, $-76$, $-42$, $-46$, $-36$, $68$, $52$, $60$, $58$, $54$, $70$, $66$, $50$, $62$, $74$, $72$, $64$, $-10$, $-84$, $-90$, $-78$, $-80$, $-92$, $-82$, $-12$, $-8$, $-86$, $-88$, $-6$, $-14$, $-4$, $96$, $110$, $94$, $106$, $100$, $24$, $98$, $108$, $-30$, $-2$, $-34$, $-26$, $-38$, $-48$, $-40$, $-28$, $-32)$\\

\noindent$11a133$ : $3$ | $(4$, $10$, $16$, $14$, $2$, $20$, $8$, $6$, $22$, $12$, $18)$ | $(50$, $56$, $44$, $42$, $58$, $48$, $52$, $54$, $46$, $60$, $26$, $62$, $-64$, $-68$, $-6$, $90$, $78$, $96$, $76$, $92$, $28$, $-30$, $-40$, $-72$, $-34$, $-36$, $-74$, $-38$, $-32$, $-70$, $-66$, $84$, $86$, $82$, $88$, $80$, $94$, $-14$, $-2$, $-10$, $-18$, $-22$, $-24$, $-20$, $-8$, $-4$, $-16$, $-12)$\\

\noindent$11a134$ : $3$ | $(4$, $10$, $16$, $14$, $2$, $20$, $8$, $18$, $6$, $22$, $12)$ | $(46$, $40$, $106$, $96$, $34$, $100$, $102$, $36$, $94$, $108$, $42$, $44$, $48$, $38$, $104$, $98$, $-80$, $-62$, $-86$, $-74$, $-68$, $-114$, $-70$, $-72$, $-112$, $-110$, $-18$, $-6$, $-30$, $-12$, $128$, $122$, $54$, $134$, $116$, $138$, $118$, $132$, $56$, $124$, $126$, $58$, $130$, $120$, $52$, $50$, $-88$, $-64$, $-78$, $-82$, $-60$, $-84$, $-76$, $-66$, $-90$, $92$, $136$, $-22$, $-2$, $-26$, $-16$, $-8$, $-32$, $-10$, $-14$, $-28$, $-4$, $-20$, $-24)$\\

\noindent$11a135$ : $3$ | $(4$, $10$, $16$, $14$, $12$, $2$, $18$, $20$, $22$, $8$, $6)$ | $(64$, $24$, $30$, $124$, $118$, $120$, $122$, $116$, $126$, $106$, $-92$, $-88$, $-82$, $-98$, $-80$, $-90$, $22$, $66$, $62$, $26$, $28$, $60$, $68$, $20$, $70$, $58$, $74$, $-12$, $-8$, $-16$, $-96$, $-84$, $-86$, $-94$, $-18$, $-6$, $-14$, $-10$, $56$, $72$, $100$, $112$, $130$, $110$, $102$, $104$, $108$, $128$, $114$, $-40$, $-4$, $-32$, $-48$, $-46$, $-34$, $-2$, $-38$, $-42$, $-52$, $-76$, $-54$, $-78$, $-50$, $-44$, $-36)$\\

\noindent$11a136$ : $3$ | $(4$, $10$, $16$, $14$, $18$, $2$, $8$, $20$, $22$, $12$, $6)$ | $(80$, $250$, $256$, $274$, $232$, $286$, $290$, $236$, $270$, $260$, $246$, $84$, $244$, $262$, $268$, $238$, $292$, $284$, $230$, $276$, $254$, $252$, $278$, $228$, $282$, $294$, $240$, $266$, $264$, $242$, $82$, $248$, $258$, $272$, $234$, $288$, $-12$, $-34$, $-50$, $-224$, $-222$, $-200$, $136$, $174$, $120$, $166$, $144$, $74$, $150$, $160$, $126$, $180$, $130$, $156$, $154$, $78$, $140$, $170$, $-56$, $-28$, $-18$, $-66$, $-6$, $-40$, $-44$, $-2$, $-62$, $-22$, $-24$, $-60$, $-52$, $-32$, $-14$, $-70$, $-10$, $-36$, $-48$, $-226$, $-46$, $-38$, $-8$, $-68$, $-16$, $-30$, $-54$, $-58$, $-26$, $-20$, $-64$, $-4$, $-42$, $128$, $158$, $152$, $76$, $142$, $168$, $118$, $172$, $138$, $134$, $176$, $122$, $164$, $146$, $72$, $148$, $162$, $124$, $178$, $132$, $296$, $280$, $-108$, $-186$, $-98$, $-214$, $-208$, $-92$, $-192$, $-114$, $-198$, $-86$, $-202$, $-220$, $-104$, $-182$, $-102$, $-218$, $-204$, $-88$, $-196$, $-116$, $-194$, $-90$, $-206$, $-216$, $-100$, $-184$, $-106$, $-110$, $-188$, $-96$, $-212$, $-210$, $-94$, $-190$, $-112)$\\

\noindent$11a137$ : $3$ | $(4$, $10$, $16$, $14$, $18$, $2$, $20$, $6$, $22$, $12$, $8)$ | $(42$, $50$, $80$, $48$, $44$, $40$, $52$, $30$, $36$, $56$, $34$, $32$, $54$, $38$, $-94$, $-60$, $-88$, $-58$, $-92$, $-64$, $-74$, $-76$, $-66$, $-68$, $-78$, $-72$, $-62$, $-90$, $108$, $86$, $112$, $104$, $120$, $98$, $118$, $114$, $102$, $122$, $100$, $116$, $-70$, $-28$, $-16$, $-22$, $110$, $106$, $84$, $82$, $46$, $-8$, $-2$, $-12$, $-26$, $-18$, $-20$, $-24$, $-14$, $-4$, $-6$, $-96$, $-10)$\\

\noindent$11a138$ : $3$ | $(4$, $10$, $16$, $14$, $18$, $2$, $20$, $22$, $12$, $8$, $6)$ | $(38$, $84$, $90$, $128$, $124$, $94$, $122$, $130$, $88$, $86$, $132$, $120$, $92$, $126$, $-8$, $-14$, $-22$, $-110$, $-106$, $-104$, $-112$, $-20$, $-16$, $-6$, $-26$, $-10$, $-12$, $-24$, $-4$, $-18$, $44$, $32$, $54$, $50$, $28$, $48$, $56$, $34$, $42$, $40$, $36$, $-108$, $-102$, $-60$, $-100$, $-72$, $-66$, $52$, $30$, $46$, $58$, $114$, $136$, $116$, $118$, $134$, $-78$, $-2$, $-82$, $-74$, $-64$, $-96$, $-68$, $-70$, $-98$, $-62$, $-76$, $-80)$\\

\noindent$11a139$ : $3$ | $(4$, $10$, $16$, $18$, $2$, $20$, $22$, $8$, $6$, $14$, $12)$ | $(84$, $100$, $92$, $76$, $-128$, $-108$, $-118$, $-138$, $-120$, $-106$, $-126$, $-6$, $-60$, $-44$, $-64$, $-2$, $-122$, $-104$, $-124$, $-4$, $-62$, $96$, $80$, $12$, $18$, $34$, $36$, $20$, $10$, $78$, $94$, $98$, $82$, $14$, $16$, $32$, $38$, $22$, $-130$, $-110$, $-116$, $-136$, $-70$, $-50$, $-54$, $-74$, $-132$, $-112$, $-114$, $-134$, $-72$, $-52$, $88$, $142$, $148$, $26$, $42$, $28$, $150$, $140$, $86$, $102$, $90$, $144$, $146$, $24$, $40$, $30$, $152$, $-68$, $-48$, $-56$, $-8$, $-58$, $-46$, $-66)$\\

\noindent$11a140$ : $2$ | $(4$, $10$, $16$, $18$, $2$, $22$, $20$, $8$, $6$, $14$, $12)$ | $(112$, $78$, $84$, $118$, $106$, $72$, $90$, $124$, $100$, $66$, $96$, $128$, $94$, $68$, $102$, $122$, $88$, $74$, $108$, $116$, $82$, $80$, $114$, $110$, $76$, $86$, $120$, $104$, $70$, $92$, $126$, $98$, $-170$, $-132$, $-178$, $-162$, $-140$, $-186$, $-154$, $-148$, $-192$, $-146$, $-156$, $-184$, $-138$, $-164$, $-176$, $-130$, $-172$, $-168$, $-134$, $-180$, $-160$, $-142$, $-188$, $-152$, $-150$, $-190$, $-144$, $-158$, $-182$, $-136$, $-166$, $-174$, $240$, $206$, $212$, $246$, $234$, $200$, $218$, $252$, $228$, $194$, $224$, $256$, $222$, $196$, $230$, $250$, $216$, $202$, $236$, $244$, $210$, $208$, $242$, $238$, $204$, $214$, $248$, $232$, $198$, $220$, $254$, $226$, $-42$, $-4$, $-50$, $-34$, $-12$, $-58$, $-26$, $-20$, $-64$, $-18$, $-28$, $-56$, $-10$, $-36$, $-48$, $-2$, $-44$, $-40$, $-6$, $-52$, $-32$, $-14$, $-60$, $-24$, $-22$, $-62$, $-16$, $-30$, $-54$, $-8$, $-38$, $-46)$\\

\noindent$11a141$ : $3$ | $(4$, $10$, $16$, $18$, $12$, $2$, $8$, $20$, $22$, $14$, $6)$ | $(82$, $70$, $86$, $66$, $88$, $68$, $84$, $72$, $80$, $112$, $78$, $74$, $116$, $-58$, $-18$, $-20$, $-56$, $-24$, $-14$, $-4$, $-10$, $-8$, $-6$, $-12$, $-2$, $-16$, $-22$, $76$, $114$, $118$, $48$, $40$, $44$, $-90$, $-62$, $-94$, $-108$, $-102$, $-100$, $-110$, $-96$, $-60$, $-92$, $-64$, $42$, $46$, $38$, $50$, $28$, $32$, $54$, $34$, $26$, $36$, $52$, $30$, $-98$, $-104$, $-106)$\\

\noindent$11a142$ : $3$ | $(4$, $10$, $16$, $18$, $12$, $2$, $20$, $22$, $6$, $8$, $14)$ | $(114$, $84$, $116$, $94$, $104$, $-158$, $-142$, $-126$, $-130$, $-146$, $-6$, $-72$, $-56$, $-54$, $-70$, $-4$, $-2$, $-68$, $-52$, $-58$, $-74$, $-8$, $-144$, $-128$, $16$, $110$, $88$, $120$, $90$, $108$, $14$, $46$, $18$, $112$, $86$, $118$, $92$, $106$, $12$, $44$, $20$, $32$, $-150$, $-134$, $-122$, $-138$, $-154$, $-80$, $-64$, $-48$, $-62$, $-78$, $-156$, $-140$, $-124$, $-132$, $-148$, $-82$, $-152$, $-136$, $26$, $38$, $164$, $100$, $98$, $166$, $36$, $28$, $24$, $40$, $162$, $102$, $96$, $168$, $34$, $30$, $22$, $42$, $160$, $-10$, $-76$, $-60$, $-50$, $-66)$\\

\noindent$11a143$ : $3$ | $(4$, $10$, $16$, $18$, $12$, $2$, $20$, $22$, $8$, $6$, $14)$ | $(28$, $30$, $26$, $38$, $36$, $24$, $32$, $42$, $46$, $40$, $34$, $-52$, $-64$, $-82$, $-84$, $-62$, $-50$, $-54$, $-66$, $-56$, $-48$, $-60$, $-58$, $76$, $70$, $92$, $90$, $68$, $74$, $78$, $72$, $94$, $88$, $-8$, $-20$, $-14$, $-80$, $-18$, $-16$, $44$, $86$, $96$, $-2$, $-6$, $-10$, $-22$, $-12$, $-4)$\\

\noindent$11a144$ : $2$ | $(4$, $10$, $16$, $18$, $12$, $2$, $22$, $20$, $6$, $8$, $14)$ | $(102$, $128$, $76$, $136$, $94$, $110$, $120$, $84$, $144$, $86$, $118$, $112$, $92$, $138$, $78$, $126$, $104$, $100$, $130$, $74$, $134$, $96$, $108$, $122$, $82$, $142$, $88$, $116$, $114$, $90$, $140$, $80$, $124$, $106$, $98$, $132$, $-200$, $-166$, $-156$, $-190$, $-210$, $-176$, $-146$, $-180$, $-214$, $-186$, $-152$, $-170$, $-204$, $-196$, $-162$, $-160$, $-194$, $-206$, $-172$, $-150$, $-184$, $-216$, $-182$, $-148$, $-174$, $-208$, $-192$, $-158$, $-164$, $-198$, $-202$, $-168$, $-154$, $-188$, $-212$, $-178$, $246$, $272$, $220$, $280$, $238$, $254$, $264$, $228$, $288$, $230$, $262$, $256$, $236$, $282$, $222$, $270$, $248$, $244$, $274$, $218$, $278$, $240$, $252$, $266$, $226$, $286$, $232$, $260$, $258$, $234$, $284$, $224$, $268$, $250$, $242$, $276$, $-56$, $-22$, $-12$, $-46$, $-66$, $-32$, $-2$, $-36$, $-70$, $-42$, $-8$, $-26$, $-60$, $-52$, $-18$, $-16$, $-50$, $-62$, $-28$, $-6$, $-40$, $-72$, $-38$, $-4$, $-30$, $-64$, $-48$, $-14$, $-20$, $-54$, $-58$, $-24$, $-10$, $-44$, $-68$, $-34)$\\

\noindent$11a145$ : $2$ | $(4$, $10$, $16$, $18$, $12$, $2$, $22$, $20$, $8$, $6$, $14)$ | $(116$, $146$, $86$, $154$, $108$, $124$, $138$, $94$, $162$, $100$, $132$, $130$, $102$, $160$, $92$, $140$, $122$, $110$, $152$, $84$, $148$, $114$, $118$, $144$, $88$, $156$, $106$, $126$, $136$, $96$, $164$, $98$, $134$, $128$, $104$, $158$, $90$, $142$, $120$, $112$, $150$, $-186$, $-230$, $-220$, $-176$, $-196$, $-240$, $-210$, $-166$, $-206$, $-244$, $-200$, $-172$, $-216$, $-234$, $-190$, $-182$, $-226$, $-224$, $-180$, $-192$, $-236$, $-214$, $-170$, $-202$, $-246$, $-204$, $-168$, $-212$, $-238$, $-194$, $-178$, $-222$, $-228$, $-184$, $-188$, $-232$, $-218$, $-174$, $-198$, $-242$, $-208$, $280$, $310$, $250$, $318$, $272$, $288$, $302$, $258$, $326$, $264$, $296$, $294$, $266$, $324$, $256$, $304$, $286$, $274$, $316$, $248$, $312$, $278$, $282$, $308$, $252$, $320$, $270$, $290$, $300$, $260$, $328$, $262$, $298$, $292$, $268$, $322$, $254$, $306$, $284$, $276$, $314$, $-22$, $-66$, $-56$, $-12$, $-32$, $-76$, $-46$, $-2$, $-42$, $-80$, $-36$, $-8$, $-52$, $-70$, $-26$, $-18$, $-62$, $-60$, $-16$, $-28$, $-72$, $-50$, $-6$, $-38$, $-82$, $-40$, $-4$, $-48$, $-74$, $-30$, $-14$, $-58$, $-64$, $-20$, $-24$, $-68$, $-54$, $-10$, $-34$, $-78$, $-44)$\\

\noindent$11a146$ : $3$ | $(4$, $10$, $16$, $18$, $14$, $2$, $20$, $22$, $6$, $12$, $8)$ | $(74$, $64$, $104$, $98$, $70$, $78$, $68$, $100$, $102$, $66$, $76$, $72$, $96$, $106$, $-4$, $-18$, $-14$, $-8$, $-24$, $-94$, $-26$, $-6$, $-16$, $32$, $114$, $118$, $36$, $42$, $28$, $40$, $38$, $-92$, $-86$, $-48$, $-80$, $-52$, $-90$, $-88$, $-50$, $116$, $34$, $44$, $30$, $112$, $120$, $110$, $108$, $-60$, $-54$, $-82$, $-46$, $-84$, $-56$, $-58$, $-62$, $-2$, $-20$, $-12$, $-10$, $-22)$\\

\noindent$11a147$ : $3$ | $(4$, $10$, $16$, $20$, $2$, $22$, $18$, $8$, $12$, $6$, $14)$ | $(28$, $90$, $40$, $98$, $36$, $86$, $-104$, $-124$, $-122$, $-106$, $-138$, $-108$, $-120$, $-126$, $-128$, $-118$, $-110$, $-136$, $-60$, $-56$, $-132$, $-114$, $94$, $32$, $24$, $18$, $14$, $166$, $156$, $154$, $168$, $12$, $20$, $22$, $34$, $96$, $42$, $92$, $30$, $26$, $16$, $-52$, $-64$, $-62$, $-54$, $-130$, $-116$, $-112$, $-134$, $-58$, $38$, $88$, $84$, $140$, $172$, $150$, $160$, $162$, $148$, $174$, $142$, $82$, $144$, $176$, $146$, $164$, $158$, $152$, $170$, $-66$, $-50$, $-80$, $-2$, $-76$, $-46$, $-70$, $-8$, $-100$, $-6$, $-72$, $-44$, $-74$, $-4$, $-102$, $-10$, $-68$, $-48$, $-78)$\\

\noindent$11a148$ : $3$ | $(4$, $10$, $16$, $20$, $12$, $2$, $18$, $6$, $22$, $8$, $14)$ | $(50$, $58$, $24$, $26$, $56$, $48$, $52$, $90$, $92$, $-14$, $-62$, $-74$, $-84$, $-72$, $-64$, $-16$, $18$, $30$, $20$, $60$, $22$, $28$, $54$, $-68$, $-80$, $-78$, $-66$, $-70$, $-82$, $-76$, $102$, $88$, $94$, $110$, $108$, $96$, $86$, $100$, $104$, $112$, $106$, $98$, $-6$, $-42$, $-36$, $-32$, $-34$, $-44$, $-8$, $-4$, $-40$, $-38$, $-2$, $-10$, $-46$, $-12)$\\

\noindent$11a149$ : $3$ | $(4$, $10$, $16$, $20$, $12$, $2$, $18$, $8$, $22$, $14$, $6)$ | $(34$, $40$, $28$, $30$, $38$, $36$, $32$, $116$, $108$, $110$, $114$, $104$, $120$, $-56$, $-66$, $-96$, $-60$, $-62$, $-98$, $-64$, $-58$, $-24$, $-54$, $-16$, $-20$, $112$, $106$, $118$, $122$, $102$, $128$, $100$, $124$, $-6$, $-14$, $-22$, $-52$, $-18$, $48$, $70$, $42$, $26$, $44$, $72$, $50$, $74$, $46$, $68$, $126$, $-10$, $-2$, $-80$, $-86$, $-92$, $-76$, $-90$, $-88$, $-78$, $-94$, $-84$, $-82$, $-4$, $-8$, $-12)$\\

\noindent$11a150$ : $3$ | $(4$, $10$, $16$, $20$, $12$, $2$, $18$, $22$, $8$, $14$, $6)$ | $(38$, $48$, $28$, $82$, $78$, $32$, $44$, $42$, $34$, $76$, $84$, $74$, $136$, $-110$, $-100$, $-94$, $-116$, $-88$, $-106$, $-104$, $-90$, $-118$, $-92$, $-102$, $-108$, $-24$, $-22$, $-20$, $-8$, $80$, $30$, $46$, $40$, $36$, $50$, $26$, $52$, $-86$, $-114$, $-96$, $-98$, $-112$, $54$, $126$, $146$, $120$, $142$, $130$, $58$, $132$, $140$, $122$, $148$, $124$, $138$, $134$, $56$, $128$, $144$, $-14$, $-64$, $-2$, $-68$, $-18$, $-10$, $-60$, $-6$, $-72$, $-70$, $-4$, $-62$, $-12$, $-16$, $-66)$\\

\noindent$11a151$ : $3$ | $(4$, $10$, $16$, $20$, $14$, $2$, $18$, $22$, $12$, $6$, $8)$ | $(92$, $70$, $82$, $26$, $84$, $68$, $90$, $94$, $-56$, $-54$, $-48$, $-62$, $-40$, $150$, $128$, $136$, $158$, $138$, $126$, $148$, $152$, $130$, $134$, $156$, $144$, $122$, $142$, $96$, $140$, $124$, $146$, $154$, $132$, $-44$, $-30$, $-36$, $-58$, $-52$, $-50$, $-60$, $-38$, $-28$, $-42$, $-64$, $-46$, $-32$, $-34$, $-118$, $18$, $74$, $78$, $22$, $88$, $66$, $86$, $24$, $80$, $72$, $16$, $20$, $76$, $-100$, $-114$, $-4$, $-10$, $-108$, $-106$, $-12$, $-2$, $-116$, $-98$, $-120$, $-102$, $-112$, $-6$, $-8$, $-110$, $-104$, $-14)$\\

\noindent$11a152$ : $3$ | $(4$, $10$, $16$, $20$, $14$, $2$, $18$, $22$, $12$, $8$, $6)$  | $(60$, $24$, $64$, $32$, $40$, $38$, $30$, $66$, $26$, $58$, $62$, $-80$, $-86$, $-46$, $-50$, $-90$, $-56$, $-44$, $-52$, $-92$, $-54$, $96$, $104$, $114$, $116$, $102$, $94$, $98$, $-76$, $-84$, $-82$, $-78$, $-88$, $-48$, $28$, $36$, $42$, $34$, $120$, $110$, $108$, $122$, $106$, $112$, $118$, $100$, $-10$, $-72$, $-74$, $-8$, $-12$, $-70$, $-18$, $-2$, $-22$, $-4$, $-16$, $-68$, $-14$, $-6$, $-20)$\\

\noindent$11a153$ : $3$ | $(4$, $10$, $18$, $12$, $2$, $8$, $20$, $22$, $6$, $16$, $14)$ | $(130$, $112$, $96$, $118$, $124$, $102$, $106$, $170$, $168$, $108$, $100$, $122$, $120$, $98$, $110$, $132$, $128$, $114$, $94$, $116$, $126$, $104$, $-12$, $-20$, $-28$, $-4$, $-36$, $-38$, $-6$, $-26$, $-22$, $-10$, $-42$, $-14$, $-18$, $-30$, $-2$, $-34$, $-40$, $-8$, $-24$, $60$, $48$, $70$, $172$, $76$, $54$, $-156$, $-146$, $-84$, $-140$, $-162$, $-164$, $-138$, $-86$, $-148$, $-154$, $-92$, $-158$, $-144$, $-82$, $-142$, $-160$, $-90$, $-152$, $-150$, $-88$, $-16$, $68$, $46$, $62$, $80$, $58$, $50$, $72$, $174$, $74$, $52$, $56$, $78$, $64$, $44$, $66$, $134$, $-136$, $-166$, $-32)$\\

\noindent$11a154$ : $2$ | $(4$, $10$, $18$, $12$, $2$, $8$, $22$, $20$, $6$, $16$, $14)$ | $(104$, $86$, $122$, $68$, $126$, $82$, $108$, $100$, $90$, $118$, $72$, $130$, $78$, $112$, $96$, $94$, $114$, $76$, $132$, $74$, $116$, $92$, $98$, $110$, $80$, $128$, $70$, $120$, $88$, $102$, $106$, $84$, $124$, $-162$, $-176$, $-148$, $-190$, $-134$, $-194$, $-144$, $-180$, $-158$, $-166$, $-172$, $-152$, $-186$, $-138$, $-198$, $-140$, $-184$, $-154$, $-170$, $-168$, $-156$, $-182$, $-142$, $-196$, $-136$, $-188$, $-150$, $-174$, $-164$, $-160$, $-178$, $-146$, $-192$, $236$, $218$, $254$, $200$, $258$, $214$, $240$, $232$, $222$, $250$, $204$, $262$, $210$, $244$, $228$, $226$, $246$, $208$, $264$, $206$, $248$, $224$, $230$, $242$, $212$, $260$, $202$, $252$, $220$, $234$, $238$, $216$, $256$, $-30$, $-44$, $-16$, $-58$, $-2$, $-62$, $-12$, $-48$, $-26$, $-34$, $-40$, $-20$, $-54$, $-6$, $-66$, $-8$, $-52$, $-22$, $-38$, $-36$, $-24$, $-50$, $-10$, $-64$, $-4$, $-56$, $-18$, $-42$, $-32$, $-28$, $-46$, $-14$, $-60)$\\

\noindent$11a155$ : $3$ | $(4$, $10$, $18$, $12$, $2$, $16$, $20$, $8$, $22$, $14$, $6)$ | $(30$, $22$, $94$, $24$, $28$, $90$, $84$, $134$, $-116$, $-108$, $-102$, $-96$, $-60$, $-62$, $-98$, $-100$, $-110$, $-114$, $-14$, $-118$, $-106$, $-104$, $-120$, $-12$, $-112$, $34$, $18$, $40$, $44$, $144$, $124$, $146$, $46$, $38$, $16$, $36$, $48$, $32$, $20$, $42$, $-56$, $-66$, $-72$, $-50$, $-74$, $-64$, $-58$, $142$, $126$, $148$, $128$, $140$, $80$, $138$, $130$, $88$, $86$, $132$, $136$, $82$, $92$, $26$, $-8$, $-2$, $-76$, $-52$, $-70$, $-68$, $-54$, $-78$, $-4$, $-6$, $-122$, $-10)$\\

\noindent$11a156$ : $3$ | $(4$, $10$, $18$, $12$, $14$, $2$, $8$, $20$, $22$, $6$, $16)$ | $(64$, $56$, $8$, $-102$, $-46$, $-32$, $-34$, $-48$, $-36$, $-30$, $-44$, $-2$, $-40$, $60$, $68$, $52$, $74$, $50$, $70$, $58$, $62$, $66$, $54$, $72$, $-90$, $-76$, $-86$, $-100$, $-94$, $-80$, $-82$, $-96$, $-98$, $-84$, $-78$, $-92$, $-88$, $12$, $18$, $108$, $26$, $110$, $20$, $10$, $14$, $16$, $106$, $24$, $112$, $22$, $104$, $-6$, $-4$, $-38$, $-28$, $-42)$\\

\noindent$11a157$ : $3$ | $(4$, $10$, $18$, $12$, $16$, $2$, $20$, $8$, $22$, $6$, $14)$ | $(76$, $72$, $24$, $30$, $38$, $40$, $28$, $26$, $42$, $36$, $32$, $-102$, $-60$, $-110$, $-62$, $-104$, $-54$, $-56$, $-106$, $-64$, $-108$, $-58$, $-52$, $-2$, $74$, $78$, $70$, $88$, $68$, $80$, $112$, $84$, $-10$, $-18$, $-90$, $-100$, $-94$, $-96$, $-98$, $-92$, $-16$, $-12$, $-8$, $-20$, $34$, $44$, $118$, $48$, $120$, $46$, $116$, $114$, $86$, $66$, $82$, $-14$, $-6$, $-22$, $-4$, $-50)$\\

\noindent$11a158$ : $3$ | $(4$, $10$, $18$, $14$, $2$, $16$, $20$, $8$, $22$, $6$, $12)$ | $(18$, $92$, $90$, $20$, $16$, $50$, $52$, $-80$, $-38$, $-40$, $-32$, $56$, $46$, $64$, $48$, $54$, $58$, $44$, $62$, $104$, $60$, $-4$, $-72$, $-66$, $-76$, $-82$, $-78$, $-68$, $-70$, $-6$, $-2$, $-74$, $96$, $86$, $22$, $88$, $94$, $98$, $84$, $102$, $100$, $-12$, $-30$, $-24$, $-34$, $-42$, $-36$, $-26$, $-28$, $-10$, $-14$, $-8)$\\

\noindent$11a159$ : $2$ | $(4$, $10$, $18$, $14$, $2$, $22$, $8$, $20$, $6$, $16$, $12)$ | $(156$, $194$, $118$, $210$, $140$, $172$, $178$, $134$, $216$, $124$, $188$, $162$, $150$, $200$, $112$, $204$, $146$, $166$, $184$, $128$, $220$, $130$, $182$, $168$, $144$, $206$, $114$, $198$, $152$, $160$, $190$, $122$, $214$, $136$, $176$, $174$, $138$, $212$, $120$, $192$, $158$, $154$, $196$, $116$, $208$, $142$, $170$, $180$, $132$, $218$, $126$, $186$, $164$, $148$, $202$, $-290$, $-232$, $-314$, $-266$, $-256$, $-324$, $-242$, $-280$, $-300$, $-222$, $-304$, $-276$, $-246$, $-328$, $-252$, $-270$, $-310$, $-228$, $-294$, $-286$, $-236$, $-318$, $-262$, $-260$, $-320$, $-238$, $-284$, $-296$, $-226$, $-308$, $-272$, $-250$, $-330$, $-248$, $-274$, $-306$, $-224$, $-298$, $-282$, $-240$, $-322$, $-258$, $-264$, $-316$, $-234$, $-288$, $-292$, $-230$, $-312$, $-268$, $-254$, $-326$, $-244$, $-278$, $-302$, $376$, $414$, $338$, $430$, $360$, $392$, $398$, $354$, $436$, $344$, $408$, $382$, $370$, $420$, $332$, $424$, $366$, $386$, $404$, $348$, $440$, $350$, $402$, $388$, $364$, $426$, $334$, $418$, $372$, $380$, $410$, $342$, $434$, $356$, $396$, $394$, $358$, $432$, $340$, $412$, $378$, $374$, $416$, $336$, $428$, $362$, $390$, $400$, $352$, $438$, $346$, $406$, $384$, $368$, $422$, $-70$, $-12$, $-94$, $-46$, $-36$, $-104$, $-22$, $-60$, $-80$, $-2$, $-84$, $-56$, $-26$, $-108$, $-32$, $-50$, $-90$, $-8$, $-74$, $-66$, $-16$, $-98$, $-42$, $-40$, $-100$, $-18$, $-64$, $-76$, $-6$, $-88$, $-52$, $-30$, $-110$, $-28$, $-54$, $-86$, $-4$, $-78$, $-62$, $-20$, $-102$, $-38$, $-44$, $-96$, $-14$, $-68$, $-72$, $-10$, $-92$, $-48$, $-34$, $-106$, $-24$, $-58$, $-82)$\\

\noindent$11a160$ : $3$ | $(4$, $10$, $18$, $14$, $2$, $22$, $20$, $8$, $12$, $6$, $16)$ | $(132$, $28$, $196$, $26$, $134$, $-76$, $-44$, $-84$, $-62$, $-58$, $-88$, $-48$, $-72$, $-194$, $-68$, $-52$, $-92$, $-54$, $-66$, $-80$, $100$, $128$, $228$, $218$, $118$, $110$, $94$, $106$, $122$, $222$, $224$, $124$, $104$, $96$, $112$, $116$, $216$, $230$, $130$, $98$, $102$, $126$, $226$, $220$, $120$, $108$, $-174$, $-148$, $-188$, $-160$, $-162$, $-186$, $-146$, $-176$, $-172$, $-150$, $-4$, $-6$, $-152$, $-170$, $-178$, $-144$, $-184$, $-164$, $-158$, $-190$, $-192$, $-74$, $-46$, $-86$, $-60$, $18$, $204$, $36$, $30$, $198$, $24$, $136$, $12$, $210$, $40$, $208$, $14$, $138$, $22$, $200$, $32$, $34$, $202$, $20$, $140$, $16$, $206$, $38$, $212$, $232$, $214$, $114$, $-70$, $-50$, $-90$, $-56$, $-64$, $-82$, $-42$, $-78$, $-2$, $-8$, $-154$, $-168$, $-180$, $-142$, $-182$, $-166$, $-156$, $-10)$\\

\noindent$11a161$ : $3$ | $(4$, $10$, $18$, $14$, $12$, $2$, $8$, $20$, $22$, $6$, $16)$ | $(20$, $40$, $36$, $38$, $-4$, $-24$, $-28$, $-30$, $-26$, $-22$, $62$, $34$, $42$, $32$, $44$, $-46$, $-50$, $-54$, $-58$, $-56$, $-52$, $-48$, $14$, $16$, $12$, $18$, $10$, $60$, $8$, $-6$, $-2)$\\

\noindent$11a162$ : $3$ | $(4$, $10$, $18$, $14$, $16$, $2$, $6$, $20$, $22$, $12$, $8)$ | $(142$, $134$, $110$, $116$, $128$, $146$, $130$, $114$, $112$, $132$, $144$, $126$, $118$, $108$, $136$, $140$, $104$, $122$, $-162$, $-24$, $-8$, $-14$, $-18$, $-4$, $-28$, $-166$, $-32$, $-158$, $-160$, $-34$, $-164$, $-26$, $-6$, $-16$, $44$, $174$, $178$, $40$, $64$, $48$, $170$, $54$, $58$, $36$, $60$, $52$, $168$, $50$, $62$, $38$, $56$, $-84$, $-90$, $-96$, $-78$, $-152$, $-68$, $-150$, $-76$, $-98$, $-88$, $-86$, $-100$, $-74$, $-148$, $-70$, $-154$, $-80$, $-94$, $-92$, $-82$, $-156$, $-72$, $176$, $42$, $66$, $46$, $172$, $180$, $102$, $124$, $120$, $106$, $138$, $-30$, $-2$, $-20$, $-12$, $-10$, $-22)$\\

\noindent$11a163$ : $3$ | $(4$, $10$, $18$, $14$, $16$, $2$, $20$, $8$, $22$, $12$, $6)$ | $(50$, $100$, $40$, $60$, $58$, $42$, $98$, $48$, $52$, $66$, $38$, $62$, $56$, $44$, $96$, $46$, $54$, $64$, $-88$, $-72$, $-112$, $-78$, $-94$, $-82$, $-108$, $-106$, $-84$, $-92$, $-76$, $-114$, $-74$, $-90$, $-86$, $-70$, $102$, $134$, $140$, $116$, $144$, $130$, $126$, $148$, $120$, $136$, $138$, $118$, $146$, $128$, $-80$, $-110$, $-104$, $68$, $122$, $150$, $124$, $132$, $142$, $-8$, $-24$, $-34$, $-18$, $-2$, $-14$, $-30$, $-28$, $-12$, $-4$, $-20$, $-36$, $-22$, $-6$, $-10$, $-26$, $-32$, $-16)$\\

\noindent$11a164$ : $3$ | $(4$, $10$, $18$, $14$, $16$, $2$, $20$, $22$, $12$, $6$, $8)$ | $(36$, $30$, $42$, $116$, $118$, $40$, $32$, $34$, $38$, $120$, $114$, $110$, $124$, $-70$, $-64$, $-58$, $-86$, $-60$, $-62$, $-88$, $-56$, $-66$, $-68$, $-24$, $-72$, $-20$, $-8$, $78$, $130$, $84$, $134$, $82$, $128$, $108$, $126$, $122$, $112$, $-22$, $-6$, $-10$, $-18$, $-2$, $-14$, $132$, $80$, $52$, $76$, $46$, $26$, $48$, $74$, $50$, $28$, $44$, $-94$, $-100$, $-102$, $-92$, $-54$, $-90$, $-104$, $-98$, $-96$, $-106$, $-4$, $-12$, $-16)$\\

\noindent$11a165$ : $3$ | $(4$, $10$, $18$, $16$, $2$, $20$, $22$, $8$, $6$, $14$, $12)$ | $(38$, $48$, $54$, $50$, $52$, $-6$, $-30$, $86$, $96$, $92$, $82$, $80$, $90$, $98$, $88$, $84$, $94$, $-18$, $-26$, $-28$, $-16$, $-32$, $-20$, $-24$, $-72$, $-76$, $-74$, $-22$, $46$, $36$, $40$, $14$, $42$, $34$, $44$, $12$, $78$, $10$, $-8$, $-4$, $-56$, $-70$, $-64$, $-62$, $-2$, $-58$, $-68$, $-66$, $-60)$\\

\noindent$11a166$ : $2$ | $(4$, $10$, $18$, $16$, $2$, $22$, $20$, $8$, $6$, $14$, $12)$ | $(72$, $100$, $106$, $78$, $66$, $94$, $112$, $84$, $60$, $88$, $116$, $90$, $62$, $82$, $110$, $96$, $68$, $76$, $104$, $102$, $74$, $70$, $98$, $108$, $80$, $64$, $92$, $114$, $86$, $-154$, $-120$, $-162$, $-146$, $-128$, $-170$, $-138$, $-136$, $-172$, $-130$, $-144$, $-164$, $-122$, $-152$, $-156$, $-118$, $-160$, $-148$, $-126$, $-168$, $-140$, $-134$, $-174$, $-132$, $-142$, $-166$, $-124$, $-150$, $-158$, $188$, $216$, $222$, $194$, $182$, $210$, $228$, $200$, $176$, $204$, $232$, $206$, $178$, $198$, $226$, $212$, $184$, $192$, $220$, $218$, $190$, $186$, $214$, $224$, $196$, $180$, $208$, $230$, $202$, $-38$, $-4$, $-46$, $-30$, $-12$, $-54$, $-22$, $-20$, $-56$, $-14$, $-28$, $-48$, $-6$, $-36$, $-40$, $-2$, $-44$, $-32$, $-10$, $-52$, $-24$, $-18$, $-58$, $-16$, $-26$, $-50$, $-8$, $-34$, $-42)$\\

\noindent$11a167$ : $3$ | $(4$, $10$, $18$, $16$, $12$, $2$, $8$, $20$, $22$, $14$, $6)$ | $(24$, $34$, $28$, $30$, $32$, $26$, $84$, $80$, $74$, $90$, $-52$, $-60$, $-104$, $-64$, $-100$, $-66$, $-102$, $-62$, $-50$, $-18$, $-54$, $-58$, $-14$, $82$, $72$, $92$, $88$, $76$, $78$, $86$, $94$, $70$, $98$, $-8$, $-4$, $-122$, $-116$, $-110$, $-128$, $-108$, $-118$, $-120$, $-106$, $-126$, $-112$, $-114$, $-124$, $-2$, $-10$, $-6$, $68$, $96$, $136$, $44$, $130$, $40$, $20$, $38$, $132$, $46$, $134$, $36$, $22$, $42$, $-56$, $-16$, $-48$, $-12)$\\

\noindent$11a168$ : $3$ | $(4$, $10$, $18$, $16$, $14$, $2$, $20$, $8$, $22$, $12$, $6)$ | $(64$, $128$, $50$, $78$, $72$, $56$, $122$, $58$, $70$, $80$, $48$, $130$, $66$, $62$, $126$, $52$, $76$, $74$, $54$, $124$, $60$, $68$, $82$, $-88$, $-114$, $-104$, $-142$, $-98$, $-120$, $-94$, $-138$, $-108$, $-110$, $-136$, $-92$, $-118$, $-100$, $-144$, $-102$, $-116$, $-90$, $-86$, $132$, $172$, $174$, $150$, $186$, $160$, $164$, $182$, $146$, $178$, $168$, $156$, $190$, $154$, $170$, $176$, $148$, $184$, $162$, $-96$, $-140$, $-106$, $-112$, $-134$, $84$, $152$, $188$, $158$, $166$, $180$, $-34$, $-8$, $-18$, $-44$, $-24$, $-2$, $-28$, $-40$, $-14$, $-12$, $-38$, $-30$, $-4$, $-22$, $-46$, $-20$, $-6$, $-32$, $-36$, $-10$, $-16$, $-42$, $-26)$\\

\noindent$11a169$ : $3$ | $(4$, $10$, $18$, $16$, $14$, $2$, $20$, $22$, $6$, $12$, $8)$ | $(54$, $58$, $86$, $80$, $62$, $82$, $84$, $60$, $78$, $88$, $-4$, $-14$, $-10$, $-8$, $-16$, $-76$, $-18$, $-6$, $-12$, $24$, $28$, $34$, $20$, $32$, $30$, $56$, $-72$, $-70$, $-74$, $-68$, $-40$, $26$, $36$, $22$, $92$, $96$, $94$, $90$, $-48$, $-42$, $-64$, $-38$, $-66$, $-44$, $-46$, $-50$, $-2$, $-52)$\\

\noindent$11a170$ : $3$ | $(4$, $10$, $18$, $20$, $2$, $16$, $6$, $8$, $22$, $14$, $12)$ | $(64$, $48$, $80$, $126$, $74$, $54$, $58$, $70$, $42$, $44$, $68$, $60$, $52$, $76$, $128$, $78$, $50$, $62$, $66$, $46$, $-140$, $-118$, $-138$, $-98$, $-104$, $-132$, $-112$, $-90$, $-88$, $-114$, $-134$, $-102$, $-100$, $-136$, $-116$, $-86$, $-92$, $-110$, $-130$, $-106$, $-96$, $176$, $164$, $188$, $150$, $186$, $166$, $174$, $178$, $158$, $194$, $156$, $180$, $172$, $168$, $184$, $152$, $190$, $162$, $-18$, $-40$, $-84$, $-94$, $-108$, $170$, $182$, $154$, $192$, $160$, $120$, $122$, $82$, $124$, $72$, $56$, $-12$, $-24$, $-34$, $-2$, $-38$, $-20$, $-16$, $-146$, $-8$, $-28$, $-30$, $-6$, $-144$, $-142$, $-4$, $-32$, $-26$, $-10$, $-148$, $-14$, $-22$, $-36)$\\

\noindent$11a171$ : $3$ | $(4$, $10$, $18$, $20$, $2$, $22$, $6$, $8$, $12$, $14$, $16)$ | $(54$, $270$, $74$, $34$, $36$, $72$, $268$, $56$, $52$, $-116$, $-136$, $-86$, $-132$, $-120$, $-98$, $-148$, $-104$, $-192$, $-106$, $-146$, $-96$, $-122$, $-130$, $-88$, $-138$, $-114$, $-200$, $-198$, $-112$, $-140$, $-90$, $-128$, $-124$, $-94$, $-144$, $-108$, $-194$, $-102$, $-150$, $-100$, $-118$, $-134$, $276$, $180$, $300$, $166$, $160$, $294$, $186$, $282$, $152$, $286$, $190$, $290$, $156$, $170$, $304$, $176$, $272$, $280$, $184$, $296$, $162$, $164$, $298$, $182$, $278$, $274$, $178$, $302$, $168$, $158$, $292$, $188$, $284$, $-236$, $-4$, $-12$, $-244$, $-228$, $-210$, $-260$, $-212$, $-226$, $-246$, $-14$, $-2$, $-18$, $-250$, $-222$, $-216$, $-256$, $-206$, $-232$, $-240$, $-8$, $288$, $154$, $172$, $306$, $174$, $80$, $28$, $42$, $66$, $262$, $62$, $46$, $24$, $84$, $22$, $48$, $60$, $264$, $68$, $40$, $30$, $78$, $76$, $32$, $38$, $70$, $266$, $58$, $50$, $20$, $82$, $26$, $44$, $64$, $-126$, $-92$, $-142$, $-110$, $-196$, $-202$, $-252$, $-220$, $-218$, $-254$, $-204$, $-234$, $-238$, $-6$, $-10$, $-242$, $-230$, $-208$, $-258$, $-214$, $-224$, $-248$, $-16)$\\

\noindent$11a172$ : $3$ | $(4$, $10$, $18$, $22$, $14$, $2$, $20$, $8$, $6$, $12$, $16)$ | $(64$, $88$, $70$, $60$, $68$, $86$, $66$, $62$, $90$, $96$, $-78$, $-74$, $-16$, $-2$, $-12$, $-10$, $-4$, $-18$, $-76$, $22$, $72$, $24$, $36$, $20$, $34$, $26$, $108$, $28$, $32$, $104$, $-82$, $-52$, $-54$, $-84$, $-58$, $-80$, $100$, $92$, $94$, $98$, $102$, $106$, $30$, $-56$, $-50$, $-44$, $-38$, $-46$, $-48$, $-40$, $-42$, $-6$, $-8$, $-14)$\\

\noindent$11a173$ : $3$ | $(4$, $10$, $20$, $14$, $2$, $8$, $18$, $22$, $6$, $12$, $16)$ | $(180$, $202$, $158$, $214$, $92$, $274$, $70$, $282$, $84$, $260$, $262$, $82$, $284$, $72$, $272$, $94$, $216$, $98$, $218$, $96$, $270$, $74$, $286$, $80$, $264$, $258$, $86$, $280$, $68$, $276$, $90$, $88$, $278$, $-238$, $-120$, $-146$, $-106$, $-224$, $-2$, $-228$, $-110$, $-150$, $-116$, $-234$, $-242$, $-124$, $-142$, $-102$, $-220$, $252$, $248$, $162$, $198$, $184$, $176$, $206$, $154$, $210$, $172$, $188$, $194$, $166$, $244$, $256$, $266$, $78$, $288$, $76$, $268$, $254$, $246$, $164$, $196$, $186$, $174$, $208$, $-20$, $-50$, $-128$, $-138$, $-40$, $-30$, $-10$, $-60$, $-64$, $-54$, $-16$, $-24$, $-46$, $-132$, $-134$, $-44$, $-26$, $-14$, $-56$, $-66$, $-58$, $-12$, $-28$, $-42$, $-136$, $-130$, $-48$, $-22$, $-18$, $-52$, $-126$, $-140$, $-100$, $250$, $160$, $200$, $182$, $178$, $204$, $156$, $212$, $170$, $190$, $192$, $168$, $-62$, $-8$, $-32$, $-38$, $-36$, $-34$, $-6$, $-232$, $-114$, $-152$, $-112$, $-230$, $-4$, $-222$, $-104$, $-144$, $-122$, $-240$, $-236$, $-118$, $-148$, $-108$, $-226)$\\

\noindent$11a174$ : $2$ | $(4$, $12$, $14$, $16$, $18$, $2$, $22$, $6$, $20$, $8$, $10)$ | $(126$, $92$, $154$, $98$, $120$, $132$, $86$, $148$, $104$, $114$, $138$, $80$, $142$, $110$, $108$, $144$, $82$, $136$, $116$, $102$, $150$, $88$, $130$, $122$, $96$, $156$, $94$, $124$, $128$, $90$, $152$, $100$, $118$, $134$, $84$, $146$, $106$, $112$, $140$, $-184$, $-230$, $-174$, $-194$, $-220$, $-164$, $-204$, $-210$, $-158$, $-214$, $-200$, $-168$, $-224$, $-190$, $-178$, $-234$, $-180$, $-188$, $-226$, $-170$, $-198$, $-216$, $-160$, $-208$, $-206$, $-162$, $-218$, $-196$, $-172$, $-228$, $-186$, $-182$, $-232$, $-176$, $-192$, $-222$, $-166$, $-202$, $-212$, $282$, $248$, $310$, $254$, $276$, $288$, $242$, $304$, $260$, $270$, $294$, $236$, $298$, $266$, $264$, $300$, $238$, $292$, $272$, $258$, $306$, $244$, $286$, $278$, $252$, $312$, $250$, $280$, $284$, $246$, $308$, $256$, $274$, $290$, $240$, $302$, $262$, $268$, $296$, $-28$, $-74$, $-18$, $-38$, $-64$, $-8$, $-48$, $-54$, $-2$, $-58$, $-44$, $-12$, $-68$, $-34$, $-22$, $-78$, $-24$, $-32$, $-70$, $-14$, $-42$, $-60$, $-4$, $-52$, $-50$, $-6$, $-62$, $-40$, $-16$, $-72$, $-30$, $-26$, $-76$, $-20$, $-36$, $-66$, $-10$, $-46$, $-56)$\\

\noindent$11a175$ : $2$ | $(4$, $12$, $14$, $16$, $18$, $2$, $22$, $6$, $20$, $10$, $8)$ | $(168$, $122$, $204$, $132$, $158$, $178$, $112$, $194$, $142$, $148$, $188$, $106$, $184$, $152$, $138$, $198$, $116$, $174$, $162$, $128$, $208$, $126$, $164$, $172$, $118$, $200$, $136$, $154$, $182$, $108$, $190$, $146$, $144$, $192$, $110$, $180$, $156$, $134$, $202$, $120$, $170$, $166$, $124$, $206$, $130$, $160$, $176$, $114$, $196$, $140$, $150$, $186$, $-272$, $-226$, $-308$, $-236$, $-262$, $-282$, $-216$, $-298$, $-246$, $-252$, $-292$, $-210$, $-288$, $-256$, $-242$, $-302$, $-220$, $-278$, $-266$, $-232$, $-312$, $-230$, $-268$, $-276$, $-222$, $-304$, $-240$, $-258$, $-286$, $-212$, $-294$, $-250$, $-248$, $-296$, $-214$, $-284$, $-260$, $-238$, $-306$, $-224$, $-274$, $-270$, $-228$, $-310$, $-234$, $-264$, $-280$, $-218$, $-300$, $-244$, $-254$, $-290$, $376$, $330$, $412$, $340$, $366$, $386$, $320$, $402$, $350$, $356$, $396$, $314$, $392$, $360$, $346$, $406$, $324$, $382$, $370$, $336$, $416$, $334$, $372$, $380$, $326$, $408$, $344$, $362$, $390$, $316$, $398$, $354$, $352$, $400$, $318$, $388$, $364$, $342$, $410$, $328$, $378$, $374$, $332$, $414$, $338$, $368$, $384$, $322$, $404$, $348$, $358$, $394$, $-64$, $-18$, $-100$, $-28$, $-54$, $-74$, $-8$, $-90$, $-38$, $-44$, $-84$, $-2$, $-80$, $-48$, $-34$, $-94$, $-12$, $-70$, $-58$, $-24$, $-104$, $-22$, $-60$, $-68$, $-14$, $-96$, $-32$, $-50$, $-78$, $-4$, $-86$, $-42$, $-40$, $-88$, $-6$, $-76$, $-52$, $-30$, $-98$, $-16$, $-66$, $-62$, $-20$, $-102$, $-26$, $-56$, $-72$, $-10$, $-92$, $-36$, $-46$, $-82)$\\

\noindent$11a176$ : $2$ | $(4$, $12$, $14$, $16$, $18$, $2$, $22$, $20$, $8$, $6$, $10)$  | $(178$, $128$, $214$, $142$, $164$, $192$, $114$, $200$, $156$, $150$, $206$, $120$, $186$, $170$, $136$, $220$, $134$, $172$, $184$, $122$, $208$, $148$, $158$, $198$, $112$, $194$, $162$, $144$, $212$, $126$, $180$, $176$, $130$, $216$, $140$, $166$, $190$, $116$, $202$, $154$, $152$, $204$, $118$, $188$, $168$, $138$, $218$, $132$, $174$, $182$, $124$, $210$, $146$, $160$, $196$, $-300$, $-238$, $-264$, $-326$, $-274$, $-228$, $-290$, $-310$, $-248$, $-254$, $-316$, $-284$, $-222$, $-280$, $-320$, $-258$, $-244$, $-306$, $-294$, $-232$, $-270$, $-330$, $-268$, $-234$, $-296$, $-304$, $-242$, $-260$, $-322$, $-278$, $-224$, $-286$, $-314$, $-252$, $-250$, $-312$, $-288$, $-226$, $-276$, $-324$, $-262$, $-240$, $-302$, $-298$, $-236$, $-266$, $-328$, $-272$, $-230$, $-292$, $-308$, $-246$, $-256$, $-318$, $-282$, $398$, $348$, $434$, $362$, $384$, $412$, $334$, $420$, $376$, $370$, $426$, $340$, $406$, $390$, $356$, $440$, $354$, $392$, $404$, $342$, $428$, $368$, $378$, $418$, $332$, $414$, $382$, $364$, $432$, $346$, $400$, $396$, $350$, $436$, $360$, $386$, $410$, $336$, $422$, $374$, $372$, $424$, $338$, $408$, $388$, $358$, $438$, $352$, $394$, $402$, $344$, $430$, $366$, $380$, $416$, $-80$, $-18$, $-44$, $-106$, $-54$, $-8$, $-70$, $-90$, $-28$, $-34$, $-96$, $-64$, $-2$, $-60$, $-100$, $-38$, $-24$, $-86$, $-74$, $-12$, $-50$, $-110$, $-48$, $-14$, $-76$, $-84$, $-22$, $-40$, $-102$, $-58$, $-4$, $-66$, $-94$, $-32$, $-30$, $-92$, $-68$, $-6$, $-56$, $-104$, $-42$, $-20$, $-82$, $-78$, $-16$, $-46$, $-108$, $-52$, $-10$, $-72$, $-88$, $-26$, $-36$, $-98$, $-62)$\\

\noindent$11a177$ : $2$ | $(4$, $12$, $14$, $16$, $18$, $2$, $22$, $20$, $8$, $10$, $6)$ | $(156$, $110$, $184$, $128$, $138$, $174$, $100$, $166$, $146$, $120$, $192$, $118$, $148$, $164$, $102$, $176$, $136$, $130$, $182$, $108$, $158$, $154$, $112$, $186$, $126$, $140$, $172$, $98$, $168$, $144$, $122$, $190$, $116$, $150$, $162$, $104$, $178$, $134$, $132$, $180$, $106$, $160$, $152$, $114$, $188$, $124$, $142$, $170$, $-268$, $-226$, $-200$, $-242$, $-284$, $-252$, $-210$, $-216$, $-258$, $-278$, $-236$, $-194$, $-232$, $-274$, $-262$, $-220$, $-206$, $-248$, $-288$, $-246$, $-204$, $-222$, $-264$, $-272$, $-230$, $-196$, $-238$, $-280$, $-256$, $-214$, $-212$, $-254$, $-282$, $-240$, $-198$, $-228$, $-270$, $-266$, $-224$, $-202$, $-244$, $-286$, $-250$, $-208$, $-218$, $-260$, $-276$, $-234$, $348$, $302$, $376$, $320$, $330$, $366$, $292$, $358$, $338$, $312$, $384$, $310$, $340$, $356$, $294$, $368$, $328$, $322$, $374$, $300$, $350$, $346$, $304$, $378$, $318$, $332$, $364$, $290$, $360$, $336$, $314$, $382$, $308$, $342$, $354$, $296$, $370$, $326$, $324$, $372$, $298$, $352$, $344$, $306$, $380$, $316$, $334$, $362$, $-76$, $-34$, $-8$, $-50$, $-92$, $-60$, $-18$, $-24$, $-66$, $-86$, $-44$, $-2$, $-40$, $-82$, $-70$, $-28$, $-14$, $-56$, $-96$, $-54$, $-12$, $-30$, $-72$, $-80$, $-38$, $-4$, $-46$, $-88$, $-64$, $-22$, $-20$, $-62$, $-90$, $-48$, $-6$, $-36$, $-78$, $-74$, $-32$, $-10$, $-52$, $-94$, $-58$, $-16$, $-26$, $-68$, $-84$, $-42)$\\

\noindent$11a178$ : $2$ | $(4$, $12$, $14$, $16$, $18$, $2$, $22$, $20$, $10$, $8$, $6)$ | $(198$, $140$, $234$, $162$, $176$, $220$, $126$, $212$, $184$, $154$, $242$, $148$, $190$, $206$, $132$, $226$, $170$, $168$, $228$, $134$, $204$, $192$, $146$, $240$, $156$, $182$, $214$, $124$, $218$, $178$, $160$, $236$, $142$, $196$, $200$, $138$, $232$, $164$, $174$, $222$, $128$, $210$, $186$, $152$, $244$, $150$, $188$, $208$, $130$, $224$, $172$, $166$, $230$, $136$, $202$, $194$, $144$, $238$, $158$, $180$, $216$, $-278$, $-346$, $-320$, $-252$, $-304$, $-362$, $-294$, $-262$, $-330$, $-336$, $-268$, $-288$, $-356$, $-310$, $-246$, $-314$, $-352$, $-284$, $-272$, $-340$, $-326$, $-258$, $-298$, $-366$, $-300$, $-256$, $-324$, $-342$, $-274$, $-282$, $-350$, $-316$, $-248$, $-308$, $-358$, $-290$, $-266$, $-334$, $-332$, $-264$, $-292$, $-360$, $-306$, $-250$, $-318$, $-348$, $-280$, $-276$, $-344$, $-322$, $-254$, $-302$, $-364$, $-296$, $-260$, $-328$, $-338$, $-270$, $-286$, $-354$, $-312$, $442$, $384$, $478$, $406$, $420$, $464$, $370$, $456$, $428$, $398$, $486$, $392$, $434$, $450$, $376$, $470$, $414$, $412$, $472$, $378$, $448$, $436$, $390$, $484$, $400$, $426$, $458$, $368$, $462$, $422$, $404$, $480$, $386$, $440$, $444$, $382$, $476$, $408$, $418$, $466$, $372$, $454$, $430$, $396$, $488$, $394$, $432$, $452$, $374$, $468$, $416$, $410$, $474$, $380$, $446$, $438$, $388$, $482$, $402$, $424$, $460$, $-34$, $-102$, $-76$, $-8$, $-60$, $-118$, $-50$, $-18$, $-86$, $-92$, $-24$, $-44$, $-112$, $-66$, $-2$, $-70$, $-108$, $-40$, $-28$, $-96$, $-82$, $-14$, $-54$, $-122$, $-56$, $-12$, $-80$, $-98$, $-30$, $-38$, $-106$, $-72$, $-4$, $-64$, $-114$, $-46$, $-22$, $-90$, $-88$, $-20$, $-48$, $-116$, $-62$, $-6$, $-74$, $-104$, $-36$, $-32$, $-100$, $-78$, $-10$, $-58$, $-120$, $-52$, $-16$, $-84$, $-94$, $-26$, $-42$, $-110$, $-68)$\\

\noindent$11a179$ : $2$ | $(4$, $12$, $14$, $16$, $18$, $20$, $2$, $6$, $8$, $22$, $10)$ | $(76$, $110$, $70$, $82$, $104$, $64$, $88$, $98$, $58$, $94$, $92$, $60$, $100$, $86$, $66$, $106$, $80$, $72$, $112$, $74$, $78$, $108$, $68$, $84$, $102$, $62$, $90$, $96$, $-132$, $-166$, $-126$, $-138$, $-160$, $-120$, $-144$, $-154$, $-114$, $-150$, $-148$, $-116$, $-156$, $-142$, $-122$, $-162$, $-136$, $-128$, $-168$, $-130$, $-134$, $-164$, $-124$, $-140$, $-158$, $-118$, $-146$, $-152$, $188$, $222$, $182$, $194$, $216$, $176$, $200$, $210$, $170$, $206$, $204$, $172$, $212$, $198$, $178$, $218$, $192$, $184$, $224$, $186$, $190$, $220$, $180$, $196$, $214$, $174$, $202$, $208$, $-20$, $-54$, $-14$, $-26$, $-48$, $-8$, $-32$, $-42$, $-2$, $-38$, $-36$, $-4$, $-44$, $-30$, $-10$, $-50$, $-24$, $-16$, $-56$, $-18$, $-22$, $-52$, $-12$, $-28$, $-46$, $-6$, $-34$, $-40)$\\

\noindent$11a180$ : $2$ | $(4$, $12$, $14$, $16$, $18$, $20$, $2$, $6$, $22$, $10$, $8)$ | $(120$, $170$, $106$, $134$, $156$, $92$, $148$, $142$, $98$, $162$, $128$, $112$, $176$, $114$, $126$, $164$, $100$, $140$, $150$, $90$, $154$, $136$, $104$, $168$, $122$, $118$, $172$, $108$, $132$, $158$, $94$, $146$, $144$, $96$, $160$, $130$, $110$, $174$, $116$, $124$, $166$, $102$, $138$, $152$, $-240$, $-190$, $-212$, $-262$, $-218$, $-184$, $-234$, $-246$, $-196$, $-206$, $-256$, $-224$, $-178$, $-228$, $-252$, $-202$, $-200$, $-250$, $-230$, $-180$, $-222$, $-258$, $-208$, $-194$, $-244$, $-236$, $-186$, $-216$, $-264$, $-214$, $-188$, $-238$, $-242$, $-192$, $-210$, $-260$, $-220$, $-182$, $-232$, $-248$, $-198$, $-204$, $-254$, $-226$, $296$, $346$, $282$, $310$, $332$, $268$, $324$, $318$, $274$, $338$, $304$, $288$, $352$, $290$, $302$, $340$, $276$, $316$, $326$, $266$, $330$, $312$, $280$, $344$, $298$, $294$, $348$, $284$, $308$, $334$, $270$, $322$, $320$, $272$, $336$, $306$, $286$, $350$, $292$, $300$, $342$, $278$, $314$, $328$, $-64$, $-14$, $-36$, $-86$, $-42$, $-8$, $-58$, $-70$, $-20$, $-30$, $-80$, $-48$, $-2$, $-52$, $-76$, $-26$, $-24$, $-74$, $-54$, $-4$, $-46$, $-82$, $-32$, $-18$, $-68$, $-60$, $-10$, $-40$, $-88$, $-38$, $-12$, $-62$, $-66$, $-16$, $-34$, $-84$, $-44$, $-6$, $-56$, $-72$, $-22$, $-28$, $-78$, $-50)$\\

\noindent$11a181$ : $3$ | $(4$, $12$, $14$, $16$, $18$, $20$, $2$, $8$, $6$, $22$, $10)$ | $(40$, $56$, $46$, $90$, $94$, $50$, $52$, $36$, $44$, $58$, $42$, $38$, $54$, $48$, $92$, $-8$, $-96$, $-98$, $-62$, $-76$, $-78$, $-60$, $-82$, $-72$, $-66$, $-100$, $-64$, $-74$, $-80$, $108$, $124$, $114$, $102$, $118$, $30$, $120$, $104$, $112$, $126$, $110$, $106$, $122$, $32$, $34$, $-84$, $-70$, $-68$, $-86$, $88$, $116$, $-12$, $-26$, $-4$, $-18$, $-20$, $-2$, $-24$, $-14$, $-10$, $-28$, $-6$, $-16$, $-22)$\\

\noindent$11a182$ : $2$ | $(4$, $12$, $14$, $16$, $18$, $20$, $2$, $22$, $8$, $10$, $6)$ | $(100$, $134$, $78$, $122$, $112$, $88$, $144$, $90$, $110$, $124$, $76$, $132$, $102$, $98$, $136$, $80$, $120$, $114$, $86$, $142$, $92$, $108$, $126$, $74$, $130$, $104$, $96$, $138$, $82$, $118$, $116$, $84$, $140$, $94$, $106$, $128$, $-204$, $-178$, $-152$, $-162$, $-188$, $-214$, $-194$, $-168$, $-146$, $-172$, $-198$, $-210$, $-184$, $-158$, $-156$, $-182$, $-208$, $-200$, $-174$, $-148$, $-166$, $-192$, $-216$, $-190$, $-164$, $-150$, $-176$, $-202$, $-206$, $-180$, $-154$, $-160$, $-186$, $-212$, $-196$, $-170$, $244$, $278$, $222$, $266$, $256$, $232$, $288$, $234$, $254$, $268$, $220$, $276$, $246$, $242$, $280$, $224$, $264$, $258$, $230$, $286$, $236$, $252$, $270$, $218$, $274$, $248$, $240$, $282$, $226$, $262$, $260$, $228$, $284$, $238$, $250$, $272$, $-60$, $-34$, $-8$, $-18$, $-44$, $-70$, $-50$, $-24$, $-2$, $-28$, $-54$, $-66$, $-40$, $-14$, $-12$, $-38$, $-64$, $-56$, $-30$, $-4$, $-22$, $-48$, $-72$, $-46$, $-20$, $-6$, $-32$, $-58$, $-62$, $-36$, $-10$, $-16$, $-42$, $-68$, $-52$, $-26)$\\

\noindent$11a183$ : $2$ | $(4$, $12$, $14$, $16$, $18$, $20$, $2$, $22$, $10$, $8$, $6)$ | $(158$, $212$, $124$, $192$, $178$, $138$, $226$, $144$, $172$, $198$, $118$, $206$, $164$, $152$, $218$, $130$, $186$, $184$, $132$, $220$, $150$, $166$, $204$, $116$, $200$, $170$, $146$, $224$, $136$, $180$, $190$, $126$, $214$, $156$, $160$, $210$, $122$, $194$, $176$, $140$, $228$, $142$, $174$, $196$, $120$, $208$, $162$, $154$, $216$, $128$, $188$, $182$, $134$, $222$, $148$, $168$, $202$, $-262$, $-330$, $-288$, $-236$, $-304$, $-314$, $-246$, $-278$, $-340$, $-272$, $-252$, $-320$, $-298$, $-230$, $-294$, $-324$, $-256$, $-268$, $-336$, $-282$, $-242$, $-310$, $-308$, $-240$, $-284$, $-334$, $-266$, $-258$, $-326$, $-292$, $-232$, $-300$, $-318$, $-250$, $-274$, $-342$, $-276$, $-248$, $-316$, $-302$, $-234$, $-290$, $-328$, $-260$, $-264$, $-332$, $-286$, $-238$, $-306$, $-312$, $-244$, $-280$, $-338$, $-270$, $-254$, $-322$, $-296$, $386$, $440$, $352$, $420$, $406$, $366$, $454$, $372$, $400$, $426$, $346$, $434$, $392$, $380$, $446$, $358$, $414$, $412$, $360$, $448$, $378$, $394$, $432$, $344$, $428$, $398$, $374$, $452$, $364$, $408$, $418$, $354$, $442$, $384$, $388$, $438$, $350$, $422$, $404$, $368$, $456$, $370$, $402$, $424$, $348$, $436$, $390$, $382$, $444$, $356$, $416$, $410$, $362$, $450$, $376$, $396$, $430$, $-34$, $-102$, $-60$, $-8$, $-76$, $-86$, $-18$, $-50$, $-112$, $-44$, $-24$, $-92$, $-70$, $-2$, $-66$, $-96$, $-28$, $-40$, $-108$, $-54$, $-14$, $-82$, $-80$, $-12$, $-56$, $-106$, $-38$, $-30$, $-98$, $-64$, $-4$, $-72$, $-90$, $-22$, $-46$, $-114$, $-48$, $-20$, $-88$, $-74$, $-6$, $-62$, $-100$, $-32$, $-36$, $-104$, $-58$, $-10$, $-78$, $-84$, $-16$, $-52$, $-110$, $-42$, $-26$, $-94$, $-68)$\\

\noindent$11a184$ : $2$ | $(4$, $12$, $14$, $16$, $18$, $22$, $2$, $20$, $8$, $10$, $6)$ | $(118$, $164$, $100$, $136$, $146$, $90$, $154$, $128$, $108$, $172$, $110$, $126$, $156$, $92$, $144$, $138$, $98$, $162$, $120$, $116$, $166$, $102$, $134$, $148$, $88$, $152$, $130$, $106$, $170$, $112$, $124$, $158$, $94$, $142$, $140$, $96$, $160$, $122$, $114$, $168$, $104$, $132$, $150$, $-240$, $-202$, $-180$, $-218$, $-256$, $-224$, $-186$, $-196$, $-234$, $-246$, $-208$, $-174$, $-212$, $-250$, $-230$, $-192$, $-190$, $-228$, $-252$, $-214$, $-176$, $-206$, $-244$, $-236$, $-198$, $-184$, $-222$, $-258$, $-220$, $-182$, $-200$, $-238$, $-242$, $-204$, $-178$, $-216$, $-254$, $-226$, $-188$, $-194$, $-232$, $-248$, $-210$, $290$, $336$, $272$, $308$, $318$, $262$, $326$, $300$, $280$, $344$, $282$, $298$, $328$, $264$, $316$, $310$, $270$, $334$, $292$, $288$, $338$, $274$, $306$, $320$, $260$, $324$, $302$, $278$, $342$, $284$, $296$, $330$, $266$, $314$, $312$, $268$, $332$, $294$, $286$, $340$, $276$, $304$, $322$, $-68$, $-30$, $-8$, $-46$, $-84$, $-52$, $-14$, $-24$, $-62$, $-74$, $-36$, $-2$, $-40$, $-78$, $-58$, $-20$, $-18$, $-56$, $-80$, $-42$, $-4$, $-34$, $-72$, $-64$, $-26$, $-12$, $-50$, $-86$, $-48$, $-10$, $-28$, $-66$, $-70$, $-32$, $-6$, $-44$, $-82$, $-54$, $-16$, $-22$, $-60$, $-76$, $-38)$\\

\noindent$11a185$ : $2$ | $(4$, $12$, $14$, $16$, $18$, $22$, $2$, $20$, $10$, $8$, $6)$ | $(148$, $206$, $126$, $170$, $184$, $112$, $192$, $162$, $134$, $214$, $140$, $156$, $198$, $118$, $178$, $176$, $120$, $200$, $154$, $142$, $212$, $132$, $164$, $190$, $110$, $186$, $168$, $128$, $208$, $146$, $150$, $204$, $124$, $172$, $182$, $114$, $194$, $160$, $136$, $216$, $138$, $158$, $196$, $116$, $180$, $174$, $122$, $202$, $152$, $144$, $210$, $130$, $166$, $188$, $-246$, $-306$, $-284$, $-224$, $-268$, $-322$, $-262$, $-230$, $-290$, $-300$, $-240$, $-252$, $-312$, $-278$, $-218$, $-274$, $-316$, $-256$, $-236$, $-296$, $-294$, $-234$, $-258$, $-318$, $-272$, $-220$, $-280$, $-310$, $-250$, $-242$, $-302$, $-288$, $-228$, $-264$, $-324$, $-266$, $-226$, $-286$, $-304$, $-244$, $-248$, $-308$, $-282$, $-222$, $-270$, $-320$, $-260$, $-232$, $-292$, $-298$, $-238$, $-254$, $-314$, $-276$, $364$, $422$, $342$, $386$, $400$, $328$, $408$, $378$, $350$, $430$, $356$, $372$, $414$, $334$, $394$, $392$, $336$, $416$, $370$, $358$, $428$, $348$, $380$, $406$, $326$, $402$, $384$, $344$, $424$, $362$, $366$, $420$, $340$, $388$, $398$, $330$, $410$, $376$, $352$, $432$, $354$, $374$, $412$, $332$, $396$, $390$, $338$, $418$, $368$, $360$, $426$, $346$, $382$, $404$, $-30$, $-90$, $-68$, $-8$, $-52$, $-106$, $-46$, $-14$, $-74$, $-84$, $-24$, $-36$, $-96$, $-62$, $-2$, $-58$, $-100$, $-40$, $-20$, $-80$, $-78$, $-18$, $-42$, $-102$, $-56$, $-4$, $-64$, $-94$, $-34$, $-26$, $-86$, $-72$, $-12$, $-48$, $-108$, $-50$, $-10$, $-70$, $-88$, $-28$, $-32$, $-92$, $-66$, $-6$, $-54$, $-104$, $-44$, $-16$, $-76$, $-82$, $-22$, $-38$, $-98$, $-60)$\\

\noindent$11a186$ : $2$ | $(4$, $12$, $14$, $16$, $20$, $2$, $10$, $6$, $22$, $8$, $18)$ | $(150$, $116$, $184$, $106$, $160$, $140$, $126$, $174$, $96$, $170$, $130$, $136$, $164$, $102$, $180$, $120$, $146$, $154$, $112$, $188$, $110$, $156$, $144$, $122$, $178$, $100$, $166$, $134$, $132$, $168$, $98$, $176$, $124$, $142$, $158$, $108$, $186$, $114$, $152$, $148$, $118$, $182$, $104$, $162$, $138$, $128$, $172$, $-244$, $-210$, $-278$, $-200$, $-254$, $-234$, $-220$, $-268$, $-190$, $-264$, $-224$, $-230$, $-258$, $-196$, $-274$, $-214$, $-240$, $-248$, $-206$, $-282$, $-204$, $-250$, $-238$, $-216$, $-272$, $-194$, $-260$, $-228$, $-226$, $-262$, $-192$, $-270$, $-218$, $-236$, $-252$, $-202$, $-280$, $-208$, $-246$, $-242$, $-212$, $-276$, $-198$, $-256$, $-232$, $-222$, $-266$, $338$, $304$, $372$, $294$, $348$, $328$, $314$, $362$, $284$, $358$, $318$, $324$, $352$, $290$, $368$, $308$, $334$, $342$, $300$, $376$, $298$, $344$, $332$, $310$, $366$, $288$, $354$, $322$, $320$, $356$, $286$, $364$, $312$, $330$, $346$, $296$, $374$, $302$, $340$, $336$, $306$, $370$, $292$, $350$, $326$, $316$, $360$, $-56$, $-22$, $-90$, $-12$, $-66$, $-46$, $-32$, $-80$, $-2$, $-76$, $-36$, $-42$, $-70$, $-8$, $-86$, $-26$, $-52$, $-60$, $-18$, $-94$, $-16$, $-62$, $-50$, $-28$, $-84$, $-6$, $-72$, $-40$, $-38$, $-74$, $-4$, $-82$, $-30$, $-48$, $-64$, $-14$, $-92$, $-20$, $-58$, $-54$, $-24$, $-88$, $-10$, $-68$, $-44$, $-34$, $-78)$\\

\noindent$11a187$ : $3$ | $(4$, $12$, $14$, $16$, $20$, $2$, $22$, $8$, $6$, $10$, $18)$ | $(48$, $44$, $52$, $68$, $50$, $46$, $-2$, $-30$, $-36$, $112$, $100$, $118$, $122$, $104$, $108$, $116$, $98$, $114$, $110$, $102$, $120$, $-24$, $-72$, $-74$, $-76$, $-70$, $-26$, $-40$, $-22$, $-32$, $-34$, $-20$, $-38$, $-28$, $106$, $124$, $12$, $126$, $14$, $58$, $62$, $18$, $64$, $56$, $42$, $54$, $66$, $16$, $60$, $-80$, $-86$, $-92$, $-6$, $-8$, $-94$, $-84$, $-82$, $-96$, $-78$, $-88$, $-90$, $-4$, $-10)$\\

\noindent$11a188$ : $2$ | $(4$, $12$, $14$, $16$, $20$, $18$, $2$, $6$, $22$, $10$, $8)$ | $(90$, $128$, $80$, $100$, $118$, $70$, $110$, $108$, $72$, $120$, $98$, $82$, $130$, $88$, $92$, $126$, $78$, $102$, $116$, $68$, $112$, $106$, $74$, $122$, $96$, $84$, $132$, $86$, $94$, $124$, $76$, $104$, $114$, $-146$, $-174$, $-196$, $-168$, $-140$, $-152$, $-180$, $-190$, $-162$, $-134$, $-158$, $-186$, $-184$, $-156$, $-136$, $-164$, $-192$, $-178$, $-150$, $-142$, $-170$, $-198$, $-172$, $-144$, $-148$, $-176$, $-194$, $-166$, $-138$, $-154$, $-182$, $-188$, $-160$, $222$, $260$, $212$, $232$, $250$, $202$, $242$, $240$, $204$, $252$, $230$, $214$, $262$, $220$, $224$, $258$, $210$, $234$, $248$, $200$, $244$, $238$, $206$, $254$, $228$, $216$, $264$, $218$, $226$, $256$, $208$, $236$, $246$, $-14$, $-42$, $-64$, $-36$, $-8$, $-20$, $-48$, $-58$, $-30$, $-2$, $-26$, $-54$, $-52$, $-24$, $-4$, $-32$, $-60$, $-46$, $-18$, $-10$, $-38$, $-66$, $-40$, $-12$, $-16$, $-44$, $-62$, $-34$, $-6$, $-22$, $-50$, $-56$, $-28)$\\

\noindent$11a189$ : $3$ | $(4$, $12$, $14$, $16$, $22$, $18$, $2$, $20$, $10$, $6$, $8)$ |  $(112$, $104$, $86$, $88$, $102$, $114$, $32$, $36$, $-130$, $-136$, $-118$, $-8$, $-4$, $-122$, $-140$, $-126$, $-134$, $-132$, $-128$, $-138$, $-120$, $-6$, $24$, $146$, $152$, $144$, $26$, $42$, $22$, $148$, $150$, $20$, $40$, $28$, $142$, $30$, $38$, $18$, $34$, $-14$, $-72$, $-66$, $-54$, $-44$, $-56$, $-64$, $-74$, $-16$, $-76$, $-62$, $-58$, $-46$, $-52$, $-68$, $-70$, $-50$, $-48$, $-60$, $100$, $90$, $84$, $106$, $110$, $80$, $94$, $96$, $78$, $98$, $92$, $82$, $108$, $-124$, $-2$, $-10$, $-116$, $-12)$\\

\noindent$11a190$ : $2$ | $(4$, $12$, $14$, $18$, $16$, $2$, $22$, $20$, $8$, $6$, $10)$ | $(136$, $98$, $164$, $108$, $126$, $146$, $88$, $154$, $118$, $116$, $156$, $90$, $144$, $128$, $106$, $166$, $100$, $134$, $138$, $96$, $162$, $110$, $124$, $148$, $86$, $152$, $120$, $114$, $158$, $92$, $142$, $130$, $104$, $168$, $102$, $132$, $140$, $94$, $160$, $112$, $122$, $150$, $-186$, $-222$, $-248$, $-212$, $-176$, $-196$, $-232$, $-238$, $-202$, $-170$, $-206$, $-242$, $-228$, $-192$, $-180$, $-216$, $-252$, $-218$, $-182$, $-190$, $-226$, $-244$, $-208$, $-172$, $-200$, $-236$, $-234$, $-198$, $-174$, $-210$, $-246$, $-224$, $-188$, $-184$, $-220$, $-250$, $-214$, $-178$, $-194$, $-230$, $-240$, $-204$, $304$, $266$, $332$, $276$, $294$, $314$, $256$, $322$, $286$, $284$, $324$, $258$, $312$, $296$, $274$, $334$, $268$, $302$, $306$, $264$, $330$, $278$, $292$, $316$, $254$, $320$, $288$, $282$, $326$, $260$, $310$, $298$, $272$, $336$, $270$, $300$, $308$, $262$, $328$, $280$, $290$, $318$, $-18$, $-54$, $-80$, $-44$, $-8$, $-28$, $-64$, $-70$, $-34$, $-2$, $-38$, $-74$, $-60$, $-24$, $-12$, $-48$, $-84$, $-50$, $-14$, $-22$, $-58$, $-76$, $-40$, $-4$, $-32$, $-68$, $-66$, $-30$, $-6$, $-42$, $-78$, $-56$, $-20$, $-16$, $-52$, $-82$, $-46$, $-10$, $-26$, $-62$, $-72$, $-36)$\\

\noindent$11a191$ : $2$ | $(4$, $12$, $14$, $18$, $20$, $2$, $10$, $22$, $6$, $8$, $16)$ | $(146$, $108$, $94$, $132$, $160$, $122$, $84$, $118$, $156$, $136$, $98$, $104$, $142$, $150$, $112$, $90$, $128$, $164$, $126$, $88$, $114$, $152$, $140$, $102$, $100$, $138$, $154$, $116$, $86$, $124$, $162$, $130$, $92$, $110$, $148$, $144$, $106$, $96$, $134$, $158$, $120$, $-212$, $-186$, $-238$, $-168$, $-230$, $-194$, $-204$, $-220$, $-178$, $-246$, $-176$, $-222$, $-202$, $-196$, $-228$, $-170$, $-240$, $-184$, $-214$, $-210$, $-188$, $-236$, $-166$, $-232$, $-192$, $-206$, $-218$, $-180$, $-244$, $-174$, $-224$, $-200$, $-198$, $-226$, $-172$, $-242$, $-182$, $-216$, $-208$, $-190$, $-234$, $310$, $272$, $258$, $296$, $324$, $286$, $248$, $282$, $320$, $300$, $262$, $268$, $306$, $314$, $276$, $254$, $292$, $328$, $290$, $252$, $278$, $316$, $304$, $266$, $264$, $302$, $318$, $280$, $250$, $288$, $326$, $294$, $256$, $274$, $312$, $308$, $270$, $260$, $298$, $322$, $284$, $-48$, $-22$, $-74$, $-4$, $-66$, $-30$, $-40$, $-56$, $-14$, $-82$, $-12$, $-58$, $-38$, $-32$, $-64$, $-6$, $-76$, $-20$, $-50$, $-46$, $-24$, $-72$, $-2$, $-68$, $-28$, $-42$, $-54$, $-16$, $-80$, $-10$, $-60$, $-36$, $-34$, $-62$, $-8$, $-78$, $-18$, $-52$, $-44$, $-26$, $-70)$\\

\noindent$11a192$ : $2$ | $(4$, $12$, $14$, $18$, $20$, $2$, $10$, $22$, $8$, $6$, $16)$ | $(122$, $174$, $160$, $108$, $136$, $188$, $146$, $98$, $150$, $184$, $132$, $112$, $164$, $170$, $118$, $126$, $178$, $156$, $104$, $140$, $192$, $142$, $102$, $154$, $180$, $128$, $116$, $168$, $166$, $114$, $130$, $182$, $152$, $100$, $144$, $190$, $138$, $106$, $158$, $176$, $124$, $120$, $172$, $162$, $110$, $134$, $186$, $148$, $-248$, $-218$, $-278$, $-196$, $-270$, $-226$, $-240$, $-256$, $-210$, $-286$, $-204$, $-262$, $-234$, $-232$, $-264$, $-202$, $-284$, $-212$, $-254$, $-242$, $-224$, $-272$, $-194$, $-276$, $-220$, $-246$, $-250$, $-216$, $-280$, $-198$, $-268$, $-228$, $-238$, $-258$, $-208$, $-288$, $-206$, $-260$, $-236$, $-230$, $-266$, $-200$, $-282$, $-214$, $-252$, $-244$, $-222$, $-274$, $314$, $366$, $352$, $300$, $328$, $380$, $338$, $290$, $342$, $376$, $324$, $304$, $356$, $362$, $310$, $318$, $370$, $348$, $296$, $332$, $384$, $334$, $294$, $346$, $372$, $320$, $308$, $360$, $358$, $306$, $322$, $374$, $344$, $292$, $336$, $382$, $330$, $298$, $350$, $368$, $316$, $312$, $364$, $354$, $302$, $326$, $378$, $340$, $-56$, $-26$, $-86$, $-4$, $-78$, $-34$, $-48$, $-64$, $-18$, $-94$, $-12$, $-70$, $-42$, $-40$, $-72$, $-10$, $-92$, $-20$, $-62$, $-50$, $-32$, $-80$, $-2$, $-84$, $-28$, $-54$, $-58$, $-24$, $-88$, $-6$, $-76$, $-36$, $-46$, $-66$, $-16$, $-96$, $-14$, $-68$, $-44$, $-38$, $-74$, $-8$, $-90$, $-22$, $-60$, $-52$, $-30$, $-82)$\\

\noindent$11a193$ : $2$ | $(4$, $12$, $14$, $18$, $20$, $16$, $2$, $22$, $10$, $8$, $6)$ | $(130$, $176$, $104$, $156$, $150$, $110$, $182$, $124$, $136$, $170$, $98$, $162$, $144$, $116$, $188$, $118$, $142$, $164$, $96$, $168$, $138$, $122$, $184$, $112$, $148$, $158$, $102$, $174$, $132$, $128$, $178$, $106$, $154$, $152$, $108$, $180$, $126$, $134$, $172$, $100$, $160$, $146$, $114$, $186$, $120$, $140$, $166$, $-254$, $-196$, $-238$, $-270$, $-212$, $-222$, $-280$, $-228$, $-206$, $-264$, $-244$, $-190$, $-248$, $-260$, $-202$, $-232$, $-276$, $-218$, $-216$, $-274$, $-234$, $-200$, $-258$, $-250$, $-192$, $-242$, $-266$, $-208$, $-226$, $-282$, $-224$, $-210$, $-268$, $-240$, $-194$, $-252$, $-256$, $-198$, $-236$, $-272$, $-214$, $-220$, $-278$, $-230$, $-204$, $-262$, $-246$, $318$, $364$, $292$, $344$, $338$, $298$, $370$, $312$, $324$, $358$, $286$, $350$, $332$, $304$, $376$, $306$, $330$, $352$, $284$, $356$, $326$, $310$, $372$, $300$, $336$, $346$, $290$, $362$, $320$, $316$, $366$, $294$, $342$, $340$, $296$, $368$, $314$, $322$, $360$, $288$, $348$, $334$, $302$, $374$, $308$, $328$, $354$, $-66$, $-8$, $-50$, $-82$, $-24$, $-34$, $-92$, $-40$, $-18$, $-76$, $-56$, $-2$, $-60$, $-72$, $-14$, $-44$, $-88$, $-30$, $-28$, $-86$, $-46$, $-12$, $-70$, $-62$, $-4$, $-54$, $-78$, $-20$, $-38$, $-94$, $-36$, $-22$, $-80$, $-52$, $-6$, $-64$, $-68$, $-10$, $-48$, $-84$, $-26$, $-32$, $-90$, $-42$, $-16$, $-74$, $-58)$\\

\noindent$11a194$ : $3$ | $(4$, $12$, $14$, $20$, $16$, $18$, $2$, $22$, $10$, $6$, $8)$ | $(82$, $112$, $64$, $100$, $94$, $70$, $118$, $76$, $88$, $106$, $110$, $84$, $80$, $114$, $66$, $98$, $96$, $68$, $116$, $78$, $86$, $108$, $-20$, $-62$, $102$, $92$, $72$, $120$, $74$, $90$, $104$, $210$, $-150$, $-174$, $-156$, $-132$, $-122$, $-136$, $-160$, $-170$, $-146$, $-144$, $-168$, $-162$, $-138$, $-124$, $-130$, $-154$, $-176$, $-152$, $-128$, $-126$, $-140$, $-164$, $-166$, $-142$, $-148$, $-172$, $-158$, $-134$, $200$, $220$, $182$, $48$, $180$, $218$, $202$, $198$, $222$, $184$, $46$, $190$, $228$, $192$, $208$, $212$, $178$, $216$, $204$, $196$, $224$, $186$, $44$, $188$, $226$, $194$, $206$, $214$, $-32$, $-8$, $-50$, $-60$, $-18$, $-42$, $-22$, $-2$, $-26$, $-38$, $-14$, $-56$, $-54$, $-12$, $-36$, $-28$, $-4$, $-6$, $-30$, $-34$, $-10$, $-52$, $-58$, $-16$, $-40$, $-24)$\\

\noindent$11a195$ : $2$ | $(4$, $12$, $14$, $20$, $18$, $16$, $2$, $22$, $10$, $8$, $6)$ | $(72$, $98$, $58$, $86$, $84$, $60$, $100$, $70$, $74$, $96$, $56$, $88$, $82$, $62$, $102$, $68$, $76$, $94$, $54$, $90$, $80$, $64$, $104$, $66$, $78$, $92$, $-112$, $-128$, $-144$, $-154$, $-138$, $-122$, $-106$, $-118$, $-134$, $-150$, $-148$, $-132$, $-116$, $-108$, $-124$, $-140$, $-156$, $-142$, $-126$, $-110$, $-114$, $-130$, $-146$, $-152$, $-136$, $-120$, $176$, $202$, $162$, $190$, $188$, $164$, $204$, $174$, $178$, $200$, $160$, $192$, $186$, $166$, $206$, $172$, $180$, $198$, $158$, $194$, $184$, $168$, $208$, $170$, $182$, $196$, $-8$, $-24$, $-40$, $-50$, $-34$, $-18$, $-2$, $-14$, $-30$, $-46$, $-44$, $-28$, $-12$, $-4$, $-20$, $-36$, $-52$, $-38$, $-22$, $-6$, $-10$, $-26$, $-42$, $-48$, $-32$, $-16)$\\

\noindent$11a196$ : $3$ | $(4$, $12$, $14$, $20$, $22$, $18$, $2$, $8$, $10$, $6$, $16)$ | $(32$, $74$, $58$, $48$, $84$, $42$, $180$, $138$, $178$, $40$, $82$, $50$, $56$, $76$, $34$, $-122$, $-124$, $-214$, $-196$, $-224$, $-134$, $-112$, $-88$, $-116$, $-130$, $-220$, $-192$, $-218$, $-128$, $-118$, $-90$, $-110$, $-14$, $-16$, $-108$, $-92$, $-120$, $-126$, $-216$, $-194$, $-222$, $-132$, $-114$, $158$, $184$, $142$, $174$, $36$, $78$, $54$, $52$, $80$, $38$, $176$, $140$, $182$, $160$, $156$, $186$, $144$, $172$, $240$, $170$, $146$, $188$, $154$, $162$, $232$, $66$, $-228$, $-200$, $-210$, $-4$, $-26$, $-98$, $-102$, $-22$, $-8$, $-206$, $-204$, $-10$, $-20$, $-104$, $-96$, $-28$, $-2$, $-212$, $-198$, $-226$, $-136$, $-230$, $-202$, $-208$, $-6$, $-24$, $-100$, $150$, $166$, $236$, $70$, $62$, $44$, $86$, $46$, $60$, $72$, $238$, $168$, $148$, $190$, $152$, $164$, $234$, $68$, $64$, $-12$, $-18$, $-106$, $-94$, $-30)$\\

\noindent$11a197$ : $3$ | $(4$, $12$, $16$, $14$, $18$, $2$, $8$, $22$, $20$, $10$, $6)$ | $(48$, $42$, $22$, $38$, $52$, $104$, $84$, $102$, $50$, $40$, $-76$, $-114$, $-122$, $-108$, $-14$, $-12$, $-106$, $-120$, $-116$, $-78$, $-74$, $-112$, $-124$, $-110$, $-72$, $-80$, $-118$, $36$, $24$, $44$, $46$, $26$, $34$, $32$, $28$, $138$, $128$, $94$, $92$, $130$, $136$, $30$, $-68$, $-64$, $-56$, $-6$, $-20$, $-4$, $-58$, $-62$, $-70$, $-82$, $-66$, $86$, $100$, $126$, $96$, $90$, $132$, $134$, $88$, $98$, $-60$, $-2$, $-18$, $-8$, $-54$, $-10$, $-16)$\\

\noindent$11a198$ : $3$ | $(4$, $12$, $16$, $14$, $18$, $2$, $22$, $6$, $20$, $10$, $8)$ | $(10$, $86$, $68$, $84$, $-64$, $-62$, $-46$, $-48$, $-60$, $-38$, $-56$, $-52$, $-42$, $-4$, $-6$, $-44$, $-50$, $-58$, $76$, $94$, $88$, $70$, $82$, $98$, $80$, $72$, $90$, $92$, $74$, $78$, $96$, $134$, $132$, $28$, $-116$, $-102$, $-124$, $-108$, $-110$, $-122$, $-100$, $-118$, $-66$, $-114$, $-104$, $-126$, $-106$, $-112$, $-120$, $14$, $32$, $24$, $128$, $20$, $36$, $18$, $130$, $26$, $30$, $12$, $16$, $34$, $22$, $-54$, $-40$, $-2$, $-8)$\\

\noindent$11a199$ : $3$ | $(4$, $12$, $16$, $14$, $18$, $20$, $2$, $6$, $22$, $10$, $8)$ | $(20$, $34$, $36$, $18$, $-116$, $-104$, $-90$, $-106$, $-114$, $-118$, $-58$, $-42$, $-50$, $-66$, $-52$, $-40$, $-56$, $-120$, $-54$, $124$, $70$, $132$, $12$, $14$, $130$, $68$, $126$, $122$, $72$, $134$, $10$, $16$, $128$, $-110$, $-94$, $-100$, $-2$, $-6$, $-64$, $-48$, $-44$, $-60$, $-8$, $-62$, $-46$, $80$, $24$, $30$, $38$, $32$, $22$, $82$, $88$, $78$, $26$, $28$, $76$, $86$, $84$, $74$, $-4$, $-98$, $-96$, $-112$, $-108$, $-92$, $-102)$\\

\noindent$11a200$ : $3$ | $(4$, $12$, $16$, $14$, $20$, $2$, $10$, $6$, $22$, $8$, $18)$ | $(44$, $22$, $40$, $48$, $26$, $30$, $28$, $24$, $46$, $42$, $-64$, $-58$, $-70$, $-32$, $-72$, $74$, $38$, $50$, $-52$, $-56$, $-66$, $-62$, $-60$, $-68$, $-54$, $90$, $88$, $92$, $78$, $84$, $96$, $82$, $80$, $94$, $86$, $76$, $-12$, $-10$, $-14$, $-4$, $-20$, $-2$, $-16$, $-8$, $-36$, $-34$, $-6$, $-18)$\\

\noindent$11a201$ : $3$ | $(4$, $12$, $16$, $14$, $20$, $18$, $2$, $6$, $22$, $10$, $8)$ | $(26$, $34$, $54$, $32$, $28$, $24$, $36$, $22$, $20$, $-60$, $-62$, $-38$, $-48$, $-50$, $-40$, $-42$, $-52$, $-46$, $72$, $80$, $66$, $78$, $74$, $70$, $82$, $68$, $76$, $-44$, $-16$, $-18$, $58$, $56$, $30$, $-8$, $-2$, $-12$, $-14$, $-4$, $-6$, $-64$, $-10$)\\

\noindent$11a202$ : $3$ | $(4$, $12$, $16$, $18$, $14$, $20$, $2$, $8$, $6$, $22$, $10)$ | $(46$, $84$, $92$, $74$, $96$, $80$, $50$, $42$, $88$, $-62$, $-4$, $-70$, $-10$, $-56$, $-58$, $-8$, $-72$, $-6$, $-60$, $-54$, $-52$, $-106$, $-122$, $-120$, $-108$, $144$, $134$, $78$, $98$, $76$, $90$, $86$, $44$, $48$, $82$, $94$, $-114$, $-128$, $-100$, $-132$, $-110$, $-118$, $-124$, $-104$, $-102$, $-126$, $-116$, $-112$, $-130$, $28$, $38$, $18$, $136$, $146$, $142$, $24$, $32$, $34$, $22$, $140$, $148$, $138$, $20$, $36$, $30$, $26$, $40$, $-16$, $-64$, $-2$, $-68$, $-12$, $-14$, $-66)$\\

\noindent$11a203$ : $2$ | $(4$, $12$, $16$, $18$, $20$, $2$, $22$, $10$, $6$, $8$, $14)$ | $(114$, $92$, $70$, $76$, $98$, $120$, $108$, $86$, $64$, $82$, $104$, $124$, $102$, $80$, $66$, $88$, $110$, $118$, $96$, $74$, $72$, $94$, $116$, $112$, $90$, $68$, $78$, $100$, $122$, $106$, $84$, $-164$, $-130$, $-176$, $-152$, $-142$, $-186$, $-140$, $-154$, $-174$, $-128$, $-166$, $-162$, $-132$, $-178$, $-150$, $-144$, $-184$, $-138$, $-156$, $-172$, $-126$, $-168$, $-160$, $-134$, $-180$, $-148$, $-146$, $-182$, $-136$, $-158$, $-170$, $238$, $216$, $194$, $200$, $222$, $244$, $232$, $210$, $188$, $206$, $228$, $248$, $226$, $204$, $190$, $212$, $234$, $242$, $220$, $198$, $196$, $218$, $240$, $236$, $214$, $192$, $202$, $224$, $246$, $230$, $208$, $-40$, $-6$, $-52$, $-28$, $-18$, $-62$, $-16$, $-30$, $-50$, $-4$, $-42$, $-38$, $-8$, $-54$, $-26$, $-20$, $-60$, $-14$, $-32$, $-48$, $-2$, $-44$, $-36$, $-10$, $-56$, $-24$, $-22$, $-58$, $-12$, $-34$, $-46)$\\

\noindent$11a204$ : $2$ | $(4$, $12$, $16$, $18$, $20$, $2$, $22$, $10$, $8$, $6$, $14)$ | $(130$, $190$, $152$, $108$, $168$, $174$, $114$, $146$, $196$, $136$, $124$, $184$, $158$, $102$, $162$, $180$, $120$, $140$, $200$, $142$, $118$, $178$, $164$, $104$, $156$, $186$, $126$, $134$, $194$, $148$, $112$, $172$, $170$, $110$, $150$, $192$, $132$, $128$, $188$, $154$, $106$, $166$, $176$, $116$, $144$, $198$, $138$, $122$, $182$, $160$, $-264$, $-210$, $-284$, $-244$, $-230$, $-298$, $-224$, $-250$, $-278$, $-204$, $-270$, $-258$, $-216$, $-290$, $-238$, $-236$, $-292$, $-218$, $-256$, $-272$, $-202$, $-276$, $-252$, $-222$, $-296$, $-232$, $-242$, $-286$, $-212$, $-262$, $-266$, $-208$, $-282$, $-246$, $-228$, $-300$, $-226$, $-248$, $-280$, $-206$, $-268$, $-260$, $-214$, $-288$, $-240$, $-234$, $-294$, $-220$, $-254$, $-274$, $330$, $390$, $352$, $308$, $368$, $374$, $314$, $346$, $396$, $336$, $324$, $384$, $358$, $302$, $362$, $380$, $320$, $340$, $400$, $342$, $318$, $378$, $364$, $304$, $356$, $386$, $326$, $334$, $394$, $348$, $312$, $372$, $370$, $310$, $350$, $392$, $332$, $328$, $388$, $354$, $306$, $366$, $376$, $316$, $344$, $398$, $338$, $322$, $382$, $360$, $-64$, $-10$, $-84$, $-44$, $-30$, $-98$, $-24$, $-50$, $-78$, $-4$, $-70$, $-58$, $-16$, $-90$, $-38$, $-36$, $-92$, $-18$, $-56$, $-72$, $-2$, $-76$, $-52$, $-22$, $-96$, $-32$, $-42$, $-86$, $-12$, $-62$, $-66$, $-8$, $-82$, $-46$, $-28$, $-100$, $-26$, $-48$, $-80$, $-6$, $-68$, $-60$, $-14$, $-88$, $-40$, $-34$, $-94$, $-20$, $-54$, $-74)$\\

\noindent$11a205$ : $2$ | $(4$, $12$, $16$, $18$, $20$, $14$, $2$, $10$, $22$, $8$, $6)$ | $(130$, $152$, $108$, $174$, $94$, $166$, $116$, $144$, $138$, $122$, $160$, $100$, $180$, $102$, $158$, $124$, $136$, $146$, $114$, $168$, $92$, $172$, $110$, $150$, $132$, $128$, $154$, $106$, $176$, $96$, $164$, $118$, $142$, $140$, $120$, $162$, $98$, $178$, $104$, $156$, $126$, $134$, $148$, $112$, $170$, $-246$, $-196$, $-214$, $-264$, $-228$, $-182$, $-232$, $-260$, $-210$, $-200$, $-250$, $-242$, $-192$, $-218$, $-268$, $-224$, $-186$, $-236$, $-256$, $-206$, $-204$, $-254$, $-238$, $-188$, $-222$, $-270$, $-220$, $-190$, $-240$, $-252$, $-202$, $-208$, $-258$, $-234$, $-184$, $-226$, $-266$, $-216$, $-194$, $-244$, $-248$, $-198$, $-212$, $-262$, $-230$, $310$, $332$, $288$, $354$, $274$, $346$, $296$, $324$, $318$, $302$, $340$, $280$, $360$, $282$, $338$, $304$, $316$, $326$, $294$, $348$, $272$, $352$, $290$, $330$, $312$, $308$, $334$, $286$, $356$, $276$, $344$, $298$, $322$, $320$, $300$, $342$, $278$, $358$, $284$, $336$, $306$, $314$, $328$, $292$, $350$, $-66$, $-16$, $-34$, $-84$, $-48$, $-2$, $-52$, $-80$, $-30$, $-20$, $-70$, $-62$, $-12$, $-38$, $-88$, $-44$, $-6$, $-56$, $-76$, $-26$, $-24$, $-74$, $-58$, $-8$, $-42$, $-90$, $-40$, $-10$, $-60$, $-72$, $-22$, $-28$, $-78$, $-54$, $-4$, $-46$, $-86$, $-36$, $-14$, $-64$, $-68$, $-18$, $-32$, $-82$, $-50)$\\

\noindent$11a206$ : $2$ | $(4$, $12$, $16$, $18$, $20$, $14$, $2$, $22$, $6$, $8$, $10)$ | $(66$, $80$, $52$, $92$, $54$, $78$, $68$, $64$, $82$, $50$, $90$, $56$, $76$, $70$, $62$, $84$, $48$, $88$, $58$, $74$, $72$, $60$, $86$, $-132$, $-118$, $-104$, $-94$, $-108$, $-122$, $-136$, $-128$, $-114$, $-100$, $-98$, $-112$, $-126$, $-138$, $-124$, $-110$, $-96$, $-102$, $-116$, $-130$, $-134$, $-120$, $-106$, $158$, $172$, $144$, $184$, $146$, $170$, $160$, $156$, $174$, $142$, $182$, $148$, $168$, $162$, $154$, $176$, $140$, $180$, $150$, $166$, $164$, $152$, $178$, $-40$, $-26$, $-12$, $-2$, $-16$, $-30$, $-44$, $-36$, $-22$, $-8$, $-6$, $-20$, $-34$, $-46$, $-32$, $-18$, $-4$, $-10$, $-24$, $-38$, $-42$, $-28$, $-14)$\\

\noindent$11a207$ : $2$ | $(4$, $12$, $16$, $18$, $20$, $14$, $2$, $22$, $6$, $10$, $8)$ | $(120$, $146$, $94$, $166$, $100$, $140$, $126$, $114$, $152$, $88$, $160$, $106$, $134$, $132$, $108$, $158$, $86$, $154$, $112$, $128$, $138$, $102$, $164$, $92$, $148$, $118$, $122$, $144$, $96$, $168$, $98$, $142$, $124$, $116$, $150$, $90$, $162$, $104$, $136$, $130$, $110$, $156$, $-194$, $-246$, $-208$, $-180$, $-232$, $-222$, $-170$, $-218$, $-236$, $-184$, $-204$, $-250$, $-198$, $-190$, $-242$, $-212$, $-176$, $-228$, $-226$, $-174$, $-214$, $-240$, $-188$, $-200$, $-252$, $-202$, $-186$, $-238$, $-216$, $-172$, $-224$, $-230$, $-178$, $-210$, $-244$, $-192$, $-196$, $-248$, $-206$, $-182$, $-234$, $-220$, $288$, $314$, $262$, $334$, $268$, $308$, $294$, $282$, $320$, $256$, $328$, $274$, $302$, $300$, $276$, $326$, $254$, $322$, $280$, $296$, $306$, $270$, $332$, $260$, $316$, $286$, $290$, $312$, $264$, $336$, $266$, $310$, $292$, $284$, $318$, $258$, $330$, $272$, $304$, $298$, $278$, $324$, $-26$, $-78$, $-40$, $-12$, $-64$, $-54$, $-2$, $-50$, $-68$, $-16$, $-36$, $-82$, $-30$, $-22$, $-74$, $-44$, $-8$, $-60$, $-58$, $-6$, $-46$, $-72$, $-20$, $-32$, $-84$, $-34$, $-18$, $-70$, $-48$, $-4$, $-56$, $-62$, $-10$, $-42$, $-76$, $-24$, $-28$, $-80$, $-38$, $-14$, $-66$, $-52)$\\

\noindent$11a208$ : $2$ | $(4$, $12$, $16$, $18$, $20$, $14$, $2$, $22$, $10$, $8$, $6)$ | $(148$, $182$, $114$, $202$, $128$, $168$, $162$, $134$, $196$, $108$, $188$, $142$, $154$, $176$, $120$, $208$, $122$, $174$, $156$, $140$, $190$, $106$, $194$, $136$, $160$, $170$, $126$, $204$, $116$, $180$, $150$, $146$, $184$, $112$, $200$, $130$, $166$, $164$, $132$, $198$, $110$, $186$, $144$, $152$, $178$, $118$, $206$, $124$, $172$, $158$, $138$, $192$, $-282$, $-220$, $-258$, $-306$, $-244$, $-234$, $-296$, $-268$, $-210$, $-272$, $-292$, $-230$, $-248$, $-310$, $-254$, $-224$, $-286$, $-278$, $-216$, $-262$, $-302$, $-240$, $-238$, $-300$, $-264$, $-214$, $-276$, $-288$, $-226$, $-252$, $-312$, $-250$, $-228$, $-290$, $-274$, $-212$, $-266$, $-298$, $-236$, $-242$, $-304$, $-260$, $-218$, $-280$, $-284$, $-222$, $-256$, $-308$, $-246$, $-232$, $-294$, $-270$, $356$, $390$, $322$, $410$, $336$, $376$, $370$, $342$, $404$, $316$, $396$, $350$, $362$, $384$, $328$, $416$, $330$, $382$, $364$, $348$, $398$, $314$, $402$, $344$, $368$, $378$, $334$, $412$, $324$, $388$, $358$, $354$, $392$, $320$, $408$, $338$, $374$, $372$, $340$, $406$, $318$, $394$, $352$, $360$, $386$, $326$, $414$, $332$, $380$, $366$, $346$, $400$, $-74$, $-12$, $-50$, $-98$, $-36$, $-26$, $-88$, $-60$, $-2$, $-64$, $-84$, $-22$, $-40$, $-102$, $-46$, $-16$, $-78$, $-70$, $-8$, $-54$, $-94$, $-32$, $-30$, $-92$, $-56$, $-6$, $-68$, $-80$, $-18$, $-44$, $-104$, $-42$, $-20$, $-82$, $-66$, $-4$, $-58$, $-90$, $-28$, $-34$, $-96$, $-52$, $-10$, $-72$, $-76$, $-14$, $-48$, $-100$, $-38$, $-24$, $-86$, $-62)$\\

\noindent$11a209$ : $3$ | $(4$, $12$, $16$, $20$, $14$, $18$, $2$, $8$, $22$, $10$, $6)$ | $(64$, $22$, $28$, $124$, $126$, $26$, $24$, $128$, $122$, $30$, $-96$, $-50$, $-88$, $-52$, $-94$, $-44$, $80$, $60$, $68$, $66$, $62$, $82$, $78$, $58$, $70$, $86$, $74$, $-8$, $-16$, $-104$, $-110$, $-98$, $-100$, $-112$, $-102$, $-14$, $-10$, $-6$, $-18$, $-106$, $-108$, $-20$, $-4$, $-12$, $72$, $56$, $76$, $84$, $116$, $134$, $114$, $130$, $120$, $32$, $118$, $132$, $-38$, $-2$, $-42$, $-34$, $-46$, $-92$, $-54$, $-90$, $-48$, $-36$, $-40)$\\

\noindent$11a210$ : $2$ | $(4$, $12$, $16$, $20$, $18$, $14$, $2$, $10$, $22$, $8$, $6)$ | $(104$, $122$, $86$, $140$, $76$, $132$, $94$, $114$, $112$, $96$, $130$, $78$, $142$, $84$, $124$, $102$, $106$, $120$, $88$, $138$, $74$, $134$, $92$, $116$, $110$, $98$, $128$, $80$, $144$, $82$, $126$, $100$, $108$, $118$, $90$, $136$, $-160$, $-192$, $-210$, $-178$, $-146$, $-174$, $-206$, $-196$, $-164$, $-156$, $-188$, $-214$, $-182$, $-150$, $-170$, $-202$, $-200$, $-168$, $-152$, $-184$, $-216$, $-186$, $-154$, $-166$, $-198$, $-204$, $-172$, $-148$, $-180$, $-212$, $-190$, $-158$, $-162$, $-194$, $-208$, $-176$, $248$, $266$, $230$, $284$, $220$, $276$, $238$, $258$, $256$, $240$, $274$, $222$, $286$, $228$, $268$, $246$, $250$, $264$, $232$, $282$, $218$, $278$, $236$, $260$, $254$, $242$, $272$, $224$, $288$, $226$, $270$, $244$, $252$, $262$, $234$, $280$, $-16$, $-48$, $-66$, $-34$, $-2$, $-30$, $-62$, $-52$, $-20$, $-12$, $-44$, $-70$, $-38$, $-6$, $-26$, $-58$, $-56$, $-24$, $-8$, $-40$, $-72$, $-42$, $-10$, $-22$, $-54$, $-60$, $-28$, $-4$, $-36$, $-68$, $-46$, $-14$, $-18$, $-50$, $-64$, $-32)$\\

\noindent$11a211$ : $2$ | $(4$, $12$, $16$, $20$, $18$, $14$, $2$, $22$, $10$, $8$, $6)$ | $(94$, $116$, $72$, $128$, $82$, $106$, $104$, $84$, $126$, $70$, $118$, $92$, $96$, $114$, $74$, $130$, $80$, $108$, $102$, $86$, $124$, $68$, $120$, $90$, $98$, $112$, $76$, $132$, $78$, $110$, $100$, $88$, $122$, $-144$, $-168$, $-192$, $-182$, $-158$, $-134$, $-154$, $-178$, $-196$, $-172$, $-148$, $-140$, $-164$, $-188$, $-186$, $-162$, $-138$, $-150$, $-174$, $-198$, $-176$, $-152$, $-136$, $-160$, $-184$, $-190$, $-166$, $-142$, $-146$, $-170$, $-194$, $-180$, $-156$, $226$, $248$, $204$, $260$, $214$, $238$, $236$, $216$, $258$, $202$, $250$, $224$, $228$, $246$, $206$, $262$, $212$, $240$, $234$, $218$, $256$, $200$, $252$, $222$, $230$, $244$, $208$, $264$, $210$, $242$, $232$, $220$, $254$, $-12$, $-36$, $-60$, $-50$, $-26$, $-2$, $-22$, $-46$, $-64$, $-40$, $-16$, $-8$, $-32$, $-56$, $-54$, $-30$, $-6$, $-18$, $-42$, $-66$, $-44$, $-20$, $-4$, $-28$, $-52$, $-58$, $-34$, $-10$, $-14$, $-38$, $-62$, $-48$, $-24)$\\

\noindent$11a212$ : $3$ | $(4$, $12$, $16$, $22$, $18$, $2$, $8$, $6$, $20$, $10$, $14)$ | $(94$, $116$, $72$, $128$, $82$, $106$, $104$, $84$, $126$, $70$, $118$, $92$, $96$, $114$, $74$, $130$, $80$, $108$, $102$, $86$, $124$, $68$, $120$, $90$, $98$, $112$, $76$, $132$, $78$, $110$, $100$, $88$, $122$, $-144$, $-168$, $-192$, $-182$, $-158$, $-134$, $-154$, $-178$, $-196$, $-172$, $-148$, $-140$, $-164$, $-188$, $-186$, $-162$, $-138$, $-150$, $-174$, $-198$, $-176$, $-152$, $-136$, $-160$, $-184$, $-190$, $-166$, $-142$, $-146$, $-170$, $-194$, $-180$, $-156$, $226$, $248$, $204$, $260$, $214$, $238$, $236$, $216$, $258$, $202$, $250$, $224$, $228$, $246$, $206$, $262$, $212$, $240$, $234$, $218$, $256$, $200$, $252$, $222$, $230$, $244$, $208$, $264$, $210$, $242$, $232$, $220$, $254$, $-12$, $-36$, $-60$, $-50$, $-26$, $-2$, $-22$, $-46$, $-64$, $-40$, $-16$, $-8$, $-32$, $-56$, $-54$, $-30$, $-6$, $-18$, $-42$, $-66$, $-44$, $-20$, $-4$, $-28$, $-52$, $-58$, $-34$, $-10$, $-14$, $-38$, $-62$, $-48$, $-24)$\\

\noindent$11a213$ : $3$ | $(4$, $12$, $16$, $22$, $18$, $2$, $20$, $10$, $6$, $14$, $8)$ | $(78$, $74$, $86$, $114$, $112$, $84$, $72$, $80$, $76$, $-2$, $-14$, $-70$, $-12$, $-4$, $-62$, $-64$, $-6$, $-10$, $-68$, $-58$, $-60$, $-66$, $-8$, $34$, $26$, $20$, $122$, $18$, $120$, $108$, $118$, $88$, $116$, $110$, $82$, $-96$, $-50$, $-56$, $-54$, $-52$, $-94$, $-98$, $-48$, $-104$, $30$, $42$, $38$, $22$, $24$, $36$, $44$, $32$, $28$, $40$, $-92$, $-100$, $-46$, $-102$, $-90$, $-106$, $-16)$\\

\noindent$11a214$ : $3$ | $(4$, $12$, $16$, $22$, $20$, $18$, $2$, $6$, $14$, $10$, $8)$ | $(32$, $34$, $36$, $40$, $-20$, $-26$, $-24$, $-22$, $-6$, $38$, $42$, $56$, $44$, $58$, $60$, $62$, $-30$, $-28$, $-52$, $-50$, $-48$, $-54$, $12$, $16$, $8$, $18$, $10$, $14$, $-46$, $-4$, $-2)$\\

\noindent$11a215$ : $3$ | $(4$, $12$, $18$, $14$, $16$, $2$, $22$, $20$, $10$, $6$, $8)$ | $(40$, $82$, $52$, $78$, $68$, $74$, $156$, $142$, $140$, $158$, $72$, $70$, $160$, $138$, $144$, $154$, $-12$, $-20$, $-24$, $-102$, $-104$, $-26$, $-18$, $-14$, $-30$, $-108$, $-98$, $80$, $38$, $42$, $84$, $50$, $76$, $-112$, $-94$, $-124$, $-92$, $-114$, $-110$, $-96$, $-100$, $-106$, $-28$, $-16$, $46$, $34$, $126$, $132$, $150$, $148$, $134$, $162$, $136$, $146$, $152$, $130$, $128$, $32$, $48$, $86$, $44$, $36$, $-22$, $-10$, $-66$, $-8$, $-56$, $-2$, $-60$, $-120$, $-88$, $-118$, $-62$, $-4$, $-54$, $-6$, $-64$, $-116$, $-90$, $-122$, $-58)$\\

\noindent$11a216$ : $3$ | $(4$, $12$, $18$, $14$, $20$, $2$, $22$, $8$, $10$, $6$, $16)$ | $(28$, $48$, $62$, $40$, $36$, $118$, $-148$, $-176$, $-174$, $-150$, $-194$, $-156$, $-168$, $-182$, $-184$, $-166$, $-158$, $-192$, $-84$, $-82$, $-190$, $-160$, $-164$, $-186$, $-180$, $-170$, $-154$, $-196$, $-152$, $-172$, $-178$, $-188$, $-162$, $44$, $32$, $24$, $52$, $58$, $18$, $14$, $212$, $144$, $214$, $12$, $20$, $56$, $54$, $22$, $34$, $42$, $64$, $46$, $30$, $26$, $50$, $60$, $16$, $-78$, $-88$, $-86$, $-80$, $-112$, $-76$, $-90$, $-94$, $-72$, $-108$, $-2$, $-104$, $-68$, $-98$, $-8$, $38$, $120$, $116$, $126$, $198$, $130$, $218$, $140$, $208$, $206$, $138$, $220$, $132$, $200$, $124$, $114$, $122$, $202$, $134$, $222$, $136$, $204$, $210$, $142$, $216$, $128$, $-92$, $-74$, $-110$, $-4$, $-102$, $-66$, $-100$, $-6$, $-146$, $-10$, $-96$, $-70$, $-106)$\\

\noindent$11a217$ : $3$ | $(4$, $12$, $18$, $14$, $20$, $16$, $2$, $8$, $22$, $6$, $10)$ | $(62$, $22$, $28$, $56$, $68$, $70$, $54$, $72$, $-92$, $-88$, $-106$, $-100$, $-82$, $-98$, $-108$, $-90$, $18$, $66$, $58$, $26$, $24$, $60$, $64$, $20$, $30$, $16$, $128$, $-110$, $-96$, $-84$, $-102$, $-104$, $-86$, $-94$, $-112$, $-52$, $-50$, $-32$, $-48$, $-4$, $-40$, $122$, $134$, $116$, $76$, $74$, $114$, $132$, $124$, $80$, $120$, $136$, $118$, $78$, $126$, $130$, $-10$, $-34$, $-46$, $-2$, $-42$, $-38$, $-6$, $-14$, $-12$, $-8$, $-36$, $-44)$\\

\noindent$11a218$ : $3$ | $(4$, $12$, $18$, $16$, $20$, $14$, $2$, $8$, $22$, $10$, $6)$ | $(92$, $46$, $44$, $90$, $42$, $48$, $32$, $58$, $30$, $50$, $40$, $36$, $54$, $-70$, $-80$, $-60$, $-84$, $-66$, $-74$, $-76$, $-64$, $-88$, $-110$, $-86$, $-62$, $-78$, $-72$, $-68$, $-82$, $124$, $134$, $114$, $94$, $120$, $128$, $130$, $118$, $96$, $116$, $132$, $126$, $122$, $136$, $112$, $140$, $-108$, $-106$, $-22$, $38$, $52$, $28$, $56$, $34$, $142$, $138$, $-4$, $-10$, $-20$, $-98$, $-24$, $-104$, $-14$, $-16$, $-102$, $-26$, $-100$, $-18$, $-12$, $-2$, $-6$, $-8$)\\

\noindent$11a219$ : $3$ | $(4$, $12$, $18$, $16$, $20$, $14$, $2$, $22$, $6$, $10$, $8)$ | $(48$, $50$, $80$, $84$, $88$, $76$, $52$, $78$, $86$, $-8$, $-40$, $-4$, $-14$, $-12$, $-44$, $-42$, $-10$, $-16$, $-6$, $82$, $90$, $74$, $92$, $94$, $-46$, $-62$, $26$, $34$, $18$, $36$, $24$, $28$, $32$, $20$, $38$, $22$, $30$, $-64$, $-54$, $-60$, $-70$, $-68$, $-58$, $-56$, $-66$, $-72$, $-2)$ \\

\noindent$11a220$ : $2$ | $(4$, $12$, $18$, $20$, $14$, $2$, $10$, $22$, $8$, $6$, $16)$ | $(132$, $110$, $154$, $88$, $162$, $102$, $140$, $124$, $118$, $146$, $96$, $168$, $94$, $148$, $116$, $126$, $138$, $104$, $160$, $86$, $156$, $108$, $134$, $130$, $112$, $152$, $90$, $164$, $100$, $142$, $122$, $120$, $144$, $98$, $166$, $92$, $150$, $114$, $128$, $136$, $106$, $158$, $-230$, $-184$, $-198$, $-244$, $-216$, $-170$, $-212$, $-248$, $-202$, $-180$, $-226$, $-234$, $-188$, $-194$, $-240$, $-220$, $-174$, $-208$, $-252$, $-206$, $-176$, $-222$, $-238$, $-192$, $-190$, $-236$, $-224$, $-178$, $-204$, $-250$, $-210$, $-172$, $-218$, $-242$, $-196$, $-186$, $-232$, $-228$, $-182$, $-200$, $-246$, $-214$, $300$, $278$, $322$, $256$, $330$, $270$, $308$, $292$, $286$, $314$, $264$, $336$, $262$, $316$, $284$, $294$, $306$, $272$, $328$, $254$, $324$, $276$, $302$, $298$, $280$, $320$, $258$, $332$, $268$, $310$, $290$, $288$, $312$, $266$, $334$, $260$, $318$, $282$, $296$, $304$, $274$, $326$, $-62$, $-16$, $-30$, $-76$, $-48$, $-2$, $-44$, $-80$, $-34$, $-12$, $-58$, $-66$, $-20$, $-26$, $-72$, $-52$, $-6$, $-40$, $-84$, $-38$, $-8$, $-54$, $-70$, $-24$, $-22$, $-68$, $-56$, $-10$, $-36$, $-82$, $-42$, $-4$, $-50$, $-74$, $-28$, $-18$, $-64$, $-60$, $-14$, $-32$, $-78$, $-46)$\\

\noindent$11a221$ : $3$ | $(4$, $12$, $18$, $20$, $14$, $16$, $2$, $10$, $22$, $6$, $8)$ | $(60$, $54$, $72$, $74$, $52$, $58$, $62$, $56$, $70$, $76$, $-8$, $-4$, $-16$, $-68$, $-18$, $-6$, $22$, $82$, $88$, $84$, $24$, $30$, $20$, $28$, $26$, $-66$, $-34$, $-40$, $-64$, $-36$, $-38$, $86$, $80$, $78$, $-48$, $-42$, $-32$, $-44$, $-46$, $-50$, $-10$, $-2$, $-14$, $-12)$\\

\noindent$11a222$ : $3$ | $(4$, $12$, $18$, $20$, $14$, $16$, $2$, $10$, $22$, $8$, $6)$ | $(64$, $76$, $62$, $66$, $26$, $42$, $24$, $38$, $30$, $28$, $40$, $-54$, $-44$, $-50$, $-58$, $-48$, $-4$, $-2$, $-46$, $-56$, $-52$, $84$, $90$, $78$, $60$, $82$, $86$, $88$, $80$, $-18$, $-68$, $-74$, $-72$, $92$, $32$, $36$, $94$, $34$, $-70$, $-16$, $-20$, $-10$, $-12$, $-22$, $-14$, $-8$, $-6)$\\

\noindent$11a223$ : $3$ | $(4$, $12$, $18$, $20$, $16$, $2$, $22$, $8$, $10$, $6$, $14)$ | $(44$, $116$, $128$, $32$, $56$, $110$, $50$, $38$, $122$, $-158$, $-144$, $-4$, $-136$, $-166$, $-134$, $-130$, $-72$, $-104$, $-80$, $-86$, $-98$, $-66$, $-64$, $-96$, $-88$, $-78$, $-106$, $-74$, $-132$, $30$, $24$, $184$, $210$, $182$, $22$, $172$, $200$, $194$, $190$, $204$, $176$, $18$, $178$, $206$, $188$, $28$, $26$, $186$, $208$, $180$, $20$, $174$, $202$, $192$, $-76$, $-90$, $-94$, $-62$, $-68$, $-100$, $-84$, $-82$, $-102$, $-70$, $170$, $198$, $196$, $168$, $114$, $46$, $42$, $118$, $126$, $34$, $54$, $108$, $52$, $36$, $124$, $120$, $40$, $48$, $112$, $-92$, $-60$, $-6$, $-146$, $-156$, $-16$, $-160$, $-142$, $-2$, $-138$, $-164$, $-12$, $-152$, $-150$, $-10$, $-58$, $-8$, $-148$, $-154$, $-14$, $-162$, $-140)$\\

\noindent$11a224$ : $2$ | $(4$, $12$, $18$, $20$, $16$, $2$, $22$, $10$, $8$, $6$, $14)$ | $(150$, $96$, $134$, $166$, $112$, $118$, $172$, $128$, $102$, $156$, $144$, $90$, $140$, $160$, $106$, $124$, $176$, $122$, $108$, $162$, $138$, $92$, $146$, $154$, $100$, $130$, $170$, $116$, $114$, $168$, $132$, $98$, $152$, $148$, $94$, $136$, $164$, $110$, $120$, $174$, $126$, $104$, $158$, $142$, $-232$, $-186$, $-252$, $-212$, $-206$, $-258$, $-192$, $-226$, $-238$, $-180$, $-246$, $-218$, $-200$, $-264$, $-198$, $-220$, $-244$, $-178$, $-240$, $-224$, $-194$, $-260$, $-204$, $-214$, $-250$, $-184$, $-234$, $-230$, $-188$, $-254$, $-210$, $-208$, $-256$, $-190$, $-228$, $-236$, $-182$, $-248$, $-216$, $-202$, $-262$, $-196$, $-222$, $-242$, $328$, $274$, $312$, $344$, $290$, $296$, $350$, $306$, $280$, $334$, $322$, $268$, $318$, $338$, $284$, $302$, $354$, $300$, $286$, $340$, $316$, $270$, $324$, $332$, $278$, $308$, $348$, $294$, $292$, $346$, $310$, $276$, $330$, $326$, $272$, $314$, $342$, $288$, $298$, $352$, $304$, $282$, $336$, $320$, $266$, $-56$, $-10$, $-76$, $-36$, $-30$, $-82$, $-16$, $-50$, $-62$, $-4$, $-70$, $-42$, $-24$, $-88$, $-22$, $-44$, $-68$, $-2$, $-64$, $-48$, $-18$, $-84$, $-28$, $-38$, $-74$, $-8$, $-58$, $-54$, $-12$, $-78$, $-34$, $-32$, $-80$, $-14$, $-52$, $-60$, $-6$, $-72$, $-40$, $-26$, $-86$, $-20$, $-46$, $-66)$\\

\noindent$11a225$ : $2$ | $(4$, $12$, $18$, $20$, $16$, $14$, $2$, $10$, $22$, $6$, $8)$ | $(76$, $86$, $66$, $96$, $56$, $104$, $58$, $94$, $68$, $84$, $78$, $74$, $88$, $64$, $98$, $54$, $102$, $60$, $92$, $70$, $82$, $80$, $72$, $90$, $62$, $100$, $-146$, $-124$, $-106$, $-128$, $-150$, $-142$, $-120$, $-110$, $-132$, $-154$, $-138$, $-116$, $-114$, $-136$, $-156$, $-134$, $-112$, $-118$, $-140$, $-152$, $-130$, $-108$, $-122$, $-144$, $-148$, $-126$, $180$, $190$, $170$, $200$, $160$, $208$, $162$, $198$, $172$, $188$, $182$, $178$, $192$, $168$, $202$, $158$, $206$, $164$, $196$, $174$, $186$, $184$, $176$, $194$, $166$, $204$, $-42$, $-20$, $-2$, $-24$, $-46$, $-38$, $-16$, $-6$, $-28$, $-50$, $-34$, $-12$, $-10$, $-32$, $-52$, $-30$, $-8$, $-14$, $-36$, $-48$, $-26$, $-4$, $-18$, $-40$, $-44$, $-22)$\\

\noindent$11a226$ : $2$ | $(4$, $12$, $18$, $20$, $16$, $14$, $2$, $10$, $22$, $8$, $6)$ | $(102$, $116$, $88$, $130$, $74$, $138$, $80$, $124$, $94$, $110$, $108$, $96$, $122$, $82$, $136$, $72$, $132$, $86$, $118$, $100$, $104$, $114$, $90$, $128$, $76$, $140$, $78$, $126$, $92$, $112$, $106$, $98$, $120$, $84$, $134$, $-160$, $-200$, $-182$, $-142$, $-178$, $-204$, $-164$, $-156$, $-196$, $-186$, $-146$, $-174$, $-208$, $-168$, $-152$, $-192$, $-190$, $-150$, $-170$, $-210$, $-172$, $-148$, $-188$, $-194$, $-154$, $-166$, $-206$, $-176$, $-144$, $-184$, $-198$, $-158$, $-162$, $-202$, $-180$, $242$, $256$, $228$, $270$, $214$, $278$, $220$, $264$, $234$, $250$, $248$, $236$, $262$, $222$, $276$, $212$, $272$, $226$, $258$, $240$, $244$, $254$, $230$, $268$, $216$, $280$, $218$, $266$, $232$, $252$, $246$, $238$, $260$, $224$, $274$, $-20$, $-60$, $-42$, $-2$, $-38$, $-64$, $-24$, $-16$, $-56$, $-46$, $-6$, $-34$, $-68$, $-28$, $-12$, $-52$, $-50$, $-10$, $-30$, $-70$, $-32$, $-8$, $-48$, $-54$, $-14$, $-26$, $-66$, $-36$, $-4$, $-44$, $-58$, $-18$, $-22$, $-62$, $-40)$\\

\noindent$11a227$ : $3$ | $(4$, $12$, $18$, $22$, $16$, $2$, $8$, $20$, $10$, $6$, $14)$ | $(110$, $44$, $206$, $202$, $48$, $114$, $106$, $40$, $102$, $118$, $52$, $198$, $210$, $172$, $208$, $200$, $50$, $116$, $104$, $-160$, $-78$, $-60$, $-66$, $-72$, $-154$, $-166$, $-84$, $196$, $212$, $174$, $184$, $222$, $186$, $92$, $128$, $90$, $188$, $220$, $182$, $176$, $214$, $194$, $-22$, $-32$, $-8$, $-130$, $-4$, $-28$, $-26$, $-58$, $-80$, $-162$, $-158$, $-76$, $-62$, $-64$, $-74$, $-156$, $-164$, $-82$, $-56$, $-24$, $-30$, $-6$, $204$, $46$, $112$, $108$, $42$, $100$, $120$, $54$, $122$, $98$, $94$, $126$, $88$, $190$, $218$, $180$, $178$, $216$, $192$, $86$, $124$, $96$, $-148$, $-138$, $-16$, $-38$, $-14$, $-136$, $-2$, $-132$, $-10$, $-34$, $-20$, $-142$, $-144$, $-168$, $-152$, $-70$, $-68$, $-150$, $-170$, $-146$, $-140$, $-18$, $-36$, $-12$, $-134)$\\

\noindent$11a228$ : $3$ | $(4$, $12$, $18$, $22$, $20$, $14$, $2$, $8$, $6$, $10$, $16)$ | $(64$, $74$, $54$, $60$, $68$, $70$, $58$, $56$, $72$, $66$, $62$, $114$, $110$, $-82$, $-22$, $-6$, $-14$, $122$, $102$, $126$, $30$, $78$, $26$, $80$, $28$, $76$, $-86$, $-96$, $-34$, $-98$, $-40$, $-90$, $-92$, $-38$, $-100$, $-36$, $-94$, $-88$, $-84$, $-24$, $112$, $108$, $116$, $118$, $106$, $128$, $104$, $120$, $32$, $124$, $-10$, $-18$, $-2$, $-50$, $-44$, $-46$, $-48$, $-42$, $-52$, $-4$, $-16$, $-12$, $-8$, $-20)$\\

\noindent$11a229$ : $2$ | $(4$, $12$, $20$, $18$, $14$, $2$, $10$, $22$, $8$, $6$, $16)$ | $(110$, $92$, $128$, $74$, $136$, $84$, $118$, $102$, $100$, $120$, $82$, $138$, $76$, $126$, $94$, $108$, $112$, $90$, $130$, $72$, $134$, $86$, $116$, $104$, $98$, $122$, $80$, $140$, $78$, $124$, $96$, $106$, $114$, $88$, $132$, $-156$, $-188$, $-202$, $-170$, $-142$, $-174$, $-206$, $-184$, $-152$, $-160$, $-192$, $-198$, $-166$, $-146$, $-178$, $-210$, $-180$, $-148$, $-164$, $-196$, $-194$, $-162$, $-150$, $-182$, $-208$, $-176$, $-144$, $-168$, $-200$, $-190$, $-158$, $-154$, $-186$, $-204$, $-172$, $250$, $232$, $268$, $214$, $276$, $224$, $258$, $242$, $240$, $260$, $222$, $278$, $216$, $266$, $234$, $248$, $252$, $230$, $270$, $212$, $274$, $226$, $256$, $244$, $238$, $262$, $220$, $280$, $218$, $264$, $236$, $246$, $254$, $228$, $272$, $-16$, $-48$, $-62$, $-30$, $-2$, $-34$, $-66$, $-44$, $-12$, $-20$, $-52$, $-58$, $-26$, $-6$, $-38$, $-70$, $-40$, $-8$, $-24$, $-56$, $-54$, $-22$, $-10$, $-42$, $-68$, $-36$, $-4$, $-28$, $-60$, $-50$, $-18$, $-14$, $-46$, $-64$, $-32)$\\

\noindent$11a230$ : $2$ | $(4$, $12$, $20$, $18$, $16$, $2$, $22$, $10$, $8$, $6$, $14)$ | $(82$, $56$, $94$, $70$, $68$, $96$, $58$, $80$, $84$, $54$, $92$, $72$, $66$, $98$, $60$, $78$, $86$, $52$, $90$, $74$, $64$, $100$, $62$, $76$, $88$, $-108$, $-124$, $-140$, $-146$, $-130$, $-114$, $-102$, $-118$, $-134$, $-150$, $-136$, $-120$, $-104$, $-112$, $-128$, $-144$, $-142$, $-126$, $-110$, $-106$, $-122$, $-138$, $-148$, $-132$, $-116$, $182$, $156$, $194$, $170$, $168$, $196$, $158$, $180$, $184$, $154$, $192$, $172$, $166$, $198$, $160$, $178$, $186$, $152$, $190$, $174$, $164$, $200$, $162$, $176$, $188$, $-8$, $-24$, $-40$, $-46$, $-30$, $-14$, $-2$, $-18$, $-34$, $-50$, $-36$, $-20$, $-4$, $-12$, $-28$, $-44$, $-42$, $-26$, $-10$, $-6$, $-22$, $-38$, $-48$, $-32$, $-16)$\\

\noindent$11a231$ : $4$ | $(4$, $14$, $10$, $12$, $18$, $6$, $2$, $20$, $22$, $8$, $16)$ | $(10$, $86$, $-36$, $-48$, $-42$, $-40$, $-50$, $-38$, $4$, $82$, $-32$, $-44$, $-46$, $-34$, $-18$, $84$, $80$, $58$, $88$, $52$, $92$, $54$, $78$, $56$, $90$, $-72$, $-2$, $-76$, $-68$, $28$, $20$, $30$, $22$, $26$, $6$, $14$, $12$, $8$, $24$, $-66$, $-60$, $-62$, $-64$, $-16$, $-70$, $-74)$\\

\noindent$11a232$ : $3$ | $(4$, $14$, $10$, $16$, $22$, $18$, $2$, $20$, $6$, $12$, $8)$ | $(94$, $78$, $104$, $84$, $134$, $90$, $98$, $74$, $100$, $88$, $132$, $86$, $102$, $76$, $96$, $92$, $136$, $82$, $106$, $80$, $138$, $-32$, $-10$, $-20$, $-22$, $-8$, $-34$, $-40$, $140$, $166$, $142$, $150$, $158$, $50$, $46$, $162$, $146$, $-126$, $-62$, $-114$, $-108$, $-118$, $-66$, $-130$, $-70$, $-122$, $-58$, $-60$, $-124$, $-72$, $-128$, $-64$, $-116$, $154$, $54$, $42$, $56$, $152$, $156$, $52$, $44$, $164$, $144$, $148$, $160$, $48$, $-68$, $-120$, $-110$, $-112$, $-2$, $-28$, $-14$, $-16$, $-26$, $-4$, $-38$, $-36$, $-6$, $-24$, $-18$, $-12$, $-30)$\\

\noindent$11a233$ : $3$ | $(4$, $14$, $10$, $20$, $22$, $18$, $2$, $8$, $6$, $12$, $16)$ | $(76$, $36$, $30$, $134$, $136$, $28$, $38$, $74$, $78$, $60$, $80$, $72$, $-108$, $-104$, $-88$, $-98$, $-114$, $-96$, $-90$, $-106$, $26$, $138$, $132$, $32$, $34$, $130$, $140$, $24$, $142$, $128$, $-50$, $-46$, $-6$, $-10$, $-42$, $-54$, $-56$, $-92$, $-94$, $-58$, $-52$, $-44$, $-8$, $66$, $118$, $124$, $144$, $126$, $116$, $64$, $84$, $68$, $120$, $122$, $70$, $82$, $62$, $-48$, $-4$, $-12$, $-40$, $-14$, $-2$, $-18$, $-20$, $-112$, $-100$, $-86$, $-102$, $-110$, $-22$, $-16)$\\

\noindent$11a234$ : $2$ | $(4$, $14$, $16$, $18$, $20$, $22$, $2$, $12$, $6$, $8$, $10)$ | $(58$, $44$, $72$, $42$, $60$, $56$, $46$, $70$, $40$, $62$, $54$, $48$, $68$, $38$, $64$, $52$, $50$, $66$, $-104$, $-94$, $-84$, $-74$, $-80$, $-90$, $-100$, $-108$, $-98$, $-88$, $-78$, $-76$, $-86$, $-96$, $-106$, $-102$, $-92$, $-82$, $130$, $116$, $144$, $114$, $132$, $128$, $118$, $142$, $112$, $134$, $126$, $120$, $140$, $110$, $136$, $124$, $122$, $138$, $-32$, $-22$, $-12$, $-2$, $-8$, $-18$, $-28$, $-36$, $-26$, $-16$, $-6$, $-4$, $-14$, $-24$, $-34$, $-30$, $-20$, $-10)$\\

\noindent$11a235$ : $2$ | $(4$, $14$, $16$, $18$, $20$, $22$, $2$, $12$, $6$, $10$, $8)$ | $(112$, $86$, $138$, $80$, $118$, $106$, $92$, $132$, $74$, $124$, $100$, $98$, $126$, $72$, $130$, $94$, $104$, $120$, $78$, $136$, $88$, $110$, $114$, $84$, $140$, $82$, $116$, $108$, $90$, $134$, $76$, $122$, $102$, $96$, $128$, $-162$, $-206$, $-172$, $-152$, $-196$, $-182$, $-142$, $-186$, $-192$, $-148$, $-176$, $-202$, $-158$, $-166$, $-210$, $-168$, $-156$, $-200$, $-178$, $-146$, $-190$, $-188$, $-144$, $-180$, $-198$, $-154$, $-170$, $-208$, $-164$, $-160$, $-204$, $-174$, $-150$, $-194$, $-184$, $252$, $226$, $278$, $220$, $258$, $246$, $232$, $272$, $214$, $264$, $240$, $238$, $266$, $212$, $270$, $234$, $244$, $260$, $218$, $276$, $228$, $250$, $254$, $224$, $280$, $222$, $256$, $248$, $230$, $274$, $216$, $262$, $242$, $236$, $268$, $-22$, $-66$, $-32$, $-12$, $-56$, $-42$, $-2$, $-46$, $-52$, $-8$, $-36$, $-62$, $-18$, $-26$, $-70$, $-28$, $-16$, $-60$, $-38$, $-6$, $-50$, $-48$, $-4$, $-40$, $-58$, $-14$, $-30$, $-68$, $-24$, $-20$, $-64$, $-34$, $-10$, $-54$, $-44)$\\

\noindent$11a236$ : $2$ | $(4$, $14$, $16$, $18$, $20$, $22$, $2$, $12$, $10$, $8$, $6)$ | $(156$, $122$, $190$, $108$, $170$, $142$, $136$, $176$, $102$, $184$, $128$, $150$, $162$, $116$, $196$, $114$, $164$, $148$, $130$, $182$, $100$, $178$, $134$, $144$, $168$, $110$, $192$, $120$, $158$, $154$, $124$, $188$, $106$, $172$, $140$, $138$, $174$, $104$, $186$, $126$, $152$, $160$, $118$, $194$, $112$, $166$, $146$, $132$, $180$, $-266$, $-208$, $-242$, $-290$, $-232$, $-218$, $-276$, $-256$, $-198$, $-252$, $-280$, $-222$, $-228$, $-286$, $-246$, $-204$, $-262$, $-270$, $-212$, $-238$, $-294$, $-236$, $-214$, $-272$, $-260$, $-202$, $-248$, $-284$, $-226$, $-224$, $-282$, $-250$, $-200$, $-258$, $-274$, $-216$, $-234$, $-292$, $-240$, $-210$, $-268$, $-264$, $-206$, $-244$, $-288$, $-230$, $-220$, $-278$, $-254$, $352$, $318$, $386$, $304$, $366$, $338$, $332$, $372$, $298$, $380$, $324$, $346$, $358$, $312$, $392$, $310$, $360$, $344$, $326$, $378$, $296$, $374$, $330$, $340$, $364$, $306$, $388$, $316$, $354$, $350$, $320$, $384$, $302$, $368$, $336$, $334$, $370$, $300$, $382$, $322$, $348$, $356$, $314$, $390$, $308$, $362$, $342$, $328$, $376$, $-70$, $-12$, $-46$, $-94$, $-36$, $-22$, $-80$, $-60$, $-2$, $-56$, $-84$, $-26$, $-32$, $-90$, $-50$, $-8$, $-66$, $-74$, $-16$, $-42$, $-98$, $-40$, $-18$, $-76$, $-64$, $-6$, $-52$, $-88$, $-30$, $-28$, $-86$, $-54$, $-4$, $-62$, $-78$, $-20$, $-38$, $-96$, $-44$, $-14$, $-72$, $-68$, $-10$, $-48$, $-92$, $-34$, $-24$, $-82$, $-58)$\\

\noindent$11a237$ : $3$ | $(4$, $14$, $16$, $18$, $22$, $20$, $2$, $12$, $8$, $6$, $10)$ | $(72$, $64$, $56$, $66$, $70$, $74$, $62$, $58$, $68$, $-8$, $-2$, $-12$, $-36$, $-16$, $-6$, $-4$, $-14$, $60$, $88$, $86$, $90$, $18$, $28$, $32$, $22$, $92$, $20$, $30$, $-46$, $-82$, $-80$, $-48$, $-54$, $-44$, $-42$, $-52$, $-50$, $-40$, $-38$, $24$, $34$, $26$, $76$, $-78$, $-84$, $-10)$\\

\noindent$11a238$ : $2$ | $(4$, $14$, $16$, $20$, $18$, $22$, $2$, $12$, $10$, $8$, $6)$ | $(102$, $80$, $124$, $70$, $112$, $92$, $90$, $114$, $68$, $122$, $82$, $100$, $104$, $78$, $126$, $72$, $110$, $94$, $88$, $116$, $66$, $120$, $84$, $98$, $106$, $76$, $128$, $74$, $108$, $96$, $86$, $118$, $-140$, $-164$, $-188$, $-174$, $-150$, $-130$, $-154$, $-178$, $-184$, $-160$, $-136$, $-144$, $-168$, $-192$, $-170$, $-146$, $-134$, $-158$, $-182$, $-180$, $-156$, $-132$, $-148$, $-172$, $-190$, $-166$, $-142$, $-138$, $-162$, $-186$, $-176$, $-152$, $230$, $208$, $252$, $198$, $240$, $220$, $218$, $242$, $196$, $250$, $210$, $228$, $232$, $206$, $254$, $200$, $238$, $222$, $216$, $244$, $194$, $248$, $212$, $226$, $234$, $204$, $256$, $202$, $236$, $224$, $214$, $246$, $-12$, $-36$, $-60$, $-46$, $-22$, $-2$, $-26$, $-50$, $-56$, $-32$, $-8$, $-16$, $-40$, $-64$, $-42$, $-18$, $-6$, $-30$, $-54$, $-52$, $-28$, $-4$, $-20$, $-44$, $-62$, $-38$, $-14$, $-10$, $-34$, $-58$, $-48$, $-24)$\\

\noindent$11a239$ : $3$ | $(4$, $14$, $18$, $20$, $22$, $6$, $2$, $8$, $10$, $12$, $16)$ | $(124$, $340$, $94$, $358$, $370$, $206$, $198$, $378$, $176$, $228$, $224$, $180$, $382$, $194$, $210$, $366$, $362$, $98$, $336$, $120$, $128$, $344$, $90$, $354$, $138$, $110$, $326$, $108$, $140$, $356$, $92$, $342$, $126$, $122$, $338$, $96$, $360$, $368$, $208$, $196$, $380$, $178$, $226$, $-20$, $-152$, $-318$, $-268$, $-278$, $-308$, $-162$, $-30$, $-76$, $-10$, $-142$, $-324$, $-146$, $-14$, $-80$, $-26$, $-158$, $-312$, $-274$, $-272$, $-314$, $-156$, $-24$, $-82$, $-16$, $-148$, $-322$, $372$, $204$, $200$, $376$, $174$, $230$, $222$, $182$, $384$, $192$, $212$, $364$, $-32$, $-74$, $-8$, $-58$, $-48$, $-256$, $-290$, $-296$, $-250$, $-42$, $-64$, $-2$, $-68$, $-38$, $-246$, $-300$, $-286$, $-260$, $-52$, $-54$, $-262$, $-284$, $-302$, $-244$, $-36$, $-70$, $-4$, $-62$, $-44$, $-252$, $-294$, $-292$, $-254$, $-46$, $-60$, $-6$, $-72$, $-34$, $164$, $100$, $334$, $118$, $130$, $346$, $88$, $352$, $136$, $112$, $328$, $106$, $170$, $234$, $218$, $186$, $388$, $188$, $216$, $236$, $168$, $104$, $330$, $114$, $134$, $350$, $86$, $348$, $132$, $116$, $332$, $102$, $166$, $238$, $214$, $190$, $386$, $184$, $220$, $232$, $172$, $374$, $202$, $-144$, $-12$, $-78$, $-28$, $-160$, $-310$, $-276$, $-270$, $-316$, $-154$, $-22$, $-84$, $-18$, $-150$, $-320$, $-266$, $-280$, $-306$, $-240$, $-242$, $-304$, $-282$, $-264$, $-56$, $-50$, $-258$, $-288$, $-298$, $-248$, $-40$, $-66)$\\

\noindent$11a240$ : $3$ | $(4$, $14$, $18$, $20$, $22$, $16$, $2$, $10$, $12$, $6$, $8)$ |  $(32$, $26$, $30$, $28$, $-4$, $-44$, $-40$, $38$, $12$, $36$, $10$, $34$, $-48$, $-54$, $-56$, $-52$, $-50$, $-46$, $-42$, $66$, $68$, $64$, $70$, $62$, $60$, $72$, $58$, $8$, $-6$, $-2$, $-24$, $-20$, $-16$, $-14$, $-18$, $-22)$\\

\noindent$11a241$ : $3$ | $(4$, $14$, $18$, $20$, $22$, $16$, $2$, $10$, $12$, $8$, $6)$ | $(30$, $24$, $14$, $28$, $32$, $26$, $-36$, $-44$, $-50$, $-54$, $-48$, $-46$, $-38$, $-34$, $-42$, $-40$, $86$, $62$, $60$, $84$, $88$, $64$, $58$, $70$, $56$, $66$, $68$, $-78$, $-72$, $-4$, $-12$, $-6$, $-74$, $-80$, $-76$, $22$, $16$, $82$, $20$, $18$, $-52$, $-2$, $-10$, $-8)$ \\

\noindent$11a242$ : $2$ | $(4$, $14$, $18$, $20$, $22$, $16$, $2$, $12$, $10$, $6$, $8)$ | $(72$, $62$, $82$, $52$, $92$, $50$, $84$, $60$, $74$, $70$, $64$, $80$, $54$, $90$, $48$, $86$, $58$, $76$, $68$, $66$, $78$, $56$, $88$, $-130$, $-112$, $-94$, $-108$, $-126$, $-134$, $-116$, $-98$, $-104$, $-122$, $-138$, $-120$, $-102$, $-100$, $-118$, $-136$, $-124$, $-106$, $-96$, $-114$, $-132$, $-128$, $-110$, $164$, $154$, $174$, $144$, $184$, $142$, $176$, $152$, $166$, $162$, $156$, $172$, $146$, $182$, $140$, $178$, $150$, $168$, $160$, $158$, $170$, $148$, $180$, $-38$, $-20$, $-2$, $-16$, $-34$, $-42$, $-24$, $-6$, $-12$, $-30$, $-46$, $-28$, $-10$, $-8$, $-26$, $-44$, $-32$, $-14$, $-4$, $-22$, $-40$, $-36$, $-18)$\\

\noindent$11a243$ : $2$ | $(4$, $14$, $18$, $20$, $22$, $16$, $2$, $12$, $10$, $8$, $6)$ | $(106$, $92$, $120$, $78$, $134$, $72$, $126$, $86$, $112$, $100$, $98$, $114$, $84$, $128$, $70$, $132$, $80$, $118$, $94$, $104$, $108$, $90$, $122$, $76$, $136$, $74$, $124$, $88$, $110$, $102$, $96$, $116$, $82$, $130$, $-156$, $-196$, $-174$, $-138$, $-178$, $-192$, $-152$, $-160$, $-200$, $-170$, $-142$, $-182$, $-188$, $-148$, $-164$, $-204$, $-166$, $-146$, $-186$, $-184$, $-144$, $-168$, $-202$, $-162$, $-150$, $-190$, $-180$, $-140$, $-172$, $-198$, $-158$, $-154$, $-194$, $-176$, $242$, $228$, $256$, $214$, $270$, $208$, $262$, $222$, $248$, $236$, $234$, $250$, $220$, $264$, $206$, $268$, $216$, $254$, $230$, $240$, $244$, $226$, $258$, $212$, $272$, $210$, $260$, $224$, $246$, $238$, $232$, $252$, $218$, $266$, $-20$, $-60$, $-38$, $-2$, $-42$, $-56$, $-16$, $-24$, $-64$, $-34$, $-6$, $-46$, $-52$, $-12$, $-28$, $-68$, $-30$, $-10$, $-50$, $-48$, $-8$, $-32$, $-66$, $-26$, $-14$, $-54$, $-44$, $-4$, $-36$, $-62$, $-22$, $-18$, $-58$, $-40)$ \\

\noindent$11a244$ : $3$ | $(4$, $14$, $18$, $22$, $16$, $20$, $2$, $8$, $12$, $6$, $10)$ | $(90$, $74$, $100$, $104$, $96$, $70$, $94$, $86$, $78$, $-62$, $-38$, $-32$, $-8$, $-106$, $-6$, $102$, $156$, $130$, $154$, $148$, $136$, $162$, $138$, $146$, $28$, $26$, $144$, $140$, $160$, $134$, $150$, $152$, $132$, $158$, $142$, $-54$, $-46$, $-40$, $-60$, $-64$, $-36$, $-34$, $-66$, $-58$, $-42$, $-44$, $-56$, $-68$, $-52$, $-48$, $-126$, $-124$, $-50$, $98$, $72$, $92$, $88$, $76$, $20$, $80$, $84$, $24$, $30$, $22$, $82$, $-120$, $-16$, $-114$, $-2$, $-110$, $-12$, $-10$, $-108$, $-4$, $-116$, $-18$, $-118$, $-128$, $-122$, $-14$, $-112)$\\

\noindent$11a245$ : $3$ | $(4$, $14$, $20$, $22$, $16$, $18$, $2$, $12$, $10$, $8$, $6)$ | $(50$, $44$, $28$, $34$, $32$, $26$, $42$, $52$, $48$, $46$, $54$, $40$, $-82$, $-66$, $-76$, $-80$, $-64$, $-78$, $30$, $96$, $-84$, $-68$, $-74$, $-88$, $-72$, $-70$, $-86$, $-38$, $-10$, $-6$, $-36$, $-8$, $58$, $100$, $92$, $106$, $90$, $102$, $60$, $62$, $56$, $98$, $94$, $104$, $-18$, $-2$, $-14$, $-22$, $-24$, $-12$, $-4$, $-20$, $-16)$\\

\noindent$11a246$ : $2$ | $(4$, $14$, $20$, $22$, $18$, $16$, $2$, $12$, $10$, $8$, $6)$ | $(62$, $56$, $68$, $50$, $74$, $44$, $80$, $42$, $76$, $48$, $70$, $54$, $64$, $60$, $58$, $66$, $52$, $72$, $46$, $78$, $-108$, $-82$, $-104$, $-112$, $-86$, $-100$, $-116$, $-90$, $-96$, $-120$, $-94$, $-92$, $-118$, $-98$, $-88$, $-114$, $-102$, $-84$, $-110$, $-106$, $142$, $136$, $148$, $130$, $154$, $124$, $160$, $122$, $156$, $128$, $150$, $134$, $144$, $140$, $138$, $146$, $132$, $152$, $126$, $158$, $-28$, $-2$, $-24$, $-32$, $-6$, $-20$, $-36$, $-10$, $-16$, $-40$, $-14$, $-12$, $-38$, $-18$, $-8$, $-34$, $-22$, $-4$, $-30$, $-26)$\\

\noindent$11a247$ : $2$ | $(4$, $14$, $22$, $20$, $18$, $16$, $2$, $12$, $10$, $8$, $6)$ | $(20$, $24$, $28$, $32$, $36$, $34$, $30$, $26$, $22$, $-46$, $-44$, $-48$, $-42$, $-50$, $-40$, $-52$, $-38$, $-54$, $56$, $60$, $64$, $68$, $72$, $70$, $66$, $62$, $58$, $-10$, $-8$, $-12$, $-6$, $-14$, $-4$, $-16$, $-2$, $-18)$\\

\noindent$11a248$ : $3$ | $(6$, $8$, $10$, $18$, $16$, $20$, $22$, $4$, $2$, $14$, $12)$ | $(218$, $272$, $288$, $294$, $112$, $98$, $308$, $84$, $126$, $82$, $306$, $100$, $110$, $296$, $76$, $300$, $106$, $104$, $302$, $78$, $122$, $88$, $312$, $94$, $116$, $290$, $292$, $114$, $96$, $310$, $86$, $124$, $80$, $304$, $102$, $108$, $298$, $-172$, $-254$, $-182$, $-162$, $-128$, $-158$, $-186$, $-250$, $-2$, $-246$, $-190$, $-154$, $-132$, $-166$, $-178$, $-258$, $-176$, $-168$, $-134$, $-152$, $-192$, $-244$, $-4$, $-252$, $-184$, $-160$, $280$, $210$, $226$, $264$, $238$, $198$, $120$, $90$, $314$, $92$, $118$, $196$, $240$, $266$, $224$, $212$, $278$, $282$, $208$, $228$, $262$, $236$, $200$, $202$, $234$, $260$, $230$, $206$, $284$, $276$, $214$, $222$, $268$, $242$, $-52$, $-6$, $-28$, $-74$, $-30$, $-38$, $-66$, $-20$, $-14$, $-60$, $-44$, $-142$, $-144$, $-46$, $-58$, $-12$, $-22$, $-68$, $-36$, $-32$, $-72$, $-26$, $-8$, $-54$, $-50$, $270$, $220$, $216$, $274$, $286$, $204$, $232$, $-34$, $-70$, $-24$, $-10$, $-56$, $-48$, $-146$, $-140$, $-42$, $-62$, $-16$, $-18$, $-64$, $-40$, $-138$, $-148$, $-194$, $-150$, $-136$, $-170$, $-174$, $-256$, $-180$, $-164$, $-130$, $-156$, $-188$, $-248)$\\

\noindent$11a249$ : $3$ | $(6$, $8$, $10$, $18$, $16$, $22$, $20$, $4$, $2$, $14$, $12)$ | $(244$, $198$, $272$, $216$, $226$, $262$, $106$, $358$, $84$, $334$, $330$, $88$, $362$, $102$, $316$, $348$, $74$, $344$, $320$, $98$, $366$, $92$, $326$, $338$, $80$, $354$, $108$, $356$, $82$, $336$, $328$, $90$, $364$, $100$, $318$, $346$, $-128$, $-168$, $-8$, $-152$, $-112$, $-144$, $-184$, $-178$, $-138$, $-118$, $-158$, $-2$, $-162$, $-122$, $-134$, $-174$, $-188$, $-148$, $258$, $230$, $212$, $276$, $202$, $240$, $248$, $194$, $268$, $220$, $222$, $266$, $192$, $250$, $238$, $204$, $278$, $210$, $232$, $256$, $260$, $228$, $214$, $274$, $200$, $242$, $246$, $196$, $270$, $218$, $224$, $264$, $190$, $252$, $236$, $206$, $280$, $208$, $234$, $254$, $-288$, $-28$, $-68$, $-38$, $-298$, $-304$, $-44$, $-62$, $-22$, $-282$, $-18$, $-58$, $-48$, $-308$, $-294$, $-34$, $-72$, $-32$, $-292$, $-310$, $-50$, $-56$, $-16$, $-284$, $-24$, $-64$, $-42$, $-302$, $-300$, $-40$, $-66$, $-26$, $-286$, $-14$, $-54$, $-52$, $-312$, $-290$, $-30$, $-70$, $-36$, $-296$, $-306$, $-46$, $-60$, $-20$, $332$, $86$, $360$, $104$, $314$, $350$, $76$, $342$, $322$, $96$, $368$, $94$, $324$, $340$, $78$, $352$, $-12$, $-172$, $-132$, $-124$, $-164$, $-4$, $-156$, $-116$, $-140$, $-180$, $-182$, $-142$, $-114$, $-154$, $-6$, $-166$, $-126$, $-130$, $-170$, $-10$, $-150$, $-110$, $-146$, $-186$, $-176$, $-136$, $-120$, $-160)$\\

\noindent$11a250$ : $3$ | $(6$, $8$, $12$, $2$, $14$, $18$, $4$, $20$, $22$, $10$, $16)$ | $(42$, $18$, $162$, $186$, $192$, $168$, $-74$, $-58$, $-84$, $-90$, $-64$, $-68$, $-94$, $-6$, $-126$, $-122$, $-10$, $-98$, $-72$, $-60$, $-86$, $-88$, $-62$, $-70$, $-96$, $-8$, $-124$, $30$, $104$, $174$, $198$, $180$, $156$, $158$, $182$, $196$, $172$, $102$, $32$, $52$, $28$, $106$, $176$, $200$, $178$, $154$, $160$, $184$, $194$, $170$, $-12$, $-120$, $-128$, $-4$, $-92$, $-66$, $22$, $46$, $38$, $14$, $166$, $190$, $188$, $164$, $16$, $40$, $44$, $20$, $110$, $24$, $48$, $36$, $100$, $34$, $50$, $26$, $108$, $-136$, $-112$, $-138$, $-152$, $-82$, $-56$, $-76$, $-146$, $-144$, $-118$, $-130$, $-2$, $-134$, $-114$, $-140$, $-150$, $-80$, $-54$, $-78$, $-148$, $-142$, $-116$, $-132)$\\

\noindent$11a251$ : $3$ | $(6$, $8$, $12$, $2$, $16$, $18$, $4$, $20$, $22$, $14$, $10)$ | $(74$, $84$, $64$, $118$, $66$, $82$, $76$, $72$, $124$, $-102$, $-90$, $-6$, $-8$, $-92$, $-100$, $-16$, $-2$, $-12$, $-96$, $138$, $28$, $20$, $88$, $22$, $26$, $136$, $36$, $140$, $30$, $18$, $86$, $-106$, $-116$, $-50$, $-40$, $-56$, $-110$, $-112$, $-54$, $-38$, $-52$, $-114$, $-108$, $-104$, $128$, $120$, $68$, $80$, $78$, $70$, $122$, $126$, $130$, $132$, $32$, $142$, $34$, $134$, $24$, $-48$, $-42$, $-58$, $-60$, $-44$, $-46$, $-62$, $-4$, $-10$, $-94$, $-98$, $-14)$\\

\noindent$11a252$ : $3$ | $(6$, $8$, $12$, $2$, $18$, $16$, $4$, $20$, $22$, $14$, $10)$ | $(142$, $132$, $110$, $106$, $128$, $146$, $124$, $102$, $114$, $136$, $138$, $116$, $100$, $122$, $144$, $130$, $108$, $-6$, $-152$, $-16$, $-24$, $-160$, $-162$, $-26$, $-14$, $-150$, $-8$, $-32$, $-4$, $-154$, $-18$, $-22$, $-158$, $-164$, $46$, $42$, $64$, $186$, $60$, $38$, $50$, $196$, $194$, $52$, $36$, $58$, $188$, $66$, $44$, $-70$, $-180$, $-80$, $-88$, $-172$, $-170$, $-90$, $-78$, $-182$, $-72$, $-96$, $-68$, $-178$, $-82$, $-86$, $-174$, $-168$, $-92$, $-76$, $-184$, $-74$, $-94$, $-166$, $-176$, $-84$, $126$, $104$, $112$, $134$, $140$, $118$, $98$, $120$, $200$, $190$, $56$, $34$, $54$, $192$, $198$, $48$, $40$, $62$, $-20$, $-156$, $-2$, $-30$, $-10$, $-148$, $-12$, $-28)$\\

\noindent$11a253$ : $3$ | $(6$, $8$, $12$, $2$, $18$, $20$, $4$, $10$, $22$, $14$, $16)$ | $(46$, $26$, $140$, $154$, $174$, $168$, $148$, $146$, $166$, $176$, $156$, $138$, $-112$, $-118$, $-130$, $-16$, $-12$, $-126$, $-122$, $-8$, $-20$, $-134$, $-114$, $-116$, $-132$, $-18$, $-10$, $-124$, $36$, $98$, $28$, $48$, $44$, $90$, $92$, $42$, $50$, $30$, $100$, $34$, $54$, $38$, $96$, $-4$, $-24$, $-136$, $-22$, $-6$, $-120$, $-128$, $-14$, $32$, $52$, $40$, $94$, $88$, $160$, $180$, $162$, $142$, $152$, $172$, $170$, $150$, $144$, $164$, $178$, $158$, $-68$, $-110$, $-86$, $-60$, $-76$, $-102$, $-78$, $-58$, $-84$, $-108$, $-70$, $-66$, $-2$, $-62$, $-74$, $-104$, $-80$, $-56$, $-82$, $-106$, $-72$, $-64)$\\

\noindent$11a254$ : $3$ | $(6$, $8$, $12$, $2$, $20$, $18$, $4$, $10$, $22$, $14$, $16)$ | $(54$, $46$, $182$, $178$, $42$, $168$, $192$, $190$, $170$, $40$, $176$, $184$, $48$, $52$, $56$, $58$, $118$, $66$, $104$, $-154$, $-126$, $-144$, $-138$, $-132$, $-160$, $-162$, $-134$, $-136$, $-164$, $-78$, $-12$, $-16$, $-82$, $-88$, $-22$, $-6$, $180$, $44$, $166$, $194$, $188$, $172$, $38$, $174$, $186$, $50$, $-26$, $-2$, $-76$, $-10$, $-18$, $-84$, $-86$, $-20$, $-8$, $-74$, $-4$, $-24$, $-90$, $-80$, $-14$, $62$, $100$, $108$, $70$, $114$, $94$, $198$, $96$, $112$, $72$, $110$, $98$, $60$, $120$, $64$, $102$, $106$, $68$, $116$, $92$, $196$, $-30$, $-146$, $-124$, $-152$, $-36$, $-156$, $-128$, $-142$, $-140$, $-130$, $-158$, $-34$, $-150$, $-122$, $-148$, $-32$, $-28)$\\

\noindent$11a255$ : $3$ | $(6$, $8$, $12$, $20$, $16$, $18$, $22$, $4$, $10$, $2$, $14)$ | $(116$, $126$, $146$, $134$, $108$, $64$, $110$, $132$, $148$, $128$, $114$, $118$, $124$, $144$, $136$, $106$, $138$, $142$, $122$, $120$, $140$, $-34$, $-28$, $-6$, $-16$, $-150$, $-18$, $-4$, $-26$, $-36$, $-40$, $-86$, $180$, $172$, $48$, $54$, $166$, $186$, $164$, $56$, $46$, $174$, $178$, $182$, $170$, $50$, $52$, $168$, $184$, $162$, $58$, $44$, $176$, $-84$, $-66$, $-88$, $-96$, $-74$, $-76$, $-98$, $-160$, $-100$, $-78$, $-72$, $-94$, $-90$, $-68$, $-82$, $-104$, $-156$, $-158$, $-102$, $-80$, $-70$, $-92$, $130$, $112$, $62$, $42$, $60$, $-38$, $-24$, $-2$, $-20$, $-152$, $-14$, $-8$, $-30$, $-32$, $-10$, $-12$, $-154$, $-22)$\\

\noindent$11a256$ : $3$ | $(6$, $8$, $12$, $20$, $18$, $16$, $22$, $4$, $10$, $2$, $14)$ | $(34$, $44$, $98$, $40$, $38$, $96$, $86$, $90$, $92$, $84$, $82$, $94$, $88$, $-14$, $-18$, $-22$, $-68$, $-66$, $-58$, $-74$, $-60$, $-64$, $-70$, $-20$, $30$, $112$, $26$, $114$, $124$, $76$, $120$, $118$, $78$, $126$, $80$, $116$, $122$, $-8$, $-102$, $-2$, $-106$, $-48$, $-108$, $-54$, $-52$, $-110$, $-50$, $-56$, $-72$, $-62$, $42$, $36$, $32$, $46$, $28$, $-16$, $-24$, $-12$, $-4$, $-100$, $-6$, $-10$, $-104)$\\

\noindent$11a257$ : $3$ | $(6$, $8$, $14$, $2$, $16$, $18$, $20$, $4$, $22$, $12$, $10)$ | $(164$, $128$, $142$, $150$, $120$, $156$, $136$, $134$, $170$, $158$, $122$, $148$, $144$, $126$, $162$, $166$, $130$, $140$, $152$, $118$, $154$, $138$, $132$, $168$, $160$, $124$, $146$, $-24$, $-6$, $-36$, $-42$, $-12$, $210$, $212$, $196$, $188$, $56$, $182$, $202$, $218$, $204$, $180$, $58$, $190$, $194$, $62$, $176$, $208$, $214$, $198$, $186$, $54$, $184$, $200$, $216$, $206$, $178$, $60$, $192$, $-88$, $-70$, $-100$, $-106$, $-76$, $-82$, $-112$, $-94$, $-174$, $-92$, $-114$, $-84$, $-74$, $-104$, $-102$, $-72$, $-86$, $-116$, $-90$, $-68$, $-98$, $-108$, $-78$, $-80$, $-110$, $-96$, $-66$, $64$, $-14$, $-44$, $-34$, $-4$, $-26$, $-52$, $-22$, $-8$, $-38$, $-40$, $-10$, $-20$, $-50$, $-28$, $-2$, $-32$, $-46$, $-16$, $-172$, $-18$, $-48$, $-30)$\\

\noindent$11a258$ : $3$ | $(6$, $8$, $14$, $2$, $16$, $20$, $18$, $4$, $22$, $12$, $10)$ | $(50$, $18$, $28$, $30$, $20$, $-72$, $-62$, $-54$, $-64$, $-74$, $-70$, $-60$, $-56$, $-66$, $-76$, $-68$, $-58$, $24$, $14$, $16$, $26$, $32$, $22$, $12$, $100$, $-48$, $-38$, $92$, $82$, $84$, $94$, $52$, $90$, $80$, $86$, $96$, $98$, $88$, $78$, $-6$, $-42$, $-34$, $-44$, $-8$, $-4$, $-40$, $-36$, $-46$, $-10$, $-2)$\\

\noindent$11a259$ : $3$ | $(6$, $8$, $14$, $2$, $18$, $20$, $4$, $22$, $12$, $10$, $16)$ | $(60$, $86$, $78$, $52$, $68$, $94$, $70$, $44$, $208$, $130$, $156$, $148$, $216$, $138$, $164$, $140$, $218$, $146$, $158$, $132$, $210$, $-178$, $-120$, $-18$, $-40$, $-98$, $-112$, $-26$, $-32$, $-106$, $-104$, $-34$, $-24$, $-114$, $-96$, $-118$, $-20$, $-38$, $-100$, $-110$, $-28$, $-30$, $-108$, $-102$, $-36$, $-22$, $-116$, $152$, $212$, $134$, $160$, $144$, $220$, $142$, $162$, $136$, $214$, $150$, $154$, $128$, $206$, $-16$, $-42$, $-14$, $-202$, $-174$, $-182$, $-6$, $-194$, $-166$, $-190$, $-2$, $-186$, $-170$, $-198$, $-10$, $-12$, $-200$, $-172$, $-184$, $-4$, $-192$, $82$, $56$, $64$, $90$, $74$, $48$, $124$, $126$, $50$, $76$, $88$, $62$, $58$, $84$, $80$, $54$, $66$, $92$, $72$, $46$, $122$, $-204$, $-176$, $-180$, $-8$, $-196$, $-168$, $-188)$\\

\noindent$11a260$ : $3$ | $(6$, $8$, $14$, $2$, $20$, $18$, $4$, $22$, $12$, $10$, $16)$ | $(44$, $24$, $38$, $30$, $32$, $36$, $26$, $42$, $68$, $-66$, $-48$, $-62$, $-102$, $-52$, $-58$, $-106$, $-56$, $-54$, $-104$, $-60$, $-50$, $-100$, $-64$, $72$, $116$, $78$, $110$, $84$, $108$, $80$, $114$, $74$, $70$, $-96$, $-98$, $-94$, $-92$, $-2$, $-14$, $-88$, $-6$, $-10$, $34$, $28$, $40$, $22$, $46$, $20$, $18$, $118$, $76$, $112$, $82$, $-8$, $-86$, $-12$, $-4$, $-90$, $-16)$\\

\noindent$11a261$ : $3$ | $(6$, $8$, $14$, $16$, $4$, $18$, $20$, $2$, $22$, $12$, $10)$ | $(38$, $92$, $78$, $70$, $100$, $68$, $80$, $90$, $36$, $40$, $94$, $76$, $72$, $98$, $66$, $82$, $88$, $-126$, $-108$, $-110$, $-62$, $-58$, $-30$, $-12$, $-6$, $-24$, $74$, $96$, $64$, $84$, $86$, $-142$, $-132$, $-120$, $-102$, $-116$, $-136$, $-138$, $-114$, $-104$, $-122$, $-130$, $-144$, $-128$, $-124$, $-106$, $-112$, $-140$, $-134$, $-118$, $160$, $48$, $150$, $170$, $176$, $168$, $152$, $50$, $158$, $162$, $46$, $148$, $172$, $174$, $146$, $44$, $164$, $156$, $52$, $154$, $166$, $42$, $-60$, $-32$, $-14$, $-4$, $-22$, $-54$, $-26$, $-8$, $-10$, $-28$, $-56$, $-20$, $-2$, $-16$, $-34$, $-18)$\\

\noindent$11a262$ : $3$ | $(6$, $8$, $14$, $16$, $4$, $20$, $18$, $2$, $22$, $12$, $10)$ | $(188$, $152$, $162$, $178$, $136$, $180$, $160$, $154$, $186$, $142$, $172$, $168$, $146$, $150$, $164$, $176$, $138$, $182$, $158$, $156$, $184$, $140$, $174$, $166$, $148$, $-26$, $-190$, $-18$, $-42$, $-34$, $-10$, $-8$, $-32$, $-44$, $-20$, $-192$, $-24$, $-48$, $-28$, $-4$, $-14$, $-38$, $206$, $58$, $66$, $198$, $76$, $50$, $74$, $196$, $68$, $56$, $208$, $82$, $204$, $60$, $64$, $200$, $78$, $52$, $72$, $194$, $70$, $54$, $210$, $80$, $202$, $62$, $-112$, $-88$, $-102$, $-126$, $-120$, $-96$, $-94$, $-118$, $-130$, $-106$, $-84$, $-108$, $-132$, $-116$, $-92$, $-98$, $-122$, $-124$, $-100$, $-90$, $-114$, $-134$, $-110$, $-86$, $-104$, $-128$, $212$, $144$, $170$, $-22$, $-46$, $-30$, $-6$, $-12$, $-36$, $-40$, $-16$, $-2)$\\

\noindent$11a263$ : $4$ | $(6$, $8$, $16$, $2$, $20$, $22$, $18$, $4$, $12$, $14$, $10)$ | $(8$, $14$, $28$, $-34$, $-12$, $30$, $50$, $-60$, $-68$, $-66$, $-58$, $-64$, $-70$, $-62$, $-32$, $10$, $74$, $22$, $20$, $16$, $26$, $24$, $18$, $-38$, $-2$, $-4$, $-40$, $-44$, $-36$, $48$, $52$, $72$, $56$, $46$, $54$, $-42$, $-6)$\\

\noindent$11a264$ : $3$ | $(6$, $8$, $16$, $14$, $4$, $18$, $20$, $2$, $22$, $12$, $10)$ | $(58$, $200$, $218$, $228$, $190$, $236$, $210$, $208$, $238$, $192$, $226$, $220$, $198$, $56$, $60$, $202$, $216$, $230$, $188$, $234$, $212$, $206$, $240$, $194$, $224$, $222$, $196$, $-178$, $-150$, $-152$, $-180$, $-44$, $-16$, $-12$, $-40$, $-184$, $-34$, $-6$, $-22$, $-50$, $124$, $106$, $62$, $134$, $70$, $114$, $116$, $72$, $132$, $-30$, $-2$, $-26$, $-54$, $-46$, $-18$, $-10$, $-38$, $-186$, $-36$, $-8$, $-20$, $-48$, $-52$, $-24$, $-4$, $-32$, $-182$, $-42$, $-14$, $66$, $110$, $120$, $76$, $128$, $102$, $244$, $100$, $130$, $74$, $118$, $112$, $68$, $136$, $64$, $108$, $122$, $78$, $126$, $104$, $242$, $204$, $214$, $232$, $-92$, $-160$, $-142$, $-170$, $-82$, $-176$, $-148$, $-154$, $-98$, $-86$, $-166$, $-138$, $-164$, $-88$, $-96$, $-156$, $-146$, $-174$, $-80$, $-172$, $-144$, $-158$, $-94$, $-90$, $-162$, $-140$, $-168$, $-84$, $-28)$\\

\noindent$11a265$ : $3$ | $(6$, $8$, $16$, $14$, $4$, $20$, $18$, $2$, $22$, $12$, $10)$ | $(76$, $92$, $102$, $32$, $108$, $34$, $170$, $172$, $156$, $136$, $160$, $176$, $166$, $142$, $150$, $-124$, $-110$, $86$, $70$, $82$, $98$, $96$, $80$, $72$, $88$, $106$, $30$, $104$, $90$, $74$, $78$, $94$, $100$, $84$, $-14$, $-130$, $-24$, $-118$, $-116$, $-22$, $-132$, $-16$, $-12$, $-128$, $-26$, $-120$, $-114$, $-20$, $-134$, $-18$, $-112$, $-122$, $-28$, $-126$, $146$, $154$, $138$, $162$, $178$, $164$, $140$, $152$, $148$, $144$, $168$, $174$, $158$, $-42$, $-64$, $-4$, $-56$, $-50$, $-10$, $-8$, $-52$, $-54$, $-6$, $-66$, $-40$, $-44$, $-62$, $-2$, $-58$, $-48$, $-36$, $-68$, $-38$, $-46$, $-60)$\\

\noindent$11a266$ : $3$ | $(6$, $10$, $12$, $14$, $16$, $18$, $22$, $20$, $4$, $2$, $8)$ | $(74$, $46$, $70$, $40$, $32$, $30$, $42$, $68$, $44$, $28$, $34$, $38$, $-58$, $-52$, $-64$, $-86$, $-62$, $-54$, $-56$, $-60$, $-88$, $-66$, $-92$, $-98$, $72$, $76$, $82$, $106$, $80$, $78$, $108$, $84$, $112$, $122$, $-90$, $-100$, $-96$, $-94$, $-102$, $-24$, $-104$, $-20$, $-10$, $110$, $124$, $120$, $114$, $128$, $116$, $118$, $126$, $26$, $36$, $-22$, $-8$, $-12$, $-18$, $-2$, $-48$, $-4$, $-16$, $-14$, $-6$, $-50)$\\

\noindent$11a267$ : $3$ | $(6$, $10$, $12$, $14$, $20$, $18$, $22$, $8$, $4$, $2$, $16)$ | $(32$, $196$, $48$, $206$, $212$, $42$, $190$, $38$, $216$, $202$, $52$, $200$, $218$, $220$, $-98$, $-74$, $-72$, $-100$, $-58$, $-88$, $-84$, $-62$, $-104$, $-68$, $-78$, $-94$, $140$, $114$, $156$, $124$, $130$, $150$, $108$, $146$, $134$, $120$, $224$, $118$, $136$, $144$, $110$, $152$, $128$, $126$, $154$, $112$, $142$, $138$, $116$, $222$, $122$, $132$, $148$, $-168$, $-8$, $-184$, $-18$, $-24$, $-178$, $-2$, $-174$, $-162$, $-14$, $-188$, $-12$, $-164$, $-172$, $-4$, $-180$, $-22$, $-20$, $-182$, $-6$, $-170$, $-166$, $-10$, $-186$, $-16$, $-26$, $28$, $36$, $192$, $44$, $210$, $208$, $46$, $194$, $34$, $30$, $198$, $50$, $204$, $214$, $40$, $-86$, $-60$, $-102$, $-70$, $-76$, $-96$, $-54$, $-92$, $-80$, $-66$, $-106$, $-64$, $-82$, $-90$, $-56$, $-158$, $-160$, $-176)$\\

\noindent$11a268$ : $3$ | $(6$, $10$, $12$, $18$, $20$, $16$, $22$, $4$, $2$, $8$, $14)$ | $(86$, $98$, $108$, $96$, $78$, $94$, $106$, $100$, $88$, $84$, $-66$, $-50$, $-52$, $-2$, $-16$, $-118$, $-122$, $-12$, $-6$, $-4$, $-14$, $-120$, $32$, $142$, $134$, $128$, $136$, $140$, $34$, $42$, $30$, $144$, $132$, $130$, $146$, $28$, $40$, $36$, $138$, $-54$, $-48$, $-64$, $-68$, $-70$, $-62$, $-46$, $-56$, $-76$, $-110$, $-72$, $-60$, $-44$, $-58$, $-74$, $38$, $26$, $22$, $82$, $90$, $102$, $104$, $92$, $80$, $24$, $-112$, $-20$, $-114$, $-126$, $-8$, $-10$, $-124$, $-116$, $-18)$\\

\noindent$11a269$ : $3$ | $(6$, $10$, $12$, $20$, $16$, $18$, $22$, $8$, $4$, $2$, $14)$ | $(30$, $78$, $72$, $36$, $24$, $26$, $34$, $74$, $76$, $32$, $28$, $-92$, $-62$, $-94$, $-52$, $-54$, $-96$, $-60$, $-46$, $118$, $104$, $106$, $120$, $116$, $102$, $108$, $126$, $110$, $100$, $114$, $122$, $-10$, $-22$, $-44$, $-40$, $-48$, $-58$, $-98$, $-56$, $-50$, $-42$, $68$, $38$, $70$, $80$, $66$, $64$, $124$, $112$, $-16$, $-4$, $-88$, $-82$, $-90$, $-2$, $-18$, $-14$, $-6$, $-86$, $-84$, $-8$, $-12$, $-20)$\\

\noindent$11a270$ : $3$ | $(6$, $10$, $12$, $20$, $18$, $16$, $22$, $4$, $2$, $8$, $14)$ | $(32$, $34$, $186$, $172$, $194$, $42$, $200$, $178$, $180$, $202$, $40$, $192$, $170$, $188$, $36$, $-146$, $-80$, $-144$, $-54$, $-62$, $-136$, $-72$, $92$, $114$, $128$, $106$, $84$, $100$, $122$, $120$, $98$, $86$, $108$, $130$, $112$, $90$, $94$, $116$, $126$, $104$, $-2$, $-30$, $-164$, $-20$, $-12$, $-156$, $-154$, $-10$, $-22$, $-166$, $-28$, $-4$, $-148$, $-162$, $-18$, $-14$, $-158$, $-152$, $-8$, $-24$, $-168$, $-26$, $-6$, $-150$, $-160$, $-16$, $110$, $88$, $96$, $118$, $124$, $102$, $82$, $206$, $184$, $174$, $196$, $44$, $198$, $176$, $182$, $204$, $38$, $190$, $-58$, $-140$, $-76$, $-48$, $-68$, $-132$, $-66$, $-50$, $-78$, $-142$, $-56$, $-60$, $-138$, $-74$, $-46$, $-70$, $-134$, $-64$, $-52)$\\

\noindent$11a271$ : $3$ | $(6$, $10$, $12$, $22$, $16$, $18$, $8$, $20$, $4$, $2$, $14)$ | $(88$, $98$, $110$, $76$, $116$, $92$, $94$, $114$, $74$, $112$, $96$, $90$, $86$, $100$, $108$, $78$, $54$, $82$, $104$, $-60$, $-28$, $-4$, $-20$, $-14$, $-10$, $-118$, $-138$, $-128$, $80$, $106$, $102$, $84$, $52$, $50$, $160$, $46$, $178$, $-154$, $-140$, $-126$, $-130$, $-136$, $-120$, $-146$, $-148$, $-68$, $-158$, $-70$, $-150$, $-144$, $-122$, $-134$, $-132$, $-124$, $-142$, $-152$, $-72$, $-156$, $-66$, $190$, $166$, $40$, $172$, $184$, $186$, $170$, $38$, $168$, $188$, $182$, $174$, $42$, $164$, $192$, $162$, $44$, $176$, $180$, $48$, $-12$, $-22$, $-2$, $-26$, $-62$, $-36$, $-58$, $-30$, $-6$, $-18$, $-16$, $-8$, $-32$, $-56$, $-34$, $-64$, $-24)$\\

\noindent$11a272$ : $3$ | $(6$, $10$, $12$, $22$, $18$, $16$, $20$, $4$, $2$, $14$, $8)$ | $(22$, $80$, $18$, $72$, $66$, $62$, $76$, $-96$, $-106$, $-54$, $-50$, $-48$, $-56$, $-46$, $-42$, $-10$, $64$, $74$, $16$, $78$, $60$, $68$, $70$, $58$, $138$, $124$, $126$, $136$, $-86$, $-110$, $-92$, $-100$, $-102$, $-90$, $-112$, $-88$, $-104$, $-98$, $-94$, $-108$, $-52$, $20$, $24$, $26$, $114$, $30$, $120$, $130$, $132$, $118$, $32$, $116$, $134$, $128$, $122$, $28$, $-44$, $-8$, $-34$, $-12$, $-40$, $-2$, $-82$, $-4$, $-38$, $-14$, $-36$, $-6$, $-84)$ \\

\noindent$11a273$ : $3$ | $(6$, $10$, $14$, $16$, $20$, $4$, $18$, $2$, $22$, $12$, $8)$ | $(102$, $108$, $112$, $122$, $150$, $130$, $152$, $120$, $114$, $106$, $100$, $104$, $116$, $118$, $154$, $128$, $148$, $124$, $110$, $-4$, $-18$, $-24$, $-28$, $-14$, $-8$, $-34$, $-146$, $-32$, $-10$, $-12$, $-30$, $-144$, $-36$, $-6$, $-16$, $-26$, $44$, $62$, $54$, $52$, $64$, $38$, $66$, $50$, $56$, $60$, $46$, $70$, $42$, $160$, $-140$, $-142$, $-98$, $-138$, $-82$, $-84$, $-136$, $-96$, $-72$, $-94$, $-134$, $-86$, $-80$, $-76$, $-90$, $58$, $48$, $68$, $40$, $158$, $162$, $156$, $126$, $-78$, $-88$, $-132$, $-92$, $-74$, $-2$, $-20$, $-22)$\\

\noindent$11a274$ : $3$ | $(6$, $10$, $14$, $18$, $2$, $20$, $22$, $4$, $12$, $8$, $16)$ | $(48$, $28$, $44$, $158$, $162$, $40$, $32$, $52$, $34$, $38$, $164$, $156$, $46$, $-64$, $-138$, $-148$, $-54$, $-74$, $-72$, $-56$, $-146$, $-140$, $-62$, $-66$, $-136$, $-150$, $112$, $126$, $106$, $90$, $100$, $120$, $118$, $98$, $92$, $108$, $128$, $110$, $94$, $96$, $50$, $30$, $42$, $160$, $-8$, $-22$, $-132$, $-152$, $-134$, $-24$, $-6$, $-86$, $-10$, $-20$, $-130$, $-16$, $-14$, $-82$, $-2$, $-26$, $-4$, $-84$, $-12$, $-18$, $36$, $166$, $154$, $170$, $116$, $122$, $102$, $88$, $104$, $124$, $114$, $168$, $-78$, $-68$, $-60$, $-142$, $-144$, $-58$, $-70$, $-76$, $-80)$\\

\noindent$11a275$ : $3$ | $(6$, $10$, $14$, $20$, $2$, $18$, $4$, $22$, $12$, $8$, $16)$ | $(20$, $52$, $50$, $22$, $18$, $30$, $26$, $-42$, $-36$, $-74$, $-76$, $-34$, $-44$, $-40$, $-38$, $-46$, $56$, $94$, $88$, $60$, $90$, $92$, $58$, $54$, $-70$, $-78$, $-72$, $-68$, $-66$, $-12$, $28$, $16$, $24$, $32$, $82$, $80$, $86$, $96$, $84$, $-4$, $-48$, $-6$, $-2$, $-10$, $-62$, $-14$, $-64$, $-8)$\\

\noindent$11a276$ : $3$ | $(6$, $10$, $14$, $20$, $4$, $18$, $2$, $22$, $8$, $12$, $16)$ | $(208$, $176$, $162$, $222$, $194$, $190$, $48$, $36$, $76$, $306$, $62$, $-256$, $-282$, $-234$, $-4$, $-226$, $-274$, $-264$, $-248$, $-290$, $-242$, $-244$, $-292$, $-246$, $-266$, $-272$, $-224$, $-6$, $-236$, $-284$, $-254$, $-258$, $-280$, $-232$, $-2$, $-228$, $-276$, $-262$, $-250$, $-288$, $-240$, $-10$, $192$, $50$, $34$, $74$, $308$, $64$, $24$, $60$, $304$, $78$, $38$, $46$, $44$, $40$, $80$, $302$, $58$, $26$, $66$, $310$, $72$, $32$, $52$, $54$, $30$, $70$, $312$, $68$, $28$, $56$, $300$, $82$, $42$, $-110$, $-106$, $-14$, $-88$, $-128$, $-132$, $-142$, $-118$, $-98$, $-22$, $-96$, $-120$, $-140$, $-134$, $-126$, $-90$, $-16$, $-104$, $-112$, $-148$, $-268$, $-270$, $-146$, $-114$, $-102$, $-18$, $-92$, $-124$, $-136$, $-138$, $-122$, $-94$, $-20$, $-100$, $-116$, $-144$, $-130$, $160$, $178$, $206$, $210$, $174$, $164$, $220$, $196$, $188$, $150$, $294$, $154$, $184$, $200$, $216$, $168$, $170$, $214$, $202$, $182$, $156$, $296$, $84$, $298$, $158$, $180$, $204$, $212$, $172$, $166$, $218$, $198$, $186$, $152$, $-108$, $-12$, $-86$, $-8$, $-238$, $-286$, $-252$, $-260$, $-278$, $-230)$\\

\noindent$11a277$ : $3$ | $(6$, $10$, $16$, $12$, $4$, $18$, $20$, $22$, $2$, $8$, $14)$ | $(158$, $64$, $176$, $172$, $68$, $162$, $186$, $54$, $148$, $52$, $188$, $160$, $66$, $174$, $-6$, $-230$, $-206$, $-200$, $-224$, $-12$, $-22$, $-108$, $-106$, $-82$, $-86$, $-128$, $-146$, $-124$, $-90$, $-78$, $-102$, $-112$, $-134$, $-140$, $-118$, $-96$, $266$, $34$, $252$, $280$, $292$, $240$, $242$, $290$, $282$, $250$, $36$, $268$, $264$, $32$, $254$, $278$, $294$, $238$, $296$, $276$, $256$, $30$, $262$, $270$, $38$, $248$, $284$, $288$, $244$, $42$, $274$, $258$, $28$, $260$, $272$, $40$, $246$, $286$, $-126$, $-88$, $-80$, $-104$, $-110$, $-132$, $-142$, $-120$, $-94$, $-74$, $-98$, $-116$, $-138$, $-136$, $-114$, $-100$, $-76$, $-92$, $-122$, $-144$, $-130$, $-84$, $48$, $152$, $58$, $182$, $166$, $72$, $168$, $180$, $60$, $154$, $46$, $190$, $50$, $150$, $56$, $184$, $164$, $70$, $170$, $178$, $62$, $156$, $44$, $-214$, $-192$, $-216$, $-20$, $-14$, $-222$, $-198$, $-208$, $-232$, $-4$, $-26$, $-8$, $-228$, $-204$, $-202$, $-226$, $-10$, $-24$, $-2$, $-234$, $-210$, $-196$, $-220$, $-16$, $-18$, $-218$, $-194$, $-212$, $-236)$\\

\noindent$11a278$ : $3$ | $(6$, $10$, $16$, $12$, $4$, $18$, $22$, $20$, $2$, $8$, $14)$ | $(84$, $266$, $102$, $66$, $184$, $232$, $296$, $202$, $214$, $308$, $220$, $196$, $290$, $286$, $192$, $224$, $304$, $210$, $206$, $300$, $228$, $188$, $186$, $230$, $298$, $204$, $212$, $306$, $222$, $194$, $288$, $-26$, $-248$, $-234$, $-12$, $-48$, $-40$, $-254$, $-32$, $-56$, $-20$, $-242$, $-240$, $-18$, $-54$, $-34$, $-256$, $-38$, $-50$, $-14$, $-236$, $-246$, $-24$, $-60$, $-28$, $-250$, $-44$, $278$, $72$, $96$, $260$, $90$, $78$, $272$, $108$, $62$, $106$, $270$, $80$, $88$, $262$, $98$, $70$, $280$, $112$, $276$, $74$, $94$, $258$, $92$, $76$, $274$, $110$, $64$, $104$, $268$, $82$, $86$, $264$, $100$, $68$, $282$, $-148$, $-114$, $-150$, $-164$, $-128$, $-134$, $-170$, $-4$, $-178$, $-142$, $-120$, $-156$, $-158$, $-122$, $-140$, $-176$, $-2$, $-172$, $-136$, $-126$, $-162$, $-152$, $-116$, $-146$, $-182$, $294$, $200$, $216$, $310$, $218$, $198$, $292$, $284$, $190$, $226$, $302$, $208$, $-36$, $-52$, $-16$, $-238$, $-244$, $-22$, $-58$, $-30$, $-252$, $-42$, $-46$, $-10$, $-8$, $-166$, $-130$, $-132$, $-168$, $-6$, $-180$, $-144$, $-118$, $-154$, $-160$, $-124$, $-138$, $-174)$\\

\noindent$11a279$ : $3$ | $(6$, $10$, $16$, $12$, $18$, $2$, $20$, $22$, $4$, $8$, $14)$ | $(50$, $34$, $40$, $30$, $20$, $26$, $36$, $38$, $28$, $-80$, $-70$, $-62$, $-72$, $-82$, $-78$, $-68$, $-64$, $-74$, $-84$, $-76$, $-66$, $32$, $22$, $24$, $98$, $-48$, $-44$, $-8$, $-18$, $-10$, $-46$, $52$, $88$, $90$, $54$, $60$, $96$, $86$, $92$, $56$, $58$, $94$, $-14$, $-4$, $-42$, $-6$, $-16$, $-12$, $-2)$\\

\noindent$11a280$ : $3$ | $(6$, $10$, $16$, $12$, $18$, $2$, $22$, $20$, $4$, $8$, $14)$ | $(36$, $34$, $38$, $40$, $68$, $64$, $44$, $60$, $72$, $58$, $42$, $66$, $-10$, $-48$, $-16$, $-18$, $-50$, $-54$, $-52$, $-96$, $-82$, $-76$, $-90$, $62$, $70$, $56$, $112$, $-2$, $-84$, $-74$, $-88$, $-92$, $-78$, $-80$, $-94$, $-86$, $30$, $100$, $106$, $24$, $108$, $98$, $32$, $28$, $102$, $104$, $26$, $110$, $-6$, $-20$, $-14$, $-46$, $-12$, $-22$, $-8$, $-4)$\\

\noindent$11a281$ : $3$ | $(6$, $10$, $16$, $12$, $20$, $2$, $18$, $4$, $22$, $8$, $14)$ | $(38$, $70$, $80$, $28$, $48$, $46$, $30$, $78$, $72$, $36$, $40$, $68$, $82$, $-106$, $-92$, $-112$, $-126$, $-118$, $-98$, $-100$, $-120$, $-124$, $-110$, $-90$, $-108$, $-122$, $42$, $34$, $74$, $76$, $32$, $44$, $50$, $136$, $-102$, $-96$, $-116$, $-128$, $-114$, $-94$, $-104$, $-66$, $-52$, $-64$, $-10$, $154$, $148$, $132$, $140$, $162$, $138$, $134$, $150$, $152$, $88$, $156$, $146$, $130$, $142$, $160$, $84$, $86$, $158$, $144$, $-58$, $-16$, $-4$, $-24$, $-22$, $-2$, $-18$, $-56$, $-60$, $-14$, $-6$, $-26$, $-8$, $-12$, $-62$, $-54$, $-20)$\\

\noindent$11a282$ : $3$ | $(6$, $10$, $16$, $14$, $4$, $18$, $20$, $22$, $2$, $12$, $8)$ | $(90$, $36$, $80$, $44$, $46$, $48$, $146$, $132$, $128$, $150$, $122$, $138$, $140$, $120$, $148$, $130$, $-16$, $-112$, $-106$, $-22$, $-26$, $-102$, $-116$, $-100$, $98$, $82$, $34$, $88$, $92$, $38$, $78$, $42$, $96$, $84$, $32$, $86$, $94$, $40$, $-24$, $-104$, $-114$, $-14$, $-30$, $-18$, $-110$, $-108$, $-20$, $-28$, $-12$, $-10$, $118$, $142$, $136$, $124$, $152$, $126$, $134$, $144$, $-8$, $-64$, $-74$, $-58$, $-2$, $-54$, $-70$, $-68$, $-52$, $-4$, $-60$, $-76$, $-62$, $-6$, $-50$, $-66$, $-72$, $-56)$\\

\noindent$11a283$ : $3$ | $(6$, $10$, $16$, $14$, $18$, $2$, $20$, $22$, $4$, $12$, $8)$ | $(80$, $282$, $292$, $90$, $96$, $298$, $276$, $74$, $272$, $302$, $100$, $86$, $288$, $286$, $84$, $102$, $304$, $270$, $76$, $278$, $296$, $94$, $92$, $294$, $280$, $78$, $268$, $306$, $104$, $82$, $284$, $290$, $88$, $98$, $300$, $274$, $-126$, $-242$, $-148$, $-256$, $-112$, $-140$, $-234$, $-134$, $-118$, $-250$, $-154$, $-248$, $-120$, $-132$, $-236$, $-142$, $-262$, $308$, $266$, $310$, $192$, $228$, $184$, $158$, $176$, $220$, $200$, $202$, $218$, $174$, $160$, $186$, $230$, $190$, $164$, $170$, $214$, $206$, $196$, $224$, $180$, $-44$, $-8$, $-66$, $-22$, $-30$, $-58$, $-106$, $-56$, $-32$, $-20$, $-68$, $-10$, $-42$, $-46$, $-6$, $-64$, $-24$, $-28$, $-60$, $-2$, $-50$, $-38$, $-14$, $-72$, $-16$, $-36$, $-52$, $-110$, $-108$, $-54$, $-34$, $-18$, $-70$, $-12$, $-40$, $-48$, $-4$, $-62$, $-26$, $188$, $162$, $172$, $216$, $204$, $198$, $222$, $178$, $156$, $182$, $226$, $194$, $208$, $212$, $168$, $166$, $210$, $-258$, $-146$, $-240$, $-128$, $-124$, $-244$, $-150$, $-254$, $-114$, $-138$, $-232$, $-136$, $-116$, $-252$, $-152$, $-246$, $-122$, $-130$, $-238$, $-144$, $-260$, $-264)$\\

\noindent$11a284$ : $3$ | $(6$, $10$, $16$, $14$, $20$, $4$, $18$, $2$, $22$, $8$, $12)$ | $(82$, $94$, $134$, $102$, $136$, $92$, $84$, $80$, $96$, $132$, $100$, $76$, $88$, $-124$, $-12$, $-2$, $-16$, $-128$, $-22$, $-8$, $-6$, $-20$, $-130$, $-18$, $-4$, $-10$, $138$, $142$, $26$, $144$, $36$, $44$, $46$, $34$, $146$, $28$, $140$, $-54$, $-60$, $-68$, $-110$, $-108$, $-70$, $-58$, $-56$, $-72$, $-106$, $-112$, $-66$, $-62$, $-116$, $40$, $50$, $30$, $148$, $32$, $48$, $42$, $38$, $52$, $74$, $90$, $86$, $78$, $98$, $-64$, $-114$, $-104$, $-118$, $-120$, $-122$, $-24$, $-126$, $-14)$\\

\noindent$11a285$ : $3$ | $(6$, $10$, $16$, $14$, $20$, $18$, $2$, $22$, $12$, $4$, $8)$ | $(74$, $58$, $76$, $72$, $60$, $86$, $62$, $70$, $78$, $38$, $82$, $66$, $-44$, $-14$, $-8$, $-2$, $-4$, $-114$, $-94$, $-104$, $80$, $68$, $64$, $84$, $36$, $34$, $30$, $138$, $32$, $-52$, $-88$, $-110$, $-98$, $-100$, $-108$, $-90$, $-54$, $-56$, $-92$, $-106$, $-102$, $-96$, $-112$, $126$, $28$, $136$, $116$, $122$, $130$, $24$, $132$, $120$, $118$, $134$, $26$, $128$, $124$, $-50$, $-18$, $-40$, $-20$, $-48$, $-10$, $-12$, $-46$, $-22$, $-42$, $-16$, $-6)$\\

\noindent$11a286$ : $3$ | $(6$, $10$, $16$, $18$, $14$, $2$, $20$, $4$, $22$, $12$, $8)$ | $(68$, $58$, $52$, $162$, $168$, $46$, $64$, $62$, $48$, $166$, $164$, $50$, $60$, $66$, $44$, $170$, $160$, $54$, $56$, $158$, $172$, $-100$, $-86$, $-148$, $-80$, $-142$, $-92$, $-106$, $-94$, $-140$, $-82$, $-150$, $-84$, $-138$, $-136$, $-40$, $-32$, $-8$, $-20$, $126$, $112$, $152$, $116$, $122$, $74$, $130$, $108$, $134$, $70$, $118$, $120$, $72$, $132$, $-36$, $-4$, $-24$, $-16$, $-12$, $-28$, $-42$, $-30$, $-10$, $-18$, $-22$, $-6$, $-34$, $-38$, $156$, $174$, $154$, $110$, $128$, $76$, $124$, $114$, $-14$, $-26$, $-2$, $-96$, $-104$, $-90$, $-144$, $-78$, $-146$, $-88$, $-102$, $-98)$\\

\noindent$11a287$ : $3$ | $(6$, $10$, $16$, $18$, $20$, $2$, $8$, $4$, $22$, $12$, $14)$ |  $(102$, $28$, $56$, $40$, $44$, $52$, $32$, $106$, $34$, $50$, $46$, $38$, $58$, $-118$, $-76$, $-66$, $-108$, $-128$, $-122$, $-114$, $-72$, $-70$, $-112$, $-124$, $-126$, $-110$, $-68$, $-74$, $-116$, $-120$, $-130$, $164$, $160$, $150$, $90$, $140$, $92$, $152$, $158$, $166$, $60$, $36$, $48$, $-20$, $-10$, $-138$, $-6$, $-24$, $-82$, $-62$, $-80$, $-78$, $-64$, $162$, $148$, $88$, $142$, $94$, $154$, $156$, $96$, $144$, $86$, $146$, $98$, $100$, $104$, $30$, $54$, $42$, $-8$, $-22$, $-84$, $-18$, $-12$, $-136$, $-4$, $-26$, $-2$, $-134$, $-14$, $-16$, $-132)$\\

\noindent$11a288$ : $3$ | $(6$, $10$, $16$, $18$, $20$, $4$, $8$, $2$, $22$, $14$, $12)$ | $(106$, $86$, $150$, $94$, $98$, $112$, $100$, $92$, $148$, $88$, $104$, $108$, $-64$, $-80$, $-58$, $-6$, $-16$, $-146$, $-18$, $-4$, $-60$, $-82$, $-62$, $-2$, $-20$, $-144$, $-14$, $-8$, $96$, $152$, $84$, $156$, $158$, $166$, $50$, $36$, $34$, $52$, $164$, $160$, $154$, $-136$, $-126$, $-120$, $-70$, $-74$, $-116$, $-130$, $-132$, $-114$, $-76$, $-68$, $-122$, $-124$, $-66$, $-78$, $162$, $54$, $32$, $38$, $48$, $168$, $46$, $40$, $30$, $56$, $28$, $42$, $44$, $26$, $110$, $102$, $90$, $-72$, $-118$, $-128$, $-134$, $-138$, $-24$, $-140$, $-10$, $-12$, $-142$, $-22)$\\

\noindent$11a289$ : $3$ | $(6$, $10$, $16$, $20$, $4$, $18$, $8$, $22$, $2$, $12$, $14)$ | $(86$, $94$, $80$, $124$, $78$, $92$, $88$, $22$, $34$, $36$, $-70$, $-64$, $-58$, $-6$, $-8$, $-56$, $-66$, $-72$, $-68$, $-54$, $-10$, $-4$, $-60$, $-62$, $-2$, $-12$, $90$, $76$, $126$, $82$, $96$, $84$, $128$, $74$, $132$, $130$, $-18$, $-102$, $-108$, $-122$, $-112$, $-98$, $-116$, $-118$, $-104$, $-106$, $-120$, $-114$, $46$, $24$, $32$, $38$, $52$, $40$, $30$, $26$, $44$, $48$, $134$, $50$, $42$, $28$, $-110$, $-100$, $-16$, $-20$, $-14)$\\

\noindent$11a290$ : $3$ | $(6$, $10$, $16$, $20$, $18$, $2$, $8$, $22$, $4$, $12$, $14)$ | $(50$, $138$, $60$, $40$, $46$, $54$, $134$, $56$, $44$, $42$, $58$, $136$, $52$, $48$, $38$, $62$, $66$, $-76$, $-70$, $-82$, $-128$, $-92$, $-130$, $-84$, $-120$, $-122$, $-86$, $-132$, $-90$, $-126$, $-80$, $-72$, $-74$, $-78$, $140$, $146$, $112$, $144$, $142$, $110$, $148$, $96$, $156$, $102$, $106$, $152$, $-18$, $-8$, $-28$, $-2$, $-24$, $-12$, $-14$, $-22$, $-4$, $-32$, $64$, $36$, $68$, $98$, $158$, $100$, $108$, $150$, $94$, $154$, $104$, $-88$, $-124$, $-118$, $-116$, $-34$, $-114$, $-30$, $-6$, $-20$, $-16$, $-10$, $-26)$\\

\noindent$11a291$ : $3$ | $(6$, $10$, $16$, $22$, $4$, $18$, $20$, $2$, $8$, $12$, $14)$ | $(116$, $66$, $56$, $46$, $74$, $48$, $54$, $68$, $118$, $114$, $64$, $58$, $44$, $72$, $50$, $52$, $70$, $120$, $112$, $62$, $60$, $-108$, $-90$, $-76$, $-94$, $-132$, $-96$, $-78$, $-88$, $-106$, $-140$, $-104$, $-86$, $-80$, $-98$, $-134$, $-92$, $154$, $128$, $146$, $162$, $174$, $164$, $144$, $126$, $156$, $152$, $130$, $148$, $160$, $172$, $166$, $142$, $124$, $122$, $-138$, $-102$, $-84$, $-82$, $-100$, $-136$, $-42$, $-24$, $-6$, $-12$, $-30$, $150$, $158$, $170$, $168$, $-40$, $-22$, $-4$, $-14$, $-32$, $-28$, $-10$, $-8$, $-26$, $-34$, $-16$, $-2$, $-20$, $-38$, $-110$, $-36$, $-18)$\\

\noindent$11a292$ : $3$ | $(6$, $10$, $16$, $22$, $4$, $18$, $20$, $2$, $8$, $14$, $12)$ | $(136$, $188$, $186$, $214$, $206$, $178$, $196$, $224$, $220$, $200$, $174$, $202$, $218$, $226$, $194$, $180$, $208$, $212$, $184$, $190$, $134$, $192$, $182$, $210$, $-34$, $-8$, $-150$, $-24$, $-18$, $-144$, $-14$, $-28$, $-42$, $-40$, $-2$, $-6$, $-36$, $-46$, $-32$, $-10$, $-148$, $-22$, $-20$, $-146$, $-12$, $-30$, $-44$, $-38$, $-4$, $70$, $52$, $80$, $88$, $60$, $142$, $58$, $86$, $82$, $54$, $138$, $64$, $92$, $76$, $48$, $74$, $94$, $66$, $-154$, $-152$, $-26$, $-16$, $222$, $198$, $176$, $204$, $216$, $228$, $68$, $96$, $72$, $50$, $78$, $90$, $62$, $140$, $56$, $84$, $-114$, $-104$, $-164$, $-124$, $-130$, $-158$, $-98$, $-156$, $-132$, $-122$, $-166$, $-106$, $-112$, $-172$, $-116$, $-102$, $-162$, $-126$, $-128$, $-160$, $-100$, $-118$, $-170$, $-110$, $-108$, $-168$, $-120)$\\

\noindent$11a293$ : $3$ | $(6$, $10$, $16$, $22$, $4$, $18$, $20$, $8$, $2$, $12$, $14)$ | $(78$, $86$, $74$, $116$, $68$, $80$, $84$, $72$, $114$, $70$, $82$, $-12$, $-2$, $-88$, $-4$, $-14$, $-20$, $-10$, $-94$, $-92$, $-8$, $-18$, $-16$, $-6$, $-90$, $-96$, $36$, $40$, $48$, $28$, $26$, $50$, $38$, $-54$, $-106$, $-112$, $-102$, $-58$, $-62$, $-98$, $-108$, $-110$, $-100$, $-60$, $76$, $118$, $66$, $120$, $122$, $24$, $30$, $46$, $42$, $34$, $22$, $32$, $44$, $-104$, $-56$, $-64$, $-52)$\\

\noindent$11a294$ : $3$ | $(6$, $10$, $16$, $22$, $4$, $18$, $20$, $8$, $2$, $14$, $12)$ | $(30$, $104$, $96$, $110$, $112$, $98$, $102$, $28$, $32$, $106$, $94$, $108$, $34$, $-88$, $-52$, $-50$, $-90$, $-86$, $68$, $82$, $70$, $60$, $74$, $78$, $64$, $120$, $66$, $80$, $72$, $-8$, $-24$, $-18$, $-2$, $-14$, $-84$, $-10$, $-6$, $-22$, $-20$, $-4$, $-12$, $100$, $114$, $116$, $118$, $62$, $76$, $-46$, $-56$, $-40$, $-38$, $-54$, $-48$, $-92$, $-44$, $-58$, $-42$, $-36$, $-26$, $-16)$\\

\noindent$11a295$ : $3$ | $(6$, $10$, $16$, $22$, $18$, $2$, $20$, $8$, $4$, $12$, $14)$ | $(60$, $52$, $46$, $50$, $58$, $62$, $54$, $-36$, $-28$, $-26$, $-4$, $-6$, $48$, $74$, $80$, $88$, $86$, $78$, $76$, $84$, $90$, $82$, $22$, $-70$, $-72$, $-44$, $-38$, $-34$, $-30$, $-42$, $-40$, $-32$, $56$, $16$, $20$, $24$, $18$, $-68$, $-14$, $-66$, $-8$, $-2$, $-12$, $-64$, $-10)$\\

\noindent$11a296$ : $3$ | $(6$, $10$, $16$, $22$, $18$, $2$, $20$, $8$, $4$, $14$, $12)$ | $(62$, $160$, $60$, $184$, $170$, $50$, $174$, $180$, $56$, $164$, $190$, $188$, $166$, $54$, $178$, $176$, $52$, $168$, $186$, $192$, $162$, $58$, $182$, $172$, $-146$, $-70$, $-130$, $-76$, $-140$, $-152$, $-158$, $-156$, $84$, $120$, $100$, $94$, $126$, $90$, $104$, $116$, $-22$, $-38$, $-6$, $-44$, $-16$, $-28$, $-32$, $-12$, $-48$, $-10$, $-34$, $-26$, $-18$, $-42$, $-4$, $-2$, $-40$, $-20$, $-24$, $-36$, $-8$, $-46$, $-14$, $-30$, $92$, $102$, $118$, $82$, $114$, $106$, $88$, $124$, $96$, $98$, $122$, $86$, $108$, $112$, $194$, $110$, $-154$, $-138$, $-78$, $-132$, $-68$, $-148$, $-144$, $-72$, $-128$, $-74$, $-142$, $-150$, $-66$, $-134$, $-80$, $-136$, $-64)$\\

\noindent$11a297$ : $3$ | $(6$, $10$, $16$, $22$, $20$, $18$, $2$, $8$, $12$, $4$, $14)$ | $(70$, $192$, $172$, $206$, $84$, $56$, $178$, $186$, $64$, $76$, $198$, $166$, $200$, $78$, $62$, $184$, $180$, $58$, $82$, $204$, $170$, $194$, $72$, $68$, $190$, $174$, $-252$, $-210$, $-154$, $-234$, $-228$, $-160$, $-216$, $-246$, $-142$, $-144$, $-244$, $-218$, $-162$, $-226$, $-236$, $-152$, $-208$, $-254$, $-6$, $-40$, $-120$, $-24$, $-22$, $-122$, $-42$, $-4$, $-50$, $-130$, $-14$, $-32$, $182$, $60$, $80$, $202$, $168$, $196$, $74$, $66$, $188$, $176$, $54$, $86$, $88$, $104$, $272$, $306$, $270$, $102$, $90$, $258$, $294$, $284$, $110$, $278$, $300$, $264$, $96$, $-222$, $-240$, $-148$, $-138$, $-250$, $-212$, $-156$, $-232$, $-230$, $-158$, $-214$, $-248$, $-140$, $-146$, $-242$, $-220$, $-164$, $-224$, $-238$, $-150$, $256$, $292$, $286$, $108$, $276$, $302$, $266$, $98$, $94$, $262$, $298$, $280$, $112$, $282$, $296$, $260$, $92$, $100$, $268$, $304$, $274$, $106$, $288$, $290$, $-136$, $-8$, $-38$, $-118$, $-26$, $-20$, $-124$, $-44$, $-2$, $-48$, $-128$, $-16$, $-30$, $-114$, $-34$, $-12$, $-132$, $-52$, $-134$, $-10$, $-36$, $-116$, $-28$, $-18$, $-126$, $-46)$\\

\noindent$11a298$ : $3$ | $(6$, $10$, $18$, $2$, $16$, $20$, $22$, $4$, $8$, $14$, $12)$ | $(104$, $40$, $94$, $114$, $92$, $136$, $152$, $168$, $146$, $142$, $164$, $156$, $132$, $90$, $134$, $154$, $166$, $144$, $-18$, $-124$, $-118$, $-12$, $-30$, $-24$, $-130$, $-26$, $-28$, $-10$, $42$, $106$, $102$, $38$, $96$, $112$, $48$, $50$, $54$, $44$, $108$, $100$, $36$, $98$, $110$, $46$, $52$, $-72$, $-116$, $-126$, $-20$, $-34$, $-16$, $-122$, $-120$, $-14$, $-32$, $-22$, $-128$, $160$, $138$, $150$, $170$, $148$, $140$, $162$, $158$, $-70$, $-88$, $-74$, $-56$, $-4$, $-64$, $-82$, $-80$, $-62$, $-2$, $-58$, $-76$, $-86$, $-68$, $-8$, $-6$, $-66$, $-84$, $-78$, $-60)$\\

\noindent$11a299$ : $3$ | $(6$, $10$, $18$, $2$, $16$, $22$, $20$, $4$, $8$, $14$, $12)$ | $(38$, $50$, $108$, $48$, $36$, $40$, $116$, $52$, $110$, $46$, $34$, $42$, $114$, $54$, $112$, $44$, $-60$, $-72$, $-104$, $-102$, $-70$, $-58$, $-62$, $-74$, $-106$, $-78$, $-66$, $96$, $84$, $126$, $124$, $86$, $94$, $98$, $82$, $128$, $122$, $88$, $92$, $100$, $-12$, $-24$, $-30$, $-18$, $-6$, $-4$, $-16$, $-28$, $-26$, $-14$, $118$, $130$, $120$, $90$, $-76$, $-64$, $-56$, $-68$, $-80$, $-2$, $-8$, $-20$, $-32$, $-22$, $-10)$\\

\noindent$11a300$ : $3$ | $(6$, $10$, $18$, $14$, $2$, $20$, $22$, $8$, $12$, $4$, $16)$ | $(178$, $182$, $28$, $66$, $58$, $36$, $88$, $82$, $42$, $52$, $72$, $-240$, $-198$, $-258$, $-222$, $-216$, $-210$, $-228$, $-252$, $-192$, $-246$, $-234$, $-204$, $-262$, $-202$, $-236$, $-244$, $-194$, $-254$, $-226$, $-212$, $-214$, $-224$, $-256$, $-196$, $-242$, $-238$, $-200$, $-260$, $-206$, $-232$, $-248$, $-190$, $-250$, $-230$, $-208$, $-2$, $180$, $30$, $64$, $60$, $34$, $90$, $80$, $44$, $50$, $74$, $24$, $70$, $54$, $40$, $84$, $86$, $38$, $56$, $68$, $26$, $76$, $48$, $46$, $78$, $92$, $32$, $62$, $-100$, $-8$, $-130$, $-118$, $-20$, $-112$, $-136$, $-140$, $-108$, $-16$, $-122$, $-126$, $-12$, $-104$, $-144$, $-218$, $-220$, $-142$, $-106$, $-14$, $-124$, $170$, $270$, $154$, $284$, $290$, $160$, $264$, $164$, $184$, $176$, $276$, $148$, $278$, $296$, $94$, $294$, $280$, $150$, $274$, $174$, $186$, $166$, $266$, $158$, $288$, $286$, $156$, $268$, $168$, $188$, $172$, $272$, $152$, $282$, $292$, $162$, $-138$, $-110$, $-18$, $-120$, $-128$, $-10$, $-102$, $-146$, $-98$, $-6$, $-132$, $-116$, $-22$, $-114$, $-134$, $-4$, $-96)$\\

\noindent$11a301$ : $3$ | $(6$, $10$, $18$, $14$, $16$, $4$, $22$, $20$, $12$, $2$, $8)$ | $(84$, $138$, $94$, $74$, $50$, $78$, $90$, $134$, $88$, $80$, $48$, $72$, $46$, $-68$, $-62$, $-2$, $-58$, $-16$, $-14$, $-56$, $-4$, $-64$, $-66$, $-70$, $-128$, $-110$, $76$, $92$, $136$, $86$, $82$, $140$, $96$, $142$, $144$, $154$, $-126$, $-112$, $-108$, $-130$, $-100$, $-120$, $-118$, $-102$, $-132$, $-106$, $-114$, $-124$, $44$, $28$, $164$, $34$, $38$, $160$, $148$, $150$, $158$, $40$, $32$, $166$, $30$, $42$, $156$, $152$, $146$, $162$, $36$, $-104$, $-116$, $-122$, $-98$, $-26$, $-20$, $-10$, $-52$, $-8$, $-22$, $-24$, $-6$, $-54$, $-12$, $-18$, $-60)$\\

\noindent$11a302$ : $3$ | $(6$, $10$, $18$, $14$, $20$, $4$, $22$, $8$, $12$, $2$, $16)$ | $(82$, $108$, $88$, $76$, $48$, $78$, $86$, $106$, $84$, $80$, $110$, $90$, $114$, $-94$, $-20$, $-2$, $-16$, $-14$, $-4$, $-22$, $-92$, $-24$, $-96$, $-100$, $120$, $126$, $36$, $32$, $130$, $30$, $38$, $44$, $26$, $42$, $40$, $28$, $46$, $-98$, $-50$, $-102$, $-58$, $-56$, $-104$, $-52$, $-62$, $-68$, $112$, $116$, $124$, $122$, $118$, $128$, $34$, $-54$, $-60$, $-70$, $-66$, $-64$, $-72$, $-10$, $-8$, $-74$, $-6$, $-12$, $-18)$\\

\noindent$11a303$ : $3$ | $(6$, $10$, $18$, $22$, $16$, $4$, $8$, $20$, $12$, $2$, $14)$ | $(104$, $98$, $116$, $168$, $114$, $96$, $106$, $176$, $122$, $174$, $108$, $94$, $112$, $170$, $118$, $180$, $184$, $-140$, $-2$, $-12$, $-22$, $-92$, $-24$, $-10$, $-4$, $-30$, $-86$, $-16$, $-18$, $-88$, $-28$, $-6$, $-8$, $-26$, $-90$, $-20$, $-14$, $-84$, $-32$, $182$, $36$, $186$, $178$, $120$, $172$, $110$, $-152$, $-128$, $-164$, $-80$, $-34$, $-82$, $-162$, $-126$, $-154$, $-150$, $-130$, $-166$, $-134$, $-146$, $-158$, $56$, $72$, $40$, $46$, $66$, $62$, $50$, $78$, $100$, $102$, $76$, $52$, $60$, $68$, $44$, $42$, $70$, $58$, $54$, $74$, $38$, $48$, $64$, $-132$, $-148$, $-156$, $-124$, $-160$, $-144$, $-136$, $-138$, $-142)$\\

\noindent$11a304$ : $3$ | $(6$, $10$, $20$, $14$, $2$, $18$, $8$, $22$, $12$, $4$, $16)$ | $(48$, $32$, $146$, $142$, $36$, $52$, $44$, $92$, $108$, $106$, $166$, $164$, $104$, $110$, $94$, $-80$, $-74$, $-68$, $-86$, $-24$, $-20$, $-90$, $-64$, $-78$, $-76$, $-66$, $-88$, $-22$, $40$, $138$, $150$, $154$, $170$, $160$, $100$, $114$, $98$, $158$, $172$, $156$, $96$, $112$, $102$, $162$, $168$, $152$, $-62$, $-124$, $-4$, $-116$, $-8$, $-128$, $-58$, $-134$, $-60$, $-126$, $-6$, $144$, $34$, $50$, $46$, $30$, $148$, $140$, $38$, $54$, $42$, $136$, $-18$, $-26$, $-84$, $-70$, $-72$, $-82$, $-28$, $-16$, $-14$, $-122$, $-2$, $-118$, $-10$, $-130$, $-56$, $-132$, $-12$, $-120)$\\

\noindent$11a305$ : $3$ | $(6$, $12$, $10$, $22$, $16$, $18$, $20$, $8$, $4$, $2$, $14)$ | $(56$, $64$, $54$, $58$, $62$, $50$, $30$, $52$, $60$, $-8$, $-2$, $-4$, $-44$, $-46$, $-40$, $82$, $90$, $18$, $92$, $80$, $84$, $88$, $78$, $86$, $-42$, $-32$, $-38$, $-66$, $-68$, $-36$, $-34$, $-70$, $22$, $94$, $20$, $24$, $28$, $48$, $26$, $-72$, $-14$, $-12$, $-74$, $-76$, $-10$, $-16$, $-6)$\\

\noindent$11a306$ : $2$ | $(6$, $12$, $14$, $16$, $18$, $20$, $4$, $2$, $22$, $10$, $8)$ | $(180$, $122$, $144$, $202$, $158$, $108$, $166$, $194$, $136$, $130$, $188$, $172$, $114$, $152$, $208$, $150$, $116$, $174$, $186$, $128$, $138$, $196$, $164$, $106$, $160$, $200$, $142$, $124$, $182$, $178$, $120$, $146$, $204$, $156$, $110$, $168$, $192$, $134$, $132$, $190$, $170$, $112$, $154$, $206$, $148$, $118$, $176$, $184$, $126$, $140$, $198$, $162$, $-284$, $-226$, $-248$, $-306$, $-262$, $-212$, $-270$, $-298$, $-240$, $-234$, $-292$, $-276$, $-218$, $-256$, $-312$, $-254$, $-220$, $-278$, $-290$, $-232$, $-242$, $-300$, $-268$, $-210$, $-264$, $-304$, $-246$, $-228$, $-286$, $-282$, $-224$, $-250$, $-308$, $-260$, $-214$, $-272$, $-296$, $-238$, $-236$, $-294$, $-274$, $-216$, $-258$, $-310$, $-252$, $-222$, $-280$, $-288$, $-230$, $-244$, $-302$, $-266$, $388$, $330$, $352$, $410$, $366$, $316$, $374$, $402$, $344$, $338$, $396$, $380$, $322$, $360$, $416$, $358$, $324$, $382$, $394$, $336$, $346$, $404$, $372$, $314$, $368$, $408$, $350$, $332$, $390$, $386$, $328$, $354$, $412$, $364$, $318$, $376$, $400$, $342$, $340$, $398$, $378$, $320$, $362$, $414$, $356$, $326$, $384$, $392$, $334$, $348$, $406$, $370$, $-76$, $-18$, $-40$, $-98$, $-54$, $-4$, $-62$, $-90$, $-32$, $-26$, $-84$, $-68$, $-10$, $-48$, $-104$, $-46$, $-12$, $-70$, $-82$, $-24$, $-34$, $-92$, $-60$, $-2$, $-56$, $-96$, $-38$, $-20$, $-78$, $-74$, $-16$, $-42$, $-100$, $-52$, $-6$, $-64$, $-88$, $-30$, $-28$, $-86$, $-66$, $-8$, $-50$, $-102$, $-44$, $-14$, $-72$, $-80$, $-22$, $-36$, $-94$, $-58)$\\

\noindent$11a307$ : $2$ | $(6$, $12$, $14$, $16$, $20$, $18$, $4$, $2$, $22$, $10$, $8)$ | $(100$, $136$, $158$, $122$, $86$, $114$, $150$, $144$, $108$, $92$, $128$, $164$, $130$, $94$, $106$, $142$, $152$, $116$, $84$, $120$, $156$, $138$, $102$, $98$, $134$, $160$, $124$, $88$, $112$, $148$, $146$, $110$, $90$, $126$, $162$, $132$, $96$, $104$, $140$, $154$, $118$, $-224$, $-178$, $-196$, $-242$, $-206$, $-168$, $-214$, $-234$, $-188$, $-186$, $-232$, $-216$, $-170$, $-204$, $-244$, $-198$, $-176$, $-222$, $-226$, $-180$, $-194$, $-240$, $-208$, $-166$, $-212$, $-236$, $-190$, $-184$, $-230$, $-218$, $-172$, $-202$, $-246$, $-200$, $-174$, $-220$, $-228$, $-182$, $-192$, $-238$, $-210$, $264$, $300$, $322$, $286$, $250$, $278$, $314$, $308$, $272$, $256$, $292$, $328$, $294$, $258$, $270$, $306$, $316$, $280$, $248$, $284$, $320$, $302$, $266$, $262$, $298$, $324$, $288$, $252$, $276$, $312$, $310$, $274$, $254$, $290$, $326$, $296$, $260$, $268$, $304$, $318$, $282$, $-60$, $-14$, $-32$, $-78$, $-42$, $-4$, $-50$, $-70$, $-24$, $-22$, $-68$, $-52$, $-6$, $-40$, $-80$, $-34$, $-12$, $-58$, $-62$, $-16$, $-30$, $-76$, $-44$, $-2$, $-48$, $-72$, $-26$, $-20$, $-66$, $-54$, $-8$, $-38$, $-82$, $-36$, $-10$, $-56$, $-64$, $-18$, $-28$, $-74$, $-46)$\\

\noindent$11a308$ : $2$ | $(6$, $12$, $14$, $18$, $20$, $2$, $4$, $22$, $10$, $8$, $16)$ | $(122$, $84$, $94$, $132$, $112$, $74$, $104$, $140$, $102$, $76$, $114$, $130$, $92$, $86$, $124$, $120$, $82$, $96$, $134$, $110$, $72$, $106$, $138$, $100$, $78$, $116$, $128$, $90$, $88$, $126$, $118$, $80$, $98$, $136$, $108$, $-196$, $-166$, $-144$, $-174$, $-204$, $-188$, $-158$, $-152$, $-182$, $-210$, $-180$, $-150$, $-160$, $-190$, $-202$, $-172$, $-142$, $-168$, $-198$, $-194$, $-164$, $-146$, $-176$, $-206$, $-186$, $-156$, $-154$, $-184$, $-208$, $-178$, $-148$, $-162$, $-192$, $-200$, $-170$, $262$, $224$, $234$, $272$, $252$, $214$, $244$, $280$, $242$, $216$, $254$, $270$, $232$, $226$, $264$, $260$, $222$, $236$, $274$, $250$, $212$, $246$, $278$, $240$, $218$, $256$, $268$, $230$, $228$, $266$, $258$, $220$, $238$, $276$, $248$, $-56$, $-26$, $-4$, $-34$, $-64$, $-48$, $-18$, $-12$, $-42$, $-70$, $-40$, $-10$, $-20$, $-50$, $-62$, $-32$, $-2$, $-28$, $-58$, $-54$, $-24$, $-6$, $-36$, $-66$, $-46$, $-16$, $-14$, $-44$, $-68$, $-38$, $-8$, $-22$, $-52$, $-60$, $-30)$\\

\noindent$11a309$ : $2$ | $(6$, $12$, $14$, $18$, $20$, $4$, $2$, $22$, $10$, $8$, $16)$ | $(160$, $110$, $124$, $174$, $146$, $96$, $138$, $182$, $132$, $102$, $152$, $168$, $118$, $116$, $166$, $154$, $104$, $130$, $180$, $140$, $94$, $144$, $176$, $126$, $108$, $158$, $162$, $112$, $122$, $172$, $148$, $98$, $136$, $184$, $134$, $100$, $150$, $170$, $120$, $114$, $164$, $156$, $106$, $128$, $178$, $142$, $-210$, $-262$, $-240$, $-188$, $-232$, $-270$, $-218$, $-202$, $-254$, $-248$, $-196$, $-224$, $-276$, $-226$, $-194$, $-246$, $-256$, $-204$, $-216$, $-268$, $-234$, $-186$, $-238$, $-264$, $-212$, $-208$, $-260$, $-242$, $-190$, $-230$, $-272$, $-220$, $-200$, $-252$, $-250$, $-198$, $-222$, $-274$, $-228$, $-192$, $-244$, $-258$, $-206$, $-214$, $-266$, $-236$, $344$, $294$, $308$, $358$, $330$, $280$, $322$, $366$, $316$, $286$, $336$, $352$, $302$, $300$, $350$, $338$, $288$, $314$, $364$, $324$, $278$, $328$, $360$, $310$, $292$, $342$, $346$, $296$, $306$, $356$, $332$, $282$, $320$, $368$, $318$, $284$, $334$, $354$, $304$, $298$, $348$, $340$, $290$, $312$, $362$, $326$, $-26$, $-78$, $-56$, $-4$, $-48$, $-86$, $-34$, $-18$, $-70$, $-64$, $-12$, $-40$, $-92$, $-42$, $-10$, $-62$, $-72$, $-20$, $-32$, $-84$, $-50$, $-2$, $-54$, $-80$, $-28$, $-24$, $-76$, $-58$, $-6$, $-46$, $-88$, $-36$, $-16$, $-68$, $-66$, $-14$, $-38$, $-90$, $-44$, $-8$, $-60$, $-74$, $-22$, $-30$, $-82$, $-52)$\\

\noindent$11a310$ : $2$ | $(6$, $12$, $14$, $20$, $18$, $2$, $4$, $22$, $10$, $8$, $16)$ | $(108$, $82$, $64$, $90$, $116$, $100$, $74$, $72$, $98$, $118$, $92$, $66$, $80$, $106$, $110$, $84$, $62$, $88$, $114$, $102$, $76$, $70$, $96$, $120$, $94$, $68$, $78$, $104$, $112$, $86$, $-134$, $-162$, $-172$, $-144$, $-124$, $-152$, $-180$, $-154$, $-126$, $-142$, $-170$, $-164$, $-136$, $-132$, $-160$, $-174$, $-146$, $-122$, $-150$, $-178$, $-156$, $-128$, $-140$, $-168$, $-166$, $-138$, $-130$, $-158$, $-176$, $-148$, $228$, $202$, $184$, $210$, $236$, $220$, $194$, $192$, $218$, $238$, $212$, $186$, $200$, $226$, $230$, $204$, $182$, $208$, $234$, $222$, $196$, $190$, $216$, $240$, $214$, $188$, $198$, $224$, $232$, $206$, $-14$, $-42$, $-52$, $-24$, $-4$, $-32$, $-60$, $-34$, $-6$, $-22$, $-50$, $-44$, $-16$, $-12$, $-40$, $-54$, $-26$, $-2$, $-30$, $-58$, $-36$, $-8$, $-20$, $-48$, $-46$, $-18$, $-10$, $-38$, $-56$, $-28)$\\

\noindent$11a311$ : $2$ | $(6$, $12$, $14$, $20$, $18$, $4$, $2$, $22$, $10$, $8$, $16)$ | $(96$, $132$, $146$, $110$, $82$, $118$, $154$, $124$, $88$, $104$, $140$, $138$, $102$, $90$, $126$, $152$, $116$, $80$, $112$, $148$, $130$, $94$, $98$, $134$, $144$, $108$, $84$, $120$, $156$, $122$, $86$, $106$, $142$, $136$, $100$, $92$, $128$, $150$, $114$, $-178$, $-222$, $-204$, $-160$, $-196$, $-230$, $-186$, $-170$, $-214$, $-212$, $-168$, $-188$, $-232$, $-194$, $-162$, $-206$, $-220$, $-176$, $-180$, $-224$, $-202$, $-158$, $-198$, $-228$, $-184$, $-172$, $-216$, $-210$, $-166$, $-190$, $-234$, $-192$, $-164$, $-208$, $-218$, $-174$, $-182$, $-226$, $-200$, $252$, $288$, $302$, $266$, $238$, $274$, $310$, $280$, $244$, $260$, $296$, $294$, $258$, $246$, $282$, $308$, $272$, $236$, $268$, $304$, $286$, $250$, $254$, $290$, $300$, $264$, $240$, $276$, $312$, $278$, $242$, $262$, $298$, $292$, $256$, $248$, $284$, $306$, $270$, $-22$, $-66$, $-48$, $-4$, $-40$, $-74$, $-30$, $-14$, $-58$, $-56$, $-12$, $-32$, $-76$, $-38$, $-6$, $-50$, $-64$, $-20$, $-24$, $-68$, $-46$, $-2$, $-42$, $-72$, $-28$, $-16$, $-60$, $-54$, $-10$, $-34$, $-78$, $-36$, $-8$, $-52$, $-62$, $-18$, $-26$, $-70$, $-44)$\\

\noindent$11a312$ : $3$ | $(6$, $12$, $16$, $14$, $18$, $20$, $2$, $22$, $4$, $10$, $8)$ | $(42$, $122$, $128$, $136$, $86$, $134$, $130$, $120$, $44$, $40$, $124$, $126$, $138$, $84$, $132$, $-12$, $-6$, $-24$, $-72$, $-64$, $-96$, $-94$, $-66$, $-70$, $-26$, $-8$, $-10$, $-28$, $-14$, $-4$, $-22$, $74$, $142$, $80$, $82$, $140$, $-114$, $-20$, $-2$, $-16$, $-118$, $-110$, $-104$, $52$, $32$, $58$, $78$, $144$, $76$, $60$, $34$, $50$, $54$, $30$, $56$, $48$, $36$, $38$, $46$, $-68$, $-92$, $-98$, $-62$, $-100$, $-90$, $-108$, $-106$, $-88$, $-102$, $-112$, $-116$, $-18)$\\

\noindent$11a313$ : $3$ | $(6$, $12$, $16$, $14$, $20$, $18$, $2$, $22$, $4$, $10$, $8)$ | $(44$, $42$, $46$, $58$, $112$, $116$, $54$, $50$, $120$, $108$, $62$, $106$, $122$, $48$, $56$, $114$, $-10$, $-130$, $-18$, $-24$, $-124$, $-102$, $-100$, $-78$, $-68$, $-90$, $-92$, $-70$, $-76$, $-98$, $-84$, $144$, $38$, $158$, $136$, $152$, $32$, $150$, $138$, $160$, $40$, $142$, $146$, $36$, $156$, $134$, $154$, $34$, $148$, $140$, $162$, $-2$, $-82$, $-64$, $-86$, $-96$, $-74$, $-72$, $-94$, $-88$, $-66$, $-80$, $104$, $60$, $110$, $118$, $52$, $-14$, $-28$, $-6$, $-126$, $-22$, $-20$, $-128$, $-8$, $-30$, $-12$, $-132$, $-16$, $-26$, $-4)$\\

\noindent$11a314$ : $3$ | $(6$, $12$, $16$, $18$, $14$, $4$, $20$, $22$, $2$, $8$, $10)$ | $(92$, $112$, $84$, $100$, $104$, $148$, $106$, $98$, $86$, $114$, $90$, $94$, $110$, $152$, $-124$, $-10$, $-12$, $-122$, $-26$, $-4$, $-18$, $-116$, $-20$, $-2$, $-24$, $-120$, $-14$, $-8$, $102$, $82$, $28$, $36$, $50$, $164$, $44$, $42$, $162$, $52$, $34$, $30$, $-126$, $-132$, $-142$, $-72$, $-68$, $-138$, $-136$, $-66$, $-74$, $-144$, $-130$, $-128$, $-146$, $-76$, $-64$, $-134$, $-140$, $-70$, $88$, $96$, $108$, $150$, $154$, $156$, $56$, $158$, $38$, $48$, $166$, $46$, $40$, $160$, $54$, $32$, $-78$, $-62$, $-60$, $-80$, $-58$, $-6$, $-16$, $-118$, $-22)$\\

\noindent$11a315$ : $3$ | $(6$, $12$, $16$, $18$, $22$, $4$, $20$, $10$, $2$, $8$, $14)$ | $(50$, $62$, $82$, $68$, $44$, $130$, $42$, $70$, $80$, $60$, $52$, $48$, $64$, $84$, $66$, $46$, $54$, $58$, $78$, $72$, $-98$, $-132$, $-150$, $-88$, $-142$, $-140$, $-90$, $-106$, $-104$, $-92$, $-138$, $-144$, $-86$, $-148$, $-134$, $-96$, $-100$, $-108$, $-102$, $-94$, $-136$, $-146$, $118$, $166$, $180$, $160$, $124$, $112$, $172$, $174$, $110$, $126$, $158$, $128$, $-154$, $-38$, $-16$, $-6$, $-28$, $-26$, $-4$, $-18$, $-36$, $-156$, $-152$, $176$, $170$, $114$, $122$, $162$, $182$, $164$, $120$, $116$, $168$, $178$, $74$, $76$, $56$, $-34$, $-20$, $-2$, $-24$, $-30$, $-8$, $-14$, $-40$, $-12$, $-10$, $-32$, $-22)$\\

\noindent$11a316$ : $3$ | $(6$, $12$, $16$, $20$, $22$, $4$, $18$, $10$, $2$, $14$, $8)$ | $(168$, $122$, $148$, $188$, $142$, $128$, $174$, $162$, $116$, $114$, $160$, $176$, $130$, $140$, $186$, $150$, $86$, $154$, $182$, $136$, $134$, $180$, $156$, $110$, $120$, $166$, $170$, $124$, $146$, $190$, $144$, $126$, $172$, $164$, $118$, $112$, $158$, $178$, $132$, $138$, $184$, $152$, $-216$, $-108$, $-48$, $-70$, $-18$, $286$, $332$, $300$, $266$, $312$, $320$, $274$, $-254$, $-202$, $-230$, $-94$, $-232$, $-200$, $-252$, $-246$, $-194$, $-238$, $-100$, $-224$, $-208$, $-260$, $-210$, $-222$, $-102$, $-240$, $-192$, $-244$, $-106$, $-218$, $-214$, $-256$, $-204$, $-228$, $-96$, $-234$, $-198$, $-250$, $-248$, $-196$, $-236$, $-98$, $-226$, $-206$, $-258$, $-212$, $-220$, $-104$, $-242$, $316$, $270$, $296$, $336$, $290$, $90$, $282$, $328$, $304$, $262$, $308$, $324$, $278$, $276$, $322$, $310$, $264$, $302$, $330$, $284$, $92$, $288$, $334$, $298$, $268$, $314$, $318$, $272$, $294$, $338$, $292$, $88$, $280$, $326$, $306$, $-40$, $-78$, $-26$, $-10$, $-62$, $-56$, $-4$, $-32$, $-84$, $-34$, $-46$, $-72$, $-20$, $-16$, $-68$, $-50$, $-2$, $-54$, $-64$, $-12$, $-24$, $-76$, $-42$, $-38$, $-80$, $-28$, $-8$, $-60$, $-58$, $-6$, $-30$, $-82$, $-36$, $-44$, $-74$, $-22$, $-14$, $-66$, $-52)$\\

\noindent$11a317$ : $3$ | $(6$, $12$, $16$, $20$, $22$, $18$, $2$, $10$, $4$, $14$, $8)$ | $(58$, $18$, $30$, $34$, $22$, $24$, $36$, $28$, $-96$, $-46$, $-52$, $-102$, $-90$, $-94$, $-98$, $-48$, $-50$, $-100$, $-92$, $26$, $108$, $112$, $106$, $68$, $80$, $70$, $104$, $110$, $-42$, $-44$, $-54$, $-40$, $-16$, $-4$, $-82$, $-8$, $-12$, $-38$, $-14$, $-6$, $32$, $20$, $60$, $56$, $62$, $74$, $76$, $64$, $66$, $78$, $72$, $-10$, $-84$, $-2$, $-88$, $-86)$\\

\noindent$11a318$ : $3$ | $(6$, $12$, $16$, $22$, $18$, $2$, $20$, $4$, $8$, $14$, $10)$ | $(22$, $28$, $16$, $32$, $52$, $50$, $30$, $-80$, $-70$, $-64$, $-76$, $-74$, $-62$, $-72$, $-78$, $-82$, $94$, $86$, $90$, $20$, $24$, $26$, $18$, $88$, $-36$, $-84$, $-34$, $-46$, $-38$, $-2$, $-42$, $102$, $56$, $48$, $54$, $100$, $60$, $104$, $58$, $98$, $96$, $92$, $-10$, $-68$, $-66$, $-12$, $-8$, $-6$, $-14$, $-4$, $-44$, $-40)$\\

\noindent$11a319$ : $3$ | $(6$, $12$, $16$, $22$, $18$, $4$, $20$, $2$, $8$, $10$, $14)$ | $(100$, $84$, $178$, $176$, $82$, $168$, $186$, $92$, $-64$, $-72$, $-50$, $-2$, $-10$, $-58$, $-78$, $-56$, $-8$, $-4$, $-52$, $-74$, $-62$, $-14$, $108$, $102$, $98$, $86$, $180$, $174$, $80$, $170$, $184$, $90$, $94$, $106$, $104$, $96$, $88$, $182$, $172$, $-154$, $-132$, $-120$, $-126$, $-148$, $-160$, $-138$, $-140$, $-162$, $-146$, $-124$, $-122$, $-144$, $-164$, $-142$, $-12$, $-60$, $-76$, $-54$, $-6$, $34$, $24$, $194$, $22$, $36$, $116$, $32$, $26$, $110$, $42$, $16$, $188$, $166$, $192$, $20$, $38$, $114$, $30$, $28$, $112$, $40$, $18$, $190$, $-68$, $-46$, $-134$, $-156$, $-152$, $-130$, $-118$, $-128$, $-150$, $-158$, $-136$, $-44$, $-66$, $-70$, $-48)$\\

\noindent$11a320$ : $3$ | $(6$, $12$, $16$, $22$, $18$, $4$, $20$, $2$, $10$, $8$, $14)$ | $(82$, $16$, $86$, $24$, $36$, $26$, $84$, $-42$, $-54$, $-64$, $-60$, $-50$, $-38$, $-46$, $-56$, $-66$, $-58$, $-48$, $72$, $120$, $110$, $108$, $122$, $70$, $74$, $118$, $112$, $80$, $68$, $76$, $116$, $114$, $78$, $-94$, $-104$, $-10$, $-8$, $-102$, $-92$, $-96$, $-106$, $-44$, $-40$, $-52$, $-62$, $20$, $32$, $30$, $18$, $88$, $22$, $34$, $28$, $-2$, $-14$, $-4$, $-98$, $-90$, $-100$, $-6$, $-12)$\\

\noindent$11a321$ : $3$ | $(6$, $12$, $16$, $22$, $18$, $20$, $2$, $8$, $4$, $10$, $14)$ | $(44$, $30$, $40$, $122$, $116$, $34$, $36$, $118$, $120$, $38$, $32$, $114$, $124$, $42$, $-66$, $-104$, $-58$, $-96$, $-56$, $-102$, $-64$, $-68$, $-70$, $-72$, $-24$, $-8$, $-16$, $80$, $94$, $84$, $50$, $76$, $90$, $88$, $128$, $112$, $126$, $-26$, $-6$, $-18$, $-14$, $-10$, $-22$, $-2$, $-28$, $-4$, $-20$, $-12$, $82$, $52$, $78$, $92$, $86$, $48$, $74$, $46$, $-106$, $-60$, $-98$, $-54$, $-100$, $-62$, $-108$, $-110)$\\

\noindent$11a322$ : $3$ | $(6$, $12$, $16$, $22$, $20$, $4$, $18$, $8$, $2$, $10$, $14)$ | $(42$, $52$, $106$, $58$, $36$, $48$, $46$, $38$, $56$, $108$, $54$, $40$, $44$, $50$, $104$, $100$, $-122$, $-118$, $-82$, $-112$, $-74$, $-72$, $-114$, $-84$, $-116$, $-70$, $-76$, $-110$, $-80$, $-120$, $34$, $124$, $32$, $126$, $134$, $88$, $142$, $90$, $132$, $128$, $94$, $138$, $-12$, $-26$, $-2$, $-22$, $-16$, $-8$, $-60$, $-66$, $-64$, $-62$, $-68$, $-78$, $130$, $92$, $140$, $86$, $136$, $96$, $98$, $102$, $-30$, $-6$, $-18$, $-20$, $-4$, $-28$, $-10$, $-14$, $-24)$\\

\noindent$11a323$ : $3$ | $(6$, $12$, $16$, $22$, $20$, $4$, $18$, $10$, $2$, $14$, $8)$ | $(22$, $30$, $58$, $28$, $20$, $24$, $54$, $56$, $26$, $-42$, $-62$, $-64$, $-40$, $-44$, $-60$, $-66$, $-70$, $-68$, $-14$, $78$, $86$, $84$, $76$, $36$, $32$, $34$, $-48$, $-38$, $-46$, $74$, $82$, $88$, $80$, $72$, $52$, $-50$, $-4$, $-6$, $-16$, $-12$, $-2$, $-8$, $-18$, $-10)$\\

\noindent$11a324$ : $3$ | $(6$, $12$, $16$, $22$, $20$, $18$, $2$, $8$, $4$, $10$, $14)$ | $(94$, $182$, $160$, $186$, $90$, $98$, $178$, $164$, $190$, $86$, $-84$, $-18$, $-146$, $-8$, $-156$, $-2$, $-152$, $-12$, $-78$, $174$, $168$, $194$, $24$, $106$, $116$, $34$, $198$, $28$, $110$, $112$, $30$, $200$, $32$, $114$, $108$, $26$, $196$, $170$, $172$, $176$, $166$, $192$, $-132$, $-44$, $-72$, $-54$, $-62$, $-64$, $-52$, $-74$, $-42$, $-82$, $-16$, $-148$, $-6$, $-158$, $-4$, $-150$, $-14$, $-80$, $-40$, $-76$, $36$, $118$, $104$, $22$, $20$, $102$, $120$, $38$, $122$, $100$, $88$, $188$, $162$, $180$, $96$, $92$, $184$, $-58$, $-68$, $-48$, $-128$, $-136$, $-142$, $-124$, $-140$, $-138$, $-126$, $-50$, $-66$, $-60$, $-56$, $-70$, $-46$, $-130$, $-134$, $-144$, $-10$, $-154)$\\

\noindent$11a325$ : $3$ | $(6$, $12$, $16$, $22$, $20$, $18$, $2$, $10$, $4$, $14$, $8)$ | $(38$, $74$, $76$, $40$, $36$, $72$, $78$, $42$, $34$, $70$, $32$, $44$, $80$, $48$, $-60$, $-52$, $-68$, $-90$, $-98$, $-100$, $-92$, $-66$, $-54$, $-58$, $-62$, $124$, $130$, $110$, $108$, $128$, $126$, $106$, $112$, $132$, $122$, $-88$, $-96$, $-102$, $-94$, $-64$, $-56$, $46$, $30$, $50$, $120$, $134$, $114$, $104$, $118$, $136$, $116$, $-4$, $-12$, $-26$, $-82$, $-24$, $-14$, $-6$, $-2$, $-10$, $-18$, $-20$, $-86$, $-28$, $-84$, $-22$, $-16$, $-8)$\\

\noindent$11a326$ : $3$ | $(6$, $12$, $18$, $14$, $22$, $4$, $20$, $8$, $10$, $2$, $16)$ | $(52$, $116$, $58$, $46$, $142$, $44$, $150$, $154$, $170$, $168$, $156$, $148$, $42$, $144$, $48$, $56$, $118$, $54$, $50$, $146$, $-126$, $-68$, $-140$, $-72$, $-122$, $-80$, $-78$, $-120$, $-74$, $-84$, $-18$, $112$, $90$, $106$, $100$, $96$, $174$, $164$, $160$, $178$, $158$, $166$, $172$, $152$, $-14$, $-22$, $-38$, $-28$, $-8$, $-10$, $-26$, $-40$, $-24$, $-12$, $-86$, $-16$, $-20$, $-82$, $-76$, $162$, $176$, $94$, $102$, $104$, $92$, $114$, $60$, $110$, $88$, $108$, $98$, $-70$, $-124$, $-128$, $-66$, $-138$, $-136$, $-64$, $-130$, $-2$, $-34$, $-32$, $-4$, $-132$, $-62$, $-134$, $-6$, $-30$, $-36)$\\

\noindent$11a327$ : $3$ | $(6$, $12$, $18$, $16$, $14$, $4$, $20$, $22$, $2$, $10$, $8)$ | $(124$, $140$, $180$, $174$, $134$, $130$, $170$, $184$, $144$, $186$, $168$, $128$, $136$, $176$, $178$, $138$, $126$, $166$, $188$, $142$, $182$, $172$, $132$, $-10$, $-24$, $-32$, $-2$, $-38$, $-18$, $-16$, $-40$, $-4$, $-30$, $-26$, $-8$, $-44$, $-12$, $-22$, $-34$, $-164$, $-36$, $-20$, $-14$, $-42$, $-6$, $-28$, $66$, $50$, $76$, $82$, $56$, $60$, $86$, $72$, $46$, $70$, $88$, $62$, $54$, $80$, $78$, $-162$, $-160$, $-96$, $-150$, $-108$, $-110$, $-152$, $-94$, $-158$, $-116$, $-102$, $58$, $84$, $74$, $48$, $68$, $90$, $64$, $52$, $192$, $190$, $-120$, $-98$, $-148$, $-106$, $-112$, $-154$, $-92$, $-156$, $-114$, $-104$, $-146$, $-100$, $-118$, $-122)$\\

\noindent$11a328$ : $3$ | $(6$, $12$, $18$, $16$, $22$, $4$, $20$, $8$, $2$, $10$, $14)$ | $(42$, $120$, $144$, $134$, $130$, $148$, $116$, $46$, $38$, $124$, $140$, $138$, $126$, $36$, $48$, $114$, $50$, $-110$, $-66$, $-100$, $-82$, $-84$, $-98$, $-68$, $-112$, $-150$, $-10$, $-166$, $-6$, $-154$, $-22$, $174$, $198$, $204$, $188$, $184$, $52$, $170$, $58$, $178$, $194$, $208$, $192$, $180$, $56$, $168$, $54$, $182$, $190$, $206$, $196$, $176$, $60$, $172$, $200$, $202$, $186$, $-70$, $-96$, $-86$, $-80$, $-102$, $-64$, $-108$, $-74$, $-92$, $-90$, $-76$, $-106$, $-62$, $-104$, $-78$, $-88$, $-94$, $-72$, $128$, $136$, $142$, $122$, $40$, $44$, $118$, $146$, $132$, $-8$, $-152$, $-24$, $-20$, $-156$, $-4$, $-164$, $-12$, $-32$, $-30$, $-14$, $-162$, $-2$, $-158$, $-18$, $-26$, $-34$, $-28$, $-16$, $-160)$\\

\noindent$11a329$ : $3$ | $(6$, $12$, $18$, $22$, $14$, $4$, $20$, $8$, $2$, $10$, $16)$ | $(44$, $238$, $34$, $54$, $72$, $62$, $26$, $22$, $66$, $68$, $-168$, $-198$, $-196$, $-170$, $-224$, $-216$, $-178$, $-188$, $-206$, $-230$, $-210$, $-184$, $-182$, $-212$, $-228$, $-204$, $-190$, $-176$, $-218$, $-222$, $-172$, $-194$, $-200$, $-166$, $-202$, $-192$, $-174$, $-220$, $58$, $30$, $242$, $48$, $40$, $234$, $38$, $50$, $244$, $28$, $60$, $74$, $56$, $32$, $240$, $46$, $42$, $236$, $36$, $52$, $70$, $64$, $24$, $-122$, $-14$, $-94$, $-80$, $-104$, $-4$, $-112$, $-88$, $-86$, $-110$, $-2$, $-106$, $-82$, $-92$, $-116$, $-8$, $-100$, $-76$, $-98$, $-10$, $-118$, $-18$, $246$, $126$, $252$, $160$, $154$, $258$, $132$, $144$, $268$, $142$, $134$, $260$, $152$, $162$, $250$, $124$, $248$, $164$, $150$, $262$, $136$, $140$, $266$, $146$, $130$, $256$, $156$, $158$, $254$, $128$, $148$, $264$, $138$, $-232$, $-208$, $-186$, $-180$, $-214$, $-226$, $-20$, $-16$, $-120$, $-12$, $-96$, $-78$, $-102$, $-6$, $-114$, $-90$, $-84$, $-108)$\\

\noindent$11a330$ : $3$ | $(6$, $12$, $18$, $22$, $16$, $4$, $20$, $8$, $10$, $2$, $14)$ | $(30$, $22$, $38$, $36$, $24$, $28$, $32$, $20$, $-74$, $-72$, $-76$, $-82$, $-46$, $-44$, $-80$, $-78$, $-42$, $-48$, $-84$, $-50$, $88$, $96$, $64$, $94$, $86$, $90$, $18$, $-16$, $-52$, $-14$, $-2$, $-10$, $26$, $34$, $40$, $58$, $54$, $56$, $60$, $98$, $62$, $92$, $-6$, $-68$, $-70$, $-4$, $-8$, $-66$, $-12)$\\

\noindent$11a331$ : $3$ | $(6$, $12$, $18$, $22$, $16$, $4$, $20$, $10$, $8$, $2$, $14)$ | $(44$, $24$, $40$, $128$, $30$, $34$, $124$, $36$, $28$, $130$, $42$, $-114$, $-64$, $-50$, $-120$, $-58$, $-56$, $-122$, $-52$, $-62$, $-116$, $-112$, $-110$, $68$, $84$, $94$, $74$, $78$, $98$, $80$, $72$, $92$, $86$, $-108$, $-20$, $-6$, $-102$, $-14$, $-12$, $-100$, $-8$, $-18$, $-106$, $-2$, $-22$, $-4$, $-104$, $-16$, $-10$, $32$, $126$, $38$, $26$, $132$, $134$, $136$, $88$, $90$, $70$, $82$, $96$, $76$, $-54$, $-60$, $-118$, $-48$, $-66$, $-46)$\\

\noindent$11a332$ : $3$ | $(6$, $14$, $10$, $18$, $2$, $22$, $20$, $4$, $12$, $8$, $16)$ | $(64$, $116$, $104$, $184$, $190$, $166$, $168$, $192$, $182$, $106$, $114$, $66$, $62$, $118$, $102$, $186$, $188$, $100$, $120$, $60$, $68$, $112$, $108$, $-18$, $-6$, $-30$, $-34$, $-10$, $-14$, $-96$, $-82$, $-136$, $-134$, $-84$, $-94$, $-16$, $-8$, $-32$, $52$, $162$, $172$, $196$, $178$, $156$, $46$, $154$, $180$, $194$, $170$, $164$, $-78$, $-140$, $-130$, $-88$, $-90$, $-92$, $-86$, $-132$, $-138$, $-80$, $-98$, $-12$, $56$, $72$, $48$, $158$, $176$, $198$, $174$, $160$, $50$, $74$, $54$, $122$, $58$, $70$, $110$, $-44$, $-20$, $-4$, $-28$, $-36$, $-146$, $-124$, $-148$, $-38$, $-26$, $-2$, $-22$, $-42$, $-152$, $-128$, $-142$, $-76$, $-144$, $-126$, $-150$, $-40$, $-24)$\\

\noindent$11a333$ : $2$ | $(6$, $14$, $12$, $16$, $20$, $18$, $4$, $2$, $22$, $10$, $8)$ | $(78$, $106$, $124$, $96$, $68$, $88$, $116$, $114$, $86$, $70$, $98$, $126$, $104$, $76$, $80$, $108$, $122$, $94$, $66$, $90$, $118$, $112$, $84$, $72$, $100$, $128$, $102$, $74$, $82$, $110$, $120$, $92$, $-142$, $-170$, $-188$, $-160$, $-132$, $-152$, $-180$, $-178$, $-150$, $-134$, $-162$, $-190$, $-168$, $-140$, $-144$, $-172$, $-186$, $-158$, $-130$, $-154$, $-182$, $-176$, $-148$, $-136$, $-164$, $-192$, $-166$, $-138$, $-146$, $-174$, $-184$, $-156$, $206$, $234$, $252$, $224$, $196$, $216$, $244$, $242$, $214$, $198$, $226$, $254$, $232$, $204$, $208$, $236$, $250$, $222$, $194$, $218$, $246$, $240$, $212$, $200$, $228$, $256$, $230$, $202$, $210$, $238$, $248$, $220$, $-14$, $-42$, $-60$, $-32$, $-4$, $-24$, $-52$, $-50$, $-22$, $-6$, $-34$, $-62$, $-40$, $-12$, $-16$, $-44$, $-58$, $-30$, $-2$, $-26$, $-54$, $-48$, $-20$, $-8$, $-36$, $-64$, $-38$, $-10$, $-18$, $-46$, $-56$, $-28)$\\

\noindent$11a334$ : $2$ | $(6$, $14$, $16$, $18$, $20$, $22$, $2$, $4$, $12$, $8$, $10)$ | $(88$, $70$, $52$, $62$, $80$, $96$, $78$, $60$, $54$, $72$, $90$, $86$, $68$, $50$, $64$, $82$, $94$, $76$, $58$, $56$, $74$, $92$, $84$, $66$, $-134$, $-112$, $-102$, $-124$, $-144$, $-122$, $-100$, $-114$, $-136$, $-132$, $-110$, $-104$, $-126$, $-142$, $-120$, $-98$, $-116$, $-138$, $-130$, $-108$, $-106$, $-128$, $-140$, $-118$, $184$, $166$, $148$, $158$, $176$, $192$, $174$, $156$, $150$, $168$, $186$, $182$, $164$, $146$, $160$, $178$, $190$, $172$, $154$, $152$, $170$, $188$, $180$, $162$, $-38$, $-16$, $-6$, $-28$, $-48$, $-26$, $-4$, $-18$, $-40$, $-36$, $-14$, $-8$, $-30$, $-46$, $-24$, $-2$, $-20$, $-42$, $-34$, $-12$, $-10$, $-32$, $-44$, $-22)$\\

\noindent$11a335$ : $2$ | $(6$, $14$, $16$, $18$, $20$, $22$, $2$, $4$, $12$, $10$, $8)$ | $(132$, $98$, $84$, $118$, $146$, $112$, $78$, $104$, $138$, $126$, $92$, $90$, $124$, $140$, $106$, $76$, $110$, $144$, $120$, $86$, $96$, $130$, $134$, $100$, $82$, $116$, $148$, $114$, $80$, $102$, $136$, $128$, $94$, $88$, $122$, $142$, $108$, $-170$, $-214$, $-188$, $-152$, $-196$, $-206$, $-162$, $-178$, $-222$, $-180$, $-160$, $-204$, $-198$, $-154$, $-186$, $-216$, $-172$, $-168$, $-212$, $-190$, $-150$, $-194$, $-208$, $-164$, $-176$, $-220$, $-182$, $-158$, $-202$, $-200$, $-156$, $-184$, $-218$, $-174$, $-166$, $-210$, $-192$, $280$, $246$, $232$, $266$, $294$, $260$, $226$, $252$, $286$, $274$, $240$, $238$, $272$, $288$, $254$, $224$, $258$, $292$, $268$, $234$, $244$, $278$, $282$, $248$, $230$, $264$, $296$, $262$, $228$, $250$, $284$, $276$, $242$, $236$, $270$, $290$, $256$, $-22$, $-66$, $-40$, $-4$, $-48$, $-58$, $-14$, $-30$, $-74$, $-32$, $-12$, $-56$, $-50$, $-6$, $-38$, $-68$, $-24$, $-20$, $-64$, $-42$, $-2$, $-46$, $-60$, $-16$, $-28$, $-72$, $-34$, $-10$, $-54$, $-52$, $-8$, $-36$, $-70$, $-26$, $-18$, $-62$, $-44)$\\

\noindent$11a336$ : $2$ | $(6$, $14$, $16$, $18$, $20$, $22$, $4$, $2$, $12$, $8$, $10)$ |  $(106$, $84$, $62$, $76$, $98$, $114$, $92$, $70$, $68$, $90$, $112$, $100$, $78$, $60$, $82$, $104$, $108$, $86$, $64$, $74$, $96$, $116$, $94$, $72$, $66$, $88$, $110$, $102$, $80$, $-132$, $-164$, $-154$, $-122$, $-142$, $-174$, $-144$, $-120$, $-152$, $-166$, $-134$, $-130$, $-162$, $-156$, $-124$, $-140$, $-172$, $-146$, $-118$, $-150$, $-168$, $-136$, $-128$, $-160$, $-158$, $-126$, $-138$, $-170$, $-148$, $222$, $200$, $178$, $192$, $214$, $230$, $208$, $186$, $184$, $206$, $228$, $216$, $194$, $176$, $198$, $220$, $224$, $202$, $180$, $190$, $212$, $232$, $210$, $188$, $182$, $204$, $226$, $218$, $196$, $-16$, $-48$, $-38$, $-6$, $-26$, $-58$, $-28$, $-4$, $-36$, $-50$, $-18$, $-14$, $-46$, $-40$, $-8$, $-24$, $-56$, $-30$, $-2$, $-34$, $-52$, $-20$, $-12$, $-44$, $-42$, $-10$, $-22$, $-54$, $-32)$\\

\noindent$11a337$ : $2$ | $(6$, $14$, $16$, $18$, $20$, $22$, $4$, $2$, $12$, $10$, $8)$ | $(112$, $160$, $146$, $98$, $126$, $174$, $132$, $92$, $140$, $166$, $118$, $106$, $154$, $152$, $104$, $120$, $168$, $138$, $90$, $134$, $172$, $124$, $100$, $148$, $158$, $110$, $114$, $162$, $144$, $96$, $128$, $176$, $130$, $94$, $142$, $164$, $116$, $108$, $156$, $150$, $102$, $122$, $170$, $136$, $-202$, $-254$, $-224$, $-180$, $-232$, $-246$, $-194$, $-210$, $-262$, $-216$, $-188$, $-240$, $-238$, $-186$, $-218$, $-260$, $-208$, $-196$, $-248$, $-230$, $-178$, $-226$, $-252$, $-200$, $-204$, $-256$, $-222$, $-182$, $-234$, $-244$, $-192$, $-212$, $-264$, $-214$, $-190$, $-242$, $-236$, $-184$, $-220$, $-258$, $-206$, $-198$, $-250$, $-228$, $288$, $336$, $322$, $274$, $302$, $350$, $308$, $268$, $316$, $342$, $294$, $282$, $330$, $328$, $280$, $296$, $344$, $314$, $266$, $310$, $348$, $300$, $276$, $324$, $334$, $286$, $290$, $338$, $320$, $272$, $304$, $352$, $306$, $270$, $318$, $340$, $292$, $284$, $332$, $326$, $278$, $298$, $346$, $312$, $-26$, $-78$, $-48$, $-4$, $-56$, $-70$, $-18$, $-34$, $-86$, $-40$, $-12$, $-64$, $-62$, $-10$, $-42$, $-84$, $-32$, $-20$, $-72$, $-54$, $-2$, $-50$, $-76$, $-24$, $-28$, $-80$, $-46$, $-6$, $-58$, $-68$, $-16$, $-36$, $-88$, $-38$, $-14$, $-66$, $-60$, $-8$, $-44$, $-82$, $-30$, $-22$, $-74$, $-52)$\\

\noindent$11a338$ : $3$ | $(6$, $14$, $16$, $20$, $22$, $18$, $2$, $4$, $10$, $12$, $8)$ | $(46$, $36$, $30$, $40$, $20$, $44$, $34$, $32$, $42$, $-70$, $-24$, $50$, $48$, $38$, $-56$, $-66$, $-74$, $-64$, $-58$, $-54$, $-68$, $-72$, $-62$, $-60$, $-52$, $22$, $90$, $80$, $78$, $88$, $96$, $86$, $76$, $82$, $92$, $94$, $84$, $-14$, $-4$, $-28$, $-6$, $-16$, $-12$, $-2$, $-26$, $-8$, $-18$, $-10)$\\

\noindent$11a339$ : $2$ | $(6$, $14$, $16$, $20$, $22$, $18$, $2$, $4$, $12$, $10$, $8)$ | $(96$, $70$, $64$, $90$, $102$, $76$, $58$, $84$, $108$, $82$, $56$, $78$, $104$, $88$, $62$, $72$, $98$, $94$, $68$, $66$, $92$, $100$, $74$, $60$, $86$, $106$, $80$, $-146$, $-112$, $-138$, $-154$, $-120$, $-130$, $-162$, $-128$, $-122$, $-156$, $-136$, $-114$, $-148$, $-144$, $-110$, $-140$, $-152$, $-118$, $-132$, $-160$, $-126$, $-124$, $-158$, $-134$, $-116$, $-150$, $-142$, $204$, $178$, $172$, $198$, $210$, $184$, $166$, $192$, $216$, $190$, $164$, $186$, $212$, $196$, $170$, $180$, $206$, $202$, $176$, $174$, $200$, $208$, $182$, $168$, $194$, $214$, $188$, $-38$, $-4$, $-30$, $-46$, $-12$, $-22$, $-54$, $-20$, $-14$, $-48$, $-28$, $-6$, $-40$, $-36$, $-2$, $-32$, $-44$, $-10$, $-24$, $-52$, $-18$, $-16$, $-50$, $-26$, $-8$, $-42$, $-34)$\\

\noindent$11a340$ : $3$ | $(6$, $14$, $16$, $20$, $22$, $18$, $4$, $2$, $10$, $12$, $8)$ | $(42$, $58$, $52$, $36$, $40$, $56$, $54$, $38$, $-8$, $-66$, $46$, $62$, $48$, $18$, $44$, $60$, $50$, $-74$, $-88$, $-90$, $-82$, $-80$, $-68$, $-64$, $-72$, $-76$, $-86$, $-92$, $-84$, $-78$, $-70$, $100$, $116$, $118$, $102$, $98$, $114$, $112$, $96$, $104$, $120$, $106$, $94$, $110$, $16$, $108$, $-12$, $-4$, $-32$, $-24$, $-22$, $-30$, $-2$, $-10$, $-14$, $-6$, $-34$, $-26$, $-20$, $-28)$\\

\noindent$11a341$ : $2$ | $(6$, $14$, $16$, $20$, $22$, $18$, $4$, $2$, $12$, $10$, $8)$ | $(76$, $108$, $102$, $70$, $82$, $114$, $96$, $64$, $88$, $120$, $90$, $62$, $94$, $116$, $84$, $68$, $100$, $110$, $78$, $74$, $106$, $104$, $72$, $80$, $112$, $98$, $66$, $86$, $118$, $92$, $-162$, $-124$, $-154$, $-170$, $-132$, $-146$, $-178$, $-140$, $-138$, $-176$, $-148$, $-130$, $-168$, $-156$, $-122$, $-160$, $-164$, $-126$, $-152$, $-172$, $-134$, $-144$, $-180$, $-142$, $-136$, $-174$, $-150$, $-128$, $-166$, $-158$, $196$, $228$, $222$, $190$, $202$, $234$, $216$, $184$, $208$, $240$, $210$, $182$, $214$, $236$, $204$, $188$, $220$, $230$, $198$, $194$, $226$, $224$, $192$, $200$, $232$, $218$, $186$, $206$, $238$, $212$, $-42$, $-4$, $-34$, $-50$, $-12$, $-26$, $-58$, $-20$, $-18$, $-56$, $-28$, $-10$, $-48$, $-36$, $-2$, $-40$, $-44$, $-6$, $-32$, $-52$, $-14$, $-24$, $-60$, $-22$, $-16$, $-54$, $-30$, $-8$, $-46$, $-38)$\\

\noindent$11a342$ : $2$ | $(6$, $14$, $16$, $22$, $20$, $18$, $2$, $4$, $12$, $10$, $8)$ | $(50$, $36$, $34$, $48$, $52$, $38$, $32$, $46$, $54$, $40$, $30$, $44$, $56$, $42$, $-60$, $-68$, $-76$, $-84$, $-78$, $-70$, $-62$, $-58$, $-66$, $-74$, $-82$, $-80$, $-72$, $-64$, $106$, $92$, $90$, $104$, $108$, $94$, $88$, $102$, $110$, $96$, $86$, $100$, $112$, $98$, $-4$, $-12$, $-20$, $-28$, $-22$, $-14$, $-6$, $-2$, $-10$, $-18$, $-26$, $-24$, $-16$, $-8)$\\

\noindent$11a343$ : $2$ | $(6$, $14$, $16$, $22$, $20$, $18$, $4$, $2$, $12$, $10$, $8)$ | $(38$, $54$, $52$, $36$, $40$, $56$, $50$, $34$, $42$, $58$, $48$, $32$, $44$, $60$, $46$, $-64$, $-72$, $-80$, $-88$, $-86$, $-78$, $-70$, $-62$, $-66$, $-74$, $-82$, $-90$, $-84$, $-76$, $-68$, $98$, $114$, $112$, $96$, $100$, $116$, $110$, $94$, $102$, $118$, $108$, $92$, $104$, $120$, $106$, $-4$, $-12$, $-20$, $-28$, $-26$, $-18$, $-10$, $-2$, $-6$, $-14$, $-22$, $-30$, $-24$, $-16$, $-8)$\\

\noindent$11a344$ : $3$ | $(6$, $14$, $18$, $16$, $20$, $22$, $4$, $8$, $2$, $12$, $10)$ | $(70$, $42$, $34$, $62$, $32$, $44$, $30$, $28$, $120$, $108$, $96$, $110$, $118$, $-90$, $-78$, $-46$, $-60$, $-50$, $-22$, $-20$, $-52$, $-58$, $102$, $116$, $112$, $98$, $106$, $122$, $104$, $100$, $114$, $-48$, $-24$, $-18$, $-54$, $-56$, $-16$, $-26$, $-14$, $72$, $68$, $40$, $36$, $64$, $76$, $74$, $66$, $38$, $-84$, $-4$, $-10$, $-12$, $-2$, $-82$, $-94$, $-86$, $-6$, $-8$, $-88$, $-92$, $-80)$\\

\noindent$11a345$ : $3$ | $(6$, $14$, $18$, $16$, $22$, $20$, $4$, $8$, $2$, $12$, $10)$ | $(98$, $124$, $112$, $86$, $110$, $136$, $134$, $108$, $88$, $114$, $122$, $96$, $100$, $126$, $140$, $130$, $104$, $92$, $118$, $-150$, $-20$, $-4$, $186$, $180$, $194$, $220$, $208$, $38$, $206$, $222$, $196$, $182$, $184$, $198$, $224$, $204$, $40$, $210$, $218$, $192$, $178$, $188$, $214$, $-172$, $-76$, $-52$, $-58$, $-70$, $-166$, $-164$, $-68$, $-60$, $-50$, $-78$, $-174$, $-84$, $-170$, $-74$, $-54$, $-56$, $-72$, $-168$, $-162$, $-66$, $-62$, $-48$, $-80$, $-176$, $-82$, $-46$, $-64$, $128$, $102$, $94$, $120$, $116$, $90$, $106$, $132$, $138$, $44$, $200$, $226$, $202$, $42$, $212$, $216$, $190$, $-12$, $-158$, $-30$, $-142$, $-28$, $-160$, $-10$, $-14$, $-156$, $-32$, $-144$, $-26$, $-2$, $-22$, $-148$, $-36$, $-152$, $-18$, $-6$, $-8$, $-16$, $-154$, $-34$, $-146$, $-24)$\\

\noindent$11a346$ : $3$ | $(6$, $16$, $12$, $18$, $20$, $4$, $22$, $10$, $2$, $8$, $14)$ | $(38$, $30$, $22$, $24$, $82$, $84$, $92$, $94$, $86$, $64$, $-78$, $-68$, $-80$, $-18$, $-76$, $-70$, $-12$, $-16$, $-74$, $-72$, $-14$, $34$, $26$, $20$, $28$, $36$, $40$, $32$, $-10$, $-8$, $-4$, $-48$, $-46$, $-6$, $60$, $90$, $96$, $88$, $62$, $66$, $58$, $-56$, $-44$, $-50$, $-2$, $-54$, $-42$, $-52)$\\

\noindent$11a347$ : $3$ | $(6$, $16$, $12$, $20$, $18$, $4$, $22$, $10$, $2$, $8$, $14)$ | $(32$, $64$, $26$, $70$, $24$, $66$, $30$, $60$, $62$, $28$, $68$, $-46$, $-72$, $-50$, $-42$, $-76$, $-78$, $-82$, $-80$, $-16$, $102$, $88$, $96$, $94$, $90$, $100$, $38$, $34$, $36$, $-54$, $-40$, $-52$, $-74$, $-44$, $-48$, $92$, $98$, $86$, $104$, $84$, $58$, $-56$, $-2$, $-12$, $-20$, $-6$, $-8$, $-22$, $-10$, $-4$, $-18$, $-14)$\\

\noindent$11a348$ : $3$ | $(6$, $18$, $16$, $12$, $4$, $2$, $20$, $22$, $10$, $8$, $14)$ | $(70$, $242$, $250$, $78$, $62$, $234$, $258$, $86$, $54$, $228$, $56$, $84$, $256$, $236$, $64$, $76$, $248$, $244$, $72$, $68$, $240$, $252$, $80$, $60$, $232$, $-90$, $-124$, $-226$, $-120$, $-94$, $-200$, $-218$, $-112$, $-102$, $-208$, $-210$, $-104$, $-110$, $-216$, $-202$, $-96$, $-118$, $-224$, $-126$, $260$, $162$, $154$, $268$, $132$, $170$, $146$, $276$, $140$, $176$, $138$, $274$, $148$, $168$, $130$, $266$, $156$, $160$, $262$, $-14$, $-24$, $-186$, $-44$, $-6$, $-32$, $-178$, $-36$, $-2$, $-40$, $-182$, $-28$, $-10$, $-48$, $-190$, $-20$, $-18$, $-192$, $-50$, $-12$, $-26$, $-184$, $-42$, $-4$, $-34$, $246$, $74$, $66$, $238$, $254$, $82$, $58$, $230$, $52$, $88$, $164$, $152$, $270$, $134$, $172$, $144$, $278$, $142$, $174$, $136$, $272$, $150$, $166$, $128$, $264$, $158$, $-122$, $-92$, $-198$, $-220$, $-114$, $-100$, $-206$, $-212$, $-106$, $-108$, $-214$, $-204$, $-98$, $-116$, $-222$, $-196$, $-194$, $-16$, $-22$, $-188$, $-46$, $-8$, $-30$, $-180$, $-38)$\\

\noindent$11a349$ : $3$ | $(6$, $18$, $16$, $12$, $4$, $2$, $22$, $20$, $10$, $8$, $14)$ | $(42$, $64$, $58$, $52$, $70$, $48$, $36$, $100$, $38$, $46$, $68$, $54$, $56$, $66$, $44$, $40$, $98$, $-80$, $-102$, $-84$, $-116$, $-114$, $-86$, $-104$, $-126$, $-122$, $-108$, $-90$, $-110$, $-120$, $-128$, $-118$, $-112$, $-88$, $-106$, $-124$, $50$, $60$, $62$, $156$, $154$, $150$, $160$, $144$, $134$, $-24$, $-6$, $-74$, $-78$, $-82$, $152$, $158$, $142$, $136$, $92$, $132$, $146$, $162$, $148$, $130$, $94$, $138$, $140$, $96$, $-76$, $-4$, $-26$, $-22$, $-8$, $-72$, $-10$, $-20$, $-28$, $-2$, $-32$, $-16$, $-14$, $-34$, $-12$, $-18$, $-30)$\\

\noindent$11a350$ : $3$ | $(6$, $18$, $16$, $14$, $4$, $2$, $20$, $22$, $10$, $12$, $8)$ | $(36$, $28$, $176$, $194$, $196$, $174$, $30$, $34$, $170$, $200$, $190$, $180$, $-64$, $-86$, $-52$, $-76$, $-74$, $-54$, $-160$, $-164$, $-154$, $-20$, $-144$, $-8$, $108$, $126$, $136$, $120$, $114$, $42$, $102$, $132$, $130$, $104$, $44$, $112$, $122$, $138$, $124$, $110$, $46$, $106$, $128$, $134$, $100$, $40$, $116$, $118$, $38$, $-162$, $-152$, $-22$, $-146$, $-6$, $-10$, $-142$, $-18$, $-156$, $-166$, $-158$, $-16$, $-140$, $-12$, $-4$, $-148$, $-24$, $-150$, $-2$, $-14$, $32$, $172$, $198$, $192$, $178$, $26$, $182$, $188$, $202$, $168$, $206$, $184$, $186$, $204$, $-94$, $-56$, $-72$, $-78$, $-50$, $-84$, $-66$, $-62$, $-88$, $-98$, $-90$, $-60$, $-68$, $-82$, $-48$, $-80$, $-70$, $-58$, $-92$, $-96)$\\

\noindent$11a351$ : $3$ | $(6$, $18$, $16$, $14$, $20$, $4$, $2$, $22$, $12$, $8$, $10)$ | $(170$, $284$, $300$, $112$, $310$, $78$, $86$, $318$, $104$, $204$, $102$, $320$, $88$, $76$, $308$, $114$, $302$, $282$, $172$, $168$, $286$, $298$, $110$, $312$, $80$, $84$, $316$, $106$, $202$, $100$, $322$, $90$, $74$, $306$, $116$, $304$, $-120$, $-250$, $-212$, $-234$, $-136$, $-134$, $-236$, $-210$, $-248$, $-122$, $-36$, $-2$, $-40$, $-126$, $-244$, $-206$, $-240$, $-130$, $-140$, $-230$, $-216$, $-254$, $148$, $192$, $280$, $174$, $166$, $288$, $296$, $158$, $182$, $272$, $184$, $156$, $-44$, $-6$, $-32$, $-70$, $-68$, $-30$, $-8$, $-46$, $-268$, $-54$, $-16$, $-22$, $-60$, $-262$, $-224$, $-222$, $-260$, $-62$, $-24$, $-14$, $-52$, $-270$, $-48$, $-10$, $-28$, $-66$, $-72$, $-34$, $-4$, $-42$, $-128$, $-242$, $292$, $162$, $178$, $276$, $188$, $152$, $144$, $196$, $94$, $326$, $96$, $198$, $142$, $154$, $186$, $274$, $180$, $160$, $294$, $290$, $164$, $176$, $278$, $190$, $150$, $146$, $194$, $92$, $324$, $98$, $200$, $108$, $314$, $82$, $-50$, $-12$, $-26$, $-64$, $-258$, $-220$, $-226$, $-264$, $-58$, $-20$, $-18$, $-56$, $-266$, $-228$, $-218$, $-256$, $-118$, $-252$, $-214$, $-232$, $-138$, $-132$, $-238$, $-208$, $-246$, $-124$, $-38)$\\

\noindent$11a352$ : $3$ | $(6$, $18$, $16$, $14$, $20$, $4$, $2$, $22$, $12$, $10$, $8)$ | $(76$, $58$, $80$, $72$, $62$, $114$, $116$, $120$, $-40$, $-46$, $-8$, $-52$, $-2$, $-56$, $-4$, $-50$, $-10$, $-12$, $84$, $68$, $66$, $86$, $64$, $70$, $82$, $60$, $74$, $78$, $-100$, $-106$, $-94$, $-112$, $-88$, $-16$, $-90$, $-110$, $-96$, $-104$, $-102$, $-98$, $-108$, $-92$, $-14$, $118$, $18$, $122$, $32$, $22$, $126$, $28$, $26$, $128$, $24$, $30$, $124$, $20$, $34$, $-36$, $-44$, $-42$, $-38$, $-48$, $-6$, $-54)$\\

\noindent$11a353$ : $3$ | $(8$, $12$, $16$, $18$, $22$, $4$, $20$, $2$, $6$, $10$, $14)$ | $(30$, $18$, $42$, $154$, $148$, $36$, $-138$, $-112$, $-128$, $-122$, $-118$, $-144$, $-142$, $-116$, $-124$, $-126$, $-114$, $-140$, $-136$, $-110$, $-130$, $-120$, $20$, $28$, $32$, $16$, $40$, $152$, $150$, $38$, $14$, $34$, $26$, $22$, $104$, $102$, $24$, $-66$, $-46$, $-2$, $-50$, $-62$, $-12$, $-60$, $-52$, $-4$, $-44$, $-6$, $-54$, $-58$, $-10$, $-146$, $-8$, $-56$, $98$, $88$, $82$, $160$, $76$, $94$, $92$, $78$, $162$, $80$, $90$, $96$, $106$, $100$, $86$, $84$, $158$, $74$, $156$, $-70$, $-134$, $-108$, $-132$, $-68$, $-72$, $-64$, $-48)$\\

\noindent$11a354$ : $3$ | $(8$, $12$, $16$, $20$, $2$, $18$, $22$, $6$, $4$, $10$, $14)$ | $(56$, $26$, $58$, $54$, $64$, $72$, $-86$, $-76$, $-88$, $-10$, $-84$, $-78$, $-90$, $92$, $102$, $94$, $18$, $12$, $20$, $96$, $100$, $24$, $16$, $14$, $22$, $98$, $-40$, $-30$, $-28$, $-42$, $-50$, $-38$, $-32$, $-44$, $-48$, $-36$, $-34$, $-46$, $68$, $60$, $52$, $62$, $70$, $74$, $66$, $-4$, $-2$, $-8$, $-82$, $-80$, $-6)$\\

\noindent$11a355$ : $2$ | $(8$, $14$, $16$, $18$, $20$, $22$, $2$, $4$, $6$, $12$, $10)$ | $(82$, $68$, $54$, $48$, $62$, $76$, $88$, $74$, $60$, $46$, $56$, $70$, $84$, $80$, $66$, $52$, $50$, $64$, $78$, $86$, $72$, $58$, $-120$, $-94$, $-108$, $-132$, $-106$, $-96$, $-122$, $-118$, $-92$, $-110$, $-130$, $-104$, $-98$, $-124$, $-116$, $-90$, $-112$, $-128$, $-102$, $-100$, $-126$, $-114$, $170$, $156$, $142$, $136$, $150$, $164$, $176$, $162$, $148$, $134$, $144$, $158$, $172$, $168$, $154$, $140$, $138$, $152$, $166$, $174$, $160$, $146$, $-32$, $-6$, $-20$, $-44$, $-18$, $-8$, $-34$, $-30$, $-4$, $-22$, $-42$, $-16$, $-10$, $-36$, $-28$, $-2$, $-24$, $-40$, $-14$, $-12$, $-38$, $-26)$\\

\noindent$11a356$ : $2$ | $(8$, $14$, $16$, $18$, $20$, $22$, $2$, $6$, $4$, $12$, $10)$ | $(102$, $150$, $116$, $88$, $136$, $130$, $82$, $122$, $144$, $96$, $108$, $156$, $110$, $94$, $142$, $124$, $80$, $128$, $138$, $90$, $114$, $152$, $104$, $100$, $148$, $118$, $86$, $134$, $132$, $84$, $120$, $146$, $98$, $106$, $154$, $112$, $92$, $140$, $126$, $-212$, $-166$, $-192$, $-232$, $-186$, $-172$, $-218$, $-206$, $-160$, $-198$, $-226$, $-180$, $-178$, $-224$, $-200$, $-158$, $-204$, $-220$, $-174$, $-184$, $-230$, $-194$, $-164$, $-210$, $-214$, $-168$, $-190$, $-234$, $-188$, $-170$, $-216$, $-208$, $-162$, $-196$, $-228$, $-182$, $-176$, $-222$, $-202$, $258$, $306$, $272$, $244$, $292$, $286$, $238$, $278$, $300$, $252$, $264$, $312$, $266$, $250$, $298$, $280$, $236$, $284$, $294$, $246$, $270$, $308$, $260$, $256$, $304$, $274$, $242$, $290$, $288$, $240$, $276$, $302$, $254$, $262$, $310$, $268$, $248$, $296$, $282$, $-56$, $-10$, $-36$, $-76$, $-30$, $-16$, $-62$, $-50$, $-4$, $-42$, $-70$, $-24$, $-22$, $-68$, $-44$, $-2$, $-48$, $-64$, $-18$, $-28$, $-74$, $-38$, $-8$, $-54$, $-58$, $-12$, $-34$, $-78$, $-32$, $-14$, $-60$, $-52$, $-6$, $-40$, $-72$, $-26$, $-20$, $-66$, $-46)$\\

\noindent$11a357$ : $2$ | $(8$, $14$, $16$, $18$, $20$, $22$, $6$, $4$, $2$, $12$, $10)$ | $(154$, $100$, $134$, $174$, $120$, $114$, $168$, $140$, $94$, $148$, $160$, $106$, $128$, $180$, $126$, $108$, $162$, $146$, $92$, $142$, $166$, $112$, $122$, $176$, $132$, $102$, $156$, $152$, $98$, $136$, $172$, $118$, $116$, $170$, $138$, $96$, $150$, $158$, $104$, $130$, $178$, $124$, $110$, $164$, $144$, $-244$, $-190$, $-224$, $-264$, $-210$, $-204$, $-258$, $-230$, $-184$, $-238$, $-250$, $-196$, $-218$, $-270$, $-216$, $-198$, $-252$, $-236$, $-182$, $-232$, $-256$, $-202$, $-212$, $-266$, $-222$, $-192$, $-246$, $-242$, $-188$, $-226$, $-262$, $-208$, $-206$, $-260$, $-228$, $-186$, $-240$, $-248$, $-194$, $-220$, $-268$, $-214$, $-200$, $-254$, $-234$, $334$, $280$, $314$, $354$, $300$, $294$, $348$, $320$, $274$, $328$, $340$, $286$, $308$, $360$, $306$, $288$, $342$, $326$, $272$, $322$, $346$, $292$, $302$, $356$, $312$, $282$, $336$, $332$, $278$, $316$, $352$, $298$, $296$, $350$, $318$, $276$, $330$, $338$, $284$, $310$, $358$, $304$, $290$, $344$, $324$, $-64$, $-10$, $-44$, $-84$, $-30$, $-24$, $-78$, $-50$, $-4$, $-58$, $-70$, $-16$, $-38$, $-90$, $-36$, $-18$, $-72$, $-56$, $-2$, $-52$, $-76$, $-22$, $-32$, $-86$, $-42$, $-12$, $-66$, $-62$, $-8$, $-46$, $-82$, $-28$, $-26$, $-80$, $-48$, $-6$, $-60$, $-68$, $-14$, $-40$, $-88$, $-34$, $-20$, $-74$, $-54)$\\

\noindent$11a358$ : $2$ | $(8$, $14$, $16$, $18$, $22$, $20$, $2$, $4$, $6$, $12$, $10)$ | $(36$, $48$, $60$, $50$, $38$, $34$, $46$, $58$, $52$, $40$, $32$, $44$, $56$, $54$, $42$, $-86$, $-76$, $-66$, $-64$, $-74$, $-84$, $-88$, $-78$, $-68$, $-62$, $-72$, $-82$, $-90$, $-80$, $-70$, $96$, $108$, $120$, $110$, $98$, $94$, $106$, $118$, $112$, $100$, $92$, $104$, $116$, $114$, $102$, $-26$, $-16$, $-6$, $-4$, $-14$, $-24$, $-28$, $-18$, $-8$, $-2$, $-12$, $-22$, $-30$, $-20$, $-10)$\\

\noindent$11a359$ : $2$ | $(8$, $14$, $16$, $18$, $22$, $20$, $2$, $6$, $4$, $12$, $10)$ | $(62$, $82$, $102$, $88$, $68$, $56$, $76$, $96$, $94$, $74$, $54$, $70$, $90$, $100$, $80$, $60$, $64$, $84$, $104$, $86$, $66$, $58$, $78$, $98$, $92$, $72$, $-120$, $-152$, $-130$, $-110$, $-142$, $-140$, $-108$, $-132$, $-150$, $-118$, $-122$, $-154$, $-128$, $-112$, $-144$, $-138$, $-106$, $-134$, $-148$, $-116$, $-124$, $-156$, $-126$, $-114$, $-146$, $-136$, $166$, $186$, $206$, $192$, $172$, $160$, $180$, $200$, $198$, $178$, $158$, $174$, $194$, $204$, $184$, $164$, $168$, $188$, $208$, $190$, $170$, $162$, $182$, $202$, $196$, $176$, $-16$, $-48$, $-26$, $-6$, $-38$, $-36$, $-4$, $-28$, $-46$, $-14$, $-18$, $-50$, $-24$, $-8$, $-40$, $-34$, $-2$, $-30$, $-44$, $-12$, $-20$, $-52$, $-22$, $-10$, $-42$, $-32)$\\

\noindent$11a360$ : $2$ | $(8$, $14$, $16$, $18$, $22$, $20$, $6$, $4$, $2$, $12$, $10)$ | $(96$, $62$, $84$, $108$, $74$, $72$, $106$, $86$, $60$, $94$, $98$, $64$, $82$, $110$, $76$, $70$, $104$, $88$, $58$, $92$, $100$, $66$, $80$, $112$, $78$, $68$, $102$, $90$, $-122$, $-142$, $-162$, $-156$, $-136$, $-116$, $-128$, $-148$, $-168$, $-150$, $-130$, $-114$, $-134$, $-154$, $-164$, $-144$, $-124$, $-120$, $-140$, $-160$, $-158$, $-138$, $-118$, $-126$, $-146$, $-166$, $-152$, $-132$, $208$, $174$, $196$, $220$, $186$, $184$, $218$, $198$, $172$, $206$, $210$, $176$, $194$, $222$, $188$, $182$, $216$, $200$, $170$, $204$, $212$, $178$, $192$, $224$, $190$, $180$, $214$, $202$, $-10$, $-30$, $-50$, $-44$, $-24$, $-4$, $-16$, $-36$, $-56$, $-38$, $-18$, $-2$, $-22$, $-42$, $-52$, $-32$, $-12$, $-8$, $-28$, $-48$, $-46$, $-26$, $-6$, $-14$, $-34$, $-54$, $-40$, $-20)$\\

\noindent$11a361$ : $3$ | $(8$, $14$, $18$, $16$, $20$, $22$, $4$, $2$, $6$, $12$, $10)$ | $(38$, $44$, $26$, $50$, $24$, $46$, $40$, $36$, $42$, $28$, $48$, $-66$, $-52$, $-62$, $-70$, $-72$, $30$, $94$, $-58$, $-56$, $-34$, $-60$, $-54$, $-68$, $-64$, $84$, $78$, $90$, $96$, $92$, $76$, $86$, $82$, $80$, $88$, $74$, $-32$, $-18$, $-4$, $-10$, $-22$, $-8$, $-6$, $-20$, $-12$, $-2$, $-16$, $-14)$\\

\noindent$11a362$ : $3$ | $(8$, $14$, $18$, $16$, $22$, $20$, $4$, $2$, $6$, $12$, $10)$ | $(28$, $30$, $32$, $36$, $-38$, $-6$, $-40$, $-42$, $14$, $10$, $8$, $12$, $48$, $50$, $-26$, $-24$, $-18$, $-22$, $-20$, $34$, $46$, $44$, $-16$, $-4$, $-2)$\\

\noindent$11a363$ : $2$ | $(8$, $16$, $14$, $18$, $22$, $20$, $6$, $4$, $2$, $12$, $10)$ | $(40$, $52$, $64$, $62$, $50$, $38$, $42$, $54$, $66$, $60$, $48$, $36$, $44$, $56$, $68$, $58$, $46$, $-74$, $-86$, $-98$, $-96$, $-84$, $-72$, $-76$, $-88$, $-100$, $-94$, $-82$, $-70$, $-78$, $-90$, $-102$, $-92$, $-80$, $108$, $120$, $132$, $130$, $118$, $106$, $110$, $122$, $134$, $128$, $116$, $104$, $112$, $124$, $136$, $126$, $114$, $-6$, $-18$, $-30$, $-28$, $-16$, $-4$, $-8$, $-20$, $-32$, $-26$, $-14$, $-2$, $-10$, $-22$, $-34$, $-24$, $-12)$\\

\noindent$11a364$ : $2$ | $(10$, $14$, $16$, $18$, $20$, $22$, $2$, $4$, $6$, $8$, $12)$ | $(32$, $48$, $34$, $30$, $46$, $36$, $28$, $44$, $38$, $26$, $42$, $40$, $-70$, $-64$, $-58$, $-52$, $-50$, $-56$, $-62$, $-68$, $-72$, $-66$, $-60$, $-54$, $80$, $96$, $82$, $78$, $94$, $84$, $76$, $92$, $86$, $74$, $90$, $88$, $-22$, $-16$, $-10$, $-4$, $-2$, $-8$, $-14$, $-20$, $-24$, $-18$, $-12$, $-6)$\\

\noindent$11a365$ : $2$ | $(10$, $14$, $16$, $18$, $20$, $22$, $2$, $4$, $8$, $6$, $12)$ | $(66$, $98$, $72$, $60$, $92$, $78$, $54$, $86$, $84$, $52$, $80$, $90$, $58$, $74$, $96$, $64$, $68$, $100$, $70$, $62$, $94$, $76$, $56$, $88$, $82$, $-116$, $-148$, $-122$, $-110$, $-142$, $-128$, $-104$, $-136$, $-134$, $-102$, $-130$, $-140$, $-108$, $-124$, $-146$, $-114$, $-118$, $-150$, $-120$, $-112$, $-144$, $-126$, $-106$, $-138$, $-132$, $166$, $198$, $172$, $160$, $192$, $178$, $154$, $186$, $184$, $152$, $180$, $190$, $158$, $174$, $196$, $164$, $168$, $200$, $170$, $162$, $194$, $176$, $156$, $188$, $182$, $-16$, $-48$, $-22$, $-10$, $-42$, $-28$, $-4$, $-36$, $-34$, $-2$, $-30$, $-40$, $-8$, $-24$, $-46$, $-14$, $-18$, $-50$, $-20$, $-12$, $-44$, $-26$, $-6$, $-38$, $-32)$\\

\noindent$11a366$ : $3$ | $(10$, $14$, $16$, $18$, $20$, $22$, $4$, $2$, $8$, $6$, $12)$ | $(60$, $174$, $160$, $166$, $168$, $66$, $54$, $180$, $154$, $182$, $186$, $158$, $176$, $58$, $62$, $172$, $162$, $164$, $170$, $64$, $56$, $178$, $156$, $184$, $-48$, $-20$, $-92$, $-130$, $-116$, $-144$, $-142$, $-114$, $-132$, $-90$, $-22$, $-46$, $-18$, $-16$, $-44$, $-24$, $-4$, $-32$, $-36$, $-8$, $106$, $188$, $-14$, $-42$, $-26$, $-2$, $-30$, $-38$, $-10$, $-6$, $-34$, $84$, $102$, $192$, $96$, $78$, $70$, $50$, $72$, $76$, $94$, $190$, $104$, $86$, $108$, $82$, $100$, $194$, $98$, $80$, $68$, $52$, $74$, $-126$, $-120$, $-148$, $-138$, $-110$, $-136$, $-88$, $-134$, $-112$, $-140$, $-146$, $-118$, $-128$, $-152$, $-124$, $-122$, $-150$, $-12$, $-40$, $-28)$\\

\noindent$11a367$ : $2$ | $(12$, $14$, $16$, $18$, $20$, $22$, $2$, $4$, $6$, $8$, $10)$ | $(20$, $18$, $16$, $14$, $12$, $-30$, $-28$, $-26$, $-24$, $-22$, $40$, $38$, $36$, $34$, $32$, $-10$, $-8$, $-6$, $-4$, $-2)$\\

\noindent$11n1$ : $3$ | $(4$, $8$, $10$, $14$, $2$, $-16$, $-20$, $6$, $-22$, $-12$, $-18)$ | $(18$, $28$, $24$, $70$, $68$, $22$, $30$, $20$, $-60$, $-58$, $-32$, $-54$, $-64$, $-62$, $-56$, $40$, $72$, $26$, $-4$, $-12$, $-14$, $-2$, $-34$, $-6$, $-10$, $-8$, $48$, $38$, $42$, $74$, $44$, $36$, $46$, $50$, $-52$, $-66$, $-16)$\\

\noindent$11n2$ : $3$ | $(4$, $8$, $10$, $-14$, $2$, $16$, $20$, $-6$, $22$, $12$, $18)$ | $(28$, $36$, $78$, $48$, $76$, $38$, $86$, $40$, $74$, $46$, $80$, $34$, $30$, $-100$, $-98$, $-66$, $-64$, $-96$, $-52$, $-88$, $-56$, $-70$, $-60$, $-92$, $128$, $140$, $116$, $120$, $136$, $132$, $124$, $144$, $112$, $146$, $122$, $134$, $-58$, $-90$, $-50$, $-94$, $-62$, $-68$, $-54$, $118$, $138$, $130$, $126$, $142$, $114$, $148$, $84$, $42$, $72$, $44$, $82$, $32$, $-24$, $-2$, $-12$, $-14$, $-102$, $-108$, $-20$, $-6$, $-8$, $-18$, $-106$, $-104$, $-16$, $-10$, $-4$, $-22$, $-110$, $-26)$\\

\noindent$11n3$ : $3$ | $(4$, $8$, $10$, $-14$, $2$, $-16$, $-20$, $-6$, $-22$, $-12$, $-18)$ | $(28$, $44$, $102$, $38$, $34$, $32$, $40$, $100$, $42$, $30$, $26$, $46$, $-48$, $-86$, $-84$, $-96$, $-76$, $-94$, $-92$, $-78$, $-98$, $-82$, $-88$, $62$, $104$, $36$, $-8$, $-12$, $-18$, $-2$, $-22$, $-24$, $-4$, $-16$, $-14$, $-6$, $-50$, $-10$, $54$, $70$, $66$, $58$, $106$, $60$, $64$, $72$, $52$, $74$, $56$, $68$, $-80$, $-90$, $-20)$\\

\noindent$11n4$ : $3$ | $(4$, $8$, $10$, $14$, $2$, $18$, $6$, $-20$, $12$, $-22$, $-16)$ | $(22$, $32$, $38$, $34$, $20$, $24$, $30$, $-64$, $-76$, $-70$, $-42$, $-44$, $-68$, $-78$, $-66$, $-46$, $-40$, $-72$, $-74$, $60$, $84$, $98$, $82$, $62$, $16$, $28$, $26$, $18$, $36$, $-6$, $-4$, $90$, $92$, $80$, $96$, $86$, $58$, $88$, $94$, $-52$, $-12$, $-2$, $-8$, $-56$, $-48$, $-14$, $-50$, $-54$, $-10)$\\

\noindent$11n5$ : $3$ | $(4$, $8$, $10$, $14$, $2$, $-18$, $6$, $20$, $-12$, $22$, $16)$ | $(28$, $68$, $30$, $26$, $44$, $22$, $34$, $72$, $116$, $-54$, $-60$, $-88$, $-84$, $-64$, $-66$, $-82$, $-90$, $-58$, $-56$, $-16$, $-52$, $-62$, $-86$, $24$, $32$, $70$, $114$, $118$, $74$, $120$, $112$, $106$, $126$, $104$, $-80$, $-92$, $-98$, $36$, $20$, $42$, $102$, $124$, $108$, $110$, $122$, $100$, $40$, $18$, $38$, $-14$, $-50$, $-4$, $-2$, $-8$, $-46$, $-10$, $-78$, $-94$, $-96$, $-76$, $-12$, $-48$, $-6)$\\

\noindent$11n6$ : $3$ | $(4$, $8$, $10$, $14$, $2$, $-18$, $6$, $-20$, $-12$, $-22$, $-16)$ | $(40$, $38$, $26$, $112$, $68$, $86$, $70$, $110$, $28$, $36$, $42$, $78$, $-60$, $-92$, $-106$, $-12$, $-10$, $-104$, $-94$, $-62$, $-96$, $-102$, $-8$, $-14$, $-108$, $48$, $30$, $34$, $44$, $76$, $80$, $-2$, $-20$, $-56$, $-88$, $-52$, $-16$, $-6$, $-100$, $-98$, $-4$, $-18$, $-54$, $32$, $46$, $74$, $82$, $64$, $116$, $24$, $114$, $66$, $84$, $72$, $-50$, $-90$, $-58$, $-22)$\\

\noindent$11n7$ : $3$ | $(4$, $8$, $10$, $-14$, $2$, $18$, $-6$, $20$, $12$, $22$, $16)$ | $(64$, $40$, $56$, $38$, $66$, $94$, $68$, $36$, $58$, $42$, $62$, $32$, $-80$, $-12$, $-10$, $-82$, $-78$, $-70$, $-90$, $-54$, $-86$, $-74$, $60$, $34$, $26$, $30$, $104$, $114$, $-88$, $-72$, $-76$, $-84$, $-52$, $-92$, $100$, $118$, $110$, $108$, $120$, $98$, $96$, $122$, $106$, $112$, $116$, $102$, $28$, $-50$, $-24$, $-48$, $-14$, $-8$, $-2$, $-20$, $-44$, $-18$, $-4$, $-6$, $-16$, $-46$, $-22)$\\

\noindent$11n8$ : $3$ | $(4$, $8$, $10$, $-14$, $2$, $-18$, $-6$, $-20$, $-12$, $-22$, $-16)$ | $(76$, $72$, $38$, $32$, $28$, $42$, $44$, $26$, $34$, $36$, $74$, $-2$, $-58$, $-100$, $-52$, $-50$, $-98$, $-60$, $-96$, $-48$, $-54$, $-102$, $-56$, $106$, $64$, $82$, $66$, $104$, $70$, $78$, $-18$, $-6$, $-84$, $-10$, $-14$, $-88$, $-94$, $-92$, $-90$, $-16$, $-8$, $30$, $40$, $46$, $24$, $110$, $22$, $108$, $62$, $80$, $68$, $-12$, $-86$, $-4$, $-20)$\\

\noindent$11n9$ : $3$ | $(4$, $8$, $10$, $14$, $2$, $-18$, $6$, $-20$, $-22$, $-12$, $-16)$ | $(26$, $44$, $42$, $28$, $24$, $46$, $22$, $30$, $40$, $104$, $-78$, $-48$, $-74$, $-72$, $-50$, $-80$, $-84$, $-88$, $-86$, $-82$, $-52$, $-70$, $-76$, $100$, $108$, $94$, $68$, $20$, $32$, $38$, $34$, $36$, $-60$, $-18$, $-4$, $96$, $110$, $98$, $102$, $106$, $92$, $66$, $90$, $64$, $-62$, $-58$, $-16$, $-6$, $-2$, $-10$, $-12$, $-54$, $-56$, $-14$, $-8)$\\

\noindent$11n10$ : $3$ | $(4$, $8$, $10$, $-14$, $2$, $18$, $-6$, $20$, $22$, $12$, $16)$ | $(28$, $38$, $30$, $92$, $82$, $90$, $32$, $36$, $86$, $-18$, $-44$, $-10$, $-8$, $-42$, $-40$, $-78$, $-68$, $-70$, $-76$, $52$, $94$, $84$, $88$, $34$, $-14$, $-4$, $-2$, $-16$, $-46$, $-12$, $-6$, $26$, $22$, $58$, $48$, $56$, $20$, $54$, $50$, $60$, $24$, $-64$, $-74$, $-72$, $-66$, $-80$, $-62)$\\

\noindent$11n11$ : $3$ | $(4$, $8$, $10$, $-14$, $2$, $-18$, $-6$, $-20$, $-22$, $-12$, $-16)$ | $(20$, $68$, $26$, $66$, $22$, $18$, $70$, $74$, $-42$, $-32$, $-30$, $-40$, $-36$, $48$, $54$, $78$, $56$, $46$, $50$, $52$, $76$, $72$, $-14$, $-4$, $-6$, $-64$, $-8$, $-2$, $-12$, $44$, $16$, $24$, $-38$, $-28$, $-34$, $-58$, $-60$, $-62$, $-10)$\\

\noindent$11n12$ : $3$ | $(4$, $8$, $10$, $-14$, $2$, $-18$, $-6$, $-22$, $-20$, $-12$, $-16)$ | $(10$, $20$, $72$, $-24$, $-36$, $-60$, $-46$, $-62$, $-48$, $-58$, $-38$, $70$, $74$, $18$, $12$, $14$, $16$, $76$, $68$, $-4$, $-26$, $-34$, $-30$, $64$, $42$, $8$, $22$, $6$, $40$, $66$, $44$, $-32$, $-28$, $-2$, $-52$, $-54$, $-56$, $-50)$\\

\noindent$11n13$ : $3$ | $(4$, $8$, $10$, $16$, $2$, $-18$, $-20$, $6$, $-22$, $-12$, $-14)$ | $(30$, $28$, $12$, $14$, $-38$, $-34$, $-20$, $-36$, $10$, $40$, $-6$, $-16$, $-8$, $-18$, $-32$, $44$, $42$, $22$, $24$, $26$, $-4$, $-2)$\\

\noindent$11n14$ : $3$ | $(4$, $8$, $10$, $-16$, $2$, $18$, $20$, $-6$, $22$, $12$, $14)$ | $(38$, $62$, $70$, $-26$, $-22$, $-30$, $-2$, $-60$, $-4$, $-28$, $66$, $74$, $72$, $68$, $64$, $76$, $16$, $14$, $78$, $-36$, $-32$, $-24$, $8$, $20$, $10$, $42$, $6$, $18$, $12$, $40$, $-34$, $-56$, $-48$, $-46$, $-54$, $-52$, $-44$, $-50$, $-58)$\\

\noindent$11n15$ : $3$ | $(4$, $8$, $10$, $-16$, $2$, $-18$, $-20$, $-6$, $-22$, $-12$, $-14)$ | $(14$, $18$, $78$, $-26$, $-30$, $-42$, $-36$, $-38$, $-40$, $-34$, $-44$, $-32$, $76$, $80$, $82$, $74$, $46$, $70$, $86$, $68$, $84$, $72$, $-62$, $-4$, $-28$, $8$, $24$, $6$, $48$, $10$, $22$, $20$, $12$, $16$, $-2$, $-64$, $-60$, $-52$, $-54$, $-56$, $-50$, $-58$, $-66)$\\

\noindent$11n16$ : $3$ | $(4$, $8$, $10$, $16$, $2$, $-18$, $-20$, $6$, $-22$, $-14$, $-12)$ | $(34$, $62$, $66$, $52$, $56$, $70$, $58$, $50$, $64$, $-4$, $-82$, $-12$, $-10$, $-84$, $-6$, $-16$, $-36$, $94$, $92$, $20$, $28$, $86$, $26$, $22$, $90$, $-74$, $-42$, $-44$, $-72$, $-48$, $-38$, $-78$, $-76$, $-40$, $-46$, $24$, $88$, $30$, $18$, $32$, $60$, $68$, $54$, $-8$, $-14$, $-80$, $-2)$\\

\noindent$11n17$ : $3$ | $(4$, $8$, $10$, $-16$, $2$, $18$, $20$, $-6$, $22$, $14$, $12)$ | $(66$, $56$, $74$, $102$, $72$, $54$, $68$, $64$, $58$, $76$, $60$, $62$, $-50$, $-8$, $-6$, $-20$, $-18$, $-4$, $-10$, $-52$, $-12$, $-2$, $-16$, $-22$, $104$, $70$, $-84$, $-98$, $-92$, $-78$, $-48$, $-80$, $-94$, $-96$, $-82$, $-86$, $-100$, $-90$, $32$, $46$, $106$, $40$, $38$, $26$, $30$, $34$, $44$, $108$, $42$, $36$, $28$, $-88$, $-24$, $-14)$\\

\noindent$11n18$ : $3$ | $(4$, $8$, $10$, $-16$, $2$, $-18$, $-20$, $-6$, $-22$, $-14$, $-12)$ | $(18$, $26$, $10$, $28$, $-30$, $-34$, $-44$, $-50$, $-40$, $-38$, $-48$, $-46$, $-36$, $-32$, $82$, $84$, $80$, $86$, $56$, $92$, $74$, $78$, $88$, $58$, $90$, $76$, $-42$, $-2$, $-68$, $22$, $14$, $52$, $16$, $20$, $24$, $12$, $54$, $-64$, $-62$, $-72$, $-6$, $-8$, $-4$, $-70$, $-60$, $-66)$\\

\noindent$11n19$ : $3$ | $(4$, $8$, $10$, $-16$, $2$, $-18$, $-20$, $-22$, $-6$, $-12$, $-14)$ | $(22$, $26$, $38$, $36$, $-6$, $-4$, $-28$, $-2$, $24$, $12$, $14$, $-30$, $-32$, $-20$, $40$, $16$, $10$, $8$, $-34$, $-18)$\\

\noindent$11n20$ : $3$ | $(4$, $8$, $10$, $-16$, $2$, $-18$, $-20$, $-22$, $-6$, $-14$, $-12)$ | $(46$, $48$, $20$, $18$, $50$, $44$, $42$, $52$, $-28$, $-38$, $-56$, $-12$, $62$, $68$, $66$, $64$, $22$, $60$, $70$, $-6$, $-32$, $-24$, $-34$, $-58$, $-36$, $-26$, $-30$, $40$, $72$, $-2$, $-10$, $-14$, $-54$, $-16$, $-8$, $-4)$\\

\noindent$11n21$ : $3$ | $(4$, $8$, $12$, $2$, $-14$, $-18$, $6$, $-20$, $-10$, $-22$, $-16)$ | $(60$, $38$, $22$, $30$, $32$, $20$, $36$, $58$, $34$, $-54$, $-66$, $-48$, $-42$, $-44$, $-50$, $-64$, $-52$, $-56$, $-68$, $-46$, $94$, $80$, $96$, $92$, $82$, $62$, $86$, $88$, $100$, $-70$, $-72$, $-8$, $84$, $90$, $98$, $102$, $28$, $24$, $40$, $26$, $-74$, $-6$, $-10$, $-18$, $-12$, $-4$, $-76$, $-78$, $-2$, $-14$, $-16)$\\

\noindent$11n22$ : $3$ | $(4$, $8$, $-12$, $2$, $14$, $18$, $-6$, $20$, $10$, $22$, $16)$ | $(40$, $24$, $34$, $76$, $32$, $26$, $42$, $28$, $30$, $96$, $64$, $-86$, $-80$, $-44$, $-56$, $-46$, $-78$, $-88$, $-20$, $-84$, $-82$, $70$, $58$, $74$, $36$, $22$, $38$, $72$, $-16$, $-10$, $-4$, $-50$, $-52$, $-2$, $-12$, $-14$, $-18$, $-90$, $92$, $60$, $68$, $98$, $66$, $62$, $94$, $-8$, $-6$, $-48$, $-54)$\\

\noindent$11n23$ : $3$ | $(4$, $8$, $-12$, $2$, $-14$, $-18$, $-6$, $-20$, $-10$, $-22$, $-16)$ | $(72$, $38$, $86$, $40$, $74$, $102$, $118$, $116$, $104$, $130$, $128$, $106$, $114$, $120$, $-60$, $-10$, $-18$, $-22$, $-94$, $-88$, $-98$, $-26$, $-14$, $80$, $32$, $70$, $36$, $84$, $42$, $76$, $78$, $44$, $82$, $34$, $-20$, $-92$, $-90$, $-100$, $-28$, $-12$, $-16$, $-24$, $-96$, $110$, $124$, $132$, $126$, $108$, $112$, $122$, $30$, $-48$, $-66$, $-54$, $-4$, $-6$, $-56$, $-64$, $-46$, $-62$, $-58$, $-8$, $-2$, $-52$, $-68$, $-50)$\\

\noindent$11n24$ : $3$ | $(4$, $8$, $12$, $2$, $14$, $-18$, $6$, $-20$, $-22$, $-10$, $-16)$ | $(28$, $62$, $70$, $72$, $48$, $74$, $68$, $64$, $66$, $-14$, $-4$, $-8$, $-10$, $-38$, $-40$, $-12$, $-6$, $22$, $50$, $76$, $52$, $20$, $32$, $24$, $-60$, $-2$, $-16$, $18$, $30$, $26$, $-36$, $-42$, $-54$, $-56$, $-44$, $-34$, $-46$, $-58)$\\

\noindent$11n25$ : $3$ | $(4$, $8$, $12$, $2$, $-14$, $-18$, $6$, $-20$, $-22$, $-10$, $-16)$ | $(46$, $104$, $88$, $106$, $112$, $94$, $98$, $100$, $92$, $110$, $108$, $90$, $102$, $96$, $-12$, $-16$, $-24$, $-4$, $-6$, $-22$, $-18$, $-10$, $-82$, $-84$, $-8$, $-20$, $40$, $32$, $68$, $28$, $44$, $48$, $-2$, $-26$, $-14$, $30$, $42$, $50$, $38$, $34$, $66$, $114$, $64$, $36$, $-78$, $-54$, $-72$, $-60$, $-58$, $-70$, $-56$, $-80$, $-86$, $-76$, $-52$, $-74$, $-62)$\\

\noindent$11n26$ : $3$ | $(4$, $8$, $-12$, $2$, $14$, $18$, $-6$, $20$, $22$, $10$, $16)$ | $(26$, $82$, $72$, $84$, $28$, $24$, $-64$, $-54$, $-4$, $-6$, $-56$, $-62$, $-66$, $-68$, $-44$, $46$, $22$, $14$, $52$, $16$, $20$, $48$, $-60$, $-58$, $-8$, $-2$, $80$, $74$, $86$, $30$, $88$, $76$, $78$, $90$, $50$, $18$, $-40$, $-34$, $-12$, $-36$, $-38$, $-70$, $-42$, $-32$, $-10)$\\

\noindent$11n27$ : $3$ | $(4$, $8$, $-12$, $2$, $-14$, $-18$, $-6$, $-20$, $-22$, $-10$, $-16)$ | $(40$, $42$, $12$, $-56$, $-64$, $-54$, $-24$, $-26$, $-36$, $-28$, $-2$, $74$, $48$, $38$, $44$, $10$, $8$, $46$, $-60$, $-50$, $-22$, $-52$, $-62$, $-58$, $14$, $20$, $72$, $66$, $68$, $16$, $18$, $70$, $-32$, $-6$, $-34$, $-30$, $-4)$\\

\noindent$11n28$ : $3$ | $(4$, $8$, $12$, $2$, $-14$, $-18$, $6$, $-22$, $-20$, $-10$, $-16)$ | $(12$, $30$, $44$, $28$, $14$, $-38$, $-40$, $-36$, $-32$, $-4$, $-6$, $18$, $16$, $-34$, $-42$, $-24$, $46$, $48$, $52$, $54$, $50$, $20$, $-22$, $-26$, $-8$, $-2$, $-10)$\\

\noindent$11n29$ : $3$ | $(4$, $8$, $-12$, $2$, $14$, $18$, $-6$, $22$, $20$, $10$, $16)$ | $(64$, $58$, $70$, $52$, $74$, $54$, $68$, $60$, $62$, $-44$, $-14$, $-4$, $-48$, $-10$, $-8$, $-50$, $-6$, $-12$, $-46$, $-2$, $-16$, $98$, $66$, $56$, $72$, $-90$, $-76$, $-94$, $-86$, $-80$, $-42$, $-82$, $-84$, $-96$, $-78$, $-88$, $-92$, $28$, $34$, $22$, $40$, $100$, $38$, $24$, $32$, $30$, $26$, $36$, $20$, $-18)$\\

\noindent$11n30$ : $3$ | $(4$, $8$, $-12$, $2$, $-14$, $-18$, $-6$, $-22$, $-20$, $-10$, $-16)$ | $(28$, $46$, $50$, $32$, $42$, $54$, $36$, $38$, $56$, $40$, $34$, $52$, $44$, $-92$, $-72$, $-66$, $-86$, $-60$, $-80$, $-58$, $-84$, $-64$, $-90$, $-70$, $-68$, $-88$, $-62$, $-82$, $106$, $78$, $110$, $102$, $74$, $114$, $98$, $116$, $-94$, $-8$, $-20$, $-14$, $108$, $104$, $76$, $112$, $100$, $118$, $120$, $30$, $48$, $-4$, $-24$, $-10$, $-18$, $-16$, $-12$, $-22$, $-6$, $-96$, $-2$, $-26)$\\

\noindent$11n31$ : $3$ | $(4$, $8$, $12$, $2$, $-16$, $6$, $-20$, $-10$, $-22$, $-14$, $-18)$ | $(70$, $60$, $98$, $96$, $28$, $30$, $26$, $94$, $100$, $88$, $-40$, $-4$, $-48$, $-32$, $-50$, $52$, $56$, $74$, $66$, $64$, $62$, $68$, $72$, $58$, $54$, $-10$, $-12$, $-8$, $-14$, $-78$, $-84$, $-20$, $-86$, $-18$, $-82$, $-80$, $-16$, $-6$, $24$, $92$, $102$, $90$, $22$, $-38$, $-42$, $-2$, $-46$, $-34$, $-76$, $-36$, $-44)$\\

\noindent$11n32$ : $3$ | $(4$, $8$, $-12$, $2$, $16$, $-6$, $20$, $10$, $22$, $14$, $18)$ | $(40$, $78$, $148$, $124$, $132$, $140$, $84$, $142$, $130$, $126$, $146$, $80$, $136$, $-70$, $-10$, $-16$, $-64$, $-20$, $-6$, $-74$, $-76$, $-4$, $-22$, $-66$, $-14$, $-12$, $-68$, $-24$, $-72$, $-8$, $-18$, $128$, $144$, $82$, $138$, $134$, $122$, $150$, $120$, $-94$, $-100$, $-110$, $50$, $30$, $60$, $36$, $44$, $56$, $26$, $54$, $46$, $34$, $62$, $32$, $48$, $52$, $28$, $58$, $38$, $42$, $-2$, $-116$, $-90$, $-104$, $-106$, $-88$, $-114$, $-96$, $-98$, $-112$, $-86$, $-108$, $-102$, $-92$, $-118)$\\

\noindent$11n33$ : $3$ | $(4$, $8$, $-12$, $2$, $-16$, $-6$, $-20$, $-10$, $-22$, $-14$, $-18)$ | $(42$, $48$, $20$, $50$, $40$, $44$, $46$, $38$, $-54$, $-56$, $74$, $66$, $80$, $64$, $76$, $16$, $52$, $18$, $-30$, $-60$, $-24$, $-26$, $-62$, $-28$, $-22$, $-58$, $70$, $36$, $72$, $68$, $78$, $-12$, $-2$, $-8$, $-32$, $-34$, $-6$, $-4$, $-14$, $-10)$\\

\noindent$11n34$ : $3$ | $(4$, $8$, $12$, $2$, $-16$, $-18$, $6$, $-20$, $-22$, $-14$, $-10)$ | $(30$, $20$, $96$, $22$, $32$, $28$, $102$, $-56$, $-40$, $-44$, $-94$, $-48$, $-36$, $-52$, $-90$, $-92$, $-50$, $112$, $68$, $60$, $104$, $106$, $62$, $66$, $110$, $114$, $70$, $58$, $-14$, $-86$, $-76$, $-4$, $-6$, $-78$, $-84$, $-12$, $-42$, $18$, $98$, $24$, $34$, $26$, $100$, $16$, $72$, $116$, $108$, $64$, $-46$, $-38$, $-54$, $-88$, $-74$, $-2$, $-8$, $-80$, $-82$, $-10)$\\

\noindent$11n35$ : $3$ | $(4$, $8$, $-12$, $2$, $16$, $18$, $-6$, $20$, $22$, $14$, $10)$ | $(30$, $26$, $42$, $90$, $40$, $24$, $32$, $82$, $86$, $36$, $-58$, $-48$, $-94$, $-100$, $-98$, $-96$, $-50$, $-56$, $-60$, $-46$, $-92$, $116$, $120$, $76$, $70$, $126$, $112$, $124$, $72$, $74$, $122$, $114$, $28$, $-2$, $-12$, $-108$, $-20$, $-106$, $-10$, $-4$, $-64$, $-52$, $-54$, $-62$, $-44$, $118$, $78$, $68$, $128$, $66$, $80$, $88$, $38$, $22$, $34$, $84$, $-16$, $-102$, $-6$, $-8$, $-104$, $-18$, $-110$, $-14)$\\

\noindent$11n36$ : $3$ | $(4$, $8$, $-12$, $2$, $-16$, $-18$, $-6$, $-20$, $-22$, $-14$, $-10)$ | $(28$, $94$, $22$, $34$, $20$, $96$, $98$, $-84$, $-86$, $-72$, $-36$, $-76$, $-90$, $-80$, $-82$, $-88$, $-74$, $106$, $60$, $68$, $14$, $18$, $32$, $24$, $92$, $26$, $30$, $16$, $-42$, $-10$, $-4$, $-48$, $-50$, $-2$, $-12$, $58$, $104$, $108$, $62$, $66$, $112$, $100$, $56$, $102$, $110$, $64$, $-78$, $-38$, $-70$, $-40$, $-54$, $-44$, $-8$, $-6$, $-46$, $-52)$\\

\noindent$11n37$ : $3$ | $(4$, $8$, $-12$, $2$, $-18$, $-6$, $-20$, $-22$, $-10$, $-16$, $-14)$ | $(52$, $24$, $16$, $18$, $26$, $50$, $42$, $-56$, $-38$, $-58$, $-12$, $-54$, $-36$, $60$, $20$, $14$, $22$, $62$, $46$, $-6$, $-30$, $-32$, $-4$, $-8$, $-28$, $-34$, $64$, $44$, $40$, $48$, $-10$, $-2)$\\

\noindent$11n38$ : $3$ | $(4$, $8$, $-12$, $2$, $-18$, $-6$, $-22$, $-20$, $-10$, $-16$, $-14)$ | $(12$, $24$, $38$, $26$, $14$, $-32$, $-34$, $-22$, $-4$, $-6$, $40$, $-30$, $-36$, $-20$, $16$, $44$, $46$, $42$, $18$, $-28$, $-8$, $-2$, $-10)$\\

\noindent$11n39$ : $3$ | $(4$, $8$, $12$, $2$, $-18$, $-16$, $6$, $-20$, $-22$, $-14$, $-10)$ | $(36$, $50$, $128$, $134$, $44$, $42$, $136$, $126$, $52$, $34$, $70$, $92$, $78$, $-64$, $-6$, $-14$, $-56$, $-104$, $-106$, $-22$, $-116$, $-96$, $-114$, $-24$, $-108$, $-102$, $-122$, $124$, $138$, $40$, $46$, $132$, $130$, $48$, $38$, $-58$, $-12$, $-8$, $-62$, $-66$, $-4$, $-16$, $-54$, $-18$, $-2$, $-68$, $-60$, $-10$, $74$, $30$, $82$, $88$, $140$, $90$, $80$, $28$, $76$, $94$, $72$, $32$, $84$, $86$, $-120$, $-100$, $-110$, $-26$, $-112$, $-98$, $-118$, $-20)$\\

\noindent$11n40$ : $3$ | $(4$, $8$, $-12$, $2$, $18$, $16$, $-6$, $20$, $22$, $14$, $10)$ | $(40$, $22$, $76$, $26$, $36$, $86$, $84$, $34$, $28$, $78$, $-56$, $-108$, $-104$, $-52$, $-94$, $-116$, $-96$, $-50$, $-102$, $-110$, $-58$, $-60$, $-68$, $-12$, $142$, $120$, $134$, $128$, $126$, $90$, $80$, $30$, $32$, $82$, $88$, $38$, $24$, $-106$, $-54$, $-92$, $-114$, $-98$, $-48$, $-100$, $-112$, $42$, $136$, $118$, $140$, $46$, $144$, $122$, $132$, $130$, $124$, $146$, $44$, $138$, $-64$, $-16$, $-8$, $-72$, $-2$, $-20$, $-4$, $-74$, $-6$, $-18$, $-62$, $-66$, $-14$, $-10$, $-70)$\\

\noindent$11n41$ : $3$ | $(4$, $8$, $-12$, $2$, $-18$, $-16$, $-6$, $-20$, $-22$, $-14$, $-10)$ | $(98$, $54$, $90$, $84$, $48$, $104$, $46$, $82$, $92$, $56$, $96$, $138$, $140$, $170$, $162$, $148$, $-70$, $-6$, $-16$, $-60$, $-22$, $-26$, $-128$, $-106$, $-124$, $-30$, $-118$, $-112$, $-134$, $94$, $80$, $44$, $102$, $50$, $86$, $88$, $52$, $100$, $42$, $-24$, $-130$, $-108$, $-122$, $-32$, $-120$, $-110$, $-132$, $-136$, $-114$, $-116$, $-28$, $-126$, $166$, $144$, $36$, $152$, $158$, $172$, $160$, $150$, $34$, $146$, $164$, $168$, $142$, $38$, $154$, $156$, $40$, $-78$, $-62$, $-14$, $-8$, $-68$, $-72$, $-4$, $-18$, $-58$, $-20$, $-2$, $-74$, $-66$, $-10$, $-12$, $-64$, $-76)$\\

\noindent$11n42$ : $3$ | $(4$, $8$, $12$, $2$, $-18$, $-20$, $6$, $-10$, $-22$, $-14$, $-16)$ | $(38$, $58$, $114$, $52$, $44$, $106$, $108$, $46$, $50$, $112$, $60$, $36$, $40$, $56$, $116$, $54$, $42$, $-122$, $-74$, $-76$, $-120$, $-84$, $-66$, $-130$, $-128$, $-68$, $-82$, $-118$, $-78$, $-72$, $-124$, $-12$, $150$, $94$, $144$, $136$, $102$, $156$, $100$, $138$, $142$, $96$, $152$, $62$, $110$, $48$, $-16$, $-8$, $-32$, $-26$, $-2$, $-22$, $-86$, $-64$, $-132$, $-126$, $-70$, $-80$, $140$, $98$, $154$, $104$, $134$, $146$, $92$, $148$, $-88$, $-20$, $-4$, $-28$, $-30$, $-6$, $-18$, $-90$, $-14$, $-10$, $-34$, $-24)$\\

\noindent$11n43$ : $3$ | $(4$, $8$, $-12$, $2$, $18$, $20$, $-6$, $10$, $22$, $14$, $16)$ | $(24$, $94$, $34$, $96$, $22$, $26$, $92$, $32$, $-82$, $-86$, $-74$, $-38$, $-40$, $-76$, $-88$, $-80$, $-44$, $60$, $72$, $108$, $70$, $62$, $58$, $20$, $28$, $90$, $30$, $18$, $-84$, $-6$, $-4$, $-14$, $-52$, $-50$, $-12$, $-2$, $-8$, $106$, $68$, $64$, $102$, $100$, $98$, $104$, $66$, $-78$, $-42$, $-36$, $-46$, $-56$, $-16$, $-54$, $-48$, $-10)$\\

\noindent$11n44$ : $3$ | $(4$, $8$, $-12$, $2$, $-18$, $-20$, $-6$, $-10$, $-22$, $-14$, $-16)$ | $(80$, $74$, $44$, $66$, $86$, $68$, $46$, $72$, $82$, $114$, $116$, $130$, $-58$, $-2$, $-16$, $-50$, $-10$, $-8$, $-52$, $-18$, $-22$, $-100$, $-90$, $-108$, $70$, $84$, $64$, $42$, $76$, $78$, $40$, $-20$, $-102$, $-88$, $-106$, $-110$, $-92$, $-98$, $-24$, $-96$, $-94$, $-112$, $-104$, $34$, $120$, $126$, $28$, $132$, $26$, $128$, $118$, $36$, $32$, $122$, $124$, $30$, $38$, $-62$, $-54$, $-6$, $-12$, $-48$, $-14$, $-4$, $-56$, $-60)$\\

\noindent$11n45$ : $3$ | $(4$, $8$, $12$, $2$, $-20$, $-18$, $6$, $-10$, $-22$, $-14$, $-16)$ | $(76$, $30$, $38$, $84$, $40$, $28$, $74$, $120$, $108$, $64$, $110$, $118$, $72$, $-102$, $-88$, $-94$, $-18$, $-20$, $-96$, $-86$, $-100$, $-24$, $-14$, $78$, $32$, $36$, $82$, $42$, $26$, $44$, $80$, $34$, $-50$, $-4$, $-10$, $-56$, $-58$, $-90$, $-92$, $-16$, $-22$, $-98$, $114$, $68$, $104$, $122$, $106$, $66$, $112$, $116$, $70$, $-12$, $-2$, $-48$, $-62$, $-52$, $-6$, $-8$, $-54$, $-60$, $-46)$\\

\noindent$11n46$ : $3$ | $(4$, $8$, $-12$, $2$, $20$, $18$, $-6$, $10$, $22$, $14$, $16)$ | $(90$, $44$, $54$, $100$, $56$, $42$, $40$, $58$, $98$, $52$, $46$, $92$, $88$, $84$, $96$, $50$, $48$, $94$, $86$, $-26$, $-118$, $-104$, $-110$, $-64$, $-74$, $-62$, $-112$, $-102$, $-116$, $136$, $128$, $76$, $124$, $140$, $38$, $138$, $126$, $-32$, $-20$, $-4$, $-10$, $-14$, $-70$, $-68$, $-106$, $-108$, $-66$, $-72$, $-60$, $-114$, $132$, $80$, $120$, $142$, $122$, $78$, $130$, $134$, $82$, $-12$, $-2$, $-22$, $-34$, $-30$, $-18$, $-6$, $-8$, $-16$, $-28$, $-36$, $-24)$\\

\noindent$11n47$ : $3$ | $(4$, $8$, $-12$, $2$, $-20$, $-18$, $-6$, $-10$, $-22$, $-14$, $-16)$ | $(70$, $32$, $24$, $38$, $46$, $48$, $40$, $26$, $30$, $72$, $112$, $-66$, $-88$, $-56$, $-54$, $-90$, $-68$, $-20$, $-64$, $-86$, $-58$, $-52$, $-62$, $-84$, $-60$, $80$, $120$, $106$, $118$, $78$, $82$, $74$, $114$, $110$, $124$, $-92$, $-98$, $-8$, $-12$, $-102$, $-100$, $-10$, $76$, $116$, $108$, $122$, $34$, $22$, $36$, $44$, $50$, $42$, $28$, $-16$, $-4$, $-94$, $-96$, $-6$, $-14$, $-104$, $-18$, $-2)$\\

\noindent$11n48$ : $3$ | $(4$, $8$, $-12$, $-16$, $2$, $-18$, $-20$, $-22$, $-10$, $-6$, $-14)$ | $(36$, $16$, $18$, $34$, $30$, $32$, $-8$, $-20$, $-26$, $28$, $50$, $52$, $48$, $-4$, $-40$, $-44$, $-42$, $-24$, $-22$, $14$, $46$, $12$, $-10$, $-6$, $-2$, $-38)$ \\

\noindent$11n49$ : $3$ | $(4$, $8$, $-12$, $-16$, $2$, $-18$, $-22$, $-20$, $-10$, $-6$, $-14)$ | $(14$, $16$, $12$, $46$, $-42$, $-38$, $-26$, $-28$, $-40$, $10$, $44$, $8$, $30$, $52$, $-2$, $-6$, $-20$, $-22$, $32$, $50$, $48$, $34$, $-36$, $-24$, $-18$, $-4)$\\

\noindent$11n50$ : $3$ | $(4$, $8$, $-12$, $-18$, $2$, $-16$, $-20$, $-6$, $-10$, $-22$, $-14)$ | $(16$, $52$, $54$, $56$, $40$, $38$, $58$, $-26$, $-28$, $-12$, $-8$, $-4$, $34$, $60$, $36$, $22$, $-6$, $-2$, $-42$, $-48$, $50$, $18$, $14$, $20$, $-10$, $-30$, $-24$, $-32$, $-46$, $-44)$\\

\noindent$11n51$ : $3$ | $(4$, $8$, $14$, $2$, $-16$, $-18$, $6$, $-22$, $-20$, $-10$, $-12)$ | $(50$, $40$, $76$, $72$, $44$, $46$, $52$, $48$, $42$, $74$, $-6$, $-8$, $-18$, $-4$, $-38$, $-2$, $-16$, $-10$, $78$, $-68$, $-62$, $-56$, $-36$, $-54$, $-64$, $-66$, $80$, $34$, $24$, $20$, $26$, $32$, $82$, $30$, $28$, $22$, $-58$, $-60$, $-70$, $-12$, $-14)$\\

\noindent$11n52$ : $3$ | $(4$, $8$, $-14$, $2$, $16$, $18$, $-6$, $22$, $20$, $10$, $12)$ | $(66$, $88$, $94$, $26$, $98$, $100$, $28$, $92$, $90$, $30$, $102$, $96$, $-46$, $-58$, $-40$, $-52$, $-2$, $-4$, $-16$, $80$, $70$, $62$, $74$, $76$, $64$, $68$, $82$, $36$, $78$, $72$, $-10$, $-22$, $-86$, $-84$, $-20$, $-12$, $-8$, $-24$, $-6$, $-14$, $-18$, $32$, $104$, $34$, $-54$, $-38$, $-56$, $-48$, $-44$, $-60$, $-42$, $-50)$\\

\noindent$11n53$ : $3$ | $(4$, $8$, $-14$, $2$, $-16$, $-18$, $-6$, $-22$, $-20$, $-10$, $-12)$ | $(68$, $78$, $88$, $84$, $82$, $90$, $76$, $66$, $70$, $80$, $86$, $-6$, $-10$, $-18$, $-2$, $-14$, $-74$, $32$, $92$, $26$, $58$, $22$, $60$, $28$, $94$, $30$, $62$, $24$, $-8$, $-20$, $-4$, $-12$, $-36$, $-50$, $-48$, $-38$, $64$, $-34$, $-52$, $-46$, $-40$, $-56$, $-42$, $-44$, $-54$, $-72$, $-16)$\\

\noindent$11n54$ : $3$ | $(4$, $8$, $14$, $2$, $-16$, $-18$, $6$, $-22$, $-20$, $-12$, $-10)$ | $(24$, $8$, $40$, $-42$, $-102$, $-120$, $-92$, $-112$, $-110$, $-94$, $-122$, $-100$, $-98$, $-124$, $-96$, $-108$, $-114$, $-90$, $-118$, $-104$, $138$, $84$, $30$, $18$, $14$, $34$, $88$, $36$, $12$, $20$, $28$, $26$, $22$, $10$, $38$, $86$, $32$, $16$, $-48$, $-58$, $-66$, $-2$, $-106$, $-116$, $132$, $144$, $78$, $150$, $126$, $154$, $82$, $140$, $136$, $128$, $148$, $76$, $146$, $130$, $134$, $142$, $80$, $152$, $-62$, $-44$, $-72$, $-52$, $-54$, $-70$, $-6$, $-4$, $-68$, $-56$, $-50$, $-74$, $-46$, $-60$, $-64)$\\

\noindent$11n55$ : $3$ | $(4$, $8$, $-14$, $2$, $16$, $18$, $-6$, $22$, $20$, $12$, $10)$ | $(20$, $58$, $-34$, $-44$, $-48$, $-30$, $-52$, $-40$, $-38$, $-54$, $-56$, $-36$, $-42$, $-50$, $112$, $66$, $60$, $106$, $118$, $70$, $116$, $108$, $62$, $64$, $110$, $114$, $68$, $120$, $104$, $-86$, $-76$, $-98$, $-78$, $-84$, $-92$, $16$, $24$, $6$, $100$, $10$, $28$, $12$, $102$, $4$, $22$, $18$, $14$, $26$, $8$, $-46$, $-32$, $-2$, $-88$, $-74$, $-96$, $-80$, $-82$, $-94$, $-72$, $-90)$\\

\noindent$11n56$ : $3$ | $(4$, $8$, $-14$, $2$, $-16$, $-18$, $-6$, $-22$, $-20$, $-12$, $-10)$ | $(32$, $84$, $90$, $26$, $94$, $96$, $28$, $88$, $86$, $30$, $98$, $92$, $-52$, $-40$, $-58$, $-46$, $-82$, $-22$, $-10$, $76$, $70$, $62$, $66$, $80$, $34$, $72$, $74$, $36$, $78$, $68$, $-16$, $-4$, $-2$, $-18$, $-14$, $-6$, $-24$, $-8$, $-12$, $-20$, $100$, $64$, $-44$, $-60$, $-42$, $-50$, $-54$, $-38$, $-56$, $-48)$\\

\noindent$11n57$ : $3$ | $(4$, $8$, $-14$, $2$, $-16$, $-18$, $-20$, $-6$, $-10$, $-22$, $-12)$ | $(16$, $24$, $14$, $56$, $46$, $54$, $-42$, $-26$, $-38$, $-8$, $-10$, $-40$, $50$, $36$, $18$, $22$, $12$, $-6$, $58$, $48$, $52$, $34$, $20$, $-30$, $-2$, $-4$, $-32$, $-44$, $-28)$\\

\noindent$11n58$ : $3$ | $(4$, $8$, $14$, $2$, $-16$, $-18$, $-20$, $6$, $-22$, $-12$, $-10)$ | $(22$, $40$, $48$, $42$, $20$, $24$, $38$, $46$, $44$, $-62$, $-52$, $-54$, $-64$, $-14$, $-4$, $-6$, $-16$, $26$, $-66$, $-56$, $-50$, $-60$, $-34$, $-58$, $30$, $72$, $78$, $70$, $28$, $32$, $74$, $76$, $68$, $-36$, $-12$, $-2$, $-8$, $-18$, $-10)$\\

\noindent$11n59$ : $3$ | $(4$, $8$, $-14$, $2$, $16$, $18$, $20$, $-6$, $22$, $12$, $10)$ | $(82$, $102$, $100$, $80$, $92$, $110$, $90$, $144$, $88$, $108$, $94$, $78$, $98$, $104$, $84$, $-70$, $-16$, $-14$, $-8$, $-22$, $-76$, $-24$, $-6$, $-2$, $-28$, $-72$, $-18$, $-12$, $-10$, $-20$, $-74$, $-26$, $-4$, $46$, $146$, $86$, $106$, $96$, $-122$, $-136$, $-68$, $-130$, $-128$, $-116$, $-142$, $-112$, $-138$, $-120$, $-124$, $-134$, $-66$, $-132$, $-126$, $-118$, $-140$, $56$, $36$, $32$, $52$, $60$, $40$, $42$, $62$, $50$, $30$, $48$, $64$, $44$, $38$, $58$, $54$, $34$, $-114)$\\

\noindent$11n60$ : $3$ | $(4$, $8$, $-14$, $2$, $-16$, $-18$, $-20$, $-6$, $-22$, $-12$, $-10)$ | $(30$, $20$, $36$, $38$, $22$, $28$, $72$, $76$, $64$, $-46$, $-82$, $-52$, $-56$, $-54$, $-84$, $-44$, $-16$, $-48$, $-80$, $-50$, $-14$, $74$, $62$, $66$, $78$, $70$, $58$, $100$, $-98$, $-10$, $-2$, $-90$, $-92$, $-4$, $-8$, $-96$, $-86$, $-88$, $-94$, $-6$, $68$, $60$, $102$, $42$, $32$, $18$, $34$, $40$, $24$, $26$, $-12)$\\

\noindent$11n61$ : $3$ | $(4$, $8$, $-14$, $2$, $-16$, $-18$, $-22$, $-6$, $-20$, $-12$, $-10)$ | $(60$, $48$, $64$, $52$, $56$, $82$, $54$, $-8$, $-68$, $-2$, $-46$, $-4$, $-66$, $-6$, $-12$, $-10$, $-70$, $18$, $24$, $26$, $16$, $32$, $62$, $-74$, $-36$, $-40$, $-78$, $-80$, $-42$, $-34$, $-72$, $-76$, $-38$, $50$, $58$, $84$, $20$, $22$, $28$, $14$, $30$, $-44)$\\

\noindent$11n62$ : $3$ | $(4$, $8$, $14$, $2$, $-16$, $-20$, $-18$, $6$, $-22$, $-12$, $-10)$ | $(20$, $54$, $22$, $24$, $56$, $60$, $26$, $58$, $-10$, $-50$, $-52$, $-38$, $-32$, $42$, $16$, $44$, $40$, $18$, $46$, $-6$, $-12$, $-14$, $-8$, $-4$, $-2$, $62$, $48$, $-36$, $-30$, $-28$, $-34)$\\

\noindent$11n63$ : $3$ | $(4$, $8$, $-14$, $2$, $16$, $20$, $18$, $-6$, $22$, $12$, $10)$ | $(32$, $28$, $20$, $18$, $26$, $34$, $50$, $58$, $-74$, $-62$, $-72$, $-14$, $-76$, $-64$, $-70$, $-68$, $-66$, $78$, $24$, $16$, $22$, $30$, $-2$, $-10$, $-36$, $-42$, $-4$, $-8$, $-38$, $-40$, $-6$, $54$, $46$, $80$, $44$, $52$, $60$, $56$, $48$, $-12)$\\

\noindent$11n64$ : $3$ | $(4$, $8$, $-14$, $2$, $-16$, $-20$, $-18$, $-6$, $-22$, $-12$, $-10)$ | $(16$, $24$, $28$, $20$, $12$, $-44$, $-56$, $-48$, $-36$, $-42$, $-54$, $-50$, $-38$, $-40$, $-52$, $22$, $14$, $78$, $70$, $62$, $68$, $76$, $-10$, $-34$, $72$, $64$, $66$, $74$, $46$, $18$, $26$, $-4$, $-30$, $-6$, $-60$, $-2$, $-32$, $-8$, $-58)$\\

\noindent$11n65$ : $3$ | $(4$, $8$, $14$, $2$, $18$, $16$, $-20$, $6$, $10$, $-22$, $-12)$ | $(34$, $100$, $108$, $114$, $94$, $92$, $116$, $106$, $102$, $104$, $-16$, $-4$, $-84$, $-88$, $-8$, $-12$, $-56$, $-58$, $-14$, $-6$, $-86$, $26$, $70$, $68$, $28$, $40$, $24$, $72$, $118$, $74$, $22$, $38$, $30$, $-10$, $-90$, $-82$, $-2$, $-18$, $20$, $36$, $32$, $98$, $110$, $112$, $96$, $-46$, $-66$, $-48$, $-44$, $-54$, $-60$, $-76$, $-78$, $-62$, $-52$, $-42$, $-50$, $-64$, $-80)$\\

\noindent$11n66$ : $3$ | $(4$, $8$, $14$, $2$, $-18$, $-16$, $20$, $6$, $-10$, $22$, $12)$ | $(30$, $62$, $22$, $38$, $36$, $24$, $64$, $28$, $32$, $60$, $-76$, $-86$, $-46$, $-50$, $-82$, $-80$, $-52$, $-44$, $-88$, $-74$, $-72$, $90$, $110$, $100$, $98$, $112$, $92$, $40$, $-42$, $-54$, $-78$, $-84$, $-48$, $26$, $34$, $58$, $56$, $106$, $104$, $94$, $114$, $96$, $102$, $108$, $-8$, $-70$, $-18$, $-2$, $-14$, $-66$, $-12$, $-4$, $-20$, $-6$, $-10$, $-68$, $-16)$\\

\noindent$11n67$ : $3$ | $(4$, $8$, $14$, $2$, $-18$, $-16$, $-20$, $6$, $-10$, $-22$, $-12)$ | $(54$, $46$, $22$, $24$, $44$, $56$, $78$, $32$, $80$, $-12$, $-70$, $-60$, $-72$, $-10$, $-2$, $-38$, $50$, $18$, $52$, $48$, $20$, $26$, $-62$, $-68$, $-14$, $-16$, $-66$, $-64$, $-42$, $-6$, $74$, $28$, $82$, $30$, $76$, $58$, $-8$, $-4$, $-40$, $-34$, $-36)$\\

\noindent$11n68$ : $3$ | $(4$, $8$, $-14$, $2$, $18$, $16$, $20$, $-6$, $10$, $22$, $12)$ | $(20$, $16$, $24$, $42$, $46$, $40$, $22$, $-52$, $-58$, $-56$, $-54$, $-50$, $68$, $62$, $66$, $18$, $-2$, $-14$, $-4$, $-30$, $-48$, $-26$, $-28$, $64$, $70$, $36$, $34$, $72$, $32$, $38$, $44$, $-8$, $-10$, $-60$, $-6$, $-12)$\\

\noindent$11n69$ : $3$ | $(4$, $8$, $-14$, $2$, $-18$, $-16$, $-20$, $-6$, $-10$, $-22$, $-12)$ | $(48$, $80$, $54$, $74$, $56$, $78$, $50$, $46$, $-2$, $-12$, $-60$, $-6$, $-8$, $-58$, $-10$, $-4$, $-62$, $20$, $28$, $26$, $22$, $32$, $-64$, $-66$, $-68$, $-34$, $-42$, $-72$, $-38$, $76$, $52$, $82$, $16$, $84$, $18$, $30$, $24$, $-40$, $-36$, $-70$, $-44$, $-14)$\\

\noindent$11n70$ : $3$ | $(4$, $8$, $-14$, $2$, $-18$, $-20$, $-6$, $-22$, $-12$, $-10$, $-16)$ | $(12$, $56$, $50$, $58$, $62$, $-14$, $64$, $34$, $26$, $30$, $38$, $32$, $-42$, $-4$, $-6$, $-40$, $-44$, $-2$, $-8$, $60$, $52$, $54$, $10$, $36$, $28$, $-18$, $-24$, $-16$, $-48$, $-20$, $-22$, $-46)$\\

\noindent$11n71$ : $4$ | $(4$, $8$, $14$, $2$, $-20$, $16$, $6$, $18$, $12$, $-22$, $-10)$ | $(32$, $6$, $-28$, $-26$, $-24$, $-36$, $38$, $10$, $8$, $-14$, $-12$, $-16$, $18$, $30$, $-4$, $-34$, $40$, $22$, $20$, $-2)$\\

\noindent$11n72$ : $4$ | $(4$, $8$, $-14$, $2$, $20$, $16$, $-6$, $18$, $12$, $22$, $10)$ | $(24$, $28$, $-46$, $-40$, $-22$, $-2$, $26$, $6$, $38$, $36$, $34$, $-32$, $-42$, $-44$, $-20$, $8$, $10$, $-30$, $-48$, $50$, $56$, $12$, $54$, $52$, $-18$, $-16$, $-4$, $-14)$\\

\noindent$11n73$ : $4$ | $(4$, $8$, $-14$, $2$, $20$, $-16$, $-6$, $-18$, $-12$, $22$, $10)$ | $(12$, $26$, $32$, $34$, $24$, $-54$, $-22$, $42$, $30$, $28$, $40$, $-52$, $-56$, $-38$, $-48$, $-36$, $-58$, $60$, $46$, $64$, $-50$, $-8$, $-4$, $-18$, $62$, $14$, $68$, $66$, $44$, $-6$, $-16$, $-20$, $-2$, $-10)$\\

\noindent$11n74$ : $4$ | $(4$, $8$, $-14$, $2$, $-20$, $16$, $-6$, $18$, $12$, $-22$, $-10)$ | $(34$, $40$, $30$, $42$, $-2$, $-22$, $10$, $-24$, $-6$, $-4$, $-12$, $32$, $44$, $8$, $18$, $-38$, $-14$, $-26$, $16$, $20$, $-28$, $-36)$\\

\noindent$11n75$ : $4$ | $(4$, $8$, $-14$, $2$, $-20$, $-16$, $-6$, $-18$, $-12$, $-22$, $-10)$ | $(10$, $26$, $-14$, $-34$, $-32$, $-44$, $24$, $28$, $50$, $30$, $8$, $-46$, $-42$, $-40$, $-48$, $36$, $22$, $-20$, $-18$, $-16$, $4$, $52$, $54$, $6$, $38$, $-2$, $-12)$\\

\noindent$11n76$ : $4$ | $(4$, $8$, $14$, $2$, $-20$, $16$, $18$, $6$, $12$, $-22$, $-10)$ | $(16$, $34$, $8$, $-24$, $-22$, $-20$, $10$, $-30$, $-32$, $-12$, $18$, $26$, $-6$, $-4$, $-2$, $36$, $28$, $-14)$\\

\noindent$11n77$ : $4$ | $(4$, $8$, $-14$, $2$, $20$, $-16$, $-18$, $-6$, $-12$, $22$, $10)$ | $(32$, $22$, $18$, $20$, $-4$, $-38$, $50$, $46$, $10$, $-42$, $-40$, $-16$, $-36$, $34$, $24$, $-12$, $-14$, $48$, $30$, $44$, $8$, $-6$, $-2$, $-28$, $-26)$\\

\noindent$11n78$ : $4$ | $(4$, $8$, $-14$, $2$, $-20$, $-16$, $-18$, $-6$, $-12$, $-22$, $-10)$ | $(24$, $28$, $26$, $-40$, $-46$, $-42$, $-2$, $-4$, $12$, $38$, $-36$, $-18$, $-30$, $-16$, $14$, $50$, $8$, $10$, $-44$, $-32$, $48$, $22$, $20$, $-34$, $-6)$\\

\noindent$11n79$ : $3$ | $(4$, $8$, $-14$, $2$, $-20$, $-18$, $-6$, $-22$, $-12$, $-10$, $-16)$ | $(30$, $56$, $32$, $16$, $36$, $40$, $34$, $-44$, $-52$, $-2$, $-4$, $-12$, $-10$, $38$, $18$, $-54$, $-46$, $-42$, $-50$, $-48$, $24$, $60$, $62$, $22$, $26$, $58$, $64$, $20$, $-28$, $-8$, $-14$, $-6)$\\

\noindent$11n80$ : $3$ | $(4$, $8$, $14$, $10$, $2$, $-18$, $6$, $-20$, $-22$, $-12$, $-16)$ | $(18$, $56$, $58$, $20$, $16$, $24$, $-46$, $-50$, $-32$, $-34$, $-52$, $-44$, $-48$, $14$, $22$, $60$, $70$, $40$, $-4$, $-2$, $-6$, $-30$, $66$, $64$, $62$, $68$, $42$, $38$, $-36$, $-54$, $-10$, $-26$, $-12$, $-28$, $-8)$\\

\noindent$11n81$ : $4$ | $(4$, $8$, $-16$, $2$, $-20$, $-22$, $-18$, $-6$, $-12$, $-14$, $-10)$ | $(12$, $14$, $44$, $16$, $32$, $-40$, $-22$, $-42$, $-28$, $34$, $38$, $-2$, $-8$, $-18$, $26$, $10$, $-30$, $-6$, $-4$, $46$, $24$, $36$, $-20)$\\

\noindent$11n82$ : $3$ | $(4$, $8$, $-16$, $-12$, $2$, $-18$, $-20$, $-22$, $-10$, $-6$, $-14)$ | $(12$, $28$, $26$, $-4$, $-16$, $-20$, $-22$, $46$, $30$, $32$, $48$, $44$, $6$, $-24$, $-34$, $-40$, $42$, $8$, $14$, $10$, $-18$, $-2$, $-36$, $-38)$\\

\noindent$11n83$ : $3$ | $(4$, $8$, $-16$, $-12$, $2$, $-18$, $-22$, $-20$, $-10$, $-6$, $-14)$ | $(56$, $26$, $30$, $20$, $22$, $32$, $24$, $58$, $54$, $-36$, $-50$, $-34$, $-46$, $-40$, $-38$, $-48$, $68$, $52$, $72$, $64$, $74$, $80$, $76$, $66$, $70$, $-8$, $-60$, $-42$, $-44$, $78$, $18$, $28$, $-4$, $-12$, $-10$, $-6$, $-62$, $-2$, $-14$, $-16)$\\

\noindent$11n84$ : $3$ | $(4$, $8$, $-18$, $-12$, $2$, $-16$, $-20$, $-6$, $-10$, $-22$, $-14)$ | $(10$, $52$, $54$, $12$, $-40$, $-38$, $-18$, $-24$, $26$, $58$, $50$, $56$, $-46$, $-6$, $-2$, $-20$, $-22$, $-4$, $28$, $60$, $30$, $34$, $16$, $14$, $32$, $-42$, $-36$, $-44$, $-48$, $-8)$\\

\noindent$11n85$ : $3$ | $(4$, $10$, $12$, $-14$, $16$, $2$, $-18$, $20$, $22$, $-6$, $8)$ | $(24$, $26$, $22$, $36$, $20$, $28$, $50$, $32$, $-68$, $-60$, $-38$, $-56$, $-54$, $-40$, $-62$, $-66$, $-70$, $-58$, $82$, $92$, $52$, $34$, $18$, $30$, $-64$, $-72$, $-2$, $94$, $80$, $84$, $90$, $74$, $76$, $88$, $86$, $78$, $-12$, $-48$, $-14$, $-42$, $-6$, $-8$, $-44$, $-16$, $-46$, $-10$, $-4)$\\

\noindent$11n86$ : $3$ | $(4$, $10$, $12$, $-14$, $-16$, $2$, $-18$, $-20$, $-22$, $-6$, $-8)$ | $(22$, $40$, $32$, $92$, $28$, $44$, $26$, $90$, $34$, $38$, $-54$, $-56$, $-68$, $-46$, $-64$, $-60$, $-50$, $-104$, $-102$, $-52$, $-58$, $-66$, $80$, $112$, $106$, $74$, $86$, $116$, $84$, $76$, $108$, $110$, $78$, $82$, $114$, $36$, $-18$, $-4$, $-96$, $-10$, $-12$, $-98$, $-2$, $-20$, $-70$, $-48$, $-62$, $30$, $42$, $24$, $88$, $72$, $-16$, $-6$, $-94$, $-8$, $-14$, $-100)$\\

\noindent$11n87$ : $3$ | $(4$, $10$, $12$, $14$, $18$, $2$, $20$, $8$, $-22$, $6$, $-16)$ | $(24$, $48$, $52$, $42$, $50$, $-4$, $-32$, $-30$, $-6$, $-8$, $-28$, $-34$, $-36$, $70$, $44$, $54$, $46$, $72$, $68$, $10$, $18$, $-62$, $-56$, $-66$, $-38$, $-40$, $-64$, $14$, $22$, $26$, $20$, $16$, $12$, $-2$, $-60$, $-58)$\\

\noindent$11n88$ : $3$ | $(4$, $10$, $12$, $16$, $18$, $2$, $-20$, $6$, $8$, $-22$, $-14)$ | $(26$, $30$, $8$, $-22$, $-38$, $-34$, $32$, $10$, $52$, $42$, $50$, $-6$, $-16$, $-18$, $-20$, $-36$, $46$, $44$, $48$, $24$, $28$, $-2$, $-12$, $-14$, $-4$, $-40)$\\

\noindent$11n89$ : $3$ | $(4$, $10$, $12$, $-16$, $-18$, $2$, $20$, $-6$, $-8$, $22$, $14)$ | $(54$, $130$, $64$, $44$, $72$, $46$, $62$, $128$, $56$, $52$, $132$, $66$, $42$, $70$, $48$, $60$, $126$, $58$, $50$, $68$, $-118$, $-88$, $-74$, $-92$, $-122$, $-100$, $-82$, $-80$, $-98$, $-124$, $-94$, $-76$, $-86$, $-116$, $-120$, $-90$, $138$, $148$, $112$, $104$, $156$, $102$, $114$, $146$, $140$, $136$, $150$, $110$, $106$, $154$, $-20$, $-2$, $-16$, $-34$, $-32$, $-14$, $-4$, $-22$, $144$, $142$, $134$, $152$, $108$, $-96$, $-78$, $-84$, $-38$, $-28$, $-10$, $-8$, $-26$, $-40$, $-24$, $-6$, $-12$, $-30$, $-36$, $-18)$\\

\noindent$11n90$ : $3$ | $(4$, $10$, $12$, $16$, $18$, $2$, $-20$, $8$, $6$, $-22$, $-14)$ | $(16$, $38$, $18$, $14$, $40$, $58$, $-28$, $-34$, $-36$, $-30$, $-10$, $-26$, $42$, $56$, $60$, $50$, $54$, $62$, $52$, $-32$, $-44$, $24$, $22$, $12$, $20$, $-4$, $-2$, $-8$, $-46$, $-48$, $-6)$\\

\noindent$11n91$ : $3$ | $(4$, $10$, $12$, $-16$, $-18$, $2$, $20$, $-8$, $-6$, $22$, $14)$ | $(28$, $30$, $34$, $20$, $38$, $18$, $32$, $-4$, $-46$, $-40$, $-50$, $36$, $16$, $66$, $-26$, $-44$, $-42$, $-52$, $-48$, $56$, $62$, $14$, $64$, $54$, $58$, $60$, $-24$, $-10$, $-8$, $-22$, $-6$, $-12$, $-2)$\\

\noindent$11n92$ : $3$ | $(4$, $10$, $12$, $-20$, $-16$, $2$, $-18$, $-8$, $-22$, $-6$, $-14)$ | $(42$, $26$, $114$, $36$, $32$, $110$, $30$, $38$, $116$, $24$, $44$, $-46$, $-102$, $-96$, $-86$, $-108$, $-90$, $-92$, $-106$, $-84$, $-98$, $-100$, $66$, $68$, $40$, $28$, $112$, $34$, $-8$, $-14$, $-50$, $-2$, $-20$, $-22$, $-18$, $-4$, $-52$, $-12$, $-10$, $-54$, $-6$, $-16$, $60$, $74$, $78$, $56$, $80$, $72$, $62$, $118$, $64$, $70$, $82$, $58$, $76$, $-88$, $-94$, $-104$, $-48)$\\

\noindent$11n93$ : $3$ | $(4$, $10$, $14$, $-12$, $2$, $-20$, $18$, $22$, $6$, $-8$, $16)$ | $(48$, $68$, $78$, $74$, $72$, $80$, $16$, $-64$, $-62$, $-2$, $-58$, $40$, $54$, $42$, $20$, $82$, $18$, $44$, $52$, $-8$, $-56$, $-4$, $-12$, $-32$, $-30$, $-10$, $-6$, $76$, $70$, $50$, $46$, $66$, $-14$, $-34$, $-28$, $-22$, $-38$, $-24$, $-26$, $-36$, $-60)$\\

\noindent$11n94$ : $3$ | $(4$, $10$, $14$, $12$, $2$, $20$, $18$, $-22$, $6$, $8$, $-16)$ | $(68$, $72$, $86$, $74$, $94$, $96$, $98$, $92$, $76$, $88$, $70$, $-2$, $-8$, $-16$, $-78$, $-18$, $-6$, $-4$, $-20$, $-80$, $-14$, $-10$, $66$, $40$, $28$, $30$, $38$, $22$, $36$, $32$, $26$, $42$, $-12$, $-82$, $-54$, $-56$, $-84$, $-60$, $-50$, $90$, $100$, $24$, $34$, $-58$, $-52$, $-48$, $-62$, $-44$, $-64$, $-46)$\\

\noindent$11n95$ : $3$ | $(4$, $10$, $14$, $-12$, $2$, $-20$, $18$, $-22$, $6$, $-8$, $-16)$ | $(20$, $18$, $44$, $42$, $16$, $38$, $-50$, $-54$, $-28$, $-60$, $-30$, $-32$, $-6$, $72$, $62$, $46$, $40$, $14$, $-52$, $-48$, $-56$, $-26$, $-58$, $22$, $68$, $66$, $64$, $70$, $24$, $74$, $-2$, $-10$, $-36$, $-12$, $-34$, $-8$, $-4)$\\

\noindent$11n96$ : $3$ | $(4$, $10$, $-14$, $-12$, $2$, $-20$, $-18$, $-22$, $-6$, $-8$, $-16)$ | $(28$, $80$, $20$, $58$, $22$, $78$, $26$, $30$, $82$, $-70$, $-60$, $-34$, $-38$, $-64$, $-66$, $-40$, $88$, $48$, $50$, $90$, $56$, $-18$, $-4$, $-10$, $-76$, $-12$, $-2$, $-16$, $-72$, $86$, $46$, $52$, $92$, $54$, $44$, $84$, $32$, $24$, $-36$, $-62$, $-68$, $-42$, $-6$, $-8$, $-74$, $-14)$\\

\noindent$11n97$ : $3$ | $(4$, $10$, $14$, $12$, $-18$, $2$, $6$, $-20$, $-22$, $-8$, $-16)$ | $(14$, $20$, $42$, $38$, $40$, $-6$, $-52$, $-28$, $-32$, $-26$, $-54$, $12$, $62$, $34$, $66$, $64$, $-48$, $-2$, $-22$, $-58$, $-56$, $-24$, $-30$, $16$, $18$, $44$, $36$, $60$, $10$, $-8$, $-4$, $-50$, $-46)$\\

\noindent$11n98$ : $3$ | $(4$, $10$, $-14$, $-12$, $18$, $2$, $-6$, $20$, $22$, $8$, $16)$ | $(20$, $24$, $38$, $100$, $30$, $32$, $98$, $36$, $26$, $-86$, $-88$, $-90$, $-84$, $-50$, $-42$, $-78$, $-44$, $-48$, $-82$, $70$, $60$, $76$, $64$, $66$, $74$, $58$, $22$, $-2$, $-18$, $-92$, $-8$, $-12$, $-96$, $-14$, $-6$, $-4$, $-16$, $-94$, $-10$, $62$, $68$, $72$, $56$, $104$, $54$, $102$, $28$, $34$, $-46$, $-80$, $-40$, $-52)$\\

\noindent$11n99$ : $3$ | $(4$, $10$, $-14$, $-12$, $-18$, $2$, $-6$, $-20$, $-22$, $-8$, $-16)$ | $(38$, $32$, $36$, $14$, $-42$, $-50$, $-44$, $-40$, $-2$, $-4$, $34$, $12$, $10$, $16$, $60$, $18$, $-48$, $-46$, $-30$, $-28$, $54$, $52$, $58$, $20$, $56$, $-24$, $-8$, $-26$, $-22$, $-6)$\\

\noindent$11n100$ : $3$ | $(4$, $10$, $14$, $12$, $-18$, $2$, $6$, $-22$, $-20$, $-8$, $-16)$ | $(16$, $122$, $132$, $-96$, $-106$, $-84$, $-112$, $-90$, $-46$, $64$, $74$, $8$, $130$, $124$, $14$, $68$, $70$, $12$, $126$, $128$, $10$, $72$, $66$, $62$, $76$, $-6$, $-24$, $-40$, $-36$, $-28$, $-2$, $-100$, $-102$, $-88$, $-114$, $-86$, $-104$, $-98$, $-4$, $-26$, $-38$, $56$, $138$, $116$, $142$, $60$, $78$, $52$, $134$, $120$, $18$, $118$, $136$, $54$, $80$, $58$, $140$, $-32$, $-44$, $-20$, $-48$, $-92$, $-110$, $-82$, $-108$, $-94$, $-50$, $-22$, $-42$, $-34$, $-30)$\\

\noindent$11n101$ : $3$ | $(4$, $10$, $-14$, $-12$, $18$, $2$, $-6$, $22$, $20$, $8$, $16)$ | $(8$, $54$, $52$, $-40$, $-42$, $-48$, $-44$, $-46$, $12$, $36$, $34$, $10$, $14$, $-16$, $-24$, $-2$, $-6$, $-50$, $-4$, $32$, $58$, $30$, $26$, $28$, $56$, $-20$, $-38$, $-18$, $-22)$\\

\noindent$11n102$ : $3$ | $(4$, $10$, $-14$, $-12$, $-18$, $2$, $-6$, $-22$, $-20$, $-8$, $-16)$ | $(14$, $58$, $36$, $22$, $38$, $60$, $-28$, $-30$, $-10$, $-4$, $-48$, $54$, $18$, $42$, $62$, $40$, $20$, $-50$, $-24$, $-46$, $-2$, $-12$, $16$, $56$, $52$, $-34$, $-6$, $-8$, $-32$, $-26$, $-44)$\\

\noindent$11n103$ : $3$ | $(4$, $10$, $14$, $12$, $-20$, $2$, $16$, $6$, $-22$, $-8$, $-18)$ | $(54$, $84$, $94$, $90$, $88$, $50$, $58$, $82$, $56$, $52$, $86$, $92$, $-6$, $-10$, $-20$, $-2$, $-16$, $-14$, $-4$, $-22$, $-8$, $26$, $98$, $96$, $-44$, $-78$, $-68$, $-66$, $-80$, $-46$, $48$, $36$, $32$, $102$, $30$, $38$, $24$, $40$, $28$, $100$, $34$, $-64$, $-70$, $-76$, $-42$, $-74$, $-72$, $-62$, $-60$, $-12$, $-18)$\\

\noindent$11n104$ : $3$ | $(4$, $10$, $-14$, $-16$, $2$, $-20$, $-22$, $-18$, $-6$, $-8$, $-12)$ | $(100$, $34$, $120$, $44$, $124$, $30$, $128$, $40$, $116$, $38$, $36$, $118$, $42$, $126$, $-60$, $-110$, $-102$, $-52$, $-68$, $-114$, $-64$, $-56$, $-106$, $122$, $32$, $130$, $88$, $80$, $74$, $94$, $96$, $72$, $82$, $86$, $-20$, $-14$, $-26$, $-8$, $-46$, $-6$, $-24$, $-16$, $-18$, $-22$, $-4$, $-48$, $-10$, $-28$, $-12$, $-50$, $132$, $90$, $78$, $76$, $92$, $98$, $70$, $84$, $-66$, $-54$, $-104$, $-108$, $-58$, $-62$, $-112$, $-2)$\\

\noindent$11n105$ : $3$ | $(4$, $10$, $14$, $16$, $2$, $-20$, $-22$, $18$, $8$, $6$, $-12)$ | $(36$, $24$, $42$, $30$, $80$, $28$, $44$, $26$, $82$, $100$, $88$, $-68$, $-74$, $-46$, $-78$, $-64$, $-18$, $-72$, $-70$, $-20$, $-66$, $-76$, $94$, $60$, $32$, $40$, $22$, $38$, $34$, $-16$, $-10$, $-8$, $98$, $90$, $86$, $102$, $84$, $92$, $96$, $62$, $-48$, $-58$, $-2$, $-54$, $-52$, $-4$, $-14$, $-12$, $-6$, $-50$, $-56)$\\

\noindent$11n106$ : $3$ | $(4$, $10$, $-14$, $-16$, $2$, $20$, $22$, $-18$, $-8$, $-6$, $12)$ | $(32$, $80$, $36$, $70$, $74$, $40$, $76$, $28$, $30$, $78$, $38$, $72$, $-8$, $-86$, $-66$, $-52$, $-50$, $-46$, $-56$, $-62$, $100$, $26$, $106$, $94$, $92$, $108$, $90$, $96$, $104$, $24$, $102$, $98$, $88$, $-16$, $-44$, $-58$, $-60$, $-42$, $-64$, $-54$, $-48$, $34$, $68$, $-14$, $-18$, $-2$, $-82$, $-4$, $-20$, $-12$, $-10$, $-22$, $-6$, $-84)$\\

\noindent$11n107$ : $3$ | $(4$, $10$, $-14$, $-16$, $2$, $-20$, $-22$, $-18$, $-8$, $-6$, $-12)$ | $(42$, $34$, $22$, $20$, $36$, $40$, $44$, $32$, $46$, $-24$, $-50$, $48$, $76$, $64$, $78$, $70$, $-58$, $-52$, $-26$, $-28$, $-54$, $-56$, $-30$, $-60$, $-6$, $74$, $66$, $80$, $68$, $72$, $18$, $38$, $-10$, $-2$, $-14$, $-16$, $-4$, $-8$, $-62$, $-12)$\\

\noindent$11n108$ : $3$ | $(4$, $10$, $-14$, $-16$, $18$, $2$, $20$, $-6$, $22$, $12$, $8)$ | $(42$, $96$, $44$, $82$, $76$, $90$, $50$, $88$, $74$, $84$, $46$, $94$, $40$, $-64$, $-8$, $-18$, $-72$, $-16$, $-10$, $-66$, $-62$, $-98$, $-58$, $-114$, $-106$, $128$, $36$, $30$, $134$, $122$, $118$, $80$, $78$, $92$, $48$, $86$, $-110$, $-54$, $-102$, $-26$, $-100$, $-56$, $-112$, $-108$, $-52$, $-104$, $-116$, $-60$, $120$, $136$, $28$, $38$, $126$, $130$, $34$, $32$, $132$, $124$, $-24$, $-2$, $-4$, $-22$, $-68$, $-12$, $-14$, $-70$, $-20$, $-6)$\\

\noindent$11n109$ : $3$ | $(4$, $10$, $-14$, $16$, $18$, $2$, $20$, $-22$, $8$, $12$, $-6)$ | $(44$, $50$, $18$, $14$, $54$, $-24$, $-20$, $-70$, $-68$, $38$, $86$, $84$, $40$, $48$, $46$, $42$, $52$, $16$, $-6$, $-58$, $-74$, $-64$, $-62$, $-76$, $-60$, $-66$, $-72$, $-22$, $12$, $92$, $78$, $96$, $80$, $90$, $36$, $88$, $82$, $94$, $-28$, $-2$, $-32$, $-10$, $-56$, $-8$, $-34$, $-4$, $-26$, $-30)$\\

\noindent$11n110$ : $3$ | $(4$, $10$, $14$, $-16$, $-18$, $2$, $20$, $22$, $-8$, $12$, $6)$ | $(26$, $56$, $74$, $82$, $50$, $80$, $76$, $54$, $52$, $78$, $-12$, $-6$, $-40$, $-42$, $-8$, $-10$, $-18$, $-14$, $-4$, $72$, $84$, $70$, $36$, $22$, $30$, $-62$, $-68$, $-44$, $-38$, $24$, $28$, $32$, $20$, $34$, $-16$, $-2$, $-58$, $-46$, $-66$, $-64$, $-48$, $-60)$\\

\noindent$11n111$ : $3$ | $(4$, $10$, $-14$, $-16$, $-18$, $2$, $-20$, $-22$, $-8$, $-12$, $-6)$ | $(14$, $18$, $38$, $34$, $36$, $20$, $-28$, $-26$, $-24$, $-30$, $-48$, $58$, $44$, $60$, $42$, $56$, $52$, $-50$, $-32$, $-22$, $-12$, $-2$, $16$, $40$, $54$, $-8$, $-6$, $-10$, $-4$, $-46)$\\

\noindent$11n112$ : $3$ | $(4$, $10$, $14$, $18$, $2$, $20$, $8$, $-22$, $6$, $12$, $-16)$ | $(60$, $50$, $52$, $58$, $44$, $-78$, $-76$, $-8$, $-32$, $-36$, $-4$, $88$, $14$, $10$, $64$, $46$, $56$, $54$, $48$, $62$, $12$, $-26$, $-66$, $-82$, $-72$, $-22$, $-70$, $-84$, $-68$, $-24$, $-74$, $-80$, $94$, $16$, $86$, $20$, $90$, $98$, $42$, $96$, $92$, $18$, $-34$, $-6$, $-2$, $-38$, $-30$, $-28$, $-40)$\\

\noindent$11n113$ : $3$ | $(4$, $10$, $-14$, $18$, $16$, $2$, $20$, $-22$, $8$, $12$, $-6)$ | $(40$, $44$, $20$, $16$, $48$, $-62$, $-52$, $-66$, $-58$, $-56$, $86$, $70$, $88$, $14$, $50$, $12$, $46$, $18$, $-6$, $-26$, $-24$, $-22$, $-54$, $-60$, $-64$, $78$, $72$, $84$, $36$, $82$, $74$, $76$, $80$, $38$, $42$, $-2$, $-32$, $-10$, $-30$, $-28$, $-8$, $-34$, $-4$, $-68)$\\

\noindent$11n114$ : $3$ | $(4$, $10$, $14$, $18$, $16$, $2$, $-20$, $22$, $8$, $-12$, $6)$ | $(28$, $56$, $70$, $54$, $26$, $-76$, $-82$, $-88$, $-4$, $-6$, $-90$, $-80$, $-78$, $-74$, $-44$, $98$, $14$, $22$, $92$, $20$, $16$, $96$, $100$, $12$, $24$, $-36$, $-38$, $-72$, $-42$, $-32$, $-46$, $-10$, $-48$, $-34$, $-40$, $18$, $94$, $62$, $64$, $50$, $66$, $60$, $30$, $58$, $68$, $52$, $-8$, $-2$, $-86$, $-84)$\\

\noindent$11n115$ : $3$ | $(4$, $10$, $14$, $-18$, $-16$, $2$, $20$, $22$, $-8$, $12$, $6)$ | $(32$, $128$, $138$, $146$, $122$, $144$, $140$, $126$, $34$, $30$, $130$, $136$, $-106$, $-86$, $-100$, $-112$, $-114$, $-98$, $-88$, $26$, $134$, $132$, $28$, $36$, $124$, $142$, $-12$, $-6$, $-24$, $-40$, $-44$, $-20$, $-2$, $-16$, $-48$, $-8$, $-10$, $-50$, $-14$, $-4$, $-22$, $-42$, $60$, $148$, $70$, $74$, $52$, $76$, $68$, $150$, $62$, $82$, $58$, $56$, $80$, $64$, $152$, $66$, $78$, $54$, $72$, $-90$, $-96$, $-116$, $-110$, $-102$, $-84$, $-104$, $-108$, $-118$, $-94$, $-92$, $-120$, $-38$, $-46$, $-18)$\\

\noindent$11n116$ : $3$ | $(4$, $10$, $-14$, $-18$, $-16$, $2$, $20$, $-22$, $-8$, $12$, $-6)$ | $(40$, $70$, $20$, $24$, $26$, $18$, $72$, $-2$, $-36$, $-52$, $-30$, $-32$, $-50$, $-34$, $46$, $64$, $48$, $68$, $42$, $-56$, $-58$, $-12$, $-4$, $-62$, $-8$, $66$, $44$, $14$, $74$, $16$, $28$, $22$, $-6$, $-10$, $-60$, $-54$, $-38)$\\

\noindent$11n117$ : $3$ | $(4$, $10$, $-14$, $-18$, $-16$, $2$, $-20$, $-22$, $-8$, $-12$, $-6)$ | $(20$, $68$, $62$, $66$, $18$, $22$, $70$, $-26$, $-32$, $-58$, $-56$, $-34$, $72$, $42$, $46$, $40$, $74$, $-12$, $-2$, $-8$, $-48$, $-4$, $-6$, $64$, $16$, $24$, $14$, $36$, $76$, $38$, $44$, $-28$, $-30$, $-60$, $-54$, $-52$, $-50$, $-10)$\\

\noindent$11n118$ : $3$ | $(4$, $10$, $-14$, $20$, $2$, $-16$, $-18$, $-8$, $-22$, $6$, $-12)$ | $(16$, $12$, $10$, $30$, $-32$, $-24$, $-34$, $-38$, $-36$, $14$, $40$, $44$, $-4$, $-20$, $-22$, $42$, $26$, $46$, $48$, $28$, $-8$, $-6$, $-2$, $-18)$\\

\noindent$11n119$ : $3$ | $(4$, $10$, $14$, $-20$, $2$, $16$, $18$, $8$, $22$, $-6$, $12)$ | $(70$, $74$, $92$, $38$, $88$, $78$, $48$, $46$, $80$, $86$, $40$, $42$, $84$, $82$, $44$, $50$, $76$, $90$, $-118$, $-104$, $-132$, $-64$, $-126$, $-98$, $-124$, $-66$, $-4$, $-34$, $-10$, $-18$, $-26$, $144$, $94$, $72$, $-110$, $-112$, $-108$, $-114$, $-122$, $-100$, $-128$, $-62$, $-130$, $-102$, $-120$, $-116$, $-106$, $-30$, $-14$, $148$, $58$, $140$, $134$, $52$, $154$, $156$, $152$, $54$, $136$, $138$, $56$, $150$, $96$, $146$, $60$, $142$, $-22$, $-6$, $-36$, $-8$, $-20$, $-24$, $-28$, $-16$, $-12$, $-32$, $-2$, $-68)$\\

\noindent$11n120$ : $3$ | $(4$, $10$, $-14$, $-20$, $2$, $-16$, $-18$, $-8$, $-22$, $-6$, $-12)$ | $(28$, $42$, $20$, $22$, $40$, $26$, $30$, $44$, $-68$, $-60$, $-50$, $-62$, $-66$, $-54$, $-56$, $82$, $16$, $18$, $24$, $-64$, $-52$, $-58$, $-32$, $-36$, $-6$, $72$, $78$, $84$, $80$, $70$, $48$, $74$, $76$, $46$, $-34$, $-4$, $-8$, $-38$, $-10$, $-2$, $-14$, $-12)$\\

\noindent$11n121$ : $3$ | $(4$, $10$, $-16$, $-12$, $18$, $2$, $-20$, $-6$, $22$, $8$, $-14)$ | $(38$, $40$, $76$, $70$, $84$, $46$, $82$, $68$, $78$, $42$, $36$, $-58$, $-6$, $-16$, $-66$, $-14$, $-8$, $-60$, $-56$, $-90$, $-54$, $-106$, $-98$, $110$, $32$, $26$, $116$, $88$, $74$, $72$, $86$, $44$, $80$, $-102$, $-50$, $-94$, $-22$, $-92$, $-52$, $-104$, $-100$, $-48$, $-96$, $-20$, $-2$, $118$, $24$, $34$, $108$, $112$, $30$, $28$, $114$, $-62$, $-10$, $-12$, $-64$, $-18$, $-4)$\\

\noindent$11n122$ : $3$ | $(4$, $10$, $-16$, $-12$, $-18$, $2$, $-20$, $-6$, $-22$, $-8$, $-14)$ | $(38$, $16$, $14$, $20$, $36$, $18$, $-32$, $-48$, $-30$, $-34$, $-26$, $40$, $42$, $60$, $56$, $46$, $44$, $58$, $-28$, $-50$, $-52$, $-6$, $-12$, $-4$, $54$, $22$, $-24$, $-2$, $-10$, $-8)$\\

\noindent$11n123$ : $3$ | $(4$, $10$, $-16$, $-14$, $12$, $2$, $-18$, $-20$, $-22$, $-8$, $-6)$ | $(94$, $76$, $100$, $88$, $82$, $104$, $80$, $90$, $98$, $74$, $96$, $92$, $78$, $102$, $84$, $86$, $-60$, $-28$, $-10$, $-12$, $-26$, $-4$, $-18$, $-20$, $-2$, $-24$, $-14$, $-8$, $-30$, $122$, $124$, $56$, $36$, $40$, $52$, $128$, $46$, $-112$, $-118$, $-68$, $-106$, $-64$, $-58$, $-62$, $-108$, $-70$, $-116$, $-114$, $-72$, $-110$, $-120$, $-66$, $38$, $54$, $126$, $44$, $48$, $130$, $50$, $42$, $34$, $-32$, $-6$, $-16$, $-22)$\\

\noindent$11n124$ : $3$ | $(4$, $10$, $16$, $-14$, $-18$, $2$, $-8$, $20$, $22$, $-12$, $6)$ | $(52$, $32$, $46$, $58$, $60$, $66$, $62$, $56$, $48$, $34$, $50$, $54$, $64$, $-12$, $-10$, $-80$, $-74$, $94$, $104$, $88$, $26$, $30$, $-82$, $-72$, $-36$, $-76$, $-78$, $-38$, $-70$, $-84$, $-68$, $-40$, $-42$, $28$, $86$, $102$, $96$, $92$, $106$, $90$, $98$, $100$, $-8$, $-14$, $-24$, $-2$, $-20$, $-18$, $-4$, $-44$, $-6$, $-16$, $-22)$\\

\noindent$11n125$ : $3$ | $(4$, $10$, $-16$, $-18$, $14$, $2$, $20$, $22$, $-6$, $12$, $8)$ | $(46$, $76$, $72$, $70$, $78$, $44$, $48$, $50$, $42$, $-10$, $-2$, $-58$, $-6$, $54$, $20$, $84$, $22$, $56$, $24$, $52$, $-62$, $-64$, $-30$, $-32$, $-66$, $-60$, $-8$, $-4$, $74$, $68$, $80$, $18$, $82$, $-14$, $-38$, $-26$, $-36$, $-34$, $-28$, $-40$, $-16$, $-12)$\\

\noindent$11n126$ : $3$ | $(4$, $10$, $-18$, $-12$, $2$, $-16$, $-20$, $-8$, $-22$, $-14$, $-6)$ | $(36$, $14$, $12$, $34$, $32$, $-22$, $-48$, $-42$, $-40$, $-8$, $52$, $50$, $18$, $38$, $16$, $-24$, $-20$, $-46$, $-44$, $58$, $30$, $56$, $60$, $54$, $-4$, $-10$, $-6$, $-2$, $-28$, $-26)$\\

\noindent$11n127$ : $3$ | $(4$, $10$, $18$, $-12$, $-16$, $2$, $-20$, $-8$, $22$, $6$, $-14)$ | $(36$, $84$, $70$, $74$, $42$, $78$, $66$, $80$, $40$, $72$, $-6$, $-46$, $-12$, $-14$, $-48$, $-4$, $-108$, $-88$, $-90$, $-106$, $-62$, $-98$, $76$, $68$, $82$, $38$, $34$, $110$, $32$, $20$, $122$, $116$, $26$, $-102$, $-94$, $-58$, $-60$, $-96$, $-100$, $-64$, $-104$, $-92$, $-56$, $-54$, $126$, $112$, $30$, $22$, $120$, $118$, $24$, $28$, $114$, $124$, $86$, $-2$, $-50$, $-16$, $-10$, $-44$, $-8$, $-18$, $-52)$\\

\noindent$11n128$ : $3$ | $(4$, $10$, $-18$, $-14$, $-16$, $2$, $-6$, $-20$, $-22$, $-12$, $-8)$ | $(44$, $68$, $40$, $70$, $-60$, $-54$, $-62$, $-2$, $-66$, $-4$, $-6$, $72$, $10$, $12$, $74$, $18$, $16$, $76$, $14$, $20$, $-30$, $-34$, $-36$, $-28$, $-24$, $-8$, $-26$, $48$, $22$, $50$, $52$, $46$, $42$, $-32$, $-38$, $-58$, $-56$, $-64)$\\

\noindent$11n129$ : $3$ | $(4$, $10$, $-18$, $-14$, $-16$, $2$, $-20$, $-8$, $-22$, $-12$, $-6)$ | $(44$, $16$, $20$, $14$, $-48$, $-56$, $-50$, $-22$, $-54$, $-52$, $36$, $60$, $64$, $62$, $10$, $12$, $18$, $-4$, $-6$, $-32$, $-28$, $-24$, $-26$, $42$, $58$, $38$, $34$, $40$, $-30$, $-8$, $-2$, $-46)$\\

\noindent$11n130$ : $3$ | $(4$, $10$, $-18$, $-14$, $-16$, $2$, $-20$, $-22$, $-12$, $-6$, $-8)$ | $(16$, $52$, $18$, $20$, $24$, $76$, $-62$, $-56$, $-54$, $-42$, $-32$, $-34$, $-40$, $50$, $14$, $44$, $66$, $46$, $12$, $48$, $68$, $-6$, $-64$, $-10$, $-60$, $-58$, $70$, $80$, $72$, $26$, $74$, $78$, $22$, $-8$, $-4$, $-2$, $-28$, $-38$, $-36$, $-30)$\\

\noindent$11n131$ : $3$ | $(4$, $10$, $18$, $-16$, $-14$, $2$, $20$, $-8$, $22$, $12$, $6)$ | $(26$, $16$, $20$, $22$, $30$, $60$, $-34$, $-48$, $-42$, $-40$, $-50$, $-38$, $-44$, $-46$, $-36$, $-76$, $62$, $58$, $86$, $52$, $66$, $54$, $88$, $56$, $64$, $-6$, $-68$, $-2$, $-72$, $-80$, $-78$, $-74$, $-4$, $18$, $24$, $28$, $14$, $84$, $82$, $12$, $-10$, $-32$, $-8$, $-70)$\\

\noindent$11n132$ : $3$ | $(4$, $10$, $-18$, $-16$, $-14$, $2$, $-20$, $-8$, $-22$, $-12$, $-6)$ | $(54$, $70$, $48$, $72$, $56$, $52$, $68$, $50$, $-36$, $-34$, $-8$, $-2$, $-12$, $-14$, $-4$, $-6$, $74$, $58$, $76$, $28$, $24$, $80$, $22$, $30$, $-38$, $-66$, $-42$, $-44$, $-64$, $18$, $32$, $20$, $78$, $26$, $-40$, $-46$, $-62$, $-60$, $-16$, $-10)$\\

\noindent$11n133$ : $3$ | $(4$, $10$, $-18$, $-20$, $2$, $-16$, $-6$, $-8$, $-22$, $-14$, $-12)$ | $(38$, $52$, $54$, $40$, $56$, $-44$, $-50$, $-2$, $-6$, $-8$, $58$, $12$, $14$, $60$, $16$, $-4$, $-10$, $-22$, $-26$, $-28$, $18$, $32$, $34$, $20$, $36$, $-24$, $-30$, $-42$, $-46$, $-48)$\\

\noindent$11n134$ : $3$ | $(4$, $10$, $-18$, $-22$, $-14$, $2$, $-20$, $-8$, $-6$, $-12$, $-16)$ | $(48$, $22$, $26$, $16$, $18$, $28$, $20$, $-32$, $-64$, $-30$, $-60$, $-54$, $-52$, $-62$, $42$, $66$, $46$, $50$, $24$, $-4$, $-12$, $-10$, $-6$, $-58$, $-56$, $72$, $38$, $70$, $74$, $68$, $40$, $44$, $-8$, $-14$, $-2$, $-36$, $-34)$\\

\noindent$11n135$ : $3$ | $(4$, $12$, $14$, $-16$, $22$, $-18$, $2$, $-20$, $-10$, $-6$, $8)$ | $(28$, $48$, $56$, $42$, $54$, $50$, $52$, $-10$, $-4$, $-32$, $-34$, $-6$, $-8$, $-36$, $-62$, $46$, $58$, $44$, $74$, $16$, $24$, $-68$, $-60$, $-64$, $-38$, $-72$, $-70$, $-40$, $-66$, $20$, $78$, $76$, $18$, $22$, $26$, $14$, $-12$, $-2$, $-30)$\\

\noindent$11n136$ : $3$ | $(4$, $12$, $-16$, $14$, $18$, $2$, $8$, $-22$, $20$, $10$, $-6)$ | $(42$, $52$, $30$, $48$, $46$, $38$, $56$, $34$, $84$, $32$, $54$, $40$, $44$, $50$, $-72$, $-62$, $-82$, $-114$, $-110$, $-78$, $-66$, $-68$, $-76$, $-108$, $-74$, $-70$, $-64$, $-80$, $-112$, $36$, $86$, $92$, $102$, $126$, $104$, $90$, $88$, $106$, $124$, $100$, $94$, $118$, $-116$, $-6$, $-28$, $-4$, $-20$, $-14$, $-58$, $-12$, $-22$, $-2$, $-26$, $-8$, $122$, $98$, $96$, $120$, $-18$, $-16$, $-60$, $-10$, $-24)$\\

\noindent$11n137$ : $3$ | $(4$, $12$, $-16$, $-14$, $18$, $20$, $2$, $-6$, $22$, $10$, $8)$ | $(28$, $48$, $54$, $60$, $52$, $50$, $30$, $26$, $46$, $56$, $58$, $94$, $-44$, $-68$, $-66$, $-42$, $-22$, $-18$, $-4$, $-10$, $86$, $96$, $92$, $-70$, $-64$, $-40$, $-62$, $-72$, $-80$, $-74$, $84$, $38$, $88$, $98$, $90$, $36$, $82$, $32$, $24$, $34$, $-20$, $-6$, $-8$, $-12$, $-2$, $-16$, $-78$, $-76$, $-14)$\\

\noindent$11n138$ : $3$ | $(4$, $12$, $-16$, $-14$, $-18$, $-20$, $2$, $-6$, $-22$, $-10$, $-8)$ | $(52$, $90$, $62$, $92$, $50$, $54$, $88$, $60$, $94$, $48$, $-18$, $-8$, $-2$, $-12$, $-44$, $-42$, $-14$, $-4$, $-6$, $-16$, $58$, $86$, $56$, $20$, $-76$, $-66$, $-68$, $-78$, $-84$, $-82$, $-72$, $36$, $24$, $98$, $26$, $34$, $38$, $22$, $96$, $28$, $32$, $40$, $30$, $-80$, $-70$, $-64$, $-74$, $-46$, $-10)$\\

\noindent$11n139$ : $3$ | $(4$, $12$, $16$, $14$, $-20$, $-18$, $2$, $6$, $-22$, $-10$, $-8)$ | $(34$, $56$, $64$, $-44$, $-42$, $-48$, $-2$, $54$, $66$, $12$, $38$, $6$, $10$, $40$, $8$, $-24$, $-18$, $-16$, $-4$, $-46$, $60$, $30$, $62$, $58$, $32$, $36$, $14$, $-52$, $-20$, $-26$, $-28$, $-22$, $-50)$\\

\noindent$11n140$ : $3$ | $(4$, $12$, $-16$, $-14$, $20$, $18$, $2$, $-6$, $22$, $10$, $8)$ | $(28$, $30$, $72$, $40$, $66$, $36$, $26$, $32$, $70$, $42$, $68$, $34$, $-48$, $-60$, $-82$, $-58$, $-46$, $-50$, $-90$, $-86$, $-54$, $96$, $106$, $78$, $76$, $104$, $98$, $94$, $108$, $80$, $74$, $38$, $-88$, $-52$, $-44$, $-56$, $-84$, $-12$, $-24$, $-14$, $-2$, $110$, $92$, $100$, $102$, $-62$, $-6$, $-18$, $-20$, $-8$, $-64$, $-10$, $-22$, $-16$, $-4)$\\

\noindent$11n141$ : $3$ | $(4$, $12$, $-16$, $-14$, $-20$, $-18$, $2$, $-6$, $-22$, $-10$, $-8)$ | $(26$, $10$, $-48$, $-38$, $-46$, $-16$, $32$, $34$, $30$, $6$, $8$, $28$, $-42$, $-44$, $-40$, $-50$, $-52$, $-4$, $-20$, $56$, $60$, $62$, $58$, $36$, $54$, $12$, $-14$, $-18$, $-22$, $-2$, $-24)$\\

\noindent$11n142$ : $3$ | $(4$, $12$, $-16$, $-20$, $-14$, $-18$, $2$, $-8$, $-22$, $-10$, $-6)$ | $(46$, $20$, $40$, $24$, $42$, $18$, $-66$, $-52$, $-64$, $-28$, $-56$, $-60$, $74$, $80$, $68$, $14$, $16$, $44$, $22$, $-58$, $-26$, $-62$, $-54$, $-30$, $-12$, $-32$, $84$, $70$, $78$, $76$, $72$, $82$, $50$, $48$, $-2$, $-36$, $-8$, $-6$, $-38$, $-4$, $-10$, $-34)$\\

\noindent$11n143$ : $3$ | $(4$, $12$, $-16$, $22$, $-18$, $2$, $-20$, $-10$, $-6$, $-14$, $8)$ | $(14$, $22$, $36$, $38$, $20$, $16$, $-56$, $-58$, $-30$, $-28$, $-60$, $-54$, $-50$, $70$, $46$, $72$, $24$, $-52$, $-62$, $-26$, $-32$, $-2$, $-10$, $18$, $40$, $64$, $66$, $42$, $74$, $44$, $68$, $-6$, $-34$, $-4$, $-8$, $-48$, $-12)$\\

\noindent$11n144$ : $3$ | $(4$, $12$, $18$, $-14$, $16$, $2$, $22$, $-20$, $10$, $6$, $-8)$ | $(68$, $36$, $122$, $154$, $132$, $144$, $142$, $134$, $152$, $120$, $150$, $136$, $140$, $146$, $130$, $156$, $124$, $-56$, $-50$, $-28$, $-10$, $-14$, $-24$, $138$, $148$, $128$, $158$, $126$, $-4$, $-20$, $-18$, $-6$, $-32$, $-54$, $-52$, $-30$, $-8$, $-16$, $-22$, $-48$, $-26$, $-12$, $42$, $74$, $62$, $58$, $78$, $46$, $80$, $38$, $70$, $66$, $34$, $64$, $72$, $40$, $82$, $44$, $76$, $60$, $-90$, $-106$, $-104$, $-2$, $-100$, $-110$, $-86$, $-116$, $-94$, $-96$, $-114$, $-84$, $-112$, $-98$, $-92$, $-118$, $-88$, $-108$, $-102)$\\

\noindent$11n145$ : $3$ | $(4$, $14$, $-10$, $-20$, $-22$, $18$, $2$, $-8$, $-6$, $12$, $-16)$ | $(28$, $84$, $18$, $62$, $22$, $80$, $24$, $60$, $-74$, $-34$, $-38$, $-78$, $-42$, $-30$, $98$, $56$, $48$, $90$, $92$, $50$, $54$, $96$, $86$, $-4$, $-66$, $-10$, $-12$, $-68$, $-2$, $-44$, $-76$, $-36$, $20$, $82$, $26$, $58$, $46$, $88$, $94$, $52$, $-40$, $-32$, $-72$, $-16$, $-6$, $-64$, $-8$, $-14$, $-70)$\\

\noindent$11n146$ : $3$ | $(4$, $14$, $10$, $20$, $22$, $-18$, $2$, $8$, $6$, $-12$, $16)$ | $(18$, $24$, $30$, $12$, $72$, $-74$, $-44$, $-56$, $-38$, $-80$, $-96$, $-78$, $-40$, $-58$, $-42$, $-76$, $-94$, $-82$, $114$, $100$, $60$, $104$, $110$, $32$, $22$, $20$, $34$, $112$, $102$, $-88$, $-4$, $-48$, $-52$, $-36$, $-54$, $-46$, $108$, $106$, $62$, $98$, $66$, $68$, $16$, $26$, $28$, $14$, $70$, $64$, $-50$, $-2$, $-86$, $-90$, $-6$, $-10$, $-8$, $-92$, $-84)$\\

\noindent$11n147$ : $3$ | $(4$, $14$, $-10$, $-20$, $-22$, $-18$, $2$, $-8$, $-6$, $-12$, $-16)$ | $(64$, $70$, $58$, $56$, $72$, $16$, $-94$, $-92$, $-28$, $-34$, $-42$, $-6$, $118$, $106$, $50$, $76$, $52$, $62$, $66$, $68$, $60$, $54$, $74$, $14$, $-12$, $-78$, $-98$, $-88$, $-24$, $-86$, $-100$, $-80$, $-82$, $-102$, $-84$, $-26$, $-90$, $-96$, $112$, $18$, $124$, $104$, $120$, $22$, $116$, $108$, $48$, $110$, $114$, $20$, $122$, $-38$, $-2$, $-10$, $-46$, $-30$, $-32$, $-44$, $-8$, $-4$, $-40$, $-36)$\\

\noindent$11n148$ : $3$ | $(4$, $14$, $-18$, $-20$, $-22$, $-6$, $2$, $-8$, $-10$, $-12$, $-16)$ | $(86$, $30$, $94$, $106$, $108$, $96$, $32$, $84$, $36$, $100$, $112$, $102$, $90$, $-118$, $-142$, $-62$, $-132$, $-128$, $-66$, $-68$, $-24$, $-74$, $-6$, $-10$, $-78$, $-20$, $-22$, $-76$, $-8$, $46$, $152$, $162$, $168$, $146$, $88$, $28$, $92$, $104$, $110$, $98$, $34$, $-130$, $-64$, $-144$, $-116$, $-120$, $-140$, $-60$, $-134$, $-126$, $-114$, $-122$, $-138$, $-58$, $-136$, $-124$, $52$, $40$, $158$, $156$, $42$, $54$, $50$, $148$, $166$, $164$, $150$, $48$, $56$, $44$, $154$, $160$, $38$, $-18$, $-80$, $-12$, $-4$, $-72$, $-26$, $-70$, $-2$, $-14$, $-82$, $-16)$\\

\noindent$11n149$ : $3$ | $(6$, $8$, $10$, $18$, $16$, $-20$, $-22$, $4$, $2$, $-14$, $-12)$ | $(62$, $30$, $56$, $20$, $18$, $58$, $32$, $60$, $-36$, $-66$, $-38$, $-44$, $-50$, $-48$, $-42$, $-70$, $92$, $108$, $64$, $104$, $96$, $88$, $102$, $98$, $86$, $100$, $-46$, $-40$, $-68$, $-34$, $-72$, $-74$, $-8$, $106$, $94$, $90$, $110$, $22$, $54$, $28$, $24$, $52$, $26$, $-80$, $-2$, $-14$, $-12$, $-4$, $-78$, $-84$, $-82$, $-76$, $-6$, $-10$, $-16)$\\

\noindent$11n150$ : $3$ | $(6$, $8$, $10$, $18$, $16$, $-22$, $-20$, $4$, $2$, $-14$, $-12)$ | $(78$, $90$, $102$, $134$, $48$, $136$, $100$, $92$, $76$, $96$, $-72$, $-14$, $-4$, $-62$, $-64$, $-6$, $-12$, $-70$, $-18$, $-74$, $-16$, $-2$, $-108$, $-124$, $84$, $32$, $24$, $42$, $140$, $44$, $130$, $104$, $132$, $46$, $138$, $98$, $94$, $-116$, $-56$, $-58$, $-68$, $-10$, $-8$, $-66$, $-60$, $-54$, $-118$, $-114$, $-112$, $-120$, $-52$, $-128$, $28$, $88$, $80$, $36$, $20$, $38$, $40$, $22$, $34$, $82$, $86$, $30$, $26$, $-106$, $-126$, $-50$, $-122$, $-110)$\\

\noindent$11n151$ : $3$ | $(6$, $8$, $12$, $2$, $-16$, $-18$, $4$, $-20$, $-22$, $-14$, $-10)$ | $(72$, $26$, $34$, $64$, $104$, $112$, $98$, $58$, $96$, $110$, $106$, $-8$, $-90$, $-84$, $-14$, $-16$, $-82$, $-92$, $24$, $74$, $70$, $28$, $32$, $66$, $78$, $22$, $76$, $68$, $30$, $-44$, $-4$, $-36$, $-80$, $-18$, $-12$, $-86$, $-88$, $-10$, $-20$, $108$, $94$, $60$, $100$, $114$, $102$, $62$, $-6$, $-46$, $-56$, $-42$, $-2$, $-38$, $-52$, $-50$, $-48$, $-54$, $-40)$\\

\noindent$11n152$ : $3$ | $(6$, $8$, $12$, $2$, $-18$, $-20$, $4$, $-10$, $-22$, $-14$, $-16)$ | $(66$, $28$, $20$, $22$, $60$, $100$, $106$, $94$, $54$, $-40$, $-50$, $-38$, $-16$, $-42$, $-48$, $-10$, $102$, $104$, $92$, $52$, $56$, $96$, $108$, $98$, $58$, $-6$, $-80$, $-82$, $-4$, $-90$, $-36$, $-14$, $-44$, $-46$, $-12$, $-34$, $32$, $70$, $62$, $24$, $18$, $26$, $64$, $68$, $30$, $-8$, $-78$, $-84$, $-2$, $-88$, $-74$, $-72$, $-76$, $-86)$\\

\noindent$11n153$ : $3$ | $(6$, $8$, $12$, $20$, $16$, $18$, $-22$, $4$, $10$, $2$, $-14)$ | $(78$, $104$, $112$, $40$, $114$, $102$, $80$, $76$, $106$, $110$, $38$, $116$, $28$, $-98$, $-96$, $-12$, $-2$, $-16$, $-92$, $-54$, $70$, $86$, $32$, $120$, $34$, $88$, $68$, $72$, $84$, $30$, $118$, $36$, $90$, $-20$, $-6$, $-8$, $-22$, $-62$, $-48$, $-46$, $-60$, $-24$, $-10$, $-4$, $-18$, $108$, $74$, $82$, $100$, $-26$, $-58$, $-44$, $-50$, $-64$, $-66$, $-52$, $-42$, $-56$, $-94$, $-14)$\\

\noindent$11n154$ : $3$ | $(6$, $8$, $12$, $20$, $-16$, $-18$, $22$, $4$, $-10$, $2$, $14)$ | $(22$, $54$, $18$, $50$, $48$, $72$, $46$, $52$, $-66$, $-36$, $-40$, $-42$, $-10$, $-12$, $-44$, $-14$, $-8$, $84$, $78$, $76$, $86$, $74$, $-62$, $-70$, $-34$, $-68$, $-64$, $-38$, $20$, $24$, $28$, $30$, $80$, $82$, $32$, $26$, $-60$, $-2$, $-56$, $-4$, $-16$, $-6$, $-58)$\\

\noindent$11n155$ : $3$ | $(6$, $8$, $12$, $20$, $18$, $16$, $-22$, $4$, $10$, $2$, $-14)$ | $(18$, $100$, $104$, $14$, $-42$, $-52$, $-40$, $-84$, $-90$, $-92$, $-82$, $-70$, $-78$, $-4$, $102$, $16$, $20$, $66$, $26$, $60$, $54$, $56$, $22$, $68$, $24$, $58$, $-74$, $-8$, $-72$, $-76$, $-6$, $-30$, $-2$, $-80$, $96$, $108$, $10$, $12$, $106$, $98$, $94$, $110$, $62$, $28$, $64$, $112$, $-34$, $-46$, $-48$, $-36$, $-88$, $-86$, $-38$, $-50$, $-44$, $-32)$\\

\noindent$11n156$ : $3$ | $(6$, $-10$, $12$, $14$, $16$, $-18$, $22$, $20$, $-4$, $2$, $8)$ | $(46$, $72$, $42$, $74$, $76$, $80$, $-60$, $-66$, $-2$, $-70$, $-4$, $-6$, $-10$, $78$, $14$, $82$, $20$, $18$, $84$, $16$, $22$, $-32$, $-36$, $-38$, $-30$, $-12$, $-28$, $-8$, $24$, $50$, $52$, $26$, $54$, $48$, $44$, $-34$, $-40$, $-56$, $-64$, $-62$, $-58$, $-68)$\\

\noindent$11n157$ : $3$ | $(6$, $-10$, $-12$, $14$, $-16$, $-18$, $22$, $20$, $-4$, $-2$, $8)$ | $(50$, $40$, $72$, $78$, $82$, $-32$, $-26$, $-2$, $-4$, $-6$, $-64$, $52$, $48$, $42$, $56$, $44$, $46$, $54$, $20$, $16$, $-70$, $-58$, $-10$, $-60$, $-68$, $-66$, $-62$, $-8$, $80$, $12$, $84$, $76$, $74$, $86$, $14$, $18$, $-38$, $-22$, $-36$, $-28$, $-30$, $-34$, $-24)$\\

\noindent$11n158$ : $3$ | $(6$, $10$, $12$, $18$, $20$, $16$, $-22$, $4$, $2$, $8$, $-14)$ | $(14$, $26$, $46$, $48$, $24$, $16$, $-76$, $-78$, $-40$, $-32$, $-34$, $-42$, $-80$, $-74$, $-66$, $92$, $62$, $90$, $28$, $94$, $60$, $88$, $84$, $-82$, $-44$, $-36$, $-30$, $-38$, $-2$, $-10$, $-70$, $20$, $52$, $54$, $18$, $22$, $50$, $56$, $96$, $58$, $86$, $-6$, $-4$, $-8$, $-72$, $-64$, $-68$, $-12)$\\

\noindent$11n159$ : $3$ | $(6$, $10$, $12$, $-18$, $-20$, $16$, $22$, $4$, $2$, $-8$, $14)$ | $(36$, $32$, $22$, $124$, $122$, $24$, $30$, $38$, $86$, $-102$, $-108$, $-46$, $-114$, $-96$, $-116$, $-48$, $-106$, $-104$, $-50$, $-56$, $-64$, $-8$, $80$, $92$, $72$, $130$, $18$, $128$, $118$, $28$, $26$, $120$, $126$, $20$, $34$, $-2$, $-14$, $-58$, $-62$, $-10$, $-6$, $-66$, $-54$, $-52$, $-68$, $-4$, $-12$, $-60$, $76$, $40$, $84$, $88$, $70$, $90$, $82$, $42$, $78$, $94$, $74$, $-98$, $-112$, $-44$, $-110$, $-100$, $-16)$\\

\noindent$11n160$ : $3$ | $(6$, $-10$, $12$, $20$, $16$, $-18$, $22$, $8$, $-4$, $2$, $14)$ | $(24$, $52$, $48$, $20$, $44$, $56$, $60$, $54$, $46$, $-66$, $-80$, $-74$, $-72$, $-82$, $-10$, $-16$, $-4$, $50$, $22$, $26$, $84$, $28$, $-64$, $-68$, $-78$, $-76$, $-70$, $-62$, $-42$, $-40$, $86$, $30$, $96$, $94$, $32$, $88$, $100$, $90$, $34$, $92$, $98$, $58$, $-12$, $-14$, $-2$, $-36$, $-6$, $-18$, $-8$, $-38)$\\

\noindent$11n161$ : $3$ | $(6$, $10$, $12$, $20$, $18$, $16$, $-22$, $4$, $2$, $8$, $-14)$ | $(38$, $44$, $110$, $120$, $118$, $112$, $42$, $40$, $36$, $86$, $84$, $34$, $76$, $-54$, $-10$, $-6$, $-58$, $-50$, $-64$, $-66$, $-92$, $-68$, $-102$, $-20$, $106$, $124$, $114$, $116$, $122$, $108$, $46$, $104$, $126$, $90$, $80$, $30$, $-8$, $-56$, $-52$, $-12$, $-4$, $-60$, $-48$, $-62$, $-2$, $-14$, $78$, $28$, $74$, $32$, $82$, $88$, $-18$, $-22$, $-100$, $-70$, $-94$, $-26$, $-96$, $-72$, $-98$, $-24$, $-16)$\\

\noindent$11n162$ : $3$ | $(6$, $-10$, $12$, $22$, $16$, $-18$, $8$, $20$, $-4$, $2$, $14)$ | $(18$, $78$, $22$, $14$, $52$, $50$, $-30$, $-44$, $-36$, $-38$, $-42$, $-32$, $-48$, $-10$, $84$, $86$, $56$, $62$, $90$, $60$, $58$, $88$, $64$, $54$, $-68$, $-28$, $-66$, $-70$, $-4$, $-76$, $-2$, $-72$, $82$, $26$, $24$, $80$, $16$, $20$, $-40$, $-34$, $-46$, $-8$, $-12$, $-6$, $-74)$\\

\noindent$11n163$ : $3$ | $(6$, $-10$, $-12$, $22$, $-16$, $-18$, $8$, $20$, $-4$, $-2$, $14)$ | $(142$, $138$, $180$, $122$, $152$, $128$, $174$, $132$, $148$, $82$, $146$, $134$, $176$, $126$, $154$, $124$, $178$, $136$, $144$, $140$, $-2$, $-36$, $-22$, $-16$, $-170$, $-10$, $-28$, $-30$, $-8$, $-168$, $-166$, $-6$, $-32$, $-26$, $-12$, $-172$, $-14$, $-24$, $-34$, $-4$, $-164$, $184$, $46$, $74$, $68$, $52$, $56$, $64$, $78$, $42$, $188$, $40$, $80$, $62$, $58$, $50$, $70$, $72$, $48$, $60$, $-108$, $-94$, $-156$, $-90$, $-112$, $-116$, $-86$, $-160$, $-98$, $-104$, $-102$, $-100$, $-162$, $-84$, $-118$, $-110$, $-92$, $54$, $66$, $76$, $44$, $186$, $182$, $120$, $150$, $130$, $-114$, $-88$, $-158$, $-96$, $-106$, $-18$, $-20$, $-38)$\\

\noindent$11n164$ : $3$ | $(6$, $-10$, $14$, $16$, $-20$, $-4$, $-18$, $2$, $22$, $-12$, $-8)$ | $(18$, $22$, $58$, $24$, $16$, $44$, $-30$, $-26$, $-54$, $-52$, $-50$, $-56$, $38$, $64$, $62$, $60$, $20$, $-2$, $-8$, $-6$, $-4$, $-32$, $-28$, $14$, $42$, $34$, $66$, $36$, $40$, $-48$, $-12$, $-46$, $-10)$\\

\noindent$11n165$ : $3$ | $(6$, $-10$, $14$, $16$, $20$, $-4$, $18$, $2$, $22$, $12$, $8)$ | $(28$, $32$, $76$, $90$, $82$, $84$, $88$, $78$, $34$, $80$, $86$, $-12$, $-8$, $-70$, $-72$, $-74$, $-68$, $50$, $58$, $22$, $60$, $26$, $30$, $-2$, $-18$, $-16$, $-4$, $-6$, $-14$, $-20$, $-10$, $24$, $56$, $52$, $48$, $92$, $46$, $54$, $-66$, $-36$, $-42$, $-62$, $-40$, $-38$, $-64$, $-44)$\\

\noindent$11n166$ : $3$ | $(6$, $10$, $-14$, $18$, $2$, $20$, $22$, $-4$, $12$, $8$, $16)$ | $(26$, $92$, $96$, $22$, $30$, $88$, $32$, $20$, $-82$, $-80$, $-84$, $-38$, $-44$, $-46$, $-36$, $-86$, $66$, $62$, $102$, $56$, $70$, $58$, $104$, $60$, $68$, $18$, $-16$, $-6$, $-72$, $-2$, $-76$, $-10$, $-12$, $-78$, $-4$, $94$, $24$, $28$, $90$, $98$, $52$, $54$, $100$, $64$, $-34$, $-48$, $-42$, $-40$, $-50$, $-14$, $-8$, $-74)$\\

\noindent$11n167$ : $3$ | $(6$, $10$, $14$, $20$, $2$, $-18$, $4$, $22$, $-12$, $8$, $16)$ | $(82$, $34$, $48$, $42$, $40$, $50$, $32$, $84$, $80$, $36$, $46$, $44$, $38$, $78$, $86$, $-88$, $-64$, $-58$, $-94$, $-52$, $-68$, $-54$, $-96$, $-56$, $-66$, $102$, $70$, $98$, $74$, $106$, $116$, $108$, $76$, $100$, $-8$, $-24$, $-20$, $-4$, $-12$, $-92$, $-60$, $-62$, $-90$, $112$, $110$, $114$, $104$, $72$, $-22$, $-6$, $-10$, $-26$, $-18$, $-2$, $-14$, $-30$, $-28$, $-16)$\\

\noindent$11n168$ : $3$ | $(6$, $10$, $14$, $20$, $4$, $-18$, $2$, $22$, $8$, $-12$, $16)$ | $(14$, $62$, $70$, $76$, $114$, $-90$, $-98$, $-96$, $-92$, $-8$, $-88$, $-100$, $-82$, $-102$, $-86$, $-6$, $-94$, $16$, $12$, $22$, $28$, $104$, $26$, $24$, $106$, $30$, $20$, $10$, $18$, $32$, $108$, $-48$, $-44$, $-42$, $-50$, $-60$, $-34$, $-58$, $-52$, $-40$, $-46$, $64$, $68$, $78$, $112$, $116$, $74$, $72$, $118$, $110$, $80$, $66$, $-84$, $-4$, $-2$, $-38$, $-54$, $-56$, $-36)$\\

\noindent$11n169$ : $3$ | $(6$, $-10$, $-16$, $12$, $-18$, $-2$, $20$, $22$, $-4$, $-8$, $14)$ | $(12$, $10$, $14$, $32$, $-24$, $-28$, $-26$, $-42$, $-38$, $36$, $16$, $34$, $30$, $56$, $-2$, $-46$, $-44$, $-40$, $52$, $50$, $54$, $48$, $8$, $-6$, $-22$, $-18$, $-20$, $-4)$\\

\noindent$11n170$ : $3$ | $(6$, $10$, $16$, $12$, $18$, $2$, $-22$, $-20$, $4$, $8$, $-14)$ | $(38$, $60$, $64$, $72$, $68$, $-46$, $-6$, $-22$, $-2$, $-4$, $62$, $36$, $40$, $12$, $44$, $10$, $42$, $-50$, $-56$, $-58$, $-48$, $-8$, $70$, $66$, $74$, $14$, $20$, $16$, $76$, $18$, $-54$, $-52$, $-26$, $-32$, $-30$, $-28$, $-34$, $-24)$\\

\noindent$11n171$ : $3$ | $(6$, $-10$, $-16$, $12$, $-18$, $-2$, $22$, $20$, $-4$, $-8$, $14)$ | $(60$, $64$, $54$, $28$, $52$, $62$, $-2$, $-44$, $-42$, $-4$, $-10$, $-74$, $-70$, $-82$, $92$, $20$, $12$, $22$, $50$, $26$, $56$, $66$, $58$, $24$, $-72$, $-84$, $-30$, $-80$, $-68$, $-76$, $-34$, $-32$, $-78$, $16$, $96$, $88$, $86$, $98$, $14$, $18$, $94$, $90$, $-36$, $-48$, $-38$, $-8$, $-6$, $-40$, $-46)$\\

\noindent$11n172$ : $3$ | $(6$, $10$, $-16$, $12$, $20$, $2$, $18$, $-4$, $22$, $8$, $14)$ | $(60$, $18$, $20$, $58$, $32$, $54$, $24$, $-82$, $-94$, $-72$, $-92$, $-84$, $-80$, $-96$, $-100$, $-76$, $-88$, $56$, $22$, $16$, $26$, $52$, $30$, $102$, $28$, $-98$, $-78$, $-86$, $-90$, $-74$, $-46$, $-48$, $-40$, $-14$, $-38$, $-4$, $104$, $122$, $116$, $110$, $66$, $68$, $108$, $118$, $120$, $106$, $70$, $64$, $112$, $114$, $62$, $-2$, $-6$, $-36$, $-12$, $-42$, $-50$, $-44$, $-10$, $-34$, $-8)$\\

\noindent$11n173$ : $3$ | $(6$, $-10$, $16$, $12$, $20$, $-2$, $18$, $4$, $22$, $8$, $14)$ | $(96$, $48$, $36$, $34$, $50$, $98$, $94$, $46$, $38$, $86$, $88$, $40$, $44$, $92$, $100$, $52$, $-68$, $-122$, $-112$, $-58$, $-106$, $-104$, $-60$, $-114$, $-120$, $-66$, $-70$, $-20$, $126$, $76$, $132$, $142$, $138$, $82$, $84$, $90$, $42$, $-24$, $-16$, $-4$, $-8$, $-32$, $-110$, $-56$, $-108$, $-102$, $-62$, $-116$, $-118$, $-64$, $140$, $130$, $74$, $128$, $54$, $124$, $78$, $134$, $144$, $136$, $80$, $-6$, $-18$, $-22$, $-72$, $-26$, $-14$, $-2$, $-10$, $-30$, $-28$, $-12)$\\

\noindent$11n174$ : $3$ | $(6$, $-10$, $16$, $14$, $20$, $-4$, $18$, $2$, $22$, $8$, $12)$ | $(30$, $32$, $28$, $76$, $22$, $38$, $126$, $146$, $-122$, $-104$, $-100$, $-52$, $-72$, $-48$, $-96$, $-94$, $-46$, $-70$, $-54$, $-102$, $20$, $78$, $92$, $158$, $134$, $138$, $154$, $88$, $82$, $16$, $84$, $86$, $18$, $80$, $90$, $156$, $136$, $-50$, $-98$, $-106$, $-120$, $-14$, $-124$, $-12$, $-118$, $-108$, $-2$, $160$, $132$, $140$, $152$, $150$, $142$, $130$, $34$, $26$, $74$, $24$, $36$, $128$, $144$, $148$, $-62$, $-56$, $-68$, $-44$, $-6$, $-112$, $-114$, $-8$, $-42$, $-66$, $-58$, $-60$, $-64$, $-40$, $-10$, $-116$, $-110$, $-4)$\\

\noindent$11n175$ : $3$ | $(6$, $-10$, $-16$, $14$, $-20$, $-18$, $-2$, $22$, $-12$, $-4$, $-8)$ | $(74$, $36$, $28$, $82$, $68$, $80$, $30$, $34$, $76$, $72$, $38$, $26$, $-98$, $-54$, $-64$, $-46$, $-44$, $-62$, $-56$, $112$, $120$, $106$, $90$, $104$, $118$, $114$, $100$, $86$, $110$, $122$, $108$, $88$, $102$, $116$, $-50$, $-40$, $-58$, $-60$, $-42$, $-48$, $-66$, $-52$, $-22$, $-4$, $-14$, $32$, $78$, $70$, $84$, $-24$, $-6$, $-12$, $-92$, $-16$, $-2$, $-20$, $-96$, $-8$, $-10$, $-94$, $-18)$\\

\noindent$11n176$ : $3$ | $(6$, $10$, $-16$, $18$, $20$, $2$, $8$, $-4$, $22$, $12$, $14)$ | $(52$, $58$, $90$, $92$, $56$, $54$, $94$, $88$, $60$, $-36$, $-34$, $-16$, $-4$, $-8$, $-12$, $96$, $86$, $62$, $84$, $98$, $82$, $26$, $-6$, $-14$, $-32$, $-38$, $-40$, $-80$, $-42$, $-74$, $-72$, $20$, $50$, $30$, $48$, $22$, $18$, $24$, $46$, $28$, $-10$, $-2$, $-64$, $-70$, $-76$, $-44$, $-78$, $-68$, $-66)$\\

\noindent$11n177$ : $3$ | $(6$, $-10$, $16$, $18$, $20$, $-2$, $8$, $4$, $22$, $12$, $14)$ | $(88$, $82$, $72$, $40$, $70$, $84$, $90$, $86$, $68$, $120$, $-54$, $-62$, $-4$, $-10$, $-16$, $-102$, $-104$, $-18$, $-110$, $-96$, $114$, $26$, $126$, $130$, $30$, $34$, $78$, $76$, $36$, $38$, $74$, $80$, $32$, $-50$, $-94$, $-42$, $-98$, $-108$, $-20$, $-106$, $-100$, $-44$, $-92$, $-48$, $-46$, $128$, $28$, $112$, $116$, $24$, $124$, $66$, $122$, $22$, $118$, $-58$, $-8$, $-6$, $-60$, $-56$, $-52$, $-64$, $-2$, $-12$, $-14)$\\

\noindent$11n178$ : $3$ | $(6$, $-10$, $16$, $18$, $20$, $-4$, $8$, $2$, $22$, $14$, $12)$ | $(52$, $44$, $28$, $34$, $62$, $32$, $30$, $46$, $50$, $166$, $154$, $112$, $-138$, $-126$, $-66$, $-120$, $-64$, $-124$, $-140$, $-22$, $-136$, $-128$, $-68$, $-86$, $-84$, $-70$, $-130$, $-134$, $-20$, $-142$, $-122$, $160$, $118$, $106$, $148$, $100$, $56$, $40$, $24$, $38$, $58$, $102$, $146$, $104$, $60$, $36$, $26$, $42$, $54$, $-132$, $-18$, $-144$, $-14$, $-4$, $-94$, $-76$, $-78$, $-92$, $-6$, $-12$, $162$, $158$, $116$, $108$, $150$, $170$, $168$, $152$, $110$, $114$, $156$, $164$, $48$, $-16$, $-2$, $-96$, $-74$, $-80$, $-90$, $-8$, $-10$, $-88$, $-82$, $-72$, $-98)$\\

\noindent$11n179$ : $3$ | $(6$, $-10$, $-20$, $14$, $-2$, $-18$, $8$, $22$, $-12$, $-4$, $16)$ | $(18$, $38$, $40$, $16$, $20$, $26$, $22$, $-30$, $-54$, $-52$, $-28$, $-32$, $-34$, $44$, $66$, $64$, $46$, $42$, $-50$, $-56$, $-48$, $-10$, $-12$, $24$, $58$, $60$, $68$, $62$, $-4$, $-36$, $-2$, $-6$, $-14$, $-8)$\\

\noindent$11n180$ : $3$ | $(6$, $-12$, $-16$, $14$, $-18$, $-20$, $-2$, $22$, $-4$, $-10$, $-8)$ | $(14$, $12$, $42$, $-30$, $-38$, $-26$, $-58$, $-24$, $-36$, $-32$, $50$, $44$, $56$, $86$, $82$, $52$, $48$, $46$, $54$, $84$, $-28$, $-60$, $-66$, $-76$, $-64$, $-62$, $-74$, $-68$, $-2$, $88$, $80$, $20$, $78$, $16$, $10$, $40$, $8$, $18$, $-34$, $-22$, $-6$, $-72$, $-70$, $-4)$\\

\noindent$11n181$ : $3$ | $(6$, $-12$, $-16$, $14$, $-20$, $-18$, $-2$, $22$, $-4$, $-10$, $-8)$ | $(44$, $40$, $36$, $38$, $46$, $50$, $-52$, $-10$, $-54$, $-8$, $-2$, $42$, $18$, $14$, $12$, $16$, $20$, $62$, $-60$, $-34$, $-58$, $-56$, $-32$, $-26$, $-30$, $-28$, $48$, $66$, $22$, $64$, $-6$, $-4$, $-24)$\\

\noindent$11n182$ : $3$ | $(6$, $-12$, $16$, $18$, $14$, $-4$, $20$, $22$, $2$, $8$, $10)$ | $(44$, $36$, $52$, $28$, $124$, $30$, $50$, $38$, $42$, $46$, $34$, $128$, $-72$, $-54$, $-122$, $-58$, $-68$, $-112$, $-118$, $-62$, $-64$, $-116$, $-114$, $-66$, $-60$, $-120$, $88$, $82$, $78$, $92$, $134$, $100$, $136$, $94$, $76$, $84$, $86$, $-22$, $-106$, $-10$, $-8$, $-108$, $-24$, $-110$, $-6$, $-12$, $-104$, $-20$, $-2$, $-16$, $40$, $48$, $32$, $126$, $26$, $130$, $96$, $138$, $98$, $132$, $90$, $80$, $-56$, $-70$, $-74$, $-4$, $-14$, $-102$, $-18)$\\

\noindent$11n183$ : $3$ | $(6$, $-14$, $10$, $-18$, $2$, $22$, $20$, $-4$, $12$, $-8$, $16)$ | $(54$, $20$, $88$, $82$, $46$, $80$, $90$, $18$, $-74$, $-40$, $-34$, $-10$, $-8$, $-32$, $86$, $84$, $48$, $78$, $92$, $76$, $50$, $58$, $24$, $-68$, $-4$, $-14$, $-38$, $-36$, $-12$, $-6$, $-30$, $28$, $60$, $26$, $22$, $56$, $52$, $-16$, $-2$, $-70$, $-44$, $-66$, $-62$, $-64$, $-42$, $-72)$\\

\noindent$11n184$ : $3$ | $(6$, $-14$, $10$, $18$, $2$, $22$, $20$, $-4$, $12$, $8$, $16)$ | $(48$, $26$, $98$, $36$, $38$, $56$, $40$, $34$, $96$, $28$, $46$, $50$, $-70$, $-64$, $-76$, $-92$, $-80$, $-60$, $-118$, $-58$, $-82$, $-90$, $-74$, $-66$, $-68$, $-72$, $-88$, $-84$, $152$, $116$, $148$, $140$, $108$, $160$, $106$, $104$, $158$, $110$, $142$, $146$, $114$, $154$, $100$, $102$, $156$, $112$, $144$, $-78$, $-62$, $-120$, $-136$, $-122$, $-24$, $-2$, $-20$, $-126$, $-132$, $-14$, $-8$, $150$, $138$, $52$, $44$, $30$, $94$, $32$, $42$, $54$, $-86$, $-4$, $-18$, $-128$, $-130$, $-16$, $-6$, $-10$, $-12$, $-134$, $-124$, $-22)$\\

\noindent$11n185$ : $3$ | $(6$, $18$, $16$, $14$, $-20$, $4$, $2$, $22$, $12$, $-8$, $-10)$ | $(44$, $154$, $28$, $144$, $34$, $32$, $142$, $30$, $156$, $162$, $104$, $88$, $-52$, $-64$, $-74$, $-140$, $-78$, $-60$, $-56$, $-20$, $-48$, $-68$, $-70$, $164$, $106$, $86$, $90$, $102$, $160$, $158$, $100$, $92$, $84$, $108$, $166$, $110$, $82$, $94$, $98$, $-10$, $-16$, $-4$, $-120$, $-130$, $-132$, $-118$, $-6$, $-14$, $-12$, $-8$, $-116$, $-134$, $-128$, $-122$, $-2$, $-18$, $-58$, $36$, $146$, $26$, $152$, $42$, $46$, $40$, $150$, $24$, $148$, $38$, $112$, $80$, $96$, $-76$, $-62$, $-54$, $-22$, $-50$, $-66$, $-72$, $-138$, $-114$, $-136$, $-126$, $-124)$\\

\bibliographystyle{hplain}
\bibliography{bridgea1}

\end{CJK*}\end{document}